\documentclass[12pt,leqno,]{article}

\usepackage{amstext}

\usepackage{amsthm}
\usepackage{amsopn}
\usepackage{amsfonts}
\usepackage{amsmath}
\usepackage{amssymb}
\usepackage{amsfonts}
\usepackage{amscd}

\theoremstyle{definition}

\theoremstyle{remark}

\numberwithin{equation}{section}

\newcommand{\be}{\begin{equation}}
\newcommand{\ee}{\end{equation}}
\def\q{\quad}
\def\qed{\hfill $\square$}
\def\qbar{{{\scriptstyle Q}\kern-.45em{\vrule height.41em width.035em
depth-.03em}}~}
\def\Cbar{\text{\sl C\kern-.35em{\vrule height.63em width.05em
depth-.033em}}~}
\def\cbar{{{\scriptstyle C}\kern-.41em{\vrule height.42em width.035em
depth-.03em}}~}
\def\ibid{\hbox to .5truein{\hrulefill}}
\def\IH{\text{{\rm I}\kern-.13em{\rm H}}}

\def\IF{\Bbb F}
\def\Z{\mathbb Z}
\def\Qbar{\mathbb Q}

\def\IR{{\Bbb R}}

\def\pd{{\rm pd}}
\def\twoheaddown{\downarrow\kern-0.78em\raise0.25em\hbox{$\downarrow$}}
\def\headtaildown{\downarrow\kern-0.79em\raise
0.5em\hbox{$\ssize\curlyvee$}}
\def\longdownarrow{\downarrow\kern-0.73em\raise0.5em\hbox{$\vert$}}
\def\hookdownarrow{\longdownarrow\kern-.12em\raise0.5em\hbox{$^{^\cap}
$}}
\def\htda{\downarrow\kern-0.79em\raise 0.5em\hbox{$\tiny\curlyvee$}} 

\def\vs{\vskip.3cm}
\def\vsk{\vskip.6cm}

\def\noi{\noindent}

\def\pr{\prime}
\def\bs{\backslash}

\def\la{\langle}
\def\ra{\rangle}
\def\wt{\widetilde}

\def\wh{\widehat}
\def\lf{\lfloor}
\def\rf{\rfloor}
\def\ol{\overline}
\def\os{\overset}
\def\ul{\underline}
\def\us{\underset}

\def\lra{\longrightarrow}
\def\lla{\longleftarrow}
\def\Lra{\Longrightarrow}
\def\Llra{\Longleftrightarrow}
\def\Qp{{\Qbar_p}}
\def\nl{\nolimits}

\def\SD{{\cal D}}
\def\CF{{\cal F}}
\def\CL{{\cal L}}

\def\SM{{\cal M}}
\def\SO{{\cal O}}
\def\SQ{{\cal Q}}

\def\CV{{\cal V}}

\def\aug{\text{\rm aug}}

\def\cent{\text{\rm cent}}
\def\defl{\text{\rm defl}}
\def\Det{\text{\rm Det}}
\def\det{\text{\rm det}}
\def\dim{\text{\rm dim}}
\def\disc{\text{\rm disc}}
\def\Ext{\text{\rm Ext}}
\def\Gal{\text{\rm Gal}}

\def\Hom{\text{\rm Hom}}
\def\id{\text{\rm id}}
\def\im{\text{\rm im}}
\def\ind{\text{\rm ind}}
\def\infl{\text{\rm  infl}}
\def\Maps{\text{\rm Maps}}
\def\Nr{\text{\rm Nr}}
\def\nr{\text{\rm nr}}

\def\pd{\text{\rm pd}}
\def\rad{\text{\rm rad}}

\def\res{\text{\rm res}}
\def\sgn{\text{\rm sgn}}

\def\Sym{\text{\rm Sym}}
\def\St{\text{\rm St}}
\def\supp{\text{\rm supp}}
\def\Tr{\text{\rm Tr}}
\def\tr{\text{\rm tr}}

\def\v{\text{\rm v}}
\def\elra{\hbox to 2in{\rightarrowfill}}
\def\ella{\hbox to 2in{\leftarrowfill}}
\def\hrf{\hbox to 2in{\hrulefill}}
\def\hdotfill{\leaders\hbox to 1em{\hss .\hss}\hfill}

\def\bd{{\pmb d}}
\def\be{{\pmb e}}

\def\bH{{\pmb H}}

\def\fg{{\frak g}}
\def\fG{{\frak G}}
\def\fq{{\frak q}}
\def\fM{{\frak M}}
\def\fO{{\frak O}}
\def\fo{{\frak o}}

\def\fQ{{\frak Q}}
\def\fR{{\frak R}}
\def\fV{{\frak V}}

\begin{document}
\title{On Uniqueness in Equivariant Iwasawa Theory}

\vsk
\author{J\"urgen Ritter and  Alfred Weiss\footnote{We acknowledge
financial support provided by NSERC and the Faculty of Science at
the University
\newline of 
Alberta.}}

\date{} \maketitle

\begin{abstract}
This sequel to our paper `On the ``main conjecture'' of equivariant
Iwawasa theory' studies the possibility of the vanishing of
$SK_1(QG),$ when $QG$ is the total ring of fractions of the
Iwasawa algebra $\Lambda  G =\Z_p\big[[G]\big],$ with $G$ the Galois
group of the Galois extension $K/k$ of loc.~cit.  The vanishing is 
equivalent to $SK_1(D)=0$ for all  division algebras $D$ in
the Wedderburn components of $QG,$ so can be studied via the
classification $D_{r,s,\IF}$ of these $D$'s, which provides
 a crossed
product order $\Delta$ in $D$ and the valuation
$\v_\bullet$ with  valuation ring
$\Delta_\bullet\,  .$

$\Delta  $ contains a special element $\Pi,$ which generates a
maximal subfield of $D$ over its centre with $\v_\bullet (\Pi)=1$
and $\Pi \Delta  \Pi^{-1}=\Delta  .$ The induced $\Pi$-filtration on
$\Delta  $ enables a study of the reduced norm  built
on a congruence for $nr(d)\!\!\!\mod \Pi\Delta $ for $d\in \Delta
.$

Given $d\in \Delta  ^\times$ with $nr(d) =1,$ and $n\ge 1$ maximal
with $d\in 1+\Pi^n\Delta  $ (called the level of $d),$ the main
problem is to find a suitable commutator product $c\equiv
d\!\!\!\mod \Pi^{n+1}\Delta  .$ Then $c^{-1}d\equiv 1\!\!\!\mod
\Pi^{n+1}\Delta  $ hence, setting $d^\pr = c^{-1}d,$ has
$nr(d^\pr)=1$ and level $n^\pr >n.$  Repetition ends in
$[\Delta  ^\times,\Delta  ^\times]$ when the level gets sufficiently
large.

\vskip.4in
\centerline{\bf Introduction}
\end{abstract}

\setlength{\baselineskip}{20pt}
\noi
Let $p$ be an odd prime number and $K/k$ a Galois extension of
totally real number fields such that $[k:\Qbar]<\infty, \, K\supset
k_\infty  $ and $[K:k_\infty  ]< \infty  ,$ with $k_\infty  $  the $\Z_p\text{-cyclotomic}$ extension of $k.$ Assuming that
Iwasawa's $\mu\text{-invariant}$ of $K/k$ vanishes, we have proved,
in [RW9]\footnote{We prefer to here use `$p$' rather than `$l$' of
[RW9].}, a weak form of the `main conjecture' of equivariant Iwasawa
theory which, roughly speaking, reads

\vs
(1) {\it the Iwasawa $L\text{-function}$ $L_{K/k,S}$ of $K/k$ knows
the refined Iwasawa module $\mho_{K/k,S}\,,$ }

\vs
\noi
where both, $L_{K/k,S}$ and $\mho_{K/k,S}\,,$ are defined with
respect to a given finite set $S$ of primes of $k$ which contains
all archimedean primes and those which ramify in $K.$

\vs\noi
More precisely (compare [RW9, p.1]), consider
$$
\begin{aligned}
&K_1(\SQ G) \os \partial\rightarrow K_0T(\Lambda   G)\ni \mho_{K/k,S}\\
&\Det \,\downarrow\\
L_{K/k,S} \in \,\Hom^*&\big(R_p(G),(\SQ^c\Gamma_k)^\times\big)
\end{aligned}
$$
in which $G$ and $\Gamma_k$ are the Galois groups of $K/k$ and
$k_\infty  /k,$ respectively, $\SQ G$ is the total ring of fractions
of the Iwasawa algebra $\Lambda  G \os{\rm def}= \Z_p[[G]]$ and
$R_p(G)$ the character ring of the characters of
$\Qbar_p^{\;c}\text{-representations}$ of $G$ with open kernel
$(\Qbar_p ^{\;c}$ is a fixed algebraic closure of
$\Qbar_p\,;\SQ^c\Gamma_k$ is short for $\Qbar_p^{\;c} \,\otimes
_{\Qbar_p}$ $\SQ\Gamma_k\,;$ $\partial$ is taken from the
localization sequence of $K\text{-theory}).$

\vs\noi
Now, the weakened `main conjecture' of equivariant Iwasawa theory
says that there is a common preimage $\Theta \in K_1(\SQ G)$ of
$L_{K/k,S}$ and $\mho_{K/k,S}\,,$ i.e.,
$$
\Det(\Theta) = L_{K/k,S}\q\text{and}\q \partial (\Theta) =
\mho_{K/k,S}\,;
$$
however (1) is stronger by requiring $\Theta$ to be unique.  The
goal of this paper is to prepare for a 
proof that $\Det$ is
injective so that the last phrase

\vs
(2) {\it up to its uniqueness statement}

\vs\noi
in {\sc Theorem} [RW9, p.1] could be erased.

\vs\noi 
The assertion ``$L$ knows $\mho$'' only declares existence of a
bridge $\Theta.$ A stronger result apparently depends on how unique
$\Theta$ is (and, perhaps, on which question about $\mho$ we ask), 
but the dream of constructing an explicit $\Theta$ seems hopeless
without it.
So we start by investigating absolute uniqueness (while being open to
something weaker).  Such questions are also relevant in the context
of [FK], because of [Bu] and [BV].  

\vs\noi
We recall that $\SQ G$ is a semisimple algebra whose Wedderburn
components are finite dimensional central simple algebras $A$ over
fields of cohomological dimension 3 (see Proposition B1).  Since our
map Det is taking reduced norms on the individual Wedderburn
components of $\SQ G,$ and since $SK_1(A) =0$ should hold for all
such $A,$ by a conjecture of Suslin, the kernel of Det is
$SK_1(\SQ G).$  Thus absolute uniqueness of $\Theta$ amounts to
proving $SK_1(\SQ G)=0.$ 

\vs\noi
\S 1 starts by fixing the most basic notation of the paper, and then
turns to the fundamental problem of identifying those division
algebras $D$ that occur in the full rings $D_{k\times k}$ of
matrices that are Wedderburn components in the algebras $\SQ G$ with
$G$ among the Galois groups of this introduction.  This reduces the
problem of proving $SK_1(\SQ G) =0,$ for all such $G,$ to proving
that $SK_1(D)=0$ for all of the identified $D\text{'s.}$  It turns out
that all these $D\text{'s}$ arise from algebras of the form
$F\otimes _{\Qbar_p}\!\!\SQ G$ where $F$ is a finite unramified
extension of $\Qbar_p$ and $G$ a pro-$p$ group of our Galois kind.
 \S1 also prepares for our study, in \S4, of the reduced norm on $D$ by
introducing an elementary valuation theoretic function $\Phi,$ and for
some detail about $A[u]$ needed in \S7.
It ends with some further notation. 

\vs\noi
In \S2 we consider a natural order $\Delta  $ in $D$ and prepare for
the $\Pi$-filtration on $\Delta  ,$ where $\Pi$ is the prime element
$\Pi=1-\zeta  _{r+s}$ of $F(\zeta  _{r+s}),$ with $F$ the unramified
extension of $\Qbar_p$ that comes with $D$ above, and $\zeta 
_{r+s}$ the root of unity of order $p^{r+s}$ that occurs in 
Theorem 1A; moreover, $p^r$ is the order of the group
of $p$-power roots of unity in the centre of $D$ and $p^s$ the Schur
index of $D.$  The goal of \S2 was to prove that

\vs\noi
$\Delta  $ {\it is a maximal order in} $D$ {\it and a local ring
with radical generated by} $\Pi$ {\it and} $u-1.$

\centerline{{\it The canonical map} \q $SK_1(\Delta  ) \to SK_1(D)$
\q {\it is surjective.}}

\vs\noi
However $\Delta$ is not quite maximal and the canonical map  not
necessarily surjective.  
Also it seems likely that there is a
  valuation $v_\bullet:D\to \Z
\cup\{\infty\},$ so this 
investigation continues;
it becomes crucial 
after \S6.  Finally, we
prepare for congruences on the crossed product order $\Delta 
=B*{\fg},$ by studying the $\fg$-action on $B,$ with $\fg
=\,\Gal\big(F(\zeta_{r+s})/F(\zeta  _r)\big).$

\newpage\noi
In the next 
section, \S3, we prove that the canonical map
$K_1(\Delta  ) \to K_1(\Delta  /\pi_1\Delta  )$ embeds $SK_1(\Delta 
)$ into  $K_1(\Delta  /\pi_1\Delta  ),$ i.e.,
$$
SK_1(\Delta  ) \to K_1(\Delta  /\pi_1\Delta  )\q \text{\it is
injective},
$$
with $\pi_1=1 - \zeta  _1$ and $\zeta  _1$ is a  primitive $p$-th
root of unity.  We also describe the cokernel.

\vs\noi
As a consequence of the last two displayed statements, $SK_1(D) = 0$
if every $d\in\Delta  ^\times$ with $\nr(d) =1$ can be multiplied by
some product $c$ of commutators in $[D^\times,D^\times ]\cap
\Delta  ^\times$ so that $cd\in 1+ \Pi^{p^{r+s-1}}\Delta  ,$ where
`$\nr$' denotes the reduced norm $D^\times\to \,\cent (D)^\times.$

\vs\noi Then, in \S4, we study the reduced norm of $d=
 1+\Pi^n\,\sum^{p^s-1}_{j=0}\,b_ju^j,$ with level $n\ge 1$ and all $b_j\in
B.$  
The conclusion, Theorem 4G, is a congruence for
$\nr(d)-1\!\!\!\mod\,\pi^{\Phi(n)+1}A,$ with $A=\,\cent(\Delta  )$ and
$\pi = \pi_r = 1 - \zeta  _r\,.$  The other side of this congruence
is an explicit expression in relatively  few of the
coefficients $b_j\,.$  First consequences are that
 $ \nr(d) \equiv 1\!\!\!\mod\pi^{\Phi(n)}A,$
and that $\nr(d)=1$ implies that $n\ge p^r-1.$ 
Later, in Theorem~6D, it is essential for dealing with critical
levels.

\vs\noi
The goal of \S5 is learning how to use Theorem 4G.  The explicit
expressions above depend only on partial sums $\lambda  _n(b_j)$ of
some of the coordinates of the $b_j$'s in the $A$-basis
$\{\Pi^j:0\le j< p^s\}$ of $B.$  Thus, for each $d\in 1+\Pi^n\Delta 
$ with $\nr(d) =1,$ Proposition~5D gives explicit relations between
the $\lambda  _n(b_j)$'s which  need to become
commutators  $\!\mod\Pi^{n+1}\Delta  .$ The rest of
this section prepares for the next one.

\vs\noi
In \S6 we succeed in some of the simplest cases, from the Theorem~4G
perspective,
when the $\lambda 
_n(b_j)$'s have only one term, to prove that $d,$ as above, is
congruent to some $\Pi$-commutator product $\!\!\mod \Pi^{n+1}\Delta,$  
 most notably when
$n$ is the critical level $p^{r+w}-1,$ with $0\le w\le s-1.$ 
Conversely, when there is some $d\in 1+\Pi^n\Delta  ,$ with
$\nr(d)=1,$ which is not congruent to any $\Pi$-commutator product
$\!\!\!\mod\Pi^{n+1}\,\Delta, $ 
then the level $n$ of $d$ has type $t,$ i.e.
$n= p^{r+t}+p^tm-1$ with $0\le t< s-1$ and $1\le m \not\equiv 0 \mod p.$ 
This accounts, for a
particular $\Delta,$ for most levels: more precisely, the
proportion of $n$'s of type $\ge t$ is roughly $1/p^t.$

\vs
\noi
The study, in \S7 of the $n$'s of type $t$ by valuation theory
methods involving $v_\bullet$ and using 
linearizations $\!\!\!\mod \Pi^{p^{r+t}}\Delta$ of
the special commutators
$[u,1-\Pi^{p^tm}x]$
to search for solutions $x\in \Delta  $ of $[u,1-\Pi^{p^tm}x] \equiv
d\!\!\!\mod \Pi^{n+1}\Delta,$ when $d$ has $nr(d)=1$ and level $n$
of type $t,$
succeeds when $t=0 $ or $1.$ 
The type 1 case of Corollary~
7D suggests a programme for higher types, and
strengthens
the conclusion of \S6 to $2\le t< s-1.$  

\vs\noi
In \S8, we take first steps toward the programme above for type
$t\ge 2.$  Imitating the proof of Corollary~7D, we study the
solutions of the linearization
congruence of Proposition~7C(iii).
This study is built on the equivalence relation induced by a special
$f:L\to\Z,$ with $L$ indexing the left side of the congruence, and a
special $\tau  :L\to L$ describing the $f$-equivalence classes.
Carrying this out, in preparation
for the
congruences of Proposition~7C(ii) in \S9, exhausts this section.

\vs\noi
\S9 is about realizing the programme of Corollary~7D for $t\ge 2$
with the new tools of \S8.
This seems to work, 
except at the
very end, so that the reduction step to Proposition~7C(ii) is not quite
reached yet.  Appendix~C gives a brief orientation of what comes next.

\vs\noi
This paper has grown since 2014 around roughly bi-annual meetings
focused on the ``next stage'', not always well-chosen (causing
interpolations in the text), and the march of time. This seemed
necessary in order to discover what will eventually be needed. 
However, the length of the journey
required making changes which took
more and more pages and time.  

\vs\noi
J\"urgen Ritter had a tragic accident on March~9, 2021.  
The shock of this must be tempered with gratitude for the
comradery of the productive years.  Carrying
this quest with him has always made it more accessible and
interesting.

\vs\noi
Continuing the task becomes harder, but the resolve persists.  It
leaves much
 later material, especially in \S8 and \S9,
unfinished, which has now been remedied to some extent.
But  the ``next stage'' momentum has not been worked out.
 Uniqueness is still plausible.

\section{Notation,  Reduction to Pro-$p$ Groups, $\Phi$ and         
$A[u]$}\label{Sec1}

Fix primitive $p^m$-th roots of unity
$\zeta_m$ with $\zeta  _0=1$ and $\zeta  ^p_m = \zeta  _{m-1}$ for $m\ge
1.$ Let ${\fO}$ be the ring of integers in its quotient field
${\fq}({\fO}) =F,$ which is a finite unramified extension field of
$\Qbar_p$ with residue field $\IF = {\fO}/p\fO,$ and note 
 that $F$ and
$\Qbar_p(\zeta  _m)$ are linearly disjoint over $\Qbar_p\,.$ 

\vs\noi
Define, with a fixed indeterminate $T,$

\vskip.25cm
\centerline
{$
\fQ _m =\fO [\zeta    _m] \;\text{and} \; A=\fO_r[[T]], \;
B=\fO_{r+s}[[T]]\;
\text{with} \; r\ge  1, s\ge 1;\; \pi_m=1 -
\zeta    _m\; \text{for}  \;m\le r+s.
$}

\vsk
\noi
Note that $B=A[\pi_{r+s}],$ that $B/\pi_{r+s}B=A/\pi_{r}A
=\IF[[T]],$ and the cyclic Galois group $\fg = \la g\ra$ of
$\fq(B)/\fq(A)$ is generated by the special $g$ that fixes
$\fq(\fO_r[[T]])$ and sends $\zeta  _{r+s}$ to $\zeta 
^{1+p^r}_{r+s}\,.$  

\vsk\noi
Now define $D_{r,s,\IF}$ to be the cyclic crossed
product algebra (cf. [Re], (30.1))
$$
D\!=\!\big(\fq(B)/\fq(A),g,1\!+\!T\big)\!=\!
\bigoplus^{p^s-1}_{j=0}\,\fq(B)\!\cdot\! u^j\;
\text{with}\; ub\! = \;^gbu,\, \text{for all}\, b\!\in\! \fq(B),
\; \text{and} \;u^{p^s}\! =\! 1+T.
$$

\vskip.3cm\noi
{\sc Theorem 1A.} 
{\it $SK_1(\SQ G) =0$ for all of the
general profinite groups $G$ of the introduction would follow from the
vanishing of $SK_1(D_{r,s,\IF})$ for all of the $D_{r,s,\IF}$
defined above..}

\vs
\begin{proof}
By [RW4, p.167, Corollary] it suffices to show that $SK_1(\SQ
G^{\pr}) =0$ for all open $\Qbar_p$-elementary subgroups $G^\pr$ of
the $G$ above; here $\Qbar_p$-elementary means either
$\Qbar_p$-$p$-elementary or $\Qbar_p$-$q$-elementary for some prime
$q\ne p$ (cf. [RW4, beginnings of \S2, p.158, or \S3, p.164]).  By
[Lau2, Theorem 1 and 2], we have

\vs
\begin{enumerate}
\item If $G^\pr$ is $\Qbar_p-p$-elementary, then $SK_1(\SQ G^\pr)
=0$ provided that $SK_1\big(F\otimes_{\Qbar_p} \SQ (V)\big) =0$ for some
 finite unramified extension field $F/\Qbar_p$ and certain
open pro-$p$ subgroups $V$ of $G^\pr.$

\item If $G^\pr$ is $\Qbar_p$-$q$-elementary, then $SK_1(\SQ G^\pr)=0.$
\end{enumerate}

\noi
This reduces us further to showing that $SK_1(F\otimes _{\Qbar_p} \SQ
G)=0$ whenever $G$ is a pro-$p$ group of our type and $F/\Qbar_p$
is finite unramified.

\vs\noi
By [Lau1, Theorem 1(iv)], every simple component $E_{n\times n}$ of
$\SQ G,$ for $G$ pro-$p$ as above, has skew field $E$ a crossed
product
$$
E= \big(\SQ_{\Qbar_p(\eta)} (\Gamma  ^{w_\chi})/\SQ_L(\Gamma 
^{w_\chi  }) ,\sigma  _{v_\chi},\gamma 
^{w_\chi}\big)
$$
where $\Qbar_p (\eta  )/L$ is a finite cyclic extension with Galois
group $\la \sigma  _{v_\chi}\ra,$ and $\eta  $ an irreducible
constituent of $\res^H_G\chi.$  To accommodate the $F\otimes_{\Qbar_p} 
\SQ G$ of the previous paragraph, observe that if $E_{n\times n}$ is
as above, then $F\otimes _{\Qbar_p} E_{n\times n} \simeq
(F\otimes_{\Qbar_p} E)_{n\times n}$ with
$$
F\otimes _{\Qbar_p} E\simeq \big(\SQ_{F(\eta  )}(\Gamma 
^{w_\chi})/\SQ_{FL}(\Gamma  ^{w_\chi  }),\sigma  _{v_\chi  },\gamma 
^{w_\chi  }\big) \leqno{\rm(1A1)}
$$
because $F\otimes_{\Qbar_p}\Qbar_p(\eta  ),$ $F\otimes
_{\Qbar_p}L$ are the composite fields $F(\eta  ),$ $FL$ and the
Galois extension $F(\eta  )/FL$ (and thus also $\SQ_{F(\eta 
)}(\Gamma  ^{w_\chi  })/\SQ_{FL}(\Gamma  ^{w_\chi  }))$ has Galois
group $\la\sigma  _{v_\chi  }\ra.$ We must show that
$F\otimes_{\Qbar_p} E$ is isomorphic to one of the division algebras
$D_{r,s,\IF}\,.$

\vs\noi
For this consider the Iwasawa isomorphism $\Z_p[[\Gamma  ^{w_\chi 
}]]\to \Z_p[[T]]$ induced by $\gamma  ^{w_\chi  }\mapsto 1+T$ as an
identification (cf. [RW2, Lemma 1]).  We also identify
$\Qbar_p(\eta),L$ with $\Qbar_p(\zeta  _{r+s}),\Qbar_p(\zeta  _r)$
for suitable $r\ge 1, s\ge 1,$ namely with $p^{r+s}=w_\chi  \,,$
$p^r=v_\chi  \,.$  For this to be possible we need $L$ to
necessarily contain a $p^{\rm th}$ root of unity if $\eta  \ne 1:$
this follows from [CRII, 74.17] by applying $L\otimes _\Qbar $ --
to $\Qbar[H];$ if $\eta  =1,$ then $w_\chi  =1$ implies, by (1A1),
that $F\otimes _{\Qbar_p}E =\SQ_F(\Gamma  )$ with $SK_1\big(\SQ_F(\Gamma 
)\big) =0.$

\vs\noi
Furthermore, $\SQ_{F(\eta  )}(\Gamma  ^{w_\chi  }),$
$\SQ_{FL}(\Gamma  ^{w_\chi  })$ are then identified with $\fq (B),
\fq(A)$ in our earlier notation.  Thus
$$
F\otimes _{\Qbar_p} E\simeq \big(\fq(B)/\fq(A),\sigma  ,1+T)
$$
with $\sigma  $ the automorphism of $\fq(B)/\fq(A)$ corresponding
to $\sigma  _{v_\chi  } $ under our identification.   Since $\sigma 
,g$ both generate $\fg,$ we have $g=\sigma  ^{m}$ with $p
\nmid m.$

\vs\noi
Now [Re, Theorem~30.4(i)] implies that
$$
F\otimes_{\Qbar_p} E\simeq \big(\fq
(B)/\fq(A),g,(1+T)^m\big),\;\text{with}  \q p\nmid m.
\leqno{\rm (1A2)}
$$
By [La, p.98], $(1+T)^z = \sum^\infty  _{i=0}\,\begin{pmatrix}
z\\ i\end{pmatrix} T^i$ in $\Z_p[[T]],$ for all $z\in \Z_p\,.$  So
the $\Z_p$-algebra homomorphism $\alpha  :\Z_p[[T]] \to \Z_p[[T]]$
induced by $T\mapsto (1+T)^m-1,$ with inverse induced by $T\mapsto
(1+T)^{1/m}-1,$ is a $\Z_p$-algebra automorphism of $\Z_p[[T]].$ 
Applying $\fO_{r+s}\otimes_{\Z_p} $ -- to $\alpha  $ gives the
$\fO_{r+s}$-algebra automorphism $B\to B,$ still denoted by $\alpha 
,$ and then $\sum^{p^s-1}_{i=0} \,b_i u^i\mapsto
\sum^{p^s-1}_{i=0}\,\alpha  (b_i)u^i_0$ defines $\alpha 
:\big(\fq(B)/\fq(A),g,1+T\big) \to
\big(\fq(B)/\fq(A),g,(1+T)^m\big)$ sending $ub = \,^gb\,u$ to $u_0\alpha 
(b) = \,^g\alpha  (b)\,u_0\,,$ since $\alpha  (^gb) = \,^g\big(\alpha 
(b)\big),$ and $u^{p^s}=1+T$ to $u^{p^s}_0=\alpha  (1+T)=(1+T)^m.$ 
Thus, our chosen example is
$$
F\otimes_{\Qbar_p} E\simeq \big(\fq(B)/\fq(A),g,1+T) = D_{r,s,\IF}\,.
$$
It remains to verify that $D_{r,s,\IF}$ is a division ring.  For
this, note that $\fg$ is the Galois group of $\fq(B)/\fq(A),$ of
$F(\zeta  _{r+s})/F(\zeta  _r),$ and of $\Qbar_p(\zeta 
_{r+s})/\Qbar_p(\zeta  _r),$ allowing us to simplify notation by
writing $N_{\fg}, \Tr_{\fg}$ for the appropriate norm, trace.  By
[Re, 30.7] and $\fg$ cyclic of order $p^s,$ we must show that if
$(1+T)^{p^t}$ vanishes in $\hat H^0\big(\fg,\fq(B)^\times \big)$ then
$t\ge s.$

\vs\noi
Now the exact sequence
$$
0\to B^\times \to \fq(B)^\times \to P_B\to0
$$
of $\fg$-modules, with $P_B$ the free $\Z$-module on the height$-1$
prime ideals of the UFD $B,$ has $\hat H^{-1}(\fg,P_B) =0$ since
$P_B$ is a permutation $\Z [\fg]$-module.  Since $1+T \in A^\times =
(B^\times)^{\fg},$ $(1+T)^{p^t}$ vanishes in 
$\hat H^0(\fg,B^\times),$ i.e. $(1+T)^{p^t}\in N_\fg (B^\times).$

\vs\noi
Now $B^\times = \fO^\times _{r+s}\cdot\big(1+\fO_{r+s}[[T]]T\big)$
implies that $(1+T)^{p^t} = aN_\fg(1+cT)$ with $a\in
N_\fg(B^\times)\subseteq A^\times$ and $c\in\fO_{r+s}[[T]],$ hence
$$
1+p^tT \equiv a\big(1+\,\Tr_\fg(c_0)T\big)\!\!\!\mod \fO_{r+s}[[T]]T^2,
$$
with $c_0\in \fO _{r+s}$ the constant term of $c.$  Thus $a=1$ and
$p^t = \Tr_\fg(c_0)\in \Tr_\fg(\fO_{r+s})\dot= 
\pi^{s\phi(p^r)}_r\SO_r    =
p^s\SO _r$ implies $t\ge s.$  Here $\dot= $ comes from
Lemma 1B(i) below.
\end{proof}

\noi
 2.  $D_{r,s,\IF}$ with $r\ge s\ge 1$ in ([Re] (30.1),[CR11] Theorem
74.20(i)) is in the {\it particular case}: this is because
$\zeta_s\in \fq(A)^\times$ implies $\la\zeta_s\ra\in \la\zeta\ra $
hence $\la \zeta_s\ra \subseteq \la \zeta_r\ra$ implies $s\le r.$

\vs\noi
{\sc Lemma 1B.} (i): {\it  Let $F/\Qbar_p$ be finite unramified and $r,s$
any positive integers.  Then}
$$
\Tr_{F(\zeta  _{r+s})/F(\zeta  _r)}(\pi^n_{r+s}\fO_{r+s})
= \pi^{\Phi_s(n)}_r \fO_r
$$
{\it
with $\Phi_k(n) = k\phi(p^r) + \lfloor\frac{n}{p^k}\rfloor$ for $n\ge
0$ and $k\ge 0.$
}

\vs\noi
(ii): $A\subseteq B$ {\it induces} $\Tr_\fg (\Pi^n B) = \pi^{\Phi_s(n)} A.$

\vs\noi
{\it Proof. For} (i): $F/\Qp$ unramified implies $F(\zeta  _m)/\Qp(\zeta 
_m)$ unramified, for any $m$ [Ko, Lemma 4.6.3], hence $\SD_{F(\zeta 
_{r+s})/F(\zeta  _r)} = \SD_{\Qp(\zeta  _{r+s})/\Qp(\zeta 
_r)}\fO_{r+s} = \SD_{\Qp(\zeta  _{r+s})/\Qp}\SD^{-1}_{\Qp(\zeta 
_r)/\Qp}\fO_{r+s} = \pi^{s\phi(p^{r+s})}_{r+s} \fO_{r+s}$
follows from $\SD_{\Qp (\zeta  _m)/\Qp} = \pi^{m\phi(p^m)-p^{m-1}}_m
\Z_p[\zeta  _m]$ (compare [Wa, Proposition 2.1]). Thus
$$
\begin{aligned}
\Tr(\pi^n_{r+s}\fO_{r+s}) &\subseteq \pi^h_r\fO_r \Llra
\,\Tr(\pi^{n-p^sh}_{r+s}\fO_{r+s}) \subseteq \fO_r \Llra\\
\pi^{n-p^sh}_{r+s} \fO_{r+s}  &\subseteq \SD^{-1}_{F(\zeta 
_{r+s})/F(\zeta  _r)} \Llra n-p^sh \ge -s\phi(p^{r+s}) \Llra h\le
s\phi(p^r) +\frac {n}{p^s}\,.
\end{aligned}
$$

\vs\noi
{\it For} (ii): Identifying the $\fO_{r+s}[\fg]$-modules, with $\fg$
acting trivially on $A,$ by $\fO_{r+s}\otimes_{\fO_r} A\to B,$ via
$x\otimes a\mapsto xa,$ we have $\Tr_\fg(\Pi^nB) =
\Tr_\fg(\Pi^n\fO_{r+s}\otimes_{\fO_r}A) =
\Tr_\fg(\Pi^n\fO_{r+s})\otimes_{\fO_r}A=
{\pi^{\Phi_s(n)}\fO_r\otimes_{\fO_r}A} =
\pi^{\Phi_s(n)}A.$  \hfill
$\square$

\vs\noi
{\sc Remark.}This Lemma will find even more roles after Lemma~4D.

\vs\noi
One of these roles of Lemma~1B arises when we investigate 
$\nr(d) -1\!\!\!\mod\!\,\pi^{\Phi(n)+1}A$ for $d\in \Delta  $ in 
\S4, which  involves studying the function
$$
\Phi(n) = \us{k\ge 0}\min\,\Phi_k(n) \q \text{for}\q n\ge 1,
$$
with $\Phi_k(n)$ as in the statement of Lemma~1B.
This implies, for $n\ge 1,$ that $\Phi_k(n) \ge \Phi(n),$ for all
$k\ge 0.$
We sometimes
write $\Phi^{(r)},\,\Phi^{(r)}_k$ to remember that $r\ge 1$ is
normally fixed.

\vs\noi
{\sc Lemma 1C.} {\it Given $n\ge 1,$ $\Phi(n)$ is the absolute
minimum of $k\mapsto \Phi_k(n).$}

\begin{description} 
\item{(i)} {\it $k \mapsto \Phi_k(n)$ is decreasing
on $\{k\ge 0:p^{r+k}\le n\}$ and increasing on $\{k\ge 0:p^{r+k}\ge
n\}.$
 In particular, there is a smallest $k\ge 0$ at which the
minimum is attained.}

\item{(ii)} If  $\Phi_{k}(n)=\Phi_{k+1}(n),$ then $p^{r+k} -
p^k\le n <p^{r+k} + p^k.$
\end{description}

\vs
\begin{proof} {\it  For} (i): 
1.  $p^{r+k}\le n$ implies $\Phi_k(n)\ge \Phi_{k+1}(n),$ and

\noi
2.~$p^{r+k}\ge n$ implies
 $\Phi_k(n)\le \Phi_{k+1}(n).$

\vs\noi
{\it For} 1: $\lfloor \frac{n}{p^k}\rfloor - \lfloor
\frac{n}{p^{k+1}}\rfloor >\big(\frac{n}{p^k} -1\Big)
-\frac{n}{p^{k+1}}= \frac{n}{p^{k+1}}\,(p-1)-1
\ge\frac{p^{r+k}}{p^{k+1}}\, (p-1) -1 = \phi(p^r) -1,$ hence
$\lfloor \frac{n}{p^k}\rfloor
- \lfloor \frac {n}{p^{k+1}}\rfloor
>\phi(p^r) -1,$ i.e.
$\lfloor \frac{n}{p^k}\rfloor \ge \phi(p^r) +
\lfloor \frac{n}{p^{k+1}}\rfloor;$
 add $k\phi (p^r)$ to both sides.

\vs\noi
{\it For} 2: $\lfloor \frac{n}{p^k}\rfloor -\lfloor
\frac{n}{p^{k+1}}\rfloor <\frac{n}{p^k} - \big(\frac{n}{p^{k+1}}
-1\big) = \frac{n}{p^{k+1}} \,(p-1)+1 \le
\frac{p^{r+k}}{p^{k+1}}\,(p-1) +1= \phi (p^r)+1,$
 hence
$\lfloor\frac{n}{p^k}\rfloor 
- \lfloor \frac{n}{p^{k+1}}\rfloor < \phi(p^r) + 1,$
i.e. $\lfloor \frac{n}{p^{k}}\rfloor\le  \phi(p^r) + \lfloor
\frac{n}{p^{k+1}}\rfloor;$
 add $k\phi(p^r)$ to
both sides.

\vs\noi
{\it For} (ii):  To investigate when $\Phi_{k+1}(n) = \Phi_k(n),$ i.e.
$\lfloor\frac{n}{p^k}\rfloor - \lfloor \frac{n}{p^{k+1}}\rfloor =
\phi(p^r),$ from Lemma~1B,
 write $n=a+p^kb+ p^{k+1}c,$ with $0\le
a<p^k,$ $0\le b <p,$ and $c\ge 0.$  Then $(b+pc) -c = \phi(p^r)$
implies $b+(p-1)c = (p-1)p^{r-1},$ hence $p-1$ divides $b,$ so there
are two cases:

\noi
3. $b=p-1$ and $c= p^{r-1} -1,$ or 4.~$b=0$ and $c=p^{r-1}.$

\noi
Then either 
3.~$n=a+p^k(p-1)+ p^{k+1}(p^{r-1}-1) = p^{r+k}-p^k+a,$ or

\noi 4. $n=a+p^{k+1}p^{r-1} = p^{r+k}+a,$  still with $0\le a<p^k.$ 

\noi
These two cases account for all $n$ in the claimed range $p^{r+k} -
p^k \le n <p^{r+k} + p^k.$
\end{proof} 

\vs\noi
We choose the following language to clarify the above behaviour.

\vs\noi
{\sc Definition.} 
 1.{\it Define the} weight $w=w(n)$ {\it of $n\ge 1$ to be the
smallest $w\ge0$ so that 

 ${}\;n<p^{r+w} + p^w$ (thus either $w=0,$  or $w\ge 1$ with
$p^{r+w-1}+ p^{w-1}\le n <p^{r+w} + p^w).$}

\begin{description}
\item{2.}  {\it Define $n,$ of weight $w,$ to be} ordinary {\it if also
$n<p^{r+w} - p^w,$ or} transitional {\it if}
\newline $p^{r+w} - p^w \le n
<p^{r+w} + p^w.$
\end{description}

\vs\noi
{\sc Corollary 1D} (i)
{\it If $n\ge 1,$ of weight $w,$ has $\Phi_w(n) = \Phi_{w+1}(n),$ then
$n$ is transitional of weight $w.$}

\vs\noi
(ii) {\it In the exceptional case $n= p^{r+w} + p^w-1,$ $n$ is
transitional of weight $w$ and $n+1$ is ordinary of weight $w+1.$
Moreover, $\Phi(n) = \Phi_w(n)$ and $\Phi(n+1)= \Phi_{w+1}(n+1).$}

\begin{proof} {\it For} (i): The hypothesis implies
$p^{r+w} - p^w \le n < p^{r+w} + p^w,$ by  Lemma~1C(ii), which  Definition~2
calls transitional of weight $w.$

\noi
{\it For} (ii): First, check that $n= p^{r+w} + p^w -1$ has weight $w$ and
$\Phi_w(n) = w\phi(p^r) + p^r = \Phi_{w+1}(n),$ hence is
transitional of weight $w,$ by (i).  
And $n+1 = p^{r+w} +p^w$ is ordinary of weight $w+1$ since $n+1<
p^{r+(w+1)} - p^{w+1},$ i.e. $(p^r+1)p^w < (p^{r+1} - p)p^w.$

\vs\noi
Moreover, similar calculations to $\Phi_w(n) = w\phi
(p^r)+ p^r$ show that $\Phi_{w-1}(n) > \Phi_w(n) = \Phi_{w+1}(n) 
<\Phi_{w+2}(n),$ which Lemma~1C(i) interprets as $\Phi(n) =\Phi_w(n).$  Also 
$\Phi_{w+1}(n+1) = w\phi(p^r) + p^r$ and $\Phi_w(n+1)
>\Phi_{w+1}(n+1) < \Phi_{w+2}(n+1)$ implies $\Phi(n+1) =
\Phi_w(n+1).$
\end{proof}

\vs\noi
{\sc Lemma 1E.} {\it If $m\ge 1,$ then $\Phi_m(2n) >\Phi(n)$ for all
$n\ge 1.$ }

\begin{proof}
We need the identity
 $\Phi_m(2n) - \Phi_m(n) = \lfloor \frac{2n}{p^m}\rfloor - \lfloor
\frac{n}{p^m}\rfloor.$   
Now if our claim were false, there exists $n\ge 1$ with $\Phi(n)\ge
\Phi_m(2n).$
Then $\Phi(n)\ge \Phi_m(2n) = \Phi_m(n) + \lfloor
\frac{2n}{p^m}\rfloor - \lfloor \frac{n}{p^m}\rfloor$ 
implies that
$0 \ge \big(\Phi_m(n) - \Phi(n)\big) + \big(\lfloor
\frac{2n}{p^m}\rfloor - \lfloor \frac{n}{p^m}\rfloor \big),$
with both parentheses $\ge 0,$ hence
that $\Phi_m(n) = \Phi(n)$ and $\lf \frac{2n}{p^m}\rf =
\lf\frac{n}{p^m}\rf.$ 
The first of these implies  $\Phi_m(n)\,
\dot\le \;\Phi_{m+1}(n),$ with $\dot\le$ by Lemma~1C(i), and the second
that $n<p^m,$ since $\frac{2n}{p^m} -1 <\lfloor
\frac{2n}{p^m}\rfloor = \lfloor \frac{n}{p^m}\rfloor \le
\frac{n}{p^m}.$

\vs\noi
There are now two possibilities, the first being that $\Phi_m(n) =
\Phi_{m+1}(n),$ hence $n$ transitional of weight $m,$
by Corollary~1D(i).  Then $p^m> n\ge p^{r+m} - p^m$ implies $p^r<2,$
contradiction.

\vs\noi
The second possibility is that $\Phi_m(n) <\Phi_{m+1}(n),$ so $n$ is
ordinary of weight $m+1\ge 1,$ hence $p^m >n\ge p^{r+(m+1)-1} + p^{(m+1)-1} =
p^{r+m} + p^m,$ which implies $p^{r+m} <0,$ contradiction.
\end{proof}

\vs
\noi
{\sc Lemma 1F.} If $n=p^{r+w} + p^w -h$ has weight $w,$ hence $h>0,$
then $\Phi(n) = \Phi_w(n).$

\begin{proof} $\Phi_w(n) = w\phi(p^r) + p^r = \Phi_{w+1}(n)$ and 
$\Phi_{w-1}(n) > \Phi_{w}(n) = \Phi_{w+1}(n) < \Phi_{w+2}(n),$ then 
$\Phi(n) = \Phi_w(n),$ as at  the end of the proof of Corollary 1D(ii).
\end{proof}

\vs\noi
The {\it exceptional} $n\ge 1$ are of special interest; we should here
also emphasize the {\it critical} $n\ge 1$ of weight $w\ge 0, $
defined by
$ n=p^{r+w}-1.$
Their special role is clarified in the Remark ending \S4.

\vs\noi
The crossed product $A$-order $\Delta = \us{0\le j<p^s}\bigoplus Bu^j$ in $D,$ 
 require properties of the commutative
domains  $A[u]= \us{0\le j<p^s}\bigoplus Au^j$ and $B$ with 
$u^{p^s} = 1 + T\in A^\times.$

\vs\noi
{\sc Lemma 1G.} (i) {\it If $x\in \Delta$ then
$x=\os{p^s-1}{\us{j=0}\sum} \Pi^jx_j$ with unique $x_j \in A[u].$}

\begin{description}
\item  {(ii)} $\pi A[u]$ {\it is a prime ideal of} $A[u].$

\item  {(iii)}  Is  $A[u]$ {\it  integrally closed?}

\item {(iv)} {\it If $\SO$ is any regular domain and $a\in
\SO[[T]]$ has $\frac{da}{dT}\not \equiv 0\!\!\!\mod
p\SO[[T]],$ then $X^{p^s} -a$ is irreducible over
$\fq(\SO[[T]]).$ }

\item  {(v)} $\fq(A[u])/\fq(A)$ {\it is a cyclic Galois extension, 
when $r\ge s,$
and
$a[u] = \sum\limits_{0\le i <p^s} a_iu^i,$ with all $a_i\in A,$ has}
$$
N_{\fq(A[u])/\fq(A)} (a[u]) \equiv \sum_i a^{p^s}_i (1+T)^i\!\!\!\!\mod
\pi_sA.
$$
\end{description}
\noi
\begin{proof}
{\it For} (i): $x = \sum_i b_iu^i,$ from Theorem~1A,
and  $b_i = \sum_j\Pi^j a_{i,j}$ with unique 
$b_i,$ $a_{i,j}$ in $B,A$ implies that we can
take $x_j \doteq\sum_i a_{i,j} u^i\in A[u].$  It follows that
$A[u]^{p^s}\to\Delta,$ by $(x_0,\dots,x_{p^s-1})\mapsto \sum
\Pi^jx_j$ is a surjective $A$-module morphism.  These free $A$-modules have
the same rank $p^{2s}$ hence the map is also injective (e.g. by [Ma],
Theorem~2.4), giving uniqueness.

\noi
{\it For} $\dot =:$ Given $b=\us{0\le i<\infty}\sum\, c_iT^i\in B=
\fO_{r+s} \big[[T]\big],$ write each $c_i\in \fO_{r+s}$ as $c_i
=\us{0\le j<p^s}\sum\, \Pi^j d_{i,j}\,,$ with unique $d_{i,j}\in \fO_r$
(by [CF], p.23) to get $b= \us i\sum\, \big(\us j\sum\, \Pi^j
d_{i,j}\big)T^i = \us j\sum\, \Pi^j\big(\us i\sum\, d_{i,j}T^i\big),$
with $\us i\sum\; d_{i,j}T^i \in \fO_r \big[[T]\big] = A.$

\vs\noi
{\it For} (ii): 
$\frac{A}{\pi A}
=\frac{\fO_r[[T]]}{\pi_r\fO_r[[T]]} =
\IF[[T]]$ implies $\frac{A[u]}{\pi A[u]} = \frac{A}{\pi
A}\,[\,\bar u\,] = \IF[[T]]\,[\,\bar u\,] = \us{0\le
j<p^s}\bigoplus\;
\IF[[T]] \bar u^j$ (with $\bar u^{p^s} = 1+T).$
Thus $\bar u$
has minimum
polynomial $X^{p^s} - (1+T)$ over the field of fractions
$\IF((T))$ of $\IF[[T]],$ (by
replacing $\SO$ by $\IF$ in (iv)),
so
that adjoining the unique root $\bar u$ of $X^{p^s} - (1+T)$ to
$\IF((T))$ yields
 a purely inseparable extension $\IF((T)) [\bar u]$
of $\IF((T)),$ by ([Ja], Proposition~8.13),
and there
is an exact sequence of $A[u]$-modules
$$
0 \to A[u] \os \pi\lra A[u] \lra \IF[[T]]\,[\bar u] \lra 0,
\leqno{\rm (1G1)}
$$
with $\IF[[T]]\,[\bar u]$ an integral domain.

\vs\noi
{\it For} (iii): This is equivalent to $A[u]$ $(\subseteq D)$ being a
maximal $A$-order in $\fq(A[u]).$  
It follows by (iv),
that $u$ has minimal polynomial $X^{p^s}- (1+T)$ over $\fq(A),$ and
from [Re, Exercise~6.13, 6.14, p.66; 11.4, p.133] that
$\text{disc}(A[u]/A)\, A = \pm$ norm $(p^s X^{p^s-1}\vert_{X=u})A=$
norm $(p^s u^{p^s-1})A =
p^{sp^s}(1+T)^{p^s-1} A= p^{sp^s}A$
reduces us to checking that
$A[u]_\bullet$ is a maximal $A_\bullet$-order.

\vs\noi
This should be decidable by maximal order routines, since
$A_\bullet$ is a Dedekind domain (cf. proof of Lemma~2E(iii)).  But
we do not expect to need it.

\vs \noi
{\it For} (iv):
This follows from 
Theorem~9.1 of ([La2], p.297),
once we verify that there is no $b\in \fq (\SO[[T]])$ with $a =
b^p.$  For if there were, then $b$ is integral over $\SO[[T]],$ which
is integrally closed, by ([Ma], Theorem~19.4), hence $b\in\SO[[T]].$  Then the formal
derivation $\frac{d}{dT}$ on $\SO[[T]]$ has $ \frac{da}{dT} 
=\frac{d}{dT} (b^p) \dot= pb^{p-1}\,\frac{db}{dT}\equiv
0\!\!\!\mod  p\SO[[T]],$ with $\dot =$ by $\frac{d}{dT}\,(\beta  ^k)
= k\beta  ^{k-1} \,\frac{d\beta  }{dT}$ by induction on $k\ge 1.$
Contradiction.  

\vs\noi
{\it For} (v):  As in (ii), $X^{p^s} - (1+T)$ is an irreducible
polynomial over $\fq(A),$ by (iv) 
(for $\SO_r),$ with $u$ as a
root, and $\fq(A)$ contains $\zeta  _s,$ by $r\ge s.$  It follows,
from  [[La2], Theorem~6.2(ii)],
that its cyclic Galois group is generated by $\sigma  ,$ with
$\sigma ^j (u) = \zeta^j_su$  
for $0\le j < p^s,$ 
and $\zeta_s = 1-\pi_s\equiv 1\!\!\!\mod \pi_sA.$  Thus
$N_{\fq(A[u])/\fq(A)}(a[u])=\Pi_j \sigma  ^j(\sum_i a_iu^i) = \Pi_j(\sum_i
a_i \zeta^{ij}_s u^i)
\equiv (\sum_i a_i u^i)^{p^s}
\equiv \sum_i a_i^{p^s} (1+T)^i \!\!\!\mod \pi_sA.$
\end{proof}

\noi
{\sc Remark:} $\bar A = A/\pi A$ has useful $\bar A^{p^k}$-basis
(cf.~1 of Proof of Corollary~4H).

\vs\noi
Some further notation:

\indent
With $r,s,\IF$ fixed: $\zeta  =\zeta  _{r+s},$
$\Pi=\pi  _{r+s}\,,$ $\pi=\pi_r\,.$
So $\zeta_m$ instead of $\zeta_{p^m}\,.$

\indent
$a,b$ will be typical elements in $A$ and $B,$ respectively,

\indent
$A_\bullet$ is the localization of $A$ at its prime ideal $\pi_rA,$
and $\Delta_\bullet = A_\bullet \us A\otimes \Delta, $

\indent
$\nr:D\to \fq (A)$ is the reduced norm,

\indent
$[x,y] = xyx^{-1}y^{-1} $ is our convention for writing commutators.

\vs\noi
The {\it weight} $w=w(n)$ of $n\ge 1$ is defined to be the smallest
$w\ge 0$ so that $n<p^{r+w}+ p^w.$

\noi
The {\it level} $n$ of $d\in 1+\Pi \Delta$ is defined by $n= \max\{k\ge
1:d\in 1+\Pi^k\Delta\}.$

\noi
The {\it type} $t\ge 0$ of $n\ge 1$ is defined by $n=p^{r+t} +
p^tm-1$ for $m\ge 1$ with $m\not\equiv 0\!\!\!\mod p.$

\section{$SK_1(\Delta  ) \to SK_1(D)$ and $\fg$-Action
 Congruences}  \label{Sec2}

{\sc Definition.} $\Delta  = \oplus^{p^s-1}_{j=0}\,Bu^j$ {\it (with
$u^{p^s} = 1+T\in Bu^0\subseteq \Delta  )$ 
is an
$A$-order in $D$
with typical element $d=\sum^{p^s-1}_{j=0}\,
b_ju^j.$}

\vs\noi
First, $\Delta  $ is a ring since $Bu^i Bu^j= B(^{g^i}B) u^i u^j =
BBu^{i+j} \doteq Bu^{p^s l +k} \os{\cdot\cdot}= Bu^k\subseteq\Delta  ,$ with $\doteq$
by $i+j = p^s l+k$ and $0\le k<p^s,$ and $\os{\cdot\cdot}=$ by
$Bu^{p^s l} = B(1+T)^l = B,$ as $(1+T)^l \in B^\times.$  

\vs\noi
Next, $\Pi\SO_{r+s} = \SO_{r+s}\Pi$ is totally ramified over
$\SO_r,$ hence $^g(\Pi\SO_{r+s}) = \Pi\SO_{r+s}$ and
$^g(\Pi B) =\Pi B.$  Thus $Bu^j\Pi = B^{g^j}\Pi u^j =\,
^{g^j}(B\Pi)u^j = B\Pi u^j = \Pi Bu^j,$ so that
$\Delta  \Pi = 
\oplus_j(Bu^j)\Pi = \oplus_j \Pi(Bu^j) =
\Pi\Delta,  $ hence $\Delta\Pi =\Pi\Delta$
is a 2-sided ideal of $\Delta  .$

\vs\noi
{\sc Lemma 2A.}
(i) $\Delta  $ {\it is a local ring with $\rad
(\Delta  ) = (u-1)\Delta  + \Pi\Delta  $
 and $\Delta  /\rad(\Delta) = \IF.$

\noi
In particular, there is a natural isomorphism} $\Delta  ^\times/[\Delta 
^\times,\Delta  ^\times] \to K_1(\Delta  ).$

  {(ii)} $\wh\Delta_\bullet  $ {\it is a maximal $\wh A_\bullet$-order of $D.$}

\vs\noi
{\it Proof.}
{\it For} (i):  Recall that $A=\fO_r[[T]]$ is a  2-dimensional local ring with $\rad
(A) = \pi_r A + TA$ and $A/\rad(A)=\IF$ (see [NSW, Chap.V, \S3]).

\noi
To understand
$\ol\Delta   = \Delta 
/\Pi\Delta , $ via
$\ol B =\frac{B}{\Pi B} =\big(\frac{\fO_{r+s}}{\Pi \fO_{r+s}}
)[[T]] 
= 
\IF [[T]]$ from \S1
as $\IF[[T]]$-module
$$
\ol\Delta   =\bigoplus_{0\le j < p^s} \bar B \bar u^j = \bigoplus_j
\IF[[T]] \bar u^j \doteq \bigoplus_j \IF[[T]] (\bar u - \bar 1)^j
=
\IF[[T]]\,[\,\ol{u-1}\,],
$$
with $\dot=$ the change of
$\IF[[T]]$-basis  $\{\ol u\,^j:0\le j < p^s\}$ to $\{\,\ol{u-1}\,^j:
0\le j<p^s\}.$
Here $\dot=$ comes from the binomial expansion of $\bar u^j =
\big((\,\bar u -\bar 1) + 1\big)^j = \us{0\le i\le j}\sum
\,\begin{pmatrix} j\\ i\end{pmatrix} (\,\ol{u-1}\,)^i,$ and
conversely by expanding $(\bar u - \bar 1)^j;$ and also $\bar u\,^0=1$ since
$\bar u \not\equiv 0\!\!\!\mod \Pi\Delta  :$ else $1+T\equiv \bar
u\,^{p^s} \equiv \bar 1\not\equiv 0\!\!\!\mod  \Pi\Delta  ,$ contradicting that
$1+T\in \Delta  ^\times$ (as $(1+T)^{-1} = 1 +\us{k\ge 1}\sum (-1)^kT^k).$

\vs\noi
Now every $k\ge 0$ is uniquely $k= p^s i+j,$ with $i\ge 0$ and $0\le
j<p^s,$ and $(\bar u - \bar 1)^{p^s} = \bar u^{p^s} - \bar 1^s = T,$ in
$\IF[[T]]\,[\bar u-\bar 1]$ hence every element of $\IF[[T]][\ol{u-1}]$
is uniquely $\us j\sum \big(\us i\sum \alpha  _{(i,j)}
T^i\big)(\ol{u-1})^j
=\us j\sum\,\us i\sum \,\alpha  _{(i,j)} (\,\ol{u-1}\,)^{p^si+j} =
\us k\sum \alpha  _k(\,\ol{u-1}\,)^k,$ with $\ol{u-1}$ the unique root of the
irreducible polynomial $X^{p^s}-T$ over $\IF\big((T)\big),$ by
Lemma~1G(iv).
 Thus
$\big(\IF\big[[T]\big]\big)[\,\ol{u-1}\,\big]\cong
\IF\big[[\,\ol{u-1}\,]\big],$ which
implies the exact sequence of $\Delta  $-modules
$$
0 \lra \Delta  \os{\bullet\Pi}\lra \Delta  \os{-}\lra
\IF[[\,\ol{u-1}\,]] \lra 0.
\leqno{\rm (2A1)}
$$

\noi
In particular, $\ol\Delta  $ is a commutative ring, hence $u-1$
commutes with $\Delta  \!\!\!\mod \Pi\Delta  ,$ so

\noi $(u-1)\Delta  +
\Pi\Delta  = \Delta  (u-1) +\Pi\Delta  = \Delta  (u-1) +\Delta  \Pi,$
i.e. $J:= (u-1)\Delta   +\Pi\Delta  $ is a 2-sided ideal of
$\Delta.$ 
Now $J^{p^s} = \big((u-1)\Delta  + \Pi\Delta  \big)^{p^s}\subseteq 
(u-1)^{p^s}\Delta  +\Pi\Delta  = T\Delta  +\Pi\Delta  $ 
and 
$J^{2p^s} \subseteq (T\Delta  + \Pi\Delta  )^{p^s} 
\dot \subseteq \,T\Delta   +\pi_r\Delta  = (\rad\,A)\Delta,  $
with $\dot=$ by (2A3*).  Thus
$J\subseteq \,\rad\,\Delta  $ (see [Re, Exercise~6.3, p.~88] with
$R,A$ our $A,\Delta  ).$

\vs\noi
Now $\bar\Delta  \cong \IF[[\ol{u-1}]]$ and $\rad(\IF[[\ol{u-1}]]) =
\ol{u-1}\,\IF[[\ol{u-1}]],$
the unique maximal ideal of $\IF[[\ol{u-1}]].$  Also $J= (u-1) \Delta +
\Pi\Delta  $
 fits into the exact sequence $0\to J\hookrightarrow
\Delta  \to \Delta  /J\to 0,$ 
with canonical maps, hence 
$0\to \frac{J}{\Pi\Delta  } \hookrightarrow \frac{\Delta 
}{\Pi\Delta  } \to \frac{\Delta  }{J} \to 0.$
It follows that $(u-1)\Delta  $ maps onto $\frac{J}{\Pi\Delta  } = \frac
{(u-1)\Delta  +\Pi\Delta  }{\Pi\Delta  },$ which has image
$(\,\ol{u-1}\,)\,\ol\Delta  $ in $\ol\Delta  ,$ so our exact
sequence becomes $0\to (\,\ol{u-1}\,)\,\ol\Delta   \hookrightarrow \ol\Delta 
\to \Delta  /J \to 0.$  Thus
$\Delta  /J \cong \ol\Delta  /\,(\,\ol{u-1}\,)\,\ol\Delta   \cong
\IF[[\,\ol{u-1}\,]]/(\,\ol{u-1}\,)\,\IF[[\,\ol{u-1}\,]] =\IF,$ a field.
This implies that $\rad~\Delta 
=J$ (as $J\subseteq \,\rad\,\Delta  )$ and $\Delta  /\rad\,\Delta  =\IF,
$ hence that
 $\Delta  $ is a local ring, by ([Re],
Theorem~6.14).

\noi
In particular, this follows from ([CR2], Theorem~(40.31) and (40.32) or,
perhaps better, from ([Va2], Theorem~2)).  This proves (i).

\vs\noi
{\it We now prove} (ii) with the $A$-basis $\lambda  _k = \zeta  ^{k_0}
u^{k_1}, \,0\le k < p^{2s},$ with $k=k_0 + p^s k_1,$ $0\le k_0,k_1
<p^s,$ of $\Delta  $ (which is also a $\fq(A)$-basis of $D).$
This follows from $\{\Pi^j:0\le j<p^s\}$ being an $\SO_r$-basis of
$\SO_{r+s}$ (by p.23 of [CF], as in the proof of Lemma~1G(i)
 which
implies that $\{\zeta^j: 0\le j<p^s\}$ is also an $\SO_r$-basis of
$\SO_{r+s}$ (since the transition matrix is invertible by $\Pi^j=
(1-\zeta)^j=\sum\limits^j_{l=0} (-1)^l \begin{pmatrix}
j\\l\end{pmatrix}
\zeta^l),$ hence an $A$-basis of $B;$ the rest 
 follows from
$\Delta = \us{0\le j <p^s}\oplus Bu^j.$
Then
\disc$(\Delta  /A)$ is the  ideal of $A$ generated by
$\det\big(\tr(\lambda  _i\lambda  _j)\big)_{0\le i,j<p^{2s}},$ where $\tr$ is the
reduced trace of $D/\fq(A).$ Recall here that the regular trace 
$\Tr\!\!:\!D \to \fq (A),$ of right multiplication by $d,$ has $\Tr(d)=
p^s\tr(d)$ for all $d\in D.$
Below we show that

\begin{description}
\item  {1.} $\disc (\Delta  /A) = p^{sp^{2s}}A.$
\item  {2.} $\Delta  _\bullet = A_\bullet \otimes _A \Delta  $ is a maximal
$A_\bullet$-order in $D.$
\end{description}

\vs\noi
From this claim (ii) follows: see [Re, Exercise 4.13, p.66;  Theorem~11.4,
p.133], since $\Delta  $ is $A$-free.

\vs\noi
{\it Proof of} 1. Since $\fq(B)$ is a maximal subfield of $D,$ 
there is an
 isomorphism

\centerline{$
\varphi:\fq(B)\otimes_{\fq (A)} D \to M_{p^s}\big(\fq(B)\big)$}

\noi
of central simple $\fq(B)$-algebras, by ([Re], Theorem 9.3 and
(9.6a)).
Using the above $\fq(B)$-basis
$\{u^{\hat i}: i\in \Z/p^s\Z\}$ of $D,$ from Theorem~1A,  with
$i=\hat i + p^s\Z$ and $0\le
\hat i<p^s,$ note that the characteristic polynomial is the
same for right and left multiplication of $d\in D,$ by ([Re],
Theorem 9.31).

\vs\noi
Writing $d\in D$ as $d=\us{k\in\Z/p^s\Z}\sum b_ku^{\hat k},$ with
all $b_k\in \fq(B),$ then
$u^{\hat j}d =
\us{k}\sum \,u^{\hat j}b_k
u^{\hat k}=$

\centerline{$\us k\sum {}\;^{g^j}b_k u^{\hat j+\hat k} =\us
k\sum\;{} ^{g^j}b_k u^{\wh{j+k}+ p^sh_{j,k}} = \us{i\in \Z/p^s\Z}\sum
\;{}^{g^j}b_{i-j} u^{\hat i + p^s h_
{j,i-j}}=\us i\sum{}\;^{g^j}b_{i-j}(1+T)^{h_{j,i-j}}u^{\hat i},$ }

\noi
where 
$p^sh_{j,k} = \hat j +\hat k -\wh{j+k}$
uniquely determines $h_{j,k} \ge 0$ for $j,k$ in $\Z/p^s\Z,$
hence
$p^s h_{j,i-j} =\hat j -\hat i +\wh{i-j},$ 
which we summarize by
$$
u^{\hat j}d =\us {i\in\Z/p^s\Z}\sum {}\;^{g^j}
b_{i-j}(1+T)^{h_j(i)}u^{\hat i} \q {\rm with} \q
h_j(i) =\begin{cases} 1 , &\hat j >\hat i\\
0, &\hat j\le \hat i.\end{cases}.
\leqno{\rm (2A2)}
$$

\noi
Here setting $h_j(i) = h_{j,i-j} $ implies the relation
$
\hat j -\hat i + \wh{i-j} = p^s h_j(i)
$
with left side implying $-p^s <p^sh_j(i) <2p^s$ hence
$h_j(i)\in \{0,1\};$  similarly $\hat j -\hat i>0$ implies $p^s
h_j(i)>0$ hence $h_j(i) =1,$ and $\hat j -\hat i \le 0$ implies 
$p^sh_j(i)<p^s$ hence
$h_j(i) = 0.$

\vs
\noi
Now $\varphi  (1\otimes d) = \big[\,^{g^j} b_{i-j}
(1+T)^{h_j(i)}\big]_{(\hat
i,\hat j)}$ has diagonal entries ${}^{g^i}b_{i-i} (1+T)^{h_i(i)} =
{}^{g^i}b_0$ 
so 
$$
\tr (d) =\us i\sum {}^{g^i}b_0 = \Tr_\fg(b_0),
$$
with $\Tr_{\fg}$
the
Galois trace of the $\fq(B) /\fq (A)$ extension, hence 
$$
\Tr_{\fg}(\zeta  ^\kappa  ) =\begin{cases} p^s\zeta_r  ^{\kappa^*} &\text{\rm if} \q
\kappa  \equiv
0\!\!\!\mod p^s\\
0 &\text{\rm else}
\end{cases}
\leqno{\rm (2B0)}
$$
with $\kappa =p^s\kappa^*.$  For
 $\Tr_{\fg}(\zeta  ^\kappa) =\sum^{p^s-1}_{j=0} \, ^{g^j} \zeta 
^\kappa
 = \zeta 
^\kappa\sum_j \,^{g^{j}-1}\zeta  ^\kappa 
\os\star = 
\zeta  ^\kappa \sum^{p^s-1}_{j=0}
\,(\zeta  ^j_s)^\kappa = \zeta  ^\kappa \sum_j(\zeta  ^\kappa_s)^j=0,$ 
when
$\kappa\not\equiv 0\!\!\!\mod p^s.$
When  $\kappa  =p^s\kappa  ^*,$ then 
$\zeta^\kappa  =\zeta^{p^s\kappa  ^*} = \zeta  ^{\kappa  ^*}_r$ and 
$
\sum_j \zeta  ^{p^s\kappa  ^*}_s = \sum_j 1^{\kappa  ^*} = p^s.$ 

\vs\noi
Here  $\os\star = $ has $^{g-1}\zeta  =\zeta  ^{p^r} = \zeta  _s$
hence 
${}^{g^j-1}\zeta   =  {}^{(1+g+\cdots +g^{j-1})(g-1)}\zeta = {}^{(1+g+\cdots +
g^{j-1})}\zeta  _s= \zeta  ^j_s,$ since ${}^{g}\zeta  _s = \zeta  _s$
by $\zeta  _s \in A.$

\vs\noi
Note here that the values of $\tr$ are all in $A.$

\vs\noi
Turning to $\tr (\lambda  _k\lambda  _l)$ we have
$\lambda  _k\lambda  _l = \zeta^{k_0} u^{k_1} \zeta^{l_0} u^{l_1}
= \zeta^{k_0}({}^{g^{k_1}}\zeta^{l_0})u^{k_1} u^{l_1}
= \zeta^{k_0}\zeta^{(1+p^r)^{k_1} l_0} u^{k_1+l_1}
= \zeta^{k_0+(1+p^r)^{k_1}l_0} u^{k_1+l_1}.$

\noi
 The reduced traces $\tr (\lambda  _k\lambda  _l)$ of these 
pairs of $A$-bases of $\Delta  $ are non-zero only when their
$\Tr_{\fg}(b_0)$'s are, which happens only when $u^{k_1+l_1}$ is in
$A$ i.e. when $k_1+l_1\equiv 0\!\!\!\mod p^s.$  The constraints
on these imply $0\le k_1<p^{s},$ $0\le l_1 <p^s$ hence $k_1+l_1$ can
only be $0$ or $p^s$ i.e. $(k_1,l_1)$ can only be $(0,0)$ or
$(k_1,p^s-k_1)$ with $0<k_1<p^s.$  Thus $u^{k_1+l_1}$ can only be
$\eta (k_1) :=\begin{cases}1, &k_1=0\\
1+T, &0<k_1<p^s,\end{cases} $ which is in $A^\times.$
It follows that $\lambda  _k\lambda  _1$ has reduced trace
$$
\tr(\lambda  _k \lambda  _l)
= \begin{cases} \Tr_{\fg} \big(\zeta  ^{k_0+(1+p^r)^{k_1}l_0}
\,\eta(k_1)\big) &\text{\rm if}\q k_1 + l_1\equiv 0\!\!\!\mod p^s\\
0 &\text{\rm else.}\end{cases}\q\q\q\q\q\q\q\q\q\q\q\q\q\q{}
$$
Evaluating $\Tr_{\fg}$ this becomes
$$
\tr(\lambda  _k\lambda  _l) 
 =\begin{cases}
p^s\zeta_{_r}  ^{(k_0+(1+p^r)^{k_1}l_0)^\pr}\eta(k_1) &\text{\rm if} \q
k_1+l_1\;\text{\rm   and}\; k_0 +(1+p^r)^{k_1}
l_0\;\text{\rm are both}\; \equiv 0\!\!\!\mod p^s.\\
0 &\text{\rm else.}
\end{cases}
$$

\noi
Then the $p^{2s}\times p^{2s}$ matrix  $\big(\tr(\lambda  _k\lambda 
_l)\big)_{(k,l)}$ has exactly one
non-zero entry in each row (and in each column): for any entry in row
$k$ is $\tr(\lambda  _k\lambda  _l)$ for some $l,$ hence 
$\begin{array}{ll} k_0 +(1+p^r)^{k_1}l_0 &\equiv 0\!\!\!\mod p^s\\
k_1 +l_1 &\equiv 0 \!\!\!\mod p^s\end{array}$ and implying
$$
 \begin{cases} l_0 \equiv -(1+p^r)^{-k_1}k_0 \!\!\!\mod p^s &\text{and}
\q 0 \le l_0 <p^s\\
l_1\equiv -k_1 &\text{and} \q 0 \le l_1 <p^s
\end{cases}
$$
determining $l_0,l_1$ and therefore $l=l_0 + p^sl_1$ uniquely.

\noi
Thus \q\q\q 
 $\disc(\Delta  /A) = \pm\big((p^s)^{p^{2s}}\;\text{\rm times powers
of}\q \zeta  _r \; \text{\rm and}\q 1+T)A,
$ \q\q\q
proving~$1.$

\vs\noi
Before continuing with the proof of 2., for future reference we
quickly supplement the computation concerning $\Tr(\zeta  ^k)$ above
by its norm analogue:

\vs\noi
 (2A3)\q\q\q\q \q\q\q\q\q\q\q\q\q\q $N_{\fq(B)/\fq(A)}(\Pi) =
\pi,$

\vs\noi
since $\os{p^s-1}{\us {j=0}\prod}\,^{g^j}\Pi = \us
j\prod(1-\,^{(g^j-1)}\zeta\cdot\zeta) \dot =
\os{p^s-1}{\us{i=0}\prod} (1-\zeta  ^i_s\zeta  ) = 1^{p^s} -\zeta 
^{p^s}
= 1-\zeta_r$ 
with $\dot\in$  as before.  This, and Lemma~2G(iii), imply that
$\pi\Pi^{-p^s}= \us{0\le j<p^s}\prod\,^{g^j-1}\Pi \equiv 1\!\!\!\mod
\Pi^{p^r-1}\SO_{r+s},$ hence 
$$
\pi\Pi ^{-p^s} \in \SO^\times_{r+s} 
\leqno{\rm (2A3^*)}
$$

\vs\noi
{\it   Proof of 2.}  
{\it Step 1.} $\pi A_\bullet$ is a maximum ideal of  $A_\bullet$  and
$A_\bullet/\pi A_\bullet$ is a simple $A_\bullet\text{-module.}$

\vs\noi
The local ring $A=\SO_r[[T]]$ gives an exact $A\text{-module}$
sequence 
$0\to A\os {\bullet\pi}\to A\os{ \alpha } \to \IF[[T]] \to 0,$ with $\alpha  $ a ring
homorphism,
localizing to 
$0\to A_\bullet \os {\bullet\pi}\to A_\bullet
\os{\alpha }\to \IF[[T]]_\bullet \to 0.$  Here $\IF[[T]]_\bullet = S^{-1}\IF[[T]],$ with
$S$ the $\alpha  \text{-image}$ of $A\bs \pi A=\IF[[T]]\bs \{0\}$
with $\IF[[T]]_\bullet = \IF\big((T)\big),$ which is a field (of formal
Laurent series)  and $A_\bullet/\pi A_\bullet \cong \IF((T)).$  Thus
$\pi A_\bullet$ is a maximal ideal of $A_\bullet$ and $A_\bullet/\pi
A_\bullet$ a simple $A_\bullet\text{-module}.$

\noi
If $M$ is a maximal  ideal $\ne \pi A_\bullet$ then
$M+\pi A_\bullet = A_\bullet,$ so $1\in M+\pi A_\bullet$ 
implies $M$ contains $1-\pi
a$ for some $a\in A_\bullet.$ 
Multiplying $1-\pi a$ by $\us{0\le
i\le k}\sum (\pi a)^i\in A_\bullet$ implies that $1- (\pi
a)^{k+1}\in M\cap \pi^{k+1}A_\bullet,$ which is a Cauchy sequence in
the $\pi A_\bullet$-topology.  So if $A_\bullet$ were complete in
the $\pi A_\bullet$-topology, then $1\in M,$ contradiction; it
would follow that $\rad A_\bullet = \pi A_\bullet$ and $A_\bullet$
is a local ring ([Ma], Theorem~8.6 and 8.7).

\noi
{\it Step 2:} $\Pi \Delta  _\bullet$ is a maximal ideal of $\Delta 
_\bullet$ and $\Delta  _\bullet/\Pi \Delta  _\bullet$ is a simple
(left and right) $\Delta  _\bullet$-module.  

\noi
We follow the proof of Step~1 closely (as $\Pi\Delta  _\bullet$ is
an ideal).  Starting from (2A1) and
\newline $_\bullet$-localizing gives the
exact sequence $0\to \Delta  _\bullet\os{\bullet\pi}\to \Delta  _\bullet
\os{\alpha }\to
\IF[[u-1]]_\bullet \to 0,$ with $\IF[[u-1]]_\bullet = \IF((u-1)),$ just as before,
to get $\Delta  _\bullet/\Pi\Delta  _\bullet\cong \IF((u-1))$ still a
field, hence $\Pi\Delta  _\bullet$ a maximal ideal of $\Delta 
_\bullet$ and  $\Delta  _\bullet/\Pi\Delta  _\bullet$ still a simple
$\Delta  _\bullet$-module.

\noi
If $M$ is a maximal $\Delta  _\bullet$-ideal $\neq \Pi\Delta 
_\bullet$ then $M+\Pi \Delta  _\bullet = \Delta  _\bullet$ still
implies that $1-\Pi\delta   \in M$ for some $ \delta  \in \Delta 
_\bullet,$ hence $1-(\Pi\delta  )^{k+1}\in M\cap \Pi^{k+1}\Delta 
_\bullet,$  by $\Pi\Delta  _\bullet = \Delta  _\bullet \Pi.$ 
The ``$\Pi\Delta  _\bullet$-topology'' gives little trouble since
$\Pi^{p^s}\Delta  _\bullet = \pi\Delta  _\bullet,$ with $\pi$
central, reduces it to a frill.  The rest is recovery of the need for
$\Delta  _\bullet$ to be complete,  of $\rad \Delta  _\bullet =
\Pi\Delta_\bullet$ and, of $\Delta_\bullet$ ``to be a local ring''
([Re], Theorem~(6.16) and of (ii) on page~85), completing Step~2.

\vs\noi
{\it Step 3:} $\wh\Delta_\bullet$ is a maximal $\wh
A_\bullet${-order}.

\noi
Now if $L$ is a full left $\Delta  _\bullet$-lattice in $D,$ then $L=\Delta 
_\bullet\Pi^j$ for some integer $j.$  For $\bigcap_i\Delta 
_\bullet\Pi^i=0$
implies that there is a largest $j$ with $ \Delta 
_\bullet\Pi^j\supseteq L.$  Then $\Delta  _\bullet
\Pi^{j+1}\subseteq
L+\Delta  _\bullet \Pi^{j+1} \subseteq \Delta  _\bullet \Pi^j$ and
$\Delta_\bullet \Pi^j/\Delta_\bullet \Pi^{j+1}$ is simple.
If the first $\subseteq$ is an equality 
then $\Delta_\bullet \Pi^{j+1} = L + \Delta_\bullet
\Pi^{j+1}\supseteq L$ contradicts maximality of $j.$  Thus
$L+\Delta_\bullet \Pi^{j+1} = \Delta_\bullet \Pi^j$ implies $L\dot =
\Delta_\bullet \Pi^j,$
 by Nakayama's Lemma, because
$L + ({\rm rad}\, \Delta  _\bullet)\Delta  _\bullet \Pi^j = 
L+ (\Delta  _\bullet \Pi)(\Delta  _\bullet \Pi^j)=
L+\Delta  _\bullet \Pi^{j+1}$ by
 Step~2, when $\Delta  _\bullet$ is complete.

\noi
Therefore, if $L$ is an $A_\bullet$-order in $D$ containing $\Delta 
_\bullet,$ then $L= \Delta  _\bullet \Pi^j,$ when $L=L\cdot L,$ as
$1\in L,$ implies $\Delta  _\bullet \Pi^j = L= L\cdot L =\Delta 
_\bullet\Pi^j\cdot \Delta  _\bullet\Pi^j = \Delta  _\bullet
\Pi^{2j},$ hence $j=0$ and $L= \Delta  _\bullet$ proving 2, and thus
Step~3. \hfill $\square$

\vs\noi
{\sc Remark:} Let $\wh A_\bullet$ denote the completion of
$A_\bullet$ in the $\pi A_\bullet$-topology (and also $\wh \Delta  _\bullet$ in the
$\Pi\Delta  
_\bullet$-adic topology).   Then
$\wh\Delta  _\bullet$ is a local ring with maximal ideal $\rad
\wh\Delta  _\bullet = \Pi\wh\Delta  _\bullet$ and $\wh\Delta 
_\bullet/\Pi\wh\Delta  _\bullet$ a simple (left or right) $\wh\Delta
 _\bullet$-module.

\vs\noi
{\sc Corollary 2B.} (i) {\it Division algebras with finite dimensional
centre (usually) have \ul{no} valuations at all ([TW] Theorem~1.4);
when a valuation $v:D\to \Z \cup\{\infty\}$ exists, it is unique (by
loc.~cit.).  This does not prevent subfields of $D$ from having
valuations  of
commutative subrings of $D,$ but only with care. Also Dubrovin
valuation rings (cf. [MMU]) appear to be a serious obstacle to
Suslin's conjecture. }

\vs\noi
(ii) {\it If $0\ne b\in B_\bullet$  $(\subseteq \Delta  _\bullet),$ then
choose $a\in A\bs \pi A$ with $ab\in B=\SO_{r+s}\big[[T]\big],$
hence $ab=f(T) =\os\infty  {\us{i=0}\sum} b_iT^i,$ with $b_i\in
\SO_{r+s}$ for all $i.$  Then setting $v_\Pi \big(f(T)\big) =n,$ with $n$ the
largest integer with $b_i\in \Pi^n\SO_{r+s},$ for all $i,$ defines a
valuation $v_\Pi:B_\bullet\to \Z\cup\{\infty  \},$ which has
$v_\bullet/B =
v_\Pi\,.$ }

\vs\noi
(iii) {\it There is a unique valuation} $v_\pi :A[u]_\bullet\to \Z \cup
\{\infty\}$ {\it so that} $v_\bullet /{A[u]} = p^sv_\pi.$

\vs\noi
(iv) {\it If $0\ne a\in A_\bullet$ $(\subseteq \Delta_\bullet),$
choose $\alpha  \in A\bs \pi A$ with $\alpha  a\in
A=\SO_r\big[[T]\big],$ hence $\alpha  a =f (T) =\os\infty 
{\us{i=0}\sum} a_iT^i,$ with $a_i\in \SO_r.$ Setting $v^\bullet(a)=n,$
with $n$ maximal so $a_i \in \pi^n\SO_r$ for all $i,$ defines a
valuation $v^\bullet:A_\bullet \to \Z \cup \{\infty  \},$ with
$v_\bullet/A_\bullet = p^s v^\bullet.$ }

\vs\noi
(v) $v_\Pi/A_\bullet  = p^sv^\bullet$ and $v_\pi/A_\bullet
=v^\bullet.$

\vs\noi
{\it Proof.  For} (ii):  First suppose that $b\in B$ and $a=1.$  Then $b=f(T)
=\Pi^n f^\pr(T)$ with $f^\pr(T) =\us i\sum b^\pr_iT^i,$ 
so there is a smallest $m$ with $b^\pr_m\notin \Pi \SO_{r+s}.$  Then
$f(T) = \Pi^n P(T)U(T)$ with $U(T)$ a unit of
$\SO_{r+s}\big[[T]\big]$ and $P(T)$ a distinguished polynomial of
degree $m,$ by ([Wa], Theorem~7.3).  Now $P(T)= T^m + c_{m-1}T^{m-1}
+\dots +c_0,$ with all $c_i\in \Pi \SO_{r+s}$ implies that
$v_\bullet\big(U(T)\big) = 0$ and $v_\bullet \big(P(T)\big) =
v_\bullet (T^m) =0,$ hence $v_\bullet \big(f(T)\big) = v_\Pi(\Pi^n) =
n.$  Since $v_\bullet/ B_\bullet $ is a valuation, so is $v_{\Pi}.$

\vs\noi
In the general case, $a\in A\bs \pi A$ implies $v_\Pi(a) = 0.$

\vs\noi
{\it For} (iii):
{\it Step 1}: \hspace{3cm} $\Delta  \cap \Pi\Delta  _\bullet = \Pi \Delta  .$

\vs\noi
This follows from
the commutative
diagram with exact rows and 
 vertical inclusions, 
$$
\begin{matrix}
0 &\lra &\Delta &\os\Pi\lra &\Delta &\lra & \IF[[U]]
&\lra &0\\
\\
&&\downarrow &&\downarrow &&\downarrow\\
\\
0 &\lra &\Delta_\bullet &\os\Pi\lra &\Delta_\bullet &\lra
&\IF\big((U)\big) &\lra &0
\end{matrix}
$$
which follows from (2A1). 
Starting with $c$ in the middle $\Delta$ gives Step 1.

\vs\noi
{\it Step 2}: \hspace{2.75cm}  $A[u]\cap \pi\Delta  \subseteq \pi A[u].$

\noi
It  this were false there
exists $x\in A[u]\cap \pi\Delta$ with $x\notin \pi A[u].$  Then
$x=\pi d$ with $d\in \Delta$ and $d\notin A[u],$ hence the image
$\tilde
d$ of $d$ under the natural $A$-module homomorphism $\Delta \to
\Delta/A[u]$ has $\pi \tilde d=\tilde x = \wt 0$ and $\tilde d\neq
\wt
0,$
i.e. $\Delta/A[u]$ has non-trivial $\pi$-torsion.  This contradicts
that $\Delta/A[u]$ is a free $A$ module: for $A[u] \subseteq\Delta$
have $A$-bases $\{u^i:0\le i <p^s\}\subseteq \{\Pi^j u^i: 0\le i
<p^s,$ $0\le j<p^s\},$ respectively.   This settles Step~2.
$$
A[u] \cap \Pi\Delta  \subseteq \pi A[u]. \leqno{\rm (2B1)}
$$

\noi
For $A[u] \cap \Pi \Delta$ is an ideal of the
commutative ring $A[u]$ which has $(A[u] \cap \Pi
\Delta)^{p^s}\subseteq A[u] \cap (\Pi\Delta)^{p^s} = A[u] \cap
\Pi^{p^s}\Delta = A[u]\cap \pi\Delta \subseteq \pi A[u],$ by Step~2.  This
implies that $A[u] \cap \Pi\Delta \subseteq \pi A[u],$ by Lemma
1G(ii).

\vs\noi
Now if  $0\ne x\in A[u],$ then there is a largest $k\ge 0$ so that
$x=\pi^kc,$ with $c\in A[u].$
If $v_\bullet (c)>0,$ then $c\in A[u]\cap \Pi \Delta_\bullet =
A[u]\cap \Delta \cap \Pi\Delta_\bullet\dot\subseteq A[u] \cap \Pi
\Delta \ddot\subseteq \pi A[u],$ by Step~1 and (2B1). 
Then $c= \pi c_1,$ with $c_1\in A[u],$
hence $x=\pi^{k+1}c_1$ contrary to the maximality of $k.$  Thus
$v_\bullet (c) =0.$

\noi
Finally (iii) follows: defining $v_\pi(x) = v_\bullet
(x)/p^s,$ for $x\in A[u],$ gives a valuation $v_\pi:A[u]\to \Z \cup
\{\infty\},$ since
$
v_\pi(x) = v_\bullet (\pi^kc)/p^s =\big(v_\bullet (\pi^k) +
0\big)/p^s = p^sk/p^s = k\in \Z. 
$

\vs\noi
{\it For} (iv): If $a\in A,$ and $\alpha  =1,$ then $f(T) =
\pi^nf^\pr(T)$ with $f^\pr(T) =\us i\sum\, a^\pr_i T^i,$ defining a
minimal $m$ with $a^\pr_m\notin \pi \SO_r\,.$  Then $f(T) =\pi^n
P(T)U(T)$ with $U(T) \in \SO_r\big[[T]\big]^\times$ and distinguished $P(T)
\in\SO_r[T]$ of degree $m,$ by [Wa].   Now $P(T) = T^m +
c_{m-1}T^{m-1} +\dots +c_0,$ with $c_i\in \pi\SO_r,$ implies
$v_\bullet \big(f(T)\big) = v_\bullet (\pi^n) = p^sn,$ hence
$v_\bullet/  A = p^sv^\bullet.$  In general, $v^\bullet (\alpha 
)=0,$
 as in
(ii).

\vs\noi
{\it For} (v): First, $v_\Pi = v_\bullet/B_\bullet,$ from (ii),
implies $v_\Pi/ A_\bullet = (v_\bullet/B_\bullet)/  A_\bullet
= v_\bullet/ A_\bullet = p^s v^\bullet,$ by (iv).  Second, $p^s v_\Pi =
v_\bullet / A[u]_\bullet,$ from (iii), implies $p^s v_\Pi/ A_\bullet
= (v_\bullet/{A[u]_\bullet})/A_\bullet 
= v_\bullet/ A_\bullet =
p^s v^\bullet,$ by (iv). Cancel $p^s.$   \hfill $\square$

\vs\noi
{\sc Lemma 2C.} { \it Every finitely generated $\la\pi\ra$-torsion
$\Delta  $-module $M$ has $\pd_\Delta  (M)\le 2.$}

\vs\noi
\begin{proof}
Suppose, first, that $M$ is a cyclic $\Delta$-module, so $M= \Delta
 /I$ for some left ideal $I$ containing $\Pi^z$ for the unique least $z\ge 0,$ 
hence
$I\supseteq \Delta  \Pi^z:$ for 
$\pi^yM=0$ implies $I\supseteq
\pi^y\Delta   =\Delta  \Pi ^{p^sy}.$
Since $\Ext^k_\Delta  (M,X) =
\,\Ext^{k-1}_\Delta  (I,X)$ for $k\ge 2,$
we need to show that
$\Ext^{k-1}_\Delta  (I,X) =0$ for $k\ge 3$ and all $\Delta 
$-modules $X.$

\vs\noi
We do this by induction on $z\ge 0;$ the case $z=0$ holds by
$I=\Delta$   projective.  For $z\ge 1,$ the 
exact sequence $0\to I\to I +\Delta \Pi^{z-1} \to
\frac{I+\Delta\Pi^{z-1}}{I} \to 0$
and
isomorphism 
$
\frac{\Delta  \Pi^{z-1}}{I\cap \Delta  \Pi^{z-1}} \to \frac{I+\Delta
 \Pi^{z-1}}{I}
$
combine to give the exact sequence
$0 \to I\to I+ \Delta  \Pi^{z-1}\to \frac{\Delta  \Pi^{z-1}}{I\cap
\Delta  \Pi^{z-1}}\to 0 $
hence  exactness of 
$$
\Ext^{k-1}_\Delta  (I+\Delta  \Pi^{z-1},X) \to \,\Ext^{k-1}_\Delta  (I,X)
\to\,\Ext^{k}_\Delta  (\Delta  \Pi^{z-1}/I\cap \Delta  \Pi^{z-1},X).
\leqno {\rm (2C1)}
$$

\noi
The left $\Ext$ is $0,$ by the induction hypothesis, since $I+
\Delta  \Pi^{z-1}\supseteq \Delta  \Pi^{z-1},$ so it suffices to show
that the right $\Ext$ is also $0,$ for all $k\ge 3$ and $X.$

\vs\noi
For that,
$\Delta \Pi^z\subseteq I\cap \Delta \Pi^{z-1} \subseteq \Delta
\Pi^{z-1}$ gives
the exact sequence
$$
0 \to I\cap \Delta  \Pi^{z-1}/\Delta  \Pi^z \to \Delta 
\Pi^{z-1}/\Delta  \Pi^{z}\to \Delta 
\Pi^{z-1}/I\cap \Delta  \Pi^{z-1}\to 0,
\leqno{\rm (2C2)}
$$
with middle term  $\simeq \Delta  /\Delta  \Pi \simeq
\IF[[\,\ol{u-1}\,]],$ by (2A1), allowing (2C2) to be rewritten
$$
0\to \wt I \to \IF\big[[\,\ol{u-1}\,]\big] \to {\Delta 
\Pi^{z-1}}/{I\cap \Delta  \Pi^{z-1}} \to 0, 
$$
with $\wt I$ a $\Delta$-submodule  of $\IF\big[[\,\ol{u-1}\,]\big],$ which implies
exactness of
$$
\Ext^{k-1}_\Delta  \,(\wt I,X)\to\,\Ext^k_\Delta  
\Big(\frac{\Delta\Pi^{z-1}}{I\cap \Delta  \Pi^{z-1}}, X\Big)\to
\Ext^k_\Delta  \big(\IF\big[[\,\ol{u-1}\,]\big],X\big)
\leqno{\rm (2C3)}
$$
 for all $k\ge 3$ and$X.$

\vs\noi
In particular, applying $\Ext_\Delta(-,X)$ to (2A1), we find that
${\rm pd}_\Delta  \IF[[\,\ol{u-1}\,]] \le 1,$ hence that
the right
$\Ext$ in (2C3) is $0,$ so
it
suffices to show that $\Ext^{k-1}_\Delta (\wt I,X) =0$ for all $k\ge
3$ and $X.$

\vs\noi
Now $\IF[[\,\ol{u-1}\,]]$ is a principal ideal domain, by the proof of
(1G1) (or [Ma], Theorem~19.5 and 20.8), so each ideal $\wt I\ne 0$
has a generator $\tilde i\ne 0$ and variant to (2A1)
$$
0\to \Delta \os{_\bullet\Pi}\lra \Delta \os \sim\lra \wt I\lra 0
$$
with $\tilde d = \ol d\,\tilde i$ for $d\in \Delta.$ 
Here $\tilde{}$  is
surjective since $\wt I =\IF[[\,\ol{u-1}\,]]\tilde i $ and
the sequence 
exact since $\ol
d\,\tilde i =0$ implies $\ol d = \ol 0$ hence $d\in \Delta \Pi.$  
Applying $\Ext_\Delta (-,X)$ to this we get 
$$
\Ext^{k-2}_\Delta (\Delta,X) \to \Ext^{k-1}_\Delta(\wt I,X) \to
\Ext^{k-1}_\Delta(\Delta,X),
$$
for every
$\wt I,$
since $k-2 \ge 1.$  This settles
 the case of cyclic
$\Delta  $-modules $M.$

\vs\noi
Finally, we show that $pd_\Delta  M\le 2,$ for $M$ having $n$
generators, by induction on $n;$ the case $n\le 1$ is above.  
When $n\ge 2,$ let $M^\pr$ be the submodule of $M$ generated by the
first $n-1$ generators of $M,$ so $M/M^\pr$ is a cyclic module.  The
exact sequences $0\to M^\pr \to M\to M/M^\pr \to 0$ and
$\Ext^k_\Delta  (M/M^\pr,X) \to \,\Ext^k_\Delta  (M,X) \to \,
\Ext^k_\Delta  (M^\pr,X),$ for all  $k\ge 3$ and $X,$ finish the proof.  
\end{proof}

\noi
{\sc Corollary 2D.}  {\it The natural map $K_0\la\pi\ra (\Delta  ) \to
G_0\la\pi\ra(\Delta  )$ is an isomorphism.}

\noi
\begin{proof}
Here $\la \pi\ra$ is the multiplicative set $\{\pi^k:k\ge 0\}$ in $A$
and $G_0\la\pi\ra(\Delta)$ (respectively, $K_0\la\pi\ra(\Delta ))$
is the Grothendieck group of the category of finitely generated $\la
\pi\ra$-torsion $\Delta$-modules $M$ (resp. of the subcategory of
these $M$ with $\pd_\Delta(M) <\infty  ).$  These are the same by
Lemma~2C.
\end{proof}

\noi
{\sc Lemma 2E.} (i) $G_0\la \pi\ra(\Delta  ) \os\bullet\to G_0T(\Delta 
_\bullet)$ {\it  is an isomorphism.}

\noi
\begin{description}
\item {(ii)} {\it The cokernel of the natural map $K_1(\Delta 
[\pi^{-1}])\os\alpha  \to K_1(D)$ has no $\Z$-torsion (which is
non-trivial of finite order in the abelian group structure on $K_1).$}

\item  {(iii)}  {\it The composite map  $K_0\la\pi\ra(\Delta  )\to
G_0\la \pi\ra (\Delta)\to G_0T(\Delta_\bullet)$ is an isomorphism.
\newline In particular,
$K_0\la \pi\ra(\Delta)$ has no $\Z$-torsion.}

\item{(iv)} {\it The map $\beta  $  of} (2E1) {\it is injective.}

\item{(v)} {\it $G_0(\Delta  /\pi\Delta  ) \to G_0(\Delta  ) \to
G_0(\Delta  [\pi^{-1}]) \to 0$ is an exact sequence with
$G_0(\Delta/\pi\Delta  )\cong \Z.$}

\item  {(vi)} {\it  Let $S$ be a multiplicatively closed set of central
nonzerodivisors in a ring $R.$  Then
$$
K_1(R) \to K_1(S^{-1}R) \os \partial\to K_0\big(\bH\!_S(R)\big) \to
K_0(R) \to K_0(S^{-1}R)
$$
is an exact sequence, with $\bH\! _S(R)$ the exact subcategory of
$S$-torsion $R$-modules in the category $\bH (R)$ of all $R$-modules
having a resolution by finitely generated projective $R$-modules. }

\end{description}

\vs\noi
{\it Proof.} {\it For} (i): By [Sw,5.15] the displayed sequence is exact:
$$
G_0 N(\Delta)\os\iota\to G_0\la \pi\ra (\Delta  ) \os\bullet\to G_0T(\Delta 
_\bullet) \to 0.
$$
In it, $G_0N(\Delta  )$ is the Grothendieck group of the category of
finitely generated $\Delta  $-modules which are pseudo-null as
$A$-modules (hence finite, by [NSW], Remark~4 on p.269), $G_0\la
\pi\ra(\Delta)  $ that of finitely generated $\la\pi\ra$-torsion
$\Delta    $-modules, and $G_0T(\Delta  _\bullet)$ that of finitely generated
torsion $\Delta  _\bullet$-modules.

\vs\noi
Assertion (i) is the vanishing of $\iota,$ which we now prove, by
induction on the length of the composition series of our pseudo-null
(as $A$-module) $\Delta  $-module  $C$ above.  Such a $C$ of length
$1$ is a simple $\Delta  $-module, hence $C\simeq \Delta 
/\rad\,\Delta  =\IF$ by Lemma~2A, so the exact sequence
$0\to\IF[[\ol{u-1}]]\to \IF[[\ol{u-1}]]\to \IF\to 0$ of $\Delta 
$-modules implies $\iota (C) =[\IF] =0$ in $G_0\la\pi\ra(\Delta  ).$
 If $C$ has length $> 1,$ and $C_0$ a simple $\Delta  $-submodule,
then $0\to C_0\to C\to C/C_0\to 0,$ with $C/C_0$ of smaller length,
hence $\iota(C)=0.$

\noi
{\it For} (ii): The centre $A[\pi^{-1}]$ of $\Delta  [\pi^{-1}]$
is a Dedekind domain: for $A=\fO  _r[[T]]$ is a regular UFD (cf.
[Ma],
Theorem~20.8 and page~30, Theorem~19.2, 19.3 perhaps clearer), as $\fO_r$ is, with the unique maximal ideal $(\pi,T)$
(cf. [NSW] V\S3), which vanishes in the localization $A[\pi^{-1}].$

\vs\noi
It follows that the maximal
$A[\pi^{-1}]$-order $\Delta[\pi^{-1}]$
is hereditary 
 ([Re], Theorem~21.4).
With $S$ the central multiplicative set $A[\pi^{-1}]\bs 
\{0\}$ of $\Delta[\pi^{-1}],$ 
$S^{-1}(\Delta  [\pi^{-1}]) = D,$
it
follows that every finitely generated $S$-torsion
$\Delta[\pi^{-1}]$-module $M$ has $\pd_{\Delta[\pi^{-1}]}M\le 1,$
by (CRI, Proposition~4.3).

\vs\noi
In
particular, $K_0S(\Delta[\pi^{-1}]) \os\sim =
G_0S(\Delta[\pi^{-1}],$ as in the proof of Corollary 2D, and
$G_0S(\Delta[\pi^{-1}])$ is a free $\Z$-module by [CRII, 38.58] and
[CRI, 16.6], hence has no $\Z$-torsion.  The localization exact sequence
$ K_1(\Delta[\pi^{-1}]) \os \alpha \to K_1(D) \os\partial\to
K_0S(\Delta[\pi^{-1}]) $
of [CRII, 40.9] finishes the proof. 

\vs\noi
{\it For} (iii): Both maps are isomorphisms, by Corollary~2D and by
(i). In particular, $A_\bullet$ is a Dedekind domain, by the first
paragraph of the proof of (ii), so $G_0T(\Delta)$ has no
$\Z$-torsion, by [CRII,38.58] and [CRI,16.6], since the $A/P$ are
fields.

\vs\noi
{\it For} (iv): The composite isomorphism of (iii) fits into a
commutative diagram
$$
\begin{matrix}
K_0\la\pi\ra(\Delta  ) &\lra &G_0\la\pi\ra(\Delta  ) &\os\bullet\lra
&G_0T(\Delta  _\bullet)\\ \\
\beta  \downarrow &&\downarrow &&\big\Vert\\ \\
K_0T(\Delta  ) &\lra &G_0T(\Delta  ) &\os\bullet\lra &G_0T(\Delta  _\bullet),
\end{matrix}
\leqno{\rm (2E1)}
$$
by (iii). The squares commute since the maps are all $[M]\mapsto [M]$
or
$[M_\bullet].$
The
diagram implies the property of $\beta.  $

\vs\noi
For (v):
The exact sequence is from ([We], Application~6.4.1 with $G_0$ as in 
Definition~6.2), so we turn to $G_0(\Delta  /\pi\Delta  ).$
The ideal $\Delta  \Pi = \Pi\Delta  $ of $\Delta  $ has
$(\Pi\Delta  )^{p^s}= \Pi^{p^s}\Delta   = \pi\Delta  ,$ by  (2A3*).
  Thus $G_0(\Delta  /\pi\Delta  ) = G_0\big(\Delta 
/(\Pi\Delta  )^{p^s}\big) \os 1= G_0(\Delta  /\Pi\Delta  ) \os 2=
G_0(\IF[[U]]),$
with $\os 1=$ by ([We], Corollary 6.3.1), and $\os 2=$ by (2A1).

\noi
Since $\IF[[U]]$ is a principal ideal domain (by the
proof of Lemma~2C), the result follows from [We, Example~6.2.1].
Here we use the Definition~6.2 convention but soon revert to
$K_0.$

\vs\noi
For (vi):The right half of the exact sequence is Corollary~7.7.4 of
(II of [We]2.),
and the left half starts in the second paragraph of the proof of
Theorem~3.2 (III of loc.cit.).  The exactness of
$K_0\big(\bH\!_S(R)\big)$ is in Lemma~3.1.5, and description of the
homomorphism $\partial:K_1(S^{-1}R) \to K_0\big(\bH\!_S(R)\big)$ in
Corollary~3.1.1.  Background on $\bH(R)$ starts at Definition~7
(of II).

\vs\noi
Restricting our $S$ to nonzerodivisors  seems adequate, for our
current needs, to cope with Theorem~3.2 (of III in [We]2.) in exact
categories, via Example~9.1.3, without Waldhausen generality.
\hfill $\square$

\vs\noi
{\sc Theorem 2F.} {\it The natural map $SK_1(\Delta  ) \to
SK_1(D)$ is surjective.}

\vs\noi
\begin{proof}
We first show that all $\Z$-torsion in $K_1(D)$ is in the image of
$K_1(\Delta  )\to K_1(D).$  This is a consequence of Lemma~2E and
the commutative diagram below with rows localization sequences 
from [CRII, 40.9],
$$
\begin{matrix}
K_1(\Delta) &\lra &K_1(\Delta  [\pi^{-1}]) &\os{\partial^\pr}\lra K_0\la\pi\ra(\Delta  )\\
\\ 
&&\alpha \downarrow &\beta  \downarrow\\ \\
 K_1(\Delta  ) &\lra&{}\q K_1(D)\q{} &\os\partial\lra K_0T(\Delta  ).
\end{matrix}\leqno{\rm (2F1)}
$$

\vs\noi
Namely, if $x\in K_1(D)$ is $\Z$-torsion and not in the image of
$K_1(\Delta  ),$ then $\partial (x)\in K_0T(\Delta  )$ has finite
order $n\ge 2.$  But, by (ii), $x= \alpha  (x^\pr)$ for some $x^\pr
\in K_1(\Delta  [\pi^{-1}]),$ so $\beta  \partial^\pr(x^\pr) =
\partial (x)$ yields $\beta\partial^\pr(x^{\pr\, n}) =0,$ whence
$n\partial^\pr(x^\pr) =0,$ since $\beta  $ is injective by (2E1). It
follows that $\partial^\pr(x^\pr)=0,$ by (iii), and thus $\partial
(x) = \beta  (0) =0$  contradicts  $n\ge 2.$

\vs\noi
To finish the proof, recall that every $x\in SK_1(D)$ has finite
order, by [Dr.,~p.157]. Thus, by the above, it is the image
of some $y\in K_1(\Delta  ).$  Now $\nr(y)=\nr(x)=1$
confirms that $y\in SK_1(\Delta  ).$
\end{proof}

\noi
{\sc Remark:} Checking commutativity of the (2F1) square, is perhaps easier
by replacing [CRII40.9] by Lemma~2E(vi) (i.e. Theorem~3.2 on
[We, p.220] helps). Indeed (v), (vi) above were intended as
preparation for carrying this out.

\vs\noi
We next exhibit some strong congruences on $B$ and on the action of
our special generator $g$ of $\fg$ on $B,$ with notation fixed by
setting $\Delta  =\Delta  _{r,s,\IF}$ as usual.

\vs\noi
{\sc Lemma 2G.}  (i) {\it If $0\ne b\in \Pi B$ and $z\ge 1,$ then}
$$
(1-b)^{p^z} \equiv (1-b^{p^z}) + pL(b^{p^{z-1}}) + p^2 M_z(b) \!\!\!\mod
p^3B,
$$
\begin{description}
\item  {\q} {\it with} $L(b) = -\sum^{p-1}_{k=1}\,\frac{b^k}{k}\,,
M_1(b) =\sum^{p-1}_{k=1}\,H_k\,\frac{b^k}{k}\,, \,H_k
=\sum^{k-1}_{j=1}\, \frac 1j\,,$
\item  {\q} {\it and} $M_z(b) = (1-b^{p^{z-1}})^{p-1} L(b^{p^{z-2}}) +
M_1(b^{p^{z-1}})$ {\it for} $z\ge 2.$

\item  {(ii)} {\it In particular,} $v_\Pi\big((1-\Pi)^{p^z} -
(1-\Pi^{p^z})\big) = v_\Pi\big(pL(\Pi^{p^{z-1}})\big)$ {\it for}
$1\le z\le r+s.$

\item  {(iii)} $^g\Pi \equiv \Pi + (1-\Pi)\Pi^{p^r}\!\!\!\mod
p\Pi^{p^{r-1}}B.$

\item  {(iv)} $^{g^{p^z}}\Pi \equiv \Pi +
\Pi^{p^{r+z}}\!\!\!\mod \Pi^{p^{r+z}+1}B,$ {\it when}
$1\le z\le
s-1.$
\end{description}

\vs
\begin{proof}
 Note that $(1-b)^{p^z} = (1-b^{p^z})$ when $z=0,$ and that $1\le k\le
p-1$ induces
$$
-\,  \frac 1 p \,\begin{pmatrix} p\\k\end{pmatrix}\dot=\,
\frac{(-1)^k}{k}\,\prod\nolimits^{k-1}_{j=1}\big(1-\frac{p}{j}\big)
\equiv \frac{(-1)^k}{k}\, (1-pH_k) \!\!\!\mod p^2 
$$
with $\dot =$ by 
$-\,\frac 1p \,\begin{pmatrix} p\\ k\end{pmatrix} = -\,\frac
1k\,\frac{(p-1)!}{(k-1)!(p-k)!} 
= \frac{-1}{k}\,\frac{\prod\nolimits_j(p-j)}{\prod\nolimits_j j}
= \frac{(-1)^k}{k} \prod_j \,\frac{(j-p)}{j}, 
$  (here $0!=1).$  So

$$
(1-b)^p - (1-b^p) =\sum\nolimits^{p-1}_{k=1} \begin{pmatrix}p\\k\end{pmatrix}
(-b)^k\equiv
-p\sum\nolimits^{p-1}_{k=1}\,\frac{(-1)^k}{k}\,(1-pH_k)(-b)^k\!\!\!\mod p^3B,
\q\text{hence}
$$
$$
(1-b)^p \equiv (1-b^p) + pL(b) + p^2M_1(b)\!\!\!\mod p^3B.
\leqno{\rm (2G1)}
$$

\vs\noi
{\it For} (i): (2G1) starts an induction on $z\ge 1.$ 
 If $z\ge 2,$ then the induction hypothesis for $z-1$ reads
$
(1-b)^{p^{z-1}}\equiv  (1-b^{p^{z-1}})
+ [pL(b^{p^{z-2}}) + p^2M_{z-1}(b)]\!\!\mod p^3B,
$
so raising to the $p^{\rm{th}}$ power gives
$(1-b)^{p^z}\equiv \big([1-b^{p^{z-1}}] + [pL(b^{p^{z-2}}) 
+ p^2M_{z-1}(b)]\big)^p
 =\sum\limits^{p}_{k=0} \begin{pmatrix} p\\
k\end{pmatrix}[1-b^{p^{z-1}}]^{p-k}\cdot
 [pL(b^{p^{z-2}})
+p^2 M_{z-1}(b)]^k.$
Now subtracting by $(1-b^{p^{z-1}})^p,$ we have
$$
(1-b)^{p^z} - \Big(1-b^{p^{z-1}}\Big)^p \dot\equiv \sum^1_{k=1} 
\begin{pmatrix} p\\ 1\end{pmatrix} (1- b^{p^{z-1}})^{p-1}
\cdot \big(pL(b^{p^{z-2}}) + p^2M_{z-1}(b)\big)\!\!\mod p^3B,
$$
with $\dot\equiv$ because the $k^{\rm{th}}$ term is in $p^3B$ for
$k> 1,$ since $1<k<p$ has $1+k\ge 3$ and $k=p$ has $k\ge 3.$

We add this relation to the relation obtained from (2G1) by
replacing $b$ by $b^{p^{z-1}}.$  This yields
$(1-b)^{p^z} - (1-b^{p^z}) \equiv pL(b^{p^{z-1}}) +
p^2\big[(1-b^{p^{z-1}})^{p-1}L(b^{p^{z-2}}) + M_1(b^{p^{z-1}})\big]
\equiv pL(b^{p^{z-1}})+ p^2 M_z(b)\!\!\mod p^3 B,$ completing the
induction.

\vs\noi
Finally, $v_\Pi\big(L(b)\big) = v_\Pi(b)$ is clear, and
$$
v_{\Pi}\big(M_z(b)\big) = \tilde zv_\Pi(b), \;\text{\rm with}\q
\tilde z = \begin{cases} 2 &z=1\\ p^{z-2} &z\ge 2.\end{cases}
\leqno{\rm (2G2)}
$$
This holds for $z=1,$ because $H_1=0$ and $H_2 =1,$ and then for
$z\ge 2$ by definition of $M_z(b),$ since the terms have different
$v_\Pi,$ completing (i).

\vs\noi
{\it For} (ii): We first check, for $1\le z\le r+s,$ that
$$
\text{1.} \; v_{\Pi}\big(pL(\Pi^{p^{z-1}})\big) <
v_{\Pi}\big(p^2M_z(\Pi)\big),\q\text{and \q
2.}\;v_\Pi\big(p^2M_z(\Pi)\big) <v_\Pi(p^3).
$$

\vs\noi
1. holds $\Llra v_\Pi\big(L(\Pi^{p^{z-1}})\big) <
v_\Pi\big(pM_z(\Pi)\big) \Llra p^{z-1} -\tilde z <\phi(p^{r+s}),$ by
(2G2).  This is trivial for $z=1,$ and works for $z\ge 2 \Llra
\phi(p^{z-1}) <\phi(p^{r+s})\Llra z-1 <r+s.$

\vs\noi
2.  holds $\Llra v_\Pi\big(M_z(\Pi)\big) <\phi(p^{r+s})\Llra \tilde
z <\phi(p^{r+s}),$ which works similarly.

\vs\noi
The congruence of (i), with $b=\Pi, $ implies that
$$
(1-\Pi)^{p^z} - (1-\Pi^{p^z}) \equiv pL(\Pi^{p^{z-1}})\!\!\!\mod
p^2\Pi^{\tilde z},\leqno{\rm (2G3)}
$$
by 2. and (2G2), which implies (ii) by 1.

\vs\noi
{\it For} (iii):  
 $(1-\Pi)^{p^{r}}= \zeta
 ^{p^r} =\zeta  _s = 1-\pi_s\,,$ implies
$$
\pi_s \equiv \Pi^{p^r} \!\!\!\mod p\Pi^{p^{r-1}}, \leqno{\rm (2G4)}
$$
since $\pi_s =1 - (1-\Pi)^{p^r} \dot\equiv \,1- (1-\Pi^{p^r}) =
\Pi^{p^r}\!\!\!\mod p\Pi^{p^{r-1}}$ with $\dot\equiv $ by (2G3) and
(2G2).

\vs\noi
Now $^g\Pi -\Pi = (1-\,^g\zeta) - (1-\zeta  ) = \zeta  -\,^g\zeta 
=\zeta  - \zeta  ^{1+p^r}
= \zeta 
(1-\zeta  _s) = (1-\Pi)\pi_s\,\dot\equiv\, (1-\Pi)\Pi^{p^r}\!\!\!\mod
p\Pi^{p^{r-1}},$ with 
 $\dot\equiv $ by (2G4).  This proves the conguence (iii).

\vs\noi
A variant of this is given by
$$
\pi_s\equiv \Pi^{p^r} + p\sum_{0<j<p} \Pi^{jp^{r-1}}/j \!\!\!\mod
p^2 \Pi^{\wt r}.
\leqno{\rm (2G4^\#)}
$$
The fragment following (iii) begins $\pi_s = 1- (1-\Pi)^{p^r}$
by
$$
\dot\equiv 1 -\big( (1-\Pi^{p^r})
+ pL (\Pi^{p^{r-1}})
\big) \equiv 1
-1+\Pi^{p^r} +p\sum_{0<j<p} \Pi^{jp^{r-1}}/j \!\!\! \mod p^2\Pi^{\tilde r},
$$
with
$\dot\equiv $ by (2G3).

\vs\noi
{\it For} (iv):
First, we show that

\vs
\noi
3. $\v_\Pi\big(pL(\Pi^{p^{r+z-1}})\big)>
p^{r+z}$ if, and only if, $z< s.$

\vs
\noi
4. $(1-\Pi^{p^{r+z}})^p \equiv 1\!\!\mod \Pi^{p^{r+z+1}},$ when
$z\le s-1.$

\vs\noi
{\it For} 3: By (2G2), 
$\v_\Pi\big(pL(\Pi^{p^{r+z-1}})\big) = \v_\Pi(p) + \v_\Pi(\Pi^{p^{r+z-1}})
 = \phi(p^{r+s}) + p^{r+z-1} > p^{r+z}\Llra
\phi(p^{r+s}) > \phi (p^{r+z}) \Llra  r+s > r+z.$

\vs\noi
{\it For} 4: $\big(1-\Pi^{p^{r+z}}\big)^p = 1+ \os{p-1}{\us{k=1}\sum} 
\begin{pmatrix} p\\ k\end{pmatrix}\big(-\Pi^{p^{r+z}}\big)^k 
-p^{r+z+1}$ with $\v_\Pi\big(p\Pi^{p^{r+z}}\big) = \phi(p^{r+s}) +
p^{r+z} \ge p^{r+z+1}$ implies 4.

\vs\noi
Now by induction on $z\ge 0,$
$$(1+p^r)^{p^z}\equiv 1+p^{r+z}\!\!\!\mod
p^{2r+z}.
\leqno{\rm (2G5)}
$$
This is clear for $z=0.$  The induction hypothesis for $z\ge 1$ has
$(1+p^r)^{p^{z-1}} =1+ p^{r+z-1} + p^{2r+z-1}h = 1+
p^{r+z-1}(1+p^rh)$ for some $h\in\Z.$  Taking $p^{\rm th}$ powers
gives
$$
(1+p^r)^{p^z} = 1+ \sum_{1\le k\le p} \begin{pmatrix} p\\
k\end{pmatrix} [p^{r+z-1}(1+p^rh)]^k.
\leqno{\rm (2G6)}
$$
The $k=1$ term is
$p^{r+z}(1+p^rh) = p^{r+z} + p^{2r+z}h \equiv p^{r+z}\!\!\!\mod
p^{2r+z},$ as required.

\noi
The $k^{\rm th}$ term for $2\le k\le p-1$ is a multiple of 
$p[p^{r+z-1}(1+p^rh)]^2,$ which is a multiple of 
$p^{1+2(r+z-1)} = p^{2r+2z-1} \equiv 0\!\!\!\mod p^{2r+z},$ as
required since $2z-1\ge z$ by $z\ge 1.$

\noi
Finally the $k=p$ term is a multiple of $(p^{r+z-1})^p =
p^{p(r+z-1)} $ so we need $p(r+z-1)\ge 2r+z \Llra (p-2)r + (p-1)z \ge
p.$  Since $r\ge 1$ and $z\ge 1,$ this needs $p\ge 3$ (and $p=2$
doesn't work).  This proves (2G5).

\vs\noi 
In the notation of
the proof of 
(2G5), with $1\le r +z \le r+s,$
we have
$$
^{g^{p^z}} \zeta = \zeta^{(1+p^r)^{p^z}} _{r+s}
= \zeta (1-\Pi)
^{p^{r+z}(1+p^rh)}=\;
  \zeta  \big((1-\Pi)^{p^{r+z}}\big)^{(1+p^rh)} \dot\equiv \;\zeta 
(1-\Pi^{p^{r+z}})^{1+p^rh}\!\!\!\mod pL(\Pi^{p^{r+z-1}}),
$$
for some $h\in \Z,$
with $\dot\equiv$  induced by the first map of (2G3), and then
taking $(1+p^rh)$-th power and left multiplication by $\zeta$ on both
sides.
Now $\v_\Pi\big(pL(\Pi^{p^{r+z-1}})\big) \ge p^{r+z}+1,$ by 3., and
$(1-\Pi^{p^{r+z}})^p\equiv 1\!\!\!\mod \Pi^{p^{r+z+1}},$ by 4., as
$p^{r+z+1} \ge p^{r+z}+1.$  Applying this to the right end of the
above congruence yields
$$
{}^{g^{p^z}}\zeta \equiv \zeta \big(1-\Pi^{p^{r+z}}\big) \!\!\!\mod
\Pi^{p^{r+z}+1}.
\leqno{\rm (2G7)}
$$

\vs\noi
Thus
$^{g^{p^z}}\Pi -\Pi 
=\zeta  -\,^{g^{p^z}}\zeta  \dot\equiv\; \zeta 
-\zeta(1-\Pi^{p^{r+z}}) = 
\zeta\Pi^{p^{r+z}} = (1-\Pi)\Pi^{p^{r+z}}
\equiv \Pi^{p^{r+z}}
\!\!\!\mod
\Pi^{p^{r+z}+1},$
with $\dot\equiv$ by (2G7).
\end{proof}

\vs\noi
In particular, this allows us to discuss the

\vs\noi
{\sc Example.} Consider the commutator $[u,\Pi] =\,^{g-1}\Pi \in
\Delta  ^\times,$ which has reduced norm $1$ as it is in
$[D^\times,D^\times].$ 
Observe that $^{g-1}\Pi \equiv
1+\Pi^{p^r-1}\!\!\!\mod \Pi^{p^r}\Delta  ,$ by Lemma 2G(iii), and
that $\Delta  /\Pi^{p^r}\Delta  $ is a commutative ring, since $u\Pi =
[u,\Pi]\Pi u \equiv (1-\Pi^{p^r-1})\Pi u\equiv \Pi u\!\!\!\mod
\Pi^{p^r}\Delta  .$ Thus $[u,\Pi]\notin [\Delta  ^\times,\Delta 
^\times],$
as it is non-trivial  under $K_1(\Delta 
)\to K_1(\Delta  /\Pi^{p^r}\Delta  ) = (\Delta  /\Pi^{p^r}\Delta 
)^\times.$ $\;\square$

\vs\noi
This is generalized by the map $SK_1(\Delta  ) \succ\!\!\to K_1(\Delta 
/\pi_1\Delta  )$ of Theorem 3E, and suggests the

\vs\noi
{\sc Lemma 2H.}  $[\Pi,_-]$: $K_1(\Delta  )\to SK_1(\Delta  )$ 
{\it is a group homomorphism.}

\vs\noi
\begin{proof}
 Conjugation by $\Pi$ is a ring automorphism of $\Delta  ,$
hence induces a group automorphism of $K_1(\Delta  );$ the
``difference'' of this and the identity map is thus an abelian group
homomorphism which takes $x\in \Delta  ^\times$ to $[\Pi,x],$ so has image in
$SK_1(\Delta  ).$
Alternatively,
$
[\Pi,xy] = [\Pi,x]\big[x,[\Pi,y]\big] [\Pi,y] \equiv
[\Pi,x][\Pi,y]\!\!\!\mod [\Delta  ^\times,\Delta  ^\times]
$
for all $x,y\in \Delta  ^\times.$
\end{proof}

\vs\noi
{\sc Definition.} {\it  Call $[\Pi,x],$ with $x\in \Delta  ^\times,$
a}
$\Pi$-commutator, {\it and any finite product of $\Pi$-commutators
or elements of}  $[\Delta  ^\times, \Delta  ^\times],$ {\it  a} $\Pi$-commutator product.

\vs\noi
{\sc Remark.} 
This precision is helpful when building elements of $\Delta 
^\times$ level by level, especially in contrast to the Lemma~2H
blurring.
 Note that $\Pi$ is essentially the only
element of $D^\times \bs \Delta^\times$ normalizing $\Delta  $ (see Proposition B2).

\section{The Map $SK_1(\Delta  ) \to K_1(\Delta  /\pi_1\Delta 
)$}\label{Sec3}

We return briefly to the dimension-1 $p$-adic Lie groups $G$ of the
introduction.  For our purposes this means that $G=H\rtimes \Gamma ,$
with $H$ finite and $\Gamma  \simeq \Z_p.$

\vs\noi
{\sc Proposition 3A.} {\it Let $C\simeq \Z_p$ be a central open
subgroup of $G$ as above, let $\wt G= G/C,$ and $\aug:\Lambda 
^c(\Gamma)\to \Z_p{}^c$ the augmentation map sending $\Gamma  $ to $1.$ 
Then
$$
\begin{matrix} K_1(\Lambda  _{\fO}G) &\os{\defl}\lra &K_1(\fO[\wt
G])\\
\\
\Det\,\downarrow &&\Det\,\downarrow\\
\\
\Hom^F\big(R_p(G),(\Lambda  ^c\Gamma  )^\times \big) &\os
\varepsilon  \lra &\Hom^F\big(R_p(\wt
G),(\Z_p{}^c)^\times\big)\
\end{matrix}
$$
commutes, with $\varepsilon  $ induced by \aug, i.e., $\varepsilon 
(f)(\wt\chi)= \,\aug\big(f(\infl\, \wt\chi)\big).$

}

\vs\noi
\begin{proof}
It is divided into three steps, the first about the left $\Det,$ the
second about the right $\Det,$ and the third bringing them together.
Note that the finite extension $F/\Qbar_p$ need not be unramified
here, and that we write  $^-{}$ for the canonical map $G\to \Gamma ,
$ and $\Lambda  ^cG\to \Lambda  ^c\Gamma.$

\vs\noi
{\it Step} 1:
$$
\begin{matrix}
K_1(\Lambda  _{\fO}G) &\lra &K_1(\SQ_F G)\\
\\
\downarrow &&\Det \downarrow\\
\\
\Hom^F\big(R_p(G),(\Lambda  ^c\Gamma  )^\times\big) &>\!\!\lra
&\Hom^F\big(R_p(G),\SQ^c(\Gamma  )^\times \big)
\end{matrix}
$$

\vs\noi
Our left $\Det$ is defined by the left vertical arrow in the
commutative square shown,which is built on the $\Det$ at the right
that is described in the first paragraph of [RW3, \S2]: here the
superscript $F$ generalized the $*$ of our introduction (as in the
first sentence on p.34 of loc.cit.).  The horizontal maps are
induced by $\Lambda  _{\fO}\hookrightarrow \SQ_F,$ the bottom map is
clearly injective.  This requires us to show, for $y\in \Lambda 
_{\fO}G^\times,$ that the right action of $y$ on $\fV_\chi$ has
$(\Det y)(\chi  )$ in $\Lambda  ^c(\Gamma  _k)^\times,$ for every $\chi\in
R_p(G).$  Brauer induction implies that it suffices to do this for
all $\chi =\,\ind^G_{G^\pr}(\chi^\pr)$ with $\chi^\pr(1)=1.$

\vs\noi
Fix such a $\chi$ and a coset decomposition $G=\dot\bigcup_{t\in T}G^\pr t.$
 For each $g\in G$ and $t\in T$ write $tg = x_g t_g$ with $x_g\in
G^\pr$ and $t_g\in T,$ as in the proof of [RW3, Proposition~3],
which constructs a $\SQ^c(\Gamma  _k)$-basis $\{f_t:t\in T\}$ of
$\fV_\chi$ on which $g\in G$ acts by
$$
f_t g= \big(\chi^\pr(x_g)\bar x_ g)f_{t_g}\,.
$$
Choose $C_1\simeq \Z_p$ central open in $G$ and contained in $\ker
\chi^\pr,$ and write $G=\dot\bigcup_{w\in W} C_1 w,$ hence $y=\sum_{w\in
W} y_w w$ with $y_w \in \Lambda  _{\fO}(C_1).$  Observe that $t\in
T,$ $c\in C_1$ implies $x_c= c,$ $t_c=t,$ hence $f_tc = \big(\chi 
^\pr(c)\bar c\big)f_t = \bar cf_t.$  Thus
$$
\begin{aligned}
f_ty &= \sum\nl_w f_ty_w w = \sum\nl_w \bar y_w f_tw =\sum\nl_w \bar
y_w\big(\chi  ^\pr(x_w)\bar x_w\big)f_{t_w}\\
&= \sum\nl_w \chi  ^\pr(x_w)\bar y_w \bar x_w f_{t_w} =\sum\nl_{u\in T}
\Big(\sum\nl_{w:t_w=u}\chi  ^\pr(x_w)\bar y_w \bar x_w\Big) f_u,
\end{aligned}
$$
which implies
$$
(\Det \,y) (\chi  ) = \det\Big(\sum\nl_{w:t_w=u} \chi ^\pr(x_w)\bar
y_w \bar x_w\Big)_{(u,t)\in T\times T}\,.
\leqno{\rm (3A1)}
$$
The $T\times T$ subscript means: choose a well-ordering of $T$ and
use it for both rows and columns; the determinant is independent of
the choice.  The entries in this matrix are all in $\Lambda 
^c(\Gamma  ),$ hence $(\Det\,y)(\chi)\in \Lambda  ^c(\Gamma  )$
and, thus, in $\Lambda  ^c(\Gamma  )^\times$ since  $(\Det\,y)(\Det\,
y^{-1}) = 1.$

\vs\noi
{\it Step} 2:
$$
\begin{matrix}
K_1(\fO [N]) &\lra &K_1(F[N])\\
\\
\downarrow &&\Det \downarrow\\
\\
\Hom^F\big(R_p(N),(\Z_p{}^c)^\times \big) &>\!\!\lra
&\Hom^F\big(R_p(N),\Qp {}^c)^\times\big)
\end{matrix}
$$

\vs\noi
Our right $\Det$ on $K_1(\fO[\wt G])$ is analogous to that on
$K_1(\Lambda  _{\fO}G)$ in Step 1, but originates in Fr\"ohlich's Hom
description $\Det$ on $K_1(F[N])$ for finite groups $N,$ and is
defined by the left vertical arrow in the commutative square shown:
here the superscript $F$ refers only to the $G_{\Qp{}^c/F^-} $
condition (there is no $W$-twist).  
Now the Fr\"ohlich $\Det(z)(\xi  )$ is the determinant of the right
action of
$$
z\in F[N]^\times \;\text{on}\;\Hom_{\Qp{}^c[N]}(V_\xi  ,\Qp{}^c[N]),
$$
for all $\xi  \in R_[(N).$  (This formulation is chosen for
compatibility with Step~1; to compare with a more usual version see
Lemma~A1 of [GRW, p.75].)

\vs\noi
As in Step 1, we show that $z\in\fO[N]^\times$ implies that $(\Det\,
z)(\xi  ) \in (\Z_p{}^c)^\times,$ and Brauer induction allows us to
assume that $\xi  =\,\ind^N_{N^\pr}(\xi  ^\pr)$ with $\xi 
^\pr(1)=1.$

\vs\noi
Fixing $z,\xi  ,$ we write $N=\dot\bigcup_{\tau  \in T_N}N^\pr \tau  $
and $\tau  n= \kappa  _n\tau  _n$ with $\kappa  _n\in N^\pr,$ $\tau 
_n\in T_N$ for each $\tau  \in T_N,$ $n\in N,$ and $z=\sum_{n\in
N}z_n n$ with $z_n \in \fO.$  As in Step~1, it remains to show that 
$$
(\Det\,z)(\xi  ) =\, \det\Big(\sum\nl_{n\in N:\tau  _n= \nu  } \xi 
^\pr(\kappa  _n)z_n\Big)_{(\nu  ,\tau)  \in T_N\times T_N}\,.
\leqno{\rm (3A2)}
$$
For this, $V_\xi  =\Qp{}^c[N]\bigotimes _{\Qp{}^c[N^\pr]} V_{\xi 
^\pr}$ and Frobenius reciprocity imply that
$$
\Hom_{\Qp{}^c[N]} (V_\xi  ,\Qp{}^c[N]) = \bigoplus_{\tau  \in T_N}
\,\Hom_{\Qp{}^c[N^\pr]}(V_{\xi  ^\pr},\Qp{}^c[N^\pr]\tau  )
$$
has $\Qp{}^c$-basis $\{j_\tau  :\tau  \in T_N\}$ with $j_\tau 
(n_0\otimes v^\pr)= n_0e_{\xi  ^\pr \tau  },$ for all $n_0\in
\Qp{}^c[N],$ where $v^\pr$ is a $\Qp{}^c$-basis of $V_{\xi  ^\pr}$
and $e_{\xi  ^\pr} = \frac{1}{\vert  N^\pr\vert  }\,\sum_{n\in
N^\pr} \xi  ^\pr(n^{-1})n.$  It follows, for $n\in N,$ that
$$
(j_\tau  n)(n_0\otimes v^\pr) = j_\tau  (n_0\otimes v)n = n_0e_{\xi 
^\pr}\tau  n = n_0e_{\xi  ^\pr}\kappa  _n \tau  _n = \xi 
^\pr(\kappa  _n)n_0 e_{\xi  ^\pr}\tau  _n = \xi  ^\pr(\kappa 
_n)j_{\tau  _n}(n_0\otimes v^\pr),
$$
hence that $j_\tau  n = \xi  ^\pr(\kappa  _n)j_{\tau  _n}.$ 
Therefore
$$
j_\tau  z = \sum_{n\in N} z_n j_\tau  n =\sum_n \xi  ^\pr(\kappa 
_n)z_n j_{\tau  _n} = \sum_{\nu  \in T_N}\Big(\sum_{n\in N: 
\tau  _n=\nu}\xi  ^\pr(\kappa  _n)z_n\Big)j_\nu  .
$$

\vs\noi
{\it Step} 3:

\noi
Finally, we must show, for all $y\in \Lambda  _{\fO} G^\times,$ that
$\big(\varepsilon  \circ (\Det\, y)\big)(\wt\chi  ) = (\Det\,\wt y)(\wt
\chi  ),$ with $\wt y =\,\defl(y),$ for all $\wt \chi  \in R_p(\wt
G);$ by Brauer induction, we may assume that $\wt\chi  =\,\ind^{\wt
G}_{\wt G^\pr}(\wt \chi  ^\pr)$ with $\wt\chi  ^
\pr(1) =1.$

\vs\noi
Fix $y,\wt \chi  $ and write the group homomorphism $G\to \wt G$ as
$g\mapsto \wt g.$  Setting $\chi  =\,\infl(\wt\chi  ),$ we have
$\chi  (g) =\wt\chi  (\wt g),$ for all $g\in G;$ also let $G^\pr$
be the pre-image of $\wt G^\pr$ under $G\to \wt G$ (hence $\wt G^\pr
= \wt{G^\pr})$ and $\chi  ^\pr = \,\infl(\wt \chi  ^\pr),$ implying
$\ker \, \chi  ^\pr \supseteq C.$  It follows, from [CR1, 10.2],
that $\chi  =\,\ind^G_{G^\pr}\,\chi  ^\pr$ and $\chi  ^\pr(1) = 1.$ 
This implies that
$$ (\varepsilon  \circ \,\Det\,y)(\wt \chi  ) =\,\det\,\Big(\sum\nl
_{w\in W:t_w=u} \chi  ^\pr(x_w)\,\aug\,(\ol{y_w})\Big)_{(u,t)\in
T\times T}
\leqno{\rm (3A3)}
$$
on combining the definition of $\varepsilon  ,$ (3A1), and $\aug
(\ol{x_w}) = 1,$ by $\ol{x_w}\in \ol{G^\pr}\subseteq \Gamma.$
Here the $w$'s are those of Step~1, so depend on a choice of $C_1.$ 
We now choose $C_1= C,$ so have $G= \dot\bigcup_{w\in W}Cw,$ hence
applying $\sim$ implies that $w\mapsto \wt w$ is a bijection $W\to
\wt G.$  Moreover, $y=\sum_{w\in W}y_w w,$ with $y_w\in \Lambda 
_{\fO}(C)$ implies that $\wt y =\sum_{w\in W} \,\aug(\ol{y_w})\wt
w,$ with $\aug (\ol{y_w})\in\fO.$

\vs\noi
We recall that $G=\dot\bigcup_{t\in T}G^\pr t,$ that $tg = x_gt_g$
with unique $x_g\in G^\pr,$ $t_g\in T,$ for all $t,g,$ and that
$[G:G^\pr] = \vert T\vert.$  It follows that $\wt G =\dot\bigcup_{\wt t\in \wt
T}\wt G^\pr\, \wt t,$ with $\wt t\,\wt g = \wt x_g\wt t_g$ by $\wt
x_g\in \wt G^\pr,$ $\wt t_g\in \wt T$ etc.  To compute $(\Det\,\
\wt y)(\wt \chi  ), $ apply Step~2 to $N=\wt G,$ $N^\pr = \wt G^\pr,$
$\xi = \wt \chi,$ 
$\xi  ^\pr = \wt \chi  ^\pr,$ and $z=\wt y.$  Write $\theta  $ for
the bijection $W \os\sim\lra \wt G = N,$ and note that $T\os\sim\lra
T_N$ is a bijection.  Moreover,$\sum_{w\in W} z_{\theta(w)}\theta 
(w) = z=\wt y =\sum_{w\in N}\,\aug(\ol{y_w})\theta  (w)$ implies
that $z_{\theta  (w)} =\,\aug(\ol{y_w})$ for all $w\in W;$ also $wt
= x_w t_w$ implies $\wt w\wt t = \wt x_w \wt t_w,$ hence $\kappa 
_{\theta  (w)} = \wt x_w$ and $\tau  _{\theta  (w)} = \wt t_w.$ 
Therefore,
$$
\begin{aligned}
(\Det\,\wt y)(\wt\chi  ) &= \,\det\Big(\sum\nl _{\theta  (w)\in
N:\tau  _{\theta  (w)}= \wt u} \wt\chi  ^\pr(\kappa  _{\theta  (w)})
z_{\theta  (w)}\Big)_{(\wt u;,\wt t)\in T_n\times T_N}\\
&= \,\det\Big(\sum\nl_{w\in W:\wt t_w=\wt u}\wt\chi  ^\pr(\wt
x_w)\,\aug(\ol{y_w})\Big)_{(u,t)\in T\times T}\\
&= \,\det \Big(\sum\nl_{w\in W:t_w=u}\chi 
^\pr(x_w)\,\aug(\ol{y_w})\Big)_{(u,t)\in T\times T}
\end{aligned}
$$
by (3A2) and the above, which agrees with $\big(\varepsilon  \circ
(\Det\, y)\big)(\wt\zeta  )$ by (3A3).
\end{proof}

\vs\noi
{\sc Remark.} The proof of [RW3, Lemma 2] is imprecise; the assertion
of the lemma, however, is correct and follows readily from (3A1).

\vsk\noi
We write $SK_1(\Lambda  _{\fO}G)$ (respectively, $SK_1(\fO [\wt G])$)
for the kernel of the left (respectively, right) vertical map in the
diagram of Proposition 3A.  The ambiguity, with respect to the
reduced norm, that this causes is dispelled by

\vs\noi
{\sc Lemma 3B.} (i) {\it If $G= H\rtimes\Gamma  .$ then
$SK_1(\Lambda  _\fO G)$ equals the kernel of the composite map

$K_1(\Lambda  _\fO G)\to K_1(\SQ_FG) \os{\nr}\lra
\,\cent(\SQ_FG)^\times.$ }

\begin{description}
\item  {(ii)} {\it If $N$ is a finite group, then $SK_1(\fO[N])$
equals the kernel of the composite map 
\newline $K_1(\fO[N]) \to K_1(F[N])
\os{\nr}\lra \,\cent (F[N])^\times.$}
\end{description}

\begin{proof}
{\it For}  (i) The kernels of (i) are the kernels of the bottom
(respectively, top) composites in the commutative diagram
$$
\begin{matrix}
&&\cent\,(\SQ_F G)^\times\q\q\q\q{}\\
& \os{nr\q{}}\nearrow\\
K_1(\Lambda  _\fO G) \lra K_1(\SQ_F G) &&\simeq \downarrow \; f\\
&\us{\Det\q{}}\searrow\\ 
&&\Hom^F\big(R_p(G),(\SQ^c\Gamma  )^\times \big)\big).
\end{matrix}
$$

\vs\noi
For the bottom composite, this follows from the Step 1 diagram of the
proof of Proposition 3A (using $>\!\to);$ for the top it is the
composite map of the assertion (i), which is implied by the triangle.
This follows from [RW2, p.558]
with $f$ the map $z\mapsto f_z$ of the proof of [loc. cit.,
Theorem~7], i.e. $f_z(\chi  )$ is the image of $ze_\chi  $ under the
map cent $\big((\SQ^c(G)e_\chi  \big) \to \SQ^c(\Gamma  )$
induced by $\gamma  _\chi  \mapsto \gamma  ^{w_\chi  }$, and
the isomorphism is by Theorem~8 (last line of the proof).

\noi
{\it For} (ii): These are the kernels of the bottom
(respectively, top) composites in the commutative diagram
$$
\begin{matrix}
&&\big(\cent(F[N])^\times\q\q\q\q{}\\
&\os{nr\q{}}\nearrow\\
K_1(\fO[N]) \lra K_1(F[N]) &&\simeq \downarrow\\
&\us{\Det\q{}}\searrow\\
&&\Hom^F\big(R_p(N),(\Qp{}^c)^\times \big)\big).
\end{matrix}
$$

\vs\noi
For the bottom composite, use the Step~2 diagram of Proposition ~3A;
the triangle is from [CRII,52.10].
\end{proof}

\vs\noi
{\sc Proposition 3C.} {\it If $G=H\rtimes \Gamma    $ with $\Gamma   
\simeq \Z_p$ and $H$ a finite abelian $p$-group, then $SK_1\big(\Lambda
_{\fO}(G)\big) =0.$ }

\begin{proof}
{}[DdSMS] provides a convenient check of a crucial point, starting
from loc.cit., p.139, with $G_k = P_k(G),$ the latter defined by
$P_1(G)= G $ and $P_{k+1}(G) = P_k(G)^p[P_k(G),G]$ (see
Corollary~1.20).  Lemma~7.1 i) gives a map
$(\fO/p^k\fO)[G/G_k]\os\alpha  \lra \Lambda  _{\fO}G/J^k$ for every
$k\ge 0,$ where
$$
J= \,\ker(\Lambda  _\fO G\os\aug\lra \fO) + p\Lambda  _\fO
G\subseteq \,\rad(\Lambda  _\fO G).
$$
Hence composing    $\fO[G/G_k]\to (\fO/p^k\fO)[G/G_k]$ with $\alpha$
  gives $\fO[G/G_k] \to \Lambda  _\fO G/J^k$ from which we deduce the
commutative triangle
$$
\begin{matrix}
&&\os\lim{\us k\lla}\, K_1(\fO[G/G_k])\\
&\nearrow\\
K_1(\Lambda  _\fO G) &&\downarrow\\
&\searrow\\
&&\os\lim{\us k\lla}\, K_1(\Lambda  _\fO G/J^k)\end{matrix}
\q\q \begin{array}{l}
\text{with}\; \searrow
\;\text{an isomorphism by}\\ 
\text{Proposition 1.5.1 of [FK]}.\\ 
\text{It follows that}\;
K_1(\Lambda  _{\fO}G)\lra\\  
 \os\lim{\us k \lla} \, K_1(\fO[G/G_k])\;\text{is
injective.}\end{array}
$$

\vs\noi
Observe here that there exists $t\ge0$ and an open central subgroup
$C_1\simeq \Z_p $ of $G$ with $G_k= C^{p^{k-t}}_1,$ for all $k\ge
t:$ for any open central $C_0\simeq \Z_p$ contains $G_t$ for some
$t,$ hence $G_t = C_1$ for such $C_1.$  Then $G_k = C^{p^{k-t}}_1,$
for $k\ge t,$ follows by induction from $G_k= P_k(G).$  

\vs\noi
Now taking kernels in the square of Proposition 3A, with $C=
C^{p^{k-t}}_1,$ yields
$$
\begin{matrix}
SK_1(\Lambda  _\fO G) &\lra &SK_1(\fO[  G/C^{p^{k-t}}_1])\\ 
\\
\downarrow &&\downarrow &\text{for all}\; k\ge t.\\
\\
K_1(\Lambda  _\fO G) &\lra & K_1(\fO[G/C^{p^{k-t}}_1])\end{matrix}
$$

\vs\noi
We now observe that $SK_1(\fO[G/C^{p^{k-t}}_1]) =0,$ by [Ol,
Corollary 7.2], since $HC^{p^{k-t}}_1/C^{p^{k-t}}_1$ is abelian and
$G/HC^{p^{k-t}}_1$ cyclic.  But the notation $SK_1$ is ambiguous,
the first one using the $\Det$ formalism and the second the $\nr$
formalism.  This ambiguity is resolved by Lemma 3B(ii); moreover,
the $\nr$ formalism is consistent with Oliver's definition of
$SK_1(\fO [N])$ because $SK_1(F[N])=0,$ by Wang's theorem [Ol,
Theorem 2.3].  Thus, the composite map
$$
SK_1(\Lambda  _\fO G) \hookrightarrow K_1(\Lambda  _\fO
G)\dot\hookrightarrow
\lim_{\lla\atop k\ge t}
K_1(\fO [G/C^{p^{k-t}}_1])
$$
is both $0$ and injective, with $\dot\hookrightarrow$ by the first
paragraph.   So $SK_1(\Lambda  _\fO G)=0.$
\end{proof}

We now return to $\Delta = \Delta_{r,s,\IF}$ of \S1, so assume that
$F/\Qbar_p$ is unramified.

\vs\noi
{\sc Lemma 
3D.} {\it Let $G= C\rtimes \Gamma  $ be the semidirect
product of the cyclic group $C=\la c\ra$ of order $p^{r+s}$ with
$\Gamma  =\ol{\la u\ra},$ where $u$ acts on $C$ by $ucu^{-1}=c^{1+p^r};$
set $\wt G= G/\la c^{p^{r+s-1}}\ra$ and $\wt c$ the image of $c$
under $G\to \wt G.$
$$
\begin{matrix}
\Lambda  _\fO G &\os\kappa  \lra &\Delta\\
\iota \downarrow &&\downarrow\\
\Lambda  _\fO \wt G &\lra &\Delta/\Delta \pi_1\end{matrix}
$$

\vs\noi
Then the above diagram is a pull-back square in which $\Delta =
\bigoplus_{0\le j< p^s
}(\fO[\zeta  _{r+s}])[[T]]u^j,$ with $u^{p^s} = 1+T$ and $u\zeta 
_{r+s} = \zeta  ^{1+p^r}_{r+s}u,$
is the maximal
order of \S1 with $\kappa $ and $\iota$  induced by $c\mapsto \zeta 
_{r+s}$ and $c \mapsto \wt c,$ respectively. }

\vs\noi{\it Proof.}
We first need to show that 
$$
\Lambda_\fO G \simeq \bigoplus_{0\le j<p^s}
\big( (\fO[[T]])[C]\big)u^j,\;  {\rm with}\; u^{p^s}=1+T
\; {\rm and} \; uc = c^{1+p^r}u \leqno{\rm (3D1)}
$$
For $ucu^{-1}=c^{1+p^r}$ implies $u^{p^s}cu^{-p^s} =
c^{(1+p^r)^{p^s}}
\dot = c^{1+p^{r+s}}=c,$ with $\dot =$ by (2G5), hence setting
$\Gamma  _0 =\ol{\la u^{p^s}\ra}$ shows that $G_0 \os{{\rm def}}=
C\Gamma  _0$ is abelian, and that $\Lambda_\fO G_0$ is the group
ring $(\Lambda_\fO\Gamma  _0)[C].$  Now the Iwasawa isomorphism gives
$\Lambda_\fO\Gamma  _0
\simeq \fO\big[[T]\big],$ with $T= u^{p^s}-1,$ and combining then proves
that $\Lambda_\fO G = \us j\bigoplus (\Lambda_\fO G_0)u^j\simeq
\us j\bigoplus \big((\fO[[T]])[C]\big)u^j,$ as claimed.

\noi
Set $n=r+s,$ denote the $p^n\text{-th}$ cyclotomic polynomial by
$\Phi(x) =\sum_{i=0}^{p-1} x^{p^{n-1}i}$ and define
$f(x)$ by 
$f(x)
= \sum^{p-2}_{j=0} (p-1-j)x^{p^{n-1}j},$ so that $\Phi(x) = p +
(x^{p^{n-1}}-1)f(x).$  Substituting $x=c$ into the above, we arrive
at
$$
(c^{p^{n-1}} -1) \Phi(c) =0 \q\text{and} \q \Phi(c) = p +
(c^{p^{n-1}}-1)f(c), \leqno{\rm (3D2)}
$$
and noting that $c^{p^{n-1}}$ is central in $\Lambda  _\fO G$
implies that $\Phi(c)$ and $f(c)$ also are.

\vs\noi
Forming the ideals $I= (c^{p^{n-1}}-1)\Lambda  _\fO G,$
$J = \Phi(c)\Lambda  _\fO G,$ we obtain the pull-back square below,
as in [CRII, 42.3],

$$
\begin{matrix}
\Lambda  _\fO G/I\cap 
J & \lra &\Lambda  _\fO G/J &\q\text{in which} \q I\cap J = 0,
\;\text{because}\q y\in I\cap J\\
\downarrow &&\downarrow &\q\text{implies}\q \Phi(c)y=0\;\text{and}\;
(c^{p^{n-1}}-1)y=0,\\
\Lambda  _\fO G/I &\lra &\Lambda  _\fO G/I + J &\q\text{hence}\;
py=0, \;\text
{by (3D2), and so} \; y=0.\q{}\end{matrix} \leqno{\rm (3D3)}
$$

\vs\noi
Next, treating  (3D1) as an identification, implies that
$$
\Lambda  _\fO G/J = \bigoplus\nl_{0\le
j<p^s}\big((\fO[[T]])[\zeta_n]\big)
 u^j = \us j\bigoplus \fO_n[[T]]u^j =\us j\bigoplus Bu^j =
\Delta,
$$
hence
$\Lambda  _\fO G/I + J = \Delta/\Delta \cdot I = \Delta/\Delta(\zeta 
^{p^{n-1}}_n -1) = \Delta/\Delta(\zeta  _1-1) = \Delta/\Delta
\pi_1.$    

\noi
Now, with (3D1) still an identification implies that

\centerline{$
\Lambda  _\fO G/I = \bigoplus\nl_ {0\le j<p^s}
\big(\fO[[T]]\big)
[\wt C]u^j \q\text{with} \q \wt C = C/\la
c^{p^{n-1}}) = \la \wt c\ra.
$}

\noi
Applying the argument of (3D1) to $\Lambda_\fO\wt G$ with the same
$\Gamma  _0$ and $C$ replaced by $\wt C,$ then 
$\wt G_0 \os{\rm
def}= \wt C\Gamma  _0$ is still abelian, with $\Lambda_\fO\wt G_0$
still
the group ring
$(\Lambda _\fO\Gamma  _0)[\wt C]=
\big(\fO[[T]]\big)[\wt C],$ hence 

\centerline{$
\Lambda_\fO\wt G =
\us{0\le j<p^s}\bigoplus (\Lambda_\fO\wt G_0)u^j= \us
j\bigoplus (\fO[[T]])[\wt C]u^j.
$}

\noi
Combining
shows that $\Lambda_\fO G/I= \Lambda_\fO \wt G.$

\noi
The Lemma follows from these ingredients, since the maps of (3D3)
are all surjective.
\hfill $ \square $

\vs\noi
{\sc Theorem 3E.} (i)  {\it There is an exact sequence
$$
0\lra SK_1(\Delta) \lra K_1(\Delta/\pi_1\Delta) \os{\delta  _1}\lra
\,\nr\big(K_1(\Delta)\big)/\nr(1+\pi_1\Delta) \lra 0
$$
with $\delta_1$ defined by } (3E4).

\noi{(ii) {\it If $d\in 1+\Pi^{p^{r+s-1}}\Delta  $ has $nr(d)=1$ then
$d\in [\Delta  ^\times,\Delta  ^\times].$}

\noi{(iii) {\it Every $n<p^{r+s-1}$ has weight} $w\le s-1.$

\begin{proof}
{\it For} (i): Turning back to the pull-back diagram of Lemma 3D, with all maps
surjective, we pass to the exact Mayer-Vietoris sequence [CRII,
p.106].
$$
K_1(\Lambda  _\fO G) \to K_1(\Delta) \times K_1(\Lambda  _\fO \wt G)
\to K_1(\Delta/\pi_1\Delta) \os \partial\lra K_0(\Lambda  _\fO G)
\to K_0 (\Delta)\times K_0(\Lambda  _\fO \wt G).
\leqno{\rm (3E1)}
$$
Note that $K_0(\Lambda  _\fO G)  $ and $ K_0(\Lambda  _\fO \wt G)$
are both isomorphic to $\Z$  because 
 projectives over local rings are free (see [Lam,19.29]).  It follows that
$K_0(\Lambda  _\fO G) \lra K_0 (\Lambda  _\fO \wt G)$ is injective,
since it sends the generator $[\Lambda  _\fO G]$ to $[\Lambda  _\fO
\wt G\bigotimes _{\Lambda  _\fO G}\Lambda  _\fO G] = [\Lambda  _{\fO}
\wt G],$  a generator.  So $\partial =0$ and the remaining
exact sequence becomes the top row of the commutative diagram with
exact rows
$$
\begin{matrix} 
&K_1(\Lambda  _\fO G) &\lra& K_1(\Delta)\times K_1(\Lambda  _\fO \wt
G) &\lra &K_1(\Delta/\pi_1\Delta) &\lra 0\\ 
\\
&\nr\downarrow &&\nr\times \nr \downarrow &&\downarrow\\
\\
0 \lra &\cent(\SQ_FG)^\times &\os=\lra &\cent(D)^\times
\times \,\cent(\SQ_F\wt G)^\times &\lra &0\\
\end{matrix}
\leqno{\rm (3E2)}
$$

\noi
and bottom row coming from the reduced norm being maps
on Wedderburn components to their centre fields.  Applying the
snake lemma,  cut short in anticipation of $SK_1(\Lambda  _\fO G)$
and $SK_1(\Lambda  _\fO\wt G)$ vanishing, by Lemma 3B(i) and
Proposition 3C, gives the short exact sequence
$$
0\lra SK_1(\Delta) \lra K_1(\Delta/\pi_1\Delta) \os\delta\lra X\lra
0, \;\text{with}\; X=
$$
$$
\ker\Big(\frac{\cent(D)^\times \times \cent(\SQ_F\wt
G)^\times}{\nr\big(K_1(\Lambda  _\fO G)\big)} \lra
\frac{\cent(D)^\times \times \cent(\SQ_F\wt
G)^\times}{\nr\big(K_1(\Delta)\big) \times \nr\big(K_1(\Lambda  _\fO
\wt G)\big)} \Big) = \frac{\nr\big(K_1(\Delta)\big)\times
\nr\big(K_1(\Lambda  _\fO \wt G)\big)}{\nr\big(K_1(\Lambda  _\fO
G)\big)} \,,
$$

\vs\noi
from which we find the natural map $\tau  : \nr\big(K_1(\Delta)\big)
\lra X.$  So we compute $\delta  $ by lifting $x\in
K_1(\Delta/\pi_1\Delta)$ to $\hat x\in K_1(\Delta)$ and setting
$\delta  (x) = \tau \big(\nr(\hat x)\big) \in X.$  This is
summarized in the right square of the bottom two rows of the
commutative diagram with exact rows
$$
\begin{matrix}
&&&&\im(\alpha  ) &\os\nr\lra &\nr \big(\im(\alpha  )\big)\\
\\
&&&&\downarrow &&\downarrow\\
\\
0 &\lra &SK_1(\Delta) &\lra &K_1(\Delta) &\os\nr\lra
&\nr\big(K_1(\Delta)\big) &\lra 0\\
\\
&&=\downarrow &&\downarrow &&\downarrow \tau  \\
\\
0 &\lra &SK_1(\Delta) &\lra &K_1(\Delta/\pi_1\Delta) &\os\delta 
\lra &X &\lra 0
\end{matrix}
\leqno{\rm (3E3)}
$$

\vs\noi
which then induces the identity map on $SK_1(\Delta):$ the map
from the bottom $SK_1(\Delta)$ comes from the middle kernel map in (3E2),
so is the inclusion $SK_1(\Delta)\hookrightarrow K_1(\Delta)$
followed by the natural map $K_1(\Delta)\lra
K_1(\Delta/\pi_1\Delta),$ in agreement with the other one.  Applying
the snake lemma to the bottom rows adds only the surjectivity of
$\tau  $ and the isomorphism on kernels induced by
$$
K_1(\Delta, \pi_1\Delta) \os\alpha  \lra K_1(\Delta) \lra
K_1(\Delta/\pi_1\Delta)\lra 0
$$
(see [Ol], Theorem 1.13).  Comparing this to the exact sequence
$$
0\lra V\lra 1 + \pi_1\Delta \lra K_1(\Delta,\pi_1\Delta)\lra 0
$$
from ([Ol], Theorem 1.15), with $V = \la (1+xy)(1+yx)^{-1}:x\in
\pi_1\Delta,y\in \Delta\ra,$ we get the exact sequence
$$
0\lra V\lra 1+\pi_1\Delta \lra \,\im(\alpha  )\lra 0.
$$
The point is that $V\subseteq [\Delta^\times , \Delta^\times]$
[locis cit.] (or proof of Theorem 3.6(b) in [Va]), which implies
$\nr(V) =1,$ hence $\nr(1+\pi_1\Delta) = \,\nr\big(\im(\alpha 
)\big). $  Replacing $\nr\big(\im(\alpha  )\big)$ in (3E3) by
$\nr(1+\pi_1\Delta)$ now allows us to replace the original right
square by
$$
\begin{matrix}
K_1(\Delta) &\os\nr\lra &\nr\big(K_1(\Delta)\big)\\
\\
\downarrow &&\downarrow\\
\\
K_1(\Delta/\pi_1\Delta) &\os{\delta  _1}\lra
&\frac{\nr\big(K_1(\Delta)\big)}{\nr(1+\pi_1\Delta)}
\end{matrix}
\leqno{\rm (3E4)}
$$
at the cost of replacing $\delta  $ by $\delta  _1,$ which is now
{\it defined} by this square.

\noi
{\it For} (ii): $d-1 \in \Pi^{p^{r+s-1}}\Delta   = \pi_1\Delta  $
implies that $[d] \in \ker \big (SK_1(\Delta  ) \to K_1(\Delta 
/\pi_1\Delta  )\big),$ by (i), hence $d\in \ker\big(\Delta  ^\times \to
K_1(\Delta  )\big) = [\Delta  ^\times,\Delta  ^\times],$ by
Lemma~2A(i).

\noi
{\it For} (iii): This follows from $n<p^{r+s-1} + p^{s-1}$ and the
definition of $w.$
\end{proof}

\section{Reduced Norm Congruences on $\Delta$}\label{Sec4}

 This section looks for congruences $\nr(d)\equiv 1\!\!\!\mod \pi
A,$  which result from the
$\Pi$-filtration on $\Delta$ induced by the powers of its ideal
$\Pi\Delta  .$ Recall that the regular norm $\Nr: D\to \fq(A)$ has
$\Nr(d) = \nr(d)^{p^s},$ for all $d\in D$ ([Re, (9.7)], [Sc,
p.~296]).
The first step is

\vs\noi
{\bf \sc Lemma 4A.} {\it If $d\in \Delta  ^\times$ has $\nr (d) =1,$
then $d\equiv 1\!\!\! \mod \Pi\Delta  .$}

\begin {proof}
Recall that $\Nr(d)$ is the determinant of the  action of $d\in D$
on the left $\fq (A)$-module $D,$ or, equally well, for 
$d\in\Delta^\times,$
of the action of $d$ on the free left $A$-module $\Delta  ,$ 
with basis
$\{\lambda  _k: 0\le k<p^{2s}\}$ from the proof of Lemma~2A(ii).

\vs\noi
The $A$-module structure of $\Delta$ implies the $\ol B$-module
structure of $\ol\Delta,$ in the notation of (2A1), via
$\ol B=A/\pi A,$ from the proof of Lemma~2A(i).  Since $\zeta \equiv
1\!\!\!\mod \Pi B,$ $\ol \Delta$ must have $\ol B$-basis $\{\ol
u^{k_1}:0\le k_1 <p^s\};$ since 
$\fq (\,\ol\Delta \, ) = \IF\big((\,\ol{u-1}\,)\big)$ is a purely
inseparable extension of $\fq (\,\ol B\,)=\IF \big((\,\ol T\,)\big),$
because
$(\,\ol{u-1}\,)^{p^s} = \ol T$ 
and with power basis $\{(\ol{u-1})^{k_1}:0\le k_1<p^s\}$ 
 by the binomial expansion change of
basis.  Thus $\ol{Nr(d)} =\ol{nr(d)}^{p^s}$ and $nr(d)=1$ implies
$\ol 1 = \ol{nr(d)}
=  N_{\fq(\,\ol\Delta\,  )/\fq(\,\ol B\,)} (\,\ol d\,) = \ol d\,^{p^s},$ hence
$(\,\ol
d-\ol 1\,)^{p^s} = \ol 0$ and $\ol d= \ol 1.$
\end{proof}

\vs\noi
We next define the level $n$ of $d\in 1 +\Pi \Delta  $ by
$$
n=\max\{k\ge 1:d\in 1+\Pi^k\Delta  \}.
$$
The motivation is to use the filtration induced by Theorem 3E on
$SK_1(\Delta  )$ to get information on the level of solutions $d$ of
$\nr(d) = 1.$  This filtration is finite because $\pi_1\Delta  =
\Pi^{p^{r+s-1}}\Delta  .$

\vs\noi
{\sc Lemma 4B.} (i) {\it If $d\in D$ has $d=\us{k\in\Z/p^s\Z}\sum
b_k u^{\hat k},$ with all $b_k\in \fq(B),$ then our chosen (reduced) norm
is $nr(d) = \det
([^{g^j}b_{i-j}(1+T)^{h_j(i)}]_{(\hat i,\hat j)}),$ with $h_j(i) =
\begin{cases} 1, &\hat j>\hat i\\ 0, &\hat j \le \hat i.\end{cases}$
}

\noi
(ii)  {\it Suppose that $d=1+ \Pi^n\sum_{k\in\Z/p^s\Z}
b_ku^{\hat k}$ with $n\ge 1$ and all $b_k\in B.$  Here  $0\le \hat k < p^s$
and
\newline $\hat k\!\!\!\mod p^s\Z$ equals $ k.$  Let $\fG_s$ denote the symmetric
group on $\Z/p^s\Z.$  Then
$$
\nr(d) =\sum_{\sigma  \in \fG_s} \, \sgn(\sigma  )\bd(\sigma  ) \q
\text{with}
$$
$$
\bd(\sigma  ) = (1+T)^{f_+(\sigma  )}\prod\nolimits_{i\in\CF(\sigma  )} \,
^{g^i}(1+\Pi^nb_0)\cdot \prod\nolimits_{i\notin \CF(\sigma  )} \,^{g^{\sigma 
(i)}}(\Pi^n b_{i-\sigma  (i)}),
\leqno{\rm (4B1)}
$$
where $\CF(\sigma  ) = \{i\in \Z/p^s\Z:\sigma  (i)=i\}$ 
and $f_+(\sigma  ) =\us{i\in\Z/p^s\Z}\sum \,h_{\sigma  (i)}(i).$}

\vs\noi
{\it Proof.}
For (i): This follows from (2A2), in the proof of Lemma~2A(ii)1.

\noi
For (ii):  The notation here differs from (i) when it is rewritten
as 
$$
d= (1+\Pi^n b_0)u^{\hat 0} +\us{k\ne 0}\sum (\Pi^n b_k)u^{\hat
k},
$$
which parametrizes a proper subset of $\Delta  .$

Following the procedure leading to (2A2) here 
 leads to the summary 
$$
u^{\hat j}d = \;^{g^j}(1+\Pi^n b_0)u^{\hat j} + \sum_{i\ne j} \,
^{g^j}\big(\Pi^n b_{i-j}\big) (1+T)^{h_j(i)} u^{\hat i}
\q{\rm with } \q
h_j(i) =\begin{cases} 1, &\hat j >\hat i\\
0, &\hat j \le \hat i.\end{cases}
$$

\noi
Setting
$$
\bd_{i,j} =
\begin{cases}{}^{g^j}(\Pi^nb_{i-j})(1+T)^{h_j(i)},\
 &i\ne j\\
^{g^j}(1+\Pi^nb_0), &i=j\end{cases}
\leqno{\rm (4B2)}
$$
it follows that
$$
\nr(d) = \,\det [\bd_{i,j}]_{(\hat i,\hat j)} =\sum_{\sigma  \in
\fG_s}\,\sgn(\sigma  )\bd(\sigma  ) \q\text{with}\q \bd(\sigma  )
=\prod_{i\in\Z/p^s\Z} \bd_{i,\sigma  (i)}\,.
$$
Combining, we have
$$
\bd(\sigma ) =  \prod_{i\in\CF(\sigma )} \,^{g^{\sigma  (i)}}(1+\Pi^nb_0)
\prod_{i\notin \CF(\sigma  )}\,^{g^{\sigma  (i)}}(\Pi^nb_{i-\sigma 
(i)})(1+T)^{h_{\sigma  (i)}(i)}
\leqno{\rm (4B3)}
$$
and the total power of $(1+T)$ is 
$\us{i\notin \CF(\sigma  )}\sum h_{\sigma  (i)}(i)
=\us{i\in\Z/p^s\Z}\sum h_{\sigma  (i)}
(i)= f_+(\sigma  )$
since $h_{\sigma 
(i)}(i) = h_i(i)=0,$ for 
$i\in \CF(\sigma  ),$
 to arrive at (4B1).  For ease of reference we end with 
$$
h_j(i) = \begin{cases} 1, &\hat j >\hat i\\
0, &\hat j\le \hat i \end{cases}, \; {\rm and} \q
f_+(\sigma  ) = \sum\nolimits_{i\in\Z/p^s\Z}h_{\sigma 
(i)}(i).\q\q\q\q\q\q \square
$$

\vsk\noi
The formula for $\nr(d)$ in Lemma 4B takes values in $A$ for Galois
theoretic reasons.  Another explanation is in

\vs\noi
{\sc Lemma 4C.} {\it Let the special generator $g$ of $\fg$ act on
$\fG_s$ by
$^{g^j}\sigma  =\omega^j\sigma  \omega^{-j},$
with $\omega\in \fG_s$ defined
by $\omega (i) = i+1$ for all $i\in \Z/p^s\Z.$  Then, for each $d\in
1+\Pi\Delta,$}

\begin{description}
\item  {(i)} {\it $\sigma  \mapsto \bd(\sigma  )$ induces a 
function $\bd:\fG_s\to B$ which is a $\fg$-set map, and}
\item  {(ii)} $\bd({}^g\sigma  ) ={}^g\bd(\sigma  )$ and
$f_+(\,^g\sigma  ) = f_+(\sigma  ) $ {\it for all
$\sigma  \in \fG_s\,.$}
\end{description}

\vsk\noi
{\it Proof.} Starting from (4B2), we have
$$
\begin{aligned}
^g\bd_{i,j} &= \bd_{i+1,j+1}(1+T)^{\eta  (i,j)} \,, \; \text{with}\\
\eta  (i,j) &= h_j(i) - h_{j+1}(i+1), \, \text{for all} \; i,j\in
\Z/p^s\Z, \, \text{hence}
\end{aligned}
$$
$$
\bd(\,^g\sigma  ) =\,^g\bd(\sigma  ) (1+T)^{-\sum_{i\in\Z/p^s\Z}\eta  (i,\sigma 
(i))} \,, \; \text{for all}\; \sigma  ,
\leqno{\rm (4C1)}
$$
since the left side equals
$
\prod\nolimits_{i\in \Z/p^s\Z} \bd_{i,\,^g \sigma  (i)}  =
\prod\nolimits_i \bd_{i+1,
\omega\sigma  \omega^{-1}(i+1)} =
 \prod\nolimits_i \bd_{i+1,\sigma  (i)+1}
=$

\noi $\prod\nolimits_i (1+T) ^{-\eta  (i,\sigma  (i))}
\,^g\bd_{i,\sigma(i)}
=\,^g\bd(\sigma  )(1+T) ^{-\sum_i\eta  (i,\sigma  (i))}. 
$

\noi
We next need
$$
\text{1.} \q \eta  (i,j) = \begin{cases}{}\q 1, &\hat i<\hat j=p^s-1\\
{}-1, &\hat j <\hat i = p^s-1\\
{}\q 0, &\text{\rm otherwise}, \end{cases}
\q \text{\rm and \q 2.}\q \sum\nolimits_i \eta  \big(i,\sigma  (i)\big)
=0, \;\text{\rm for all} \q \sigma  \in \fG_s\,.
$$

\vs\noi
{\it For} 1: First, $\wh{i+1} - \hat i = 1-p^s\kappa(i)$ with
$\kappa(i)= \begin{cases} 1, &\hat i = p^s-1\\
0, &{\rm else}, \end{cases}$ 
because $0\le \wh{i+1} = p^s-1,$ $1\le \hat i +1\le p^s$ and 
$\wh {i+1} \equiv \hat i+1\!\!\!\mod p^s$ unless $\wh{i+1}=0$ or
$\hat i+1= p^s,$ i.e. $\hat i =p^s-1.$ It follows that  
$\wh{j+1} >\wh{i+1} \Longleftrightarrow \hat j >\hat i +
p^s\big(\kappa(j) - \kappa(i)\big),$ hence that $\eta  (i,j) =0$
(i.e.$ h_j(i) = h_{j+1}(i+1))$ unless $\kappa (j) \ne \kappa(i).$ 
This happens only in the two cases:

\vs
\text{3.} $i<j =p^s-1,$ and  \q\q\q\ 4. $j<i=p^s-1.$

\vs\noi
In case 3. $\eta  (i,j) = h_j(i) - h_{j+1}(i+1) \dot =\,
1-h_0(i+1)\ddot =\, 1-0=1,$ with $\dot=$ by $\hat j>\hat i$ and 
$\ddot =$
by $0=\hat 0 \le \wh{i+1}.$

\vs\noi
In case 4. $\eta  (i,j) \dot = \,0-h_{j+1}(i+1) = -h_{j+1}(0)\ddot=
\,-1,$ with $\dot =$ by $\hat j <\hat i$ and $\ddot =$ by $\wh
j<p^s-1$ implying $\kappa  (j)=0,$ hence $\wh {j + 1} = \hat j + 1>0
=\hat 0,$ and so $h_{j+1}(0) = 1.$

\vs\noi
{\it For} 2: The $\big(i,\sigma  (i)\big)$ with $\eta  \big(i,\sigma 
(i)\big)\ne 0$ occur in pairs $\big(\sigma  ^{-1}(p^s-1),
p^s-1\big),$ $\big(p^s-1,\sigma  (p^s-1)\big)$ in different cases
3.,4. above, hence $\sigma  (p^s-1)\ne p^s-1,$ and
$$
\eta  \big(\sigma  ^{-1}(p^s-1),p^s-1\big) + \eta  \big(p^s-1,\sigma 
(p^s-1)\big )= 1 + -1=0, \,\text{proving 2.}
$$
Now, $\bd(\,^g\sigma  ) =\,^g\bd(\sigma  )$ follows from (4C1) and 2.;
and then $\bd(\,^{g^j}\sigma  ) =\, ^{g^j}\bd(\sigma  ),$ for all
$j,$ follows by induction on $j\ge 1.$  This proves (i).

\vs\noi
Finally (ii) follows from 2, (4C1) and (4B4), by
$$
\begin{aligned}
\sum\nolimits_i\eta  \big(i,\sigma  (i)\big)
&=\sum\nolimits _i h_{\sigma  (i)}(i) -\sum\nolimits_i h_{\sigma 
(i)+1}(i+1)\\
&= \sum\nolimits _i h_{\sigma  (i)}(i)-\sum\nolimits_i h_{\sigma 
(i-1)+1}(i-1+1)
= \sum\nolimits _i h_{\sigma  (i)}(i) -\sum\nolimits _i
h_{(\,^g\sigma  )(i)}(i)\\
&= f_+(\sigma  ) - f_+(\,^g\sigma ).\q\q\q\q\q\q\q\q\q\q\q\q\q\q\q\q\q\q\q\q\q
\square
\end{aligned}
$$

\vs\noi
Fix a set $\fR$ of representatives for the $\fg$-orbits in $\fG_s,$
so Lemma 4C(i) and (ii) allow us to write
$$
\nr(d) =\sum_{\sigma  \in \fR}\,\sgn(\sigma  )\Tr_{\fg/_{\St(\sigma }
)}\big(\bd(\sigma  )\big),
 $$
where $\St(\sigma  ) = \St_\fg(\sigma  )$ is the stabilizer in $\fg$
of its action on $\sigma  ,$ and where $\Tr$ is both the permutation
trace and the Galois trace, hence taking values in $A.$  So we have
$$
\nr (d) =\sum\nolimits ^s_{k=0}\,\sum\nolimits_{\sigma  \in\fR_k}
\sgn(\sigma  )\Tr_{\fg/\fg^{p^k}}\big(\bd(\sigma  )\big) \leqno{\rm
(*)}
$$
with $\fR_k = \{\sigma  \in \fR:\St_\fg(\sigma  )= \fg^{p^k}\},$
since $\fg$ is cyclic of order $p^s.$

\noi
We write $C_{\fG_s}(\tau  )$ for the centralizer of $\tau  \in
\fG_s$ in $\fG_s\,.$

\vs\noi
{\sc Lemma 4D.} {\it If $1\ne \sigma  \in \fR_k,$ with $0\le k\le
s,$ and $d=1 + \Pi^n \sum\nolimits_{j\in\Z/p^s\Z} b_j u^{\hat j}$ with
$n\ge 1,$ then }

\begin{description}
\item  {\rm (i)} $\sigma  \in C_{\fG_s}(\omega^{p^k}),$  {\it with
$\omega$ as in
Lemma} 4C,
{\it with} $\sigma  (i+p^ky) = \sigma  (i) + p^ky$ {\it for all}
$i,y$ in $\Z/p^s\Z,$
{\it and} $\Tr_{\fg/\fg^{p^k}}\big(\bd(\sigma  )\big)\in
\pi^{\Phi_k(n)}A.$

\item  {\rm (ii)} $\Tr_{\fg/\fg^{p^k}}\big(\bd(\sigma  )\big)$  {\it is in
$\pi^{\Phi(n)+1} A$ unless $\sigma  $ is in the special case when
$\Phi_k(n) = \Phi(n)$ and $\vert  \supp(\sigma  )\vert  = p^{s-k},$
where $\supp (\sigma  ) =\{i\in\Z/p^s\Z:\sigma  (i)\ne i)\}.$}

\item{\rm (iii)} $\CF^\pr_k(\sigma  )\subseteq \sigma 
_k\big(\CF^\pr_k(\sigma  )\big),$ 
with $\CF^\pr_k(\sigma  ) = (\Z/p^s\Z) \bs \CF_k(\sigma  ).$
\end{description}

\begin{proof} Since $St_{\fg}(\sigma  ) = \fg^{p^k},$ hence
$\sigma  =\,^{g^{p^k}}\!\sigma  = \omega^{p^k} \sigma 
\omega^{-p^{k}},$ we have $\sigma  \in C_{\fG_s}(\omega^{p^k}),$ so

\noi
$ \sigma  (j+ p^ky) = (\sigma  \omega^{p^{k}y})(j) = (\omega^{p^{k}y}\sigma 
)(j) =
\sigma  (j) + p^k y,
$
for all $j,y\in \Z/p^s\Z.$
This implies the existence of a unique 
$\sigma  _k\in \fG_s$
 making the  square below commute, 
$$
\begin{matrix}
\Z/p^s\Z  &\os \sigma  \longrightarrow &\Z/p^s\Z\\
\\
-\downarrow &&-\downarrow\\
\\
\Z/p^k\Z &\os {\sigma  _k}\longrightarrow &\Z/p^k\Z
\end{matrix}
$$
 which is
well defined as
$\ol{\sigma  (j+p^ky)} = \ol{\sigma  (j) +p^ky} = \ol{\sigma  (j)}$
and surjective as $\sigma  $ is.  Thus if $i\in\Z/p^k\Z$ has
$\wt i\in\Z/p^s\Z$ with $\ol{\wt i}=i,$ then $\sigma  _k(i) =
\ol{\sigma  \big(\,\tilde i\,\big)}.$

\vs\noi
It also follows that the
 $\CF^\pr_k(\sigma  )$ 
 are disjoint unions of the cosets of $p^k\Z/p^s\Z$ in
$\Z/p^s\Z:$ for $j\in \CF(\sigma  )$ implies that $j+p^k\Z/p^s\Z
\subseteq \CF(\sigma  ),$ because $\sigma  (j+p^ky) - (j+p^ky) =
\sigma  (j)-j$ (and similarly for $\CF^\pr (\sigma  )).$

\vs\noi
Writing $\CF_k(\sigma  ),$ $\CF^\pr_k(\sigma  )$ for the
images of $\CF(\sigma  ),$ $\CF^\pr(\sigma  )$ under the natural map
$\Z/p^s\Z \os -\to \Z/p^k\Z$ we can organize the products in (4B1)
according to these special cosets 
$\wt e  + p^k\Z/p^s\Z,$ with $\wt e\in
\Z/p^s\Z$ having image $e$ in $\CF_k(\sigma  ),$ respectively
$\CF^\pr_k(\sigma  ),$  to
get
$$
\bd(\sigma  ) = (1+T)^{f_+(\sigma  )} \prod\nolimits_{e \in \CF_k(\sigma
 )}\,^{g^e_k}N_{\fg^{p^k}}(1+\Pi^nb_0)
\prod\nolimits_{e\in\CF^\pr_k(\sigma  )} \, ^{g^{\sigma  _k(e)}_k}
N_{\fg^{p^k}}(\Pi^nb\,_{\wt e-\sigma  (\,\tilde e\,)}) \leqno{\rm (4D1)}
$$  
with $g_k = g\fg^{p^k} \in \fg/\fg^{p^k},$ since the norms take
values in $B^\pr_k,$ the subring of $B$ fixed by $\fg^{p^k}.$
Note that $\wt e + p^k\Z/p^s\Z$ is independent of the choice of $\wt
e\in \Z/p^s\Z$ with $\ol{\tilde e} = e.$

\vs\noi
More precisely,\, if\, $e\in\Z/p^k\Z$\, has \,$\wt e\in\Z/p^s\Z$\, with
\,$\ol{\wt e} = e,$\, and \,$Z=p^k\Z/p^s\Z,$\, then
$e\in \CF_k(\sigma  )$ defines
 $^{g^{\tilde e}}(^{g^Z}(1+\Pi^n
b_0)) =\, ^{{g^e_k}} N_{\fg^{p^k}}
 (1+\Pi^n b_0);$ and $e\in
\CF^\pr_k(\sigma  )$ defines
$
{^{g^{\sigma  (\tilde e)}}}(^{g^Z}(\Pi^n b\,_{\tilde
e-\sigma  (\tilde e)}))
\dot = $
$^{{g^{\sigma  (\tilde e)}_k}} N_{\fg^{p^k}}
(\Pi^n b\,_{\tilde e - \sigma 
(\tilde e)})
 \,\ddot =\,  ^{\fg^{\sigma  _k(e)}_k} N_{\fg^{p^k}}(\Pi^n b\,_{\tilde
e\,-\,\sigma  (\tilde e)}),
$
\, with \,\,
$\dot = $ \, by \,$\sigma  (\tilde e +Z)\, =\, \sigma  (\tilde e)+Z$ and $\,\ddot
=\,$ by
\newline  $\ol{\sigma  (\tilde e)} =  \sigma  _k(e).$

\vs\noi
Note also that the $f_+(\sigma  )$ and $\CF_k(\sigma  )$ factors of
$\bd(\sigma  )$ in (4D1) have $\v_\Pi =0,$ so only the
$\CF^\pr_k(\sigma  )$ factor, which has $\v_\Pi >0,$ matters
below.  In particular, $\CF_k(\sigma  )\cap \CF^\pr_k(\sigma  ) =
\emptyset.$

\vs\noi
Now $N_{\fg^{p^k}}(\Pi) = \pi_{(r+s)-(s-k)},$ by the proof of  (2A3), implies that
$\bd(\sigma )\in \prod_{i\in\CF^\pr_k(\sigma  )}\,^{g^{\sigma 
_k(i)}_k} (\pi^n_{r+k}B^\pr_k) = \pi^{n\vert \CF^\pr_k(\sigma 
)\vert}_{r+k} B^\pr_k,$ hence $\Tr_{\fg/\fg^{p^k}}\big(\bd(\sigma 
)\big)\in \,\Tr _{\fg/\fg^{p^k}}(\pi^{n\vert  \CF^\pr_k(\sigma 
)\vert  }_{r+k} B^\pr_k) =
\pi^{\Phi_k(n\vert\CF^\pr_k(\sigma  )\vert)}A,$ by Lemma
1B(ii). Here $\vert  \CF^\pr_k(\sigma  )\vert  \ge 1,$ by $\sigma  \ne 1,$
implies that $\Phi_k\big(n\vert  \CF^\pr_k(\sigma  )\vert 
\big)\ge\Phi_k(n),$ which implies (i).

\vs\noi
Next, if $\vert  \CF^\pr_k(\sigma  )\vert  \ge 2$  then
$\Phi_k\big(n\vert  \CF^\pr_k(\sigma  )\vert  \big) \ge \Phi_k(2n)
\ge \Phi(n)+1,$ by Lemma 1E, hence 
$\Tr_{\fg/\fg^{p^k}}\big(\bd(\sigma  )\big)\in \pi^{\Phi(n)+1}A,$ unless
we are in the special case when $\Phi_k(n) = \Phi(n)$ and $\vert 
\CF^\pr_k(\sigma  )\vert  =1$ both hold.  This gives (ii) since
$\supp(\sigma  )=\CF^\pr(\sigma  ) $ and $\vert  \CF^\pr(\sigma 
)\vert  =p^{s-k}\vert  \CF^\pr_k(\sigma  )\vert = p^{s-k} .$

\vs\noi
Finally, for (iii), $\sigma  _k \big(\CF_k(\sigma  )\big) =
\CF_k(\sigma  )$ since $j\in \CF_k(\sigma  )$ implies $\wt j\in
\CF(\sigma  )$ hence   $\sigma  _k(j) = \ol{\wt j} = j.$  Thus
$\Z/p^k\Z = \CF_k(\sigma  )\dot\cup \CF^\pr_k(\sigma  )$
 implies $\Z/p^k\Z  =
\CF_k(\sigma  ) \cup \sigma  _k \big(\CF^\pr(\sigma  )\big),$
hence the assertion. 
\end{proof}

\vs\noi
{\sc Lemma 4E.} {\it If $n\ge 1$ and $b\in B,$ with $b\not\equiv 0\!\!\!\mod
\Pi B,$
then
$$
N_{\fg}(1+\Pi^nb)\equiv 1 + \sum^s_{m=0} \,\Tr_{\fg/\fg^{p^m}}
\big(N_{\fg^{p^m}}(\Pi^nb)\big) \!\!\!\mod\, \pi^{\Phi(n)+1}A,
$$
with $m^{\rm th}$ term in $\pi^{\Phi_m(n)}A\subseteq\pi^{\Phi(n)}A,$
for all $m.$}

\vs\noi
\begin{proof}
Consider  $\fM \os {\rm def}= \,\Maps(\Z/p^s\Z,I),$ with $I=
\{0,1\}\subseteq \Z,$
as a $\fg$-set by $^{g^j}\!\mu  (i) = \mu  \big(\omega^{-j}(i)\big) = \mu 
(i-j)$ for all  $i\in \Z/p^s\Z.$  Then 
 $h(\mu  ) =\,
^{(\sum_{j\in \Z/p^s\Z}\mu
(j)g^j)}(\Pi^nb)$
defines a $\fg$-set map $h:\fM\to B:$ for $\us j\sum
(^{g^k}\mu)(j)g^j= \us j\sum  \mu  (j-k)g^j
=\us j\sum \mu  (j)g^{j+k} = g^k\Big(\sum_j \mu 
(j)g^j\Big).
$
The usual $\fg$-set analysis, with $\fM^*$ a set of
representatives for the $\fg$-orbits on $\fM,$ implies
$$
N_{\fg}(1+\Pi^nb) =
\prod_{j \in \Z/p^s\Z} \big(1+{}^{g^j}(\Pi^nb)\big)
\dot= \sum_{\mu 
\in\fM} h(\mu  )
= \sum_{\mu  \in\fM^*} \,\Tr_{\fg/\St(\mu  )}
\big(h(\mu  )\big)
= \sum^s_{m=0}\,\sum_{\mu  \in
\fM^*_m}\,\Tr_{\fg/\fg^{p^m}}\big(h(\mu  )\big),
$$
where $\fM^*_m =\{\mu  \in \fM^*:\,\St_\fg (\mu  ) = \fg^{p^m}\},$
and both sides of $\dot =$ are indexed by the subsets of $\Z/{p^s\Z}.$

\vs\noi
Now $\mu  \in \fM^*_0$ implies $\mu  ={}^{g^{p^0}}\mu  ={}^{g}\mu  $
hence $\mu  (i) = \mu  (i-1)$ for all $i,$ so $\mu  $ is the
constant function $\ul 0,$ or  $\ul 1.$
Thus  $h(\ul 0) + h(\ul 1) = {}^0(\Pi^n b) + \us j\prod\, ^{g^j}(\Pi^n b)
 = 1 + N_\fg (\Pi^nb)$
identifies the initial terms claimed, and it remains to show that
$$
\sum\nolimits _{\mu  \in\fM^*_m}\,\Tr_{\fg/\fg^{p^m}}\big(h(\mu 
)\big) \equiv \,\Tr_{\fg/\fg^{p^m}}\big(N_{\fg^{p^m}}(\Pi^nb)\big)
\!\!\!\mod\,\pi^{\Phi(n)+1}A \;\text{for} \;1\le m\le s.
\leqno {\rm (4E1)}
$$

\noi
Fixing $m\ge 1,$ then
$\mu \in\fM^*_m$ has $St_\fg(\mu) = \fg^{p^m},$ hence $\mu$ 
 has $\mu  (j) = (\,^{g^{p^m}}\mu  )(j) =
\mu  \big(\omega^{-p^m }(j)\big) =\mu  (j-p^m)$
for all $j,$ hence $\mu  (j) = \mu  (j-p^mk)$ by induction on $k.$
Thus there is a unique
function $\mu  _m$ so that
$$
\begin{array}{ccl}
\Z/p^s\Z &\os \mu  \lra &I  \\
\\
-\downarrow &&\Vert\\
\\
\Z/p^m\Z &\os{\mu  _m}\lra &I
\end{array}
\leqno {\rm (4E2)}
$$
commutes.
Since $\Z/p^s\Z = \us{i\in\Z/p^m\Z}{\dot\bigcup}(\wt i +
p^m\Z/p^s\Z),$ 
with $\tilde i \in \Z/p^s\Z$ having
$\ol{\tilde i} = i,$ we have
$\us {j\in\Z/p^s\Z}\sum \mu  (j)g^j
= \us i\sum\,
\us{z\in p^m\Z/p^s\Z}\sum
\mu  (\tilde i + z) g^{\tilde i+z} 
\doteq \us z \sum \big(\us i\sum\mu  _m(i)g^{\tilde i}\big)g^z
= \big(\us i\sum \mu  _m(i)g^{\tilde i}\big)
\big(\us{0\le k<p^{s-m}}\sum g^{p^mk}\big),$
with 
$\doteq$ by (4E2), i.e. $\mu(\tilde i+z) =
\mu _m \Big(\,\ol{\tilde i +z}\,\Big) = \mu  _m(i).$ 

\vs\noi
Thus, writing $B^\pr_m$ for the subring of $B$
fixed by $\fg^{p^m}$ and $g_m = g\fg^{p^m}\in \fg/\fg^{p^m},$ it
follows that
$h(\mu)  
={} ^{\us i\sum \mu  _m(i)g^{\tilde i}}\!
N_{\fg^{p^m}}
(\Pi^n b)\,\dot \in\,{}^{\us i\sum \mu  _m(i)g^i_m}
(\pi^n_{r+m} B^\pr_m)
= $
the ideal product
$\us{i\in \Z/p^m\Z}\Pi {}^{\mu_m(i)g^i_m}(\pi^n_{r+m}B^\pr_m)
= \us i\Pi{}^{ \mu_m(i)}(\pi^n_{r+m} B^\pr_m)
= \pi_{r+m}^{\vert  \mu  _m\vert  n}B^\pr_m$
with $\dot\in$ by (2A3), and $\vert \mu  _m\vert
=\sum_i \mu  _m(i).$ 

\vs\noi
Note that $\vert  \mu  _m\vert  \ge 1,$ since
$\vert  \mu  _m\vert  =0$ implies $\mu  =\ul 0\in \fM^*_0,$ 
contrary to $m\ge 1,$
hence
$\Tr_{\fg/\fg^{p^m}}\big(h(\mu  )\big) \in
\,\Tr_{\fg/\fg^{p^m}}(\pi^{n\vert  \mu  _m\vert  }_{r+m} B^\pr_m) =
\pi^{\Phi_m(n\vert  \mu  _m\vert  )} A,$ by Lemma~1B(ii).  If $\vert 
\mu  _m\vert   \ge 2,$ then $\Phi_m(n\vert  \mu  _m\vert  )\ge
\Phi_m(2n) \ge \Phi(n)+1$ by Lemma 1E.  Combining, we have
restricted the sum on the left side of (4E1) to $\mu  \in \fM^*_m$ with
$\vert  \mu  _m\vert  =1.$

\vs\noi
Replacing $\mu  $ in its $\fg$-orbit (by $^{g^{p^mk}}\mu  $ for suitable
$k),$  we may assume that $\mu  (0) =1.$ Then $\vert  \mu  _m\vert 
=1$ uniquely determines $\mu  $ by $\mu(i)  =\begin{cases} 1, &i\equiv
0\!\!\!\mod\, p^m,\\ 0, &\text{otherwise},\q\end{cases} $ hence $h(\mu  ) =
N_{\fg^{p^m}}(\Pi^n b).$

\noi
Finally, $\Tr_{\fg/\fg^{p^m}}\big(N_{\fg^{p^m}}(\Pi^n b)\big) \in
\Tr_{\fg/\fg^{p^m}}(\pi^n_{r+m} B^\pr_m) =\pi^{\Phi_m'(n)}A,$ by
Lemma~1B(ii), which is in $\pi^{\Phi(n)}A.$
\end{proof}

\vs\noi
{\sc Proposition 4F.} {\it With the notation and hypothesis of
Lemma} 4D, {\it suppose that $\sigma  \in \fR_k$ is in the special
case of} (ii).  {\it Then $\sigma  $ is in the $\fg$-orbit of
$\alpha  ^{-v}_k,$ for some $0<v<p^{s-k}$ and $\alpha  _k\in \fG_s$
defined by $\alpha  _k(i) = \begin{cases} i+p^k, &\ol i = \ol 0\\ i,
&\ol i\ne \ol 0\end{cases}$ where $^-: \Z/p^s\Z \to \Z/p^k\Z$ is the natural
map.  Moreover, $\sgn (\alpha  _k) = 1$ and}
$$
\Tr_{\fg/\fg^{p^k}}\big(\bd(\alpha^{-v}_k)\big)\equiv  
(1+T)^v\;\Tr_{\fg/\fg^{p^k}}\big(N_{\fg^{p^k}}(\Pi^n b_{p^kv})\big)
\!\!\!\mod \, \pi^{\Phi(n)+1} A.
$$

\vs\noi
{\it Proof.} 
We start by checking 
that the powers $\alpha  ^v_k$ are given by
$\alpha  ^v_k(i) =\begin{cases} i + p^kv, &\ol i = \ol 0\\
i, &\ol i \ne \ol 0\end{cases}$ by induction on $v\ge 1.$ 
This implies $\ol{\alpha  ^v_k(i)} =\ol i,$ 
so the induction step is 
$$
\begin{aligned}
\alpha  ^v_k(i) = 
\alpha  _{k}\big(\alpha 
^{v-1}_k(i)\big) &=
\begin{cases} \alpha  ^{v-1}_k(i) + p^k, &\ol{\alpha  ^{v-1}_k(i)}=
\ol 0\\ 
\alpha  ^{v-1}_k(i), &\ol{\alpha  ^{v-1}_k(i)} \ne \ol 0 \end{cases}
 =
\begin{cases} \alpha  ^{v-1}_k(i) + p^k, &\ol i=\ol 0\\
\alpha  ^{v-1}_k(i), &\ol i\ne \ol 0\end{cases} \\
&=
\begin{cases} i+p^k(v-1)+ p^k, &\ol i = \ol 0\\
i, &\ol i \ne \ol 0\end{cases} =
\begin{cases} i+p^kv, &\ol i = \ol 0\\
 i, &\ol i\ne \ol 0.\end{cases}
\end{aligned}
$$

\vs\noi
If $j\in \Z/p^s\Z$
then
$(\,^{g^j}\alpha^v  _k)(x) = (\omega  ^j\alpha^v  _k\omega  ^{-j})(x) =
j+\alpha^v  _k(x-j) $
$$ 
{}\q\q\q\q\q\q\q = j +\begin{cases} x-j+p^k v, &\ol{x-j} = \ol 0\\
x-j, &\ol{x-j} \ne \ol 0\end{cases}
= \begin{cases} x+p^k v, &\ol x=\ol j\\
x, &\ol x \ne \ol j. \end{cases} 
$$

\noi
 In particular, $\,^{g^j}\alpha  _k$ is the cycle 
$\big(j,j+p^k, \dots, j+(p^{s-k}-1)p^k\big)$ with $\vert  \supp (\,^{g^j}\alpha  _k)\vert 
=p^{s-k}$ and the special coset
$j+p^k\Z/p^s\Z.$ 
It follows that $\omega  ^{p^k}$ has the cycle decomposition
$\omega  ^{p^k} = \us{j\in \Z/p^k\Z}\Pi{} (^{g^j}\alpha  _k)$
(with product by function composition), and that
$\supp(\omega  ^{p^k}) = \Z/p^s\Z.$

\vs\noi  
{\it Step} 1:
We next show that $\sigma  =\alpha  ^{-v}_k$ has $f_+(\sigma 
)=v$ when $0<v<p^{s-k}.$  

\noi
We have $\sigma  (x) =\begin{cases}
x-p^kv, &\ol x=\ol 0\\ x, &\ol x=\ol 0\end{cases},$ from the first
paragraph of the proof, and

\noi $f_+(x) = \us{x\in\Z/p^s\Z}\sum
h_{\sigma  (x)}(x) = \us x\sum \begin{cases} 1, &\wh{\sigma  (x)}
>\wh x\\ 0, &\wh{\sigma  (x)} \le \wh x\end{cases}$ from (B4B), so
we need to count the $x\in \Z/p^s\Z$ with $\wh{\sigma  (x)} - \wh x
>0.$  The $x$ with $\ol x\ne \ol 0$ have $\wh{\sigma  (x)} - \wh x
=0$ so need only count the $x\in \Z/p^s\Z$ with $\ol x = \ol 0,$ and
$\wh{\sigma  (x)} -\wh{x} >0.$  These $x$ have 
$$
\wh{\sigma  (x)} -
\wh x = \wh{x-p^kv} - \wh x =\begin{cases}
-p^kv, &\wh x\ge p^kv\\
p^k(p^{s-k}-v), &\wh x< p^kv\end{cases} .
$$
The remaining $x$ with $\wh x <p^kv$ are the $p^ky$ for which
$0<y<v,$ giving Step~1.

\vs\noi
We also use the notation, and  commutative square,
from the proof of Lemma 4D, as the language of permutations fits
well here.  We seek to classify the $\sigma  \ne 1$ in
$C_{\fG_s}(\omega^{p^k})$
with `small'  support, i.e., $\vert  \supp(\sigma  )\vert   = p^{s-k}
= \vert  p^k\Z/p^s\Z\vert  ,$ 
$\Phi_k(n) = \Phi(n),$
and the special  coset $ \wt h +
p^k\Z/p^s\Z$  of (4D1) with $\wt h\in \Z/p^s\Z$ a pre-image in $\Z/p^s\Z$ of
the unique $h\in\CF^\pr_k(\sigma  )$ by Lemma~4D(ii). 
  Then $\CF^\pr_k(\sigma  ) = \{h\}$
implies 
that $\sigma  _k(h) = h$ and
that $\wt h -\sigma  (\wt h) = p^kv,$ for some $v\not\equiv
0\!\!\!\mod p^{s-k},$ 
by the (4D1) convention.
Thus
$$
\sigma  (x) =\begin{cases} x-p^kv, &\ol x= h\\
x, &\ol x\ne h, \end{cases} 
$$
since $\sigma  (\wt h + p^ky) = \sigma  (\wt h) + p^ky = \wt h-p^kv
+ p^ky = (\wt h + p^ky) - p^kv.$  
This accounts for all $p^{s-k}$ elements in $\supp(\sigma  ),$ hence
$\sigma  (x) = x$ for all $\ol x\ne h.$

\vs\noi
Comparing
$\sigma  $ with the $^{g^j}\alpha^v  _k$
of the first paragraph of this proof
now implies that 
$\sigma  =
\,^{g^h}\alpha^{-v}  _k$
is in the $\fg$-orbit of $\alpha  ^{-v}_k,$ as
asserted. 

\vs\noi
Conversely, we  check that $\St_{\fg}(\alpha 
^{-v}_k) = \fg^{p^k}$
(as $0<v<p^{s-k}),$
which  means that
$\omega  ^{p^{k-1}}$ does not commute with $\alpha  ^{-v}_k$ 
i.e. $^{g^{p^{k-1}}}\alpha  ^{-v}_k \ne \alpha  ^{-v}_k.$ 
The same first paragraph gives
$(^{g^{p^{k-1}}}\alpha  ^{-v}_k)(0)=0,$ and
$(\alpha  ^{-v}_k)(0) = -p^kv,$ with $0\ne -p^kv.$

\vs\noi
From now on we may set $\sigma  =\alpha   ^{-v}_k,$ with
$0\le k\le s-1$ (as there is no $v$ when $k=s).$ 

\noi
Our next task is to simplify (4D1) to get the final formula of the
Proposition.

\vs\noi
With $\CF^\pr_k(\sigma  )  = \{h\}$ and $\CF_k(\sigma  ) =\Z/p^k\Z\bs
\{h\},$ applying $\Tr_{\fg/\fg^{p^k}}$ to (4D1) in the form 
$
\bd(\sigma) =(1+T)^{f_+(\sigma  )} N_\fg(1+\Pi^n b_0) \,
\frac{{}^{\fg^{h}_k}N_{\fg^{p^k}}\big(\Pi^n b\,_{\wt h -\sigma  (\wt
h)\big)}}
{{^{\fg^{h}_k}} N_{\fg^{p^k}}(1+\Pi^n b_0)}
$
yields
$$
\Tr_{\fg/\fg^{p^k}}\big(\bd(\sigma  )\big) = (1+T)^{f_+(\sigma 
)}N_\fg(1+\Pi^nb_0)\Tr_{\fg/\fg^{p^k}}\Big(\frac{^{g^h_k}N_{\fg^{p^k}}(
\Pi^nb_{p^kv})}{{^{g^{h}_k}}N_{\fg^{p^k}}(1+\Pi^n b_0)}\,\Big),
$$
because $\fg$ fixes $N_\fg (1+\Pi^n b_0),$
 $\sigma  _k(h) = h$ and $\wt h-\sigma  (\wt h) = p^kv.$  The
$g^h_k$ can be removed, since $g^h_k(\fg/\fg^{p^k}) =
\fg/\fg^{p^k},$ 
which implies
$$
\begin{aligned}
\Tr_{\fg/\fg^{p^k}} 
\big(\bd(\sigma  )\big) &- (1+T)^{f_+(\sigma  )}\Tr_{\fg/\fg^{p^k}}
\Big(\frac{N_{\fg^{p^k}}(\Pi^n b_{p^k v})}{N_{\fg^{p^k}}(1+\Pi^nb_0)}\Big)
\\
&= (1+T)^{f_+(\sigma 
)}\big(N_\fg(1+\Pi^n{b_0)}-1\big)\Tr_{\fg/\fg^{p^k}}
\Big(\frac{N_{\fg^{p^k}}(\Pi^n b_{p^k v})}
{N_{\fg^{p^k}}(1+\Pi^nb_0)}\Big).
\end{aligned}
\leqno{\rm (4F1)}
 $$ 
The methods of (4D1), in the notation of Lemma~1B(ii),  bounds
the right term of (4F1) by
$$
\Tr_{\fg/\fg^{p^k}}
\big(
\frac {N_{\fg^{p^k}}(\Pi^nb_{p^k v})}
{N_{\fg^{p^k}}(1+\Pi^{^{n}}b_0)}\,B^\pr _k \big)
\doteq \Tr_{\fg/\fg^{p^k}} \big(N_{\fg^{p^k}}(\Pi^n
b_{p^kv})B^\pr_k\big)\subseteq
\Tr_{\fg/\fg^{p^k}}(\pi^n_{r+k}B^\pr_k )
=\pi^{\Phi_k(n)}A\subseteq \pi^{\Phi(n)}A,
$$
 with
$\doteq$ since $N_{\fg^{p^k}}(1+\Pi^n b_0)$
is a unit of $B^\pr_k. $ 
Now the middle term  $N_\fg(1+\Pi^n b_0)-1 \equiv 0
\!\!\!\mod\,\pi^{\Phi(n)}A,$ by Lemma 4E,
hence the right half of (4F1) is $\subseteq (1+T)^{f_+(\sigma  )}
\pi^{2\Phi(n)}A \subseteq \pi^{\Phi(n)+1}A,$ as $\Phi(n)>0,$ when
(4F1)
 implies that
$$
\Tr_{\fg/\fg^{p^k}}\big(\bd(\sigma  )\big) \equiv (1+T)^{f_+(\sigma 
)}\Tr_{\fg/\fg^{p^k}}\Big(\frac{N_{\fg^{p^k}}(\Pi^nb_{p^kv})}
{ N_{\fg^{p^k}}(1+\Pi^nb_0)}\,\Big)\!\!\!\mod\, \pi^{\Phi(n)+1}A.
\leqno{\rm (4F2)}
$$

\vs\noi
Similarly, we know $N_{\fg^{p^k}}(\Pi^n b_{p^kv}) \in \pi^n_{r+k}
B^\pr_k,$ from (2A3), and also that 
$$
N_{\fg^{p^k}}(1+\Pi^n b_0)-1 \equiv 0
\!\!\!\mod\,\pi^{\Phi^\pr(n)}_{r+k} B^\pr_k,
\leqno{\rm (4F3)}
$$
with $\Phi^\pr =
\Phi^{(r+k)},$ by Lemma 4E and Lemma~1B(ii) with $r,s$ replaced by $r+k,$ $s-k,$
respectively. 
Next we compare (4F2) with

\vs\noi
$\Tr_{\fg/\fg^{p^k}}
\Big(\frac{N_{\fg^{p^k}}
(\Pi^n b_{p^kv})}
{N_{\fg^{p^k}}
(1+\Pi^nb_0)}\Big)
- \Tr_{\fg/\fg^{p^k}}
\big(N_{\fg^{p^k}}(\Pi^n b_{p^kv})\big) 
= \Tr_{\fg/\fg^{p^k}}\Big(\big(N_{\fg^{p^k}}(1+\Pi^n b_0)^{-1}-1\big)
N_{\fg^{p^k}}(\Pi^n b_{p^kv})\Big)$

\vs
$
{}\q \dot\in \,
\Tr_{\fg/\fg^{p^k}}\Big(\big(N_{g^{p^k}}(1+\Pi^n b_0)-1\big)
\pi^n_{r+k}B^\pr_k\Big)\subseteq
\Tr_{\fg/\fg^{p^k}}(\pi^{\Phi^\pr (n)+n}_{r+k}B^\pr_k) =
\pi^{\Phi_k(n+\Phi^\pr(n))}A,
$

\vs\noi
with $\dot\in$ by a factor $-N_{\fg^{p^k}} (1+\Pi
b_0)^{-1} \ddot\in (B^\pr_k)^\times$ inside the right $\Tr,$ and
$\ddot\in$ by (4F3).  We can now replace the right hand $\Tr$ in
(4F2) by the second $\Tr$ in the formula above by 
 the final formula
$$
\Tr_{\fg/\fg^{p^k}}\big(\bd(\sigma  )\big) \equiv (1+T)^{f_+(\sigma 
)}\Tr_{\fg/\fg^{p^k}}\big(N_{\fg^{p^k}}(\Pi^n
b_{p^kv})\big)\!\!\!\mod \pi^{\Phi(n)+1} A,
$$
provided that we can prove that $\Phi_k\big(n+\Phi^\pr(n)\big) >
\Phi(n).$

\vs\noi
For this, we write $\Phi^\pr(n) = \Phi^\pr_{w^\pr}(n)$ with $w^\pr$
the minimal weight of $n$ with respect to $\Phi^\pr$
(cf Lemma~1C).
Then
$
\Phi_k\big(n+\Phi^\pr_{w^\pr}(n)\big)
= k\phi(p^r) + \Big\lfloor\frac{n+\Phi^\pr_{w^\pr}(n)}{p^k}\Big\rfloor
= k\phi(p^r) + \Big\lfloor \frac{n+w^\pr \phi(p^{r+k})+\lfloor
n/p^{w^\pr}\rfloor}{p^k}\Big\rfloor
= \Phi_k(n) + w^\pr\phi(p^r) + \Big(\Big\lfloor
\frac{n+\lfloor n/p^{w^\pr}\rfloor}{p^k}\Big \rfloor -
 \lfloor n/p^k\rfloor\Big)
\ge \Phi_k(n) + w^\pr\phi(p^r)
\ge \Phi(n).
$

\vs\noi
A counterexample would need $w^\pr =0,$ where $\Phi^\pr_{w^\pr}(n) =
n$  implies
$\Phi_k\big(n+\Phi^\pr(n)\big) = \Phi_k(2n) > \Phi(n)$ by Lemma 1E,
unless $k=0.$ 
When $w^\pr =0$ and $k=0,$
$\Phi_0\big(n+\Phi^\pr_0(n)\big) = \Phi_0(2n)> \Phi_0(n) \ge
\Phi(n),$ for $n\ge 1,$ so there is no counter example.

\noi
Moreover, $\sgn(\alpha  _k)=1$ since $\alpha  _k$ has $p$-power
order, with $p\ne 2.$ \hfill $\square$

\vs\noi
{\sc Theorem 4G.} {\it If $d=1+\Pi^n\sum_{j\in\Z/p^s\Z} b_j u^{\hat j}$ with
$1\le n <p^{r+s-1},$ then
$$
\nr(d)-1\equiv \sum\nolimits^*_k \,\sum_{0\le v<p^{s-k}} \,(1+T)^v\;
\Tr_{\fg/\fg^{p^k}}
\big(N_{\fg^{p^k}}(\Pi^nb_{p^kv})\big)\!\!\!\mod\, \pi^{\Phi(n)+1} A
$$
with each term of the sum in $\pi^{\Phi (n)}A.$  Here $\sum^*_k$
means that $k$ is restricted to have $\Phi_k(n) = \Phi(n),$ i.e., to
$k=w$ when $n$ is ordinary of weight $w,$ and to $k\in \{w,w+1\}$
when $n$ is transitional of weight $w.$ }

\begin{proof}
Note that we may assume that $n<p^{r+s-1},$ by Theorem 3E(ii).

\noi
We show that most of the terms in the sum $(*)$ (prior to Lemma 
4D)
are in $\pi^{\Phi(n)+1}A.$

\noi
We start with $\sigma  =1,$ when $\bd(1) = N_\fg (1+\Pi^n b_0),$ by
Lemma 4B(ii) with $\CF(1) =\Z/p^s\Z$ and $f_+(1) =0,$ when Lemma 4E gives
$$
\bd(1)-1 \equiv \sum\nolimits ^*_k\,\Tr_{\fg/\fg^{p^k}}
\big(N_{\fg^{p^k}}(\Pi^n b_0)\big)\!\!\! \mod \pi^{\Phi(n)+1} A, 
\leqno{\rm (4G1)}
$$
since $\Phi_k(n)>\Phi(n)$ for all of the other $k.$

\noi
Lemma 4D(i) assures us that all terms in $\sum^*$ are in
$\pi^{\Phi(n)}A,$ as claimed, and then even in $\pi^{\Phi(n)+1}A$
unless $\sigma  $ is in the special case of Lemma~4D(ii) when  $\Phi_k(n) = \Phi(n)$
and $\vert  \supp(\sigma  )\vert  = p^{s-k}.$  Then Proposition 4F
yields the special case terms in $\fR_k\bs\{1\}$ as
$$
\Tr_{\fg/\fg^{p^k}}\big(\bd(\sigma  )\big) \equiv
(1+T)^v\,\;\Tr_{\fg/\fg^{p^k}}\big(N_{\fg^{p^k}}(\Pi^nb_{p^kv})\big)
\!\!\!\mod \pi^{\Phi(n)+1}A, \leqno{\rm (4G2)}
$$
for $0<v<p^{s-k}.$  Combining (4G1) and (4G2) we finally get
$$
\nr(d)-1 \equiv
\sum\nolimits^*_k\Big(\Tr_{\fg/\fg^{p^k}}\big(N_{\fg^{p^k}}(\Pi^n
b_0)\big )+\sum_{1\le v<p^{s-k}}
(1+T)^v\;\Tr_{\fg/\fg^{p^k}}\big(N_{\fg^{p^k}}(\Pi^n
b_{p^kv})\big)\Big)
\!\!\!\mod \pi^{\Phi(n)+1}A,
$$
as required.
\end{proof}

\vs\noi
{\sc Corollary 4H.} {\it If $d\in \Delta^\times$ has $\nr(d)=1$ (or
$\equiv 1\!\!\!\mod  \pi^{\Phi(n)+1}A),$ then $d$ has level $n\ge
p^r-1.$}

\vs\noi
\begin{proof} We first show that
\begin{description}
\item  {1.} $\ol A =  A/\pi A$ has $\ol A^{p^k}$-basis $\{(1+\ol
T)^i: 0\le i <p^k\}$ for all $k\ge 0.$

\item  {2.} If $b\in B$ and $a\in A$ satisfy $b\equiv a\!\!\!\mod
\Pi B,$ then $N_\fg(b)\equiv a^{p^s}\!\!\!\mod \pi A.$
\end{description}

\vs\noi
{\it For} 1: $\ol A=\IF\big[[\,\ol T\,]\big]$ has basis $\{\,\ol T^i
:0\le i <p^k\}$ over $\ol A\,^{p^k} =\IF\big[[\,\ol
T\,^{p^k}]\big].$  Passing from $\{\,\ol T^i:0\le i <p^k\}$ to
$\{(1+\ol T\,)^i: 0\le i <p^k\}$ is straightforward by the binomial
expansion.

\noi
{\it For} 2: This follows from $^{g^i}b\equiv a\!\!\!\mod \Pi B,$ for
all $i\in \Z/p^s\Z,$ and $A\cap \Pi B= \pi A.$ 

\vs\noi
Returning to $d\in
\Delta  ^\times$ with $\nr(d)=1,$ recall that $d\in 1 +\Pi\Delta  $
by Lemma 4A, and let $n\ge 1$ be the level of $d.$  Suppose, to the
contrary, that $n<p^r-1,$ i.e., $n$ is ordinary of weight $0,$ by
definition,so Theorem 4G, with $\Nr(d)=1,$ yields
$$
0\equiv \sum_{0\le v < p^s} (1+T)^v N_{\fg}(\Pi^n d_v) = \pi^n \sum_v
(1+T)^v N_{\fg}(d_v)\! \!\!\!\mod \pi^{n+1}A.
$$
Writing $ d_v\equiv a_v\!\!\!\mod \Pi B,$
with $a_v\in A,$
cancelling $\pi^n$ and passing to $\ol A,$ this becomes
$$
\sum_v \,\ol a^{p^s}_v (1+\ol T\,)^v = \sum_v \ol{N_{\fg}(d_v)}
(1+\ol T)^v = \ol 0,
$$
by 2.  Hence $\ol a^{p^s}_v = \ol 0$ and $\ol a_v=\ol 0$
for all $v,$ by $1.$ and
 $\ol A$ a  domain. Thus $d_v\equiv 0\!\!\! \mod \Pi
B,$ for all $v,$ so $d\in 1+\Pi^{n+1}\Delta  ,$ contrary to $d$
having level $n.$
\end{proof}

\vs\noi
{\sc Remark.}  Fix $s=1.$  Then $d\in \Delta  ^\times$ with $\nr(d)=1$ must
have level $n\ge p^r-1,$ by Corollary 4H, and must be in $[\Delta 
^\times, \Delta  ^\times]$ when $n\ge p^r,$ by Theorem 3E.  The
Example before Lemma~2H exhibits such an element $d$ in the first
critical level $n= p^r -1,$ which is not in $[\Delta  ^\times,
\Delta  ^\times ],$ yet we know that $SK_1(D) =0$ by Wang's theorem 
 [Pi, Theorem, p.312].  Our resolution of this, in  Theorem
6D(ii), is
of quite different nature.

\section{Interpreting Theorem 4G}\label{Sec5}

We investigate the terms of the reduced norm 1 congruence in more
detail.

\vs\noi
{\sc Lemma 5A.}
{\it If $n\ge 1$ and $0\le k\le s,$ then $\pi^{-\Phi_k(n)}_r
\Tr_{\fg/\fg^{p^k}} (\pi^n_{r+k}) \equiv (-1)^k\!\!\! \mod \pi_r A.$}

\vs\noi
\begin{proof}
The $k=0$ case just says $\Phi_0(n)=n,$ so we may assume that $k\ge
1.$  Now $\Phi_k(n)= k\phi(p^r) + \lfloor \frac{n}{p^k}\rfloor,$
so, writing $m=n-p^k \lfloor\frac{n}{p^k}\rfloor $ and $M=\lfloor
\frac{n}{p^k}\rfloor,$ we have $0\le m<p^k$ and $M\ge 0.$  Then
$$
\begin{aligned}
\pi^{-\Phi_k(n)}\Tr(\pi^n_{r+k}) &=
\pi^{-(k\phi(p^r)+M)}\Tr(\pi^{m+p^kM}_{r+k})\\
&=
\pi^{-k\phi(p^r)}\Tr\big(\pi^m_{r+k}(\pi^{-1}\pi^{p^k}_{r+k})^M\big)
\os 1\equiv \pi^{-k\phi(p^r)}\Tr(\pi^m_{r+k})\\
&\os 2\equiv \pi^{-k\phi(p^r)}p^k \os 3\equiv (-1)^k\!\!\!\mod \pi A.
\end{aligned}
$$

\vs\noi
{\it For} $\os 1\equiv:$ Letting $B^\pr$ denote the fixed ring of
$\fg^{p^k},$ $1-\pi= \zeta  _r = \zeta  ^{p^k}_{r+k} = (1-\pi_{r+k})^{p^k}
\equiv 1 - \pi^{p^k}_{r+k}\!\!\!\mod pB^\pr,$ hence 
$\pi^{-1}\pi^{p^k}_{r+k}- 1\in \pi^{-1}pB^\pr\subseteq \pi B^\pr,$
hence
$(\pi^{-1}\pi^{p^k}_{r+k})^M\equiv 1\!\!\!\mod \pi
B^\pr.$   Thus
$$
\begin{aligned}
\pi^{-k\phi(p^r)}\Tr\big(\pi^m_{r+k}(\pi^{-1}\pi^{p^k}_{r+k})^M\big)
 &-\pi^{-k\phi(p^r)}\Tr(\pi^m_{r+k})
= \pi^{-k\phi(p^r)}
\Tr\big(\pi^m_{r+k}\big((\pi^{-1}\pi^{p^k}_{r+k})^M-1\big)\big)\\
\in
\pi^{-k\phi(p^r)}\Tr(\pi^m_{r+k}\pi B^\pr)
&= \pi^{1-k\phi(p^r)}\Tr(\pi^m_{r+k}B^\pr)
\dot= \pi^{1-k\phi(p^r)+\Phi_k(m)}A= \pi^{1+\lfloor\frac
{m}{p^k}\rfloor}
A=\pi A,
\end{aligned}
$$
with $\dot =$ by Lemma~1B(i), which implies
$\pi^{-k\phi(p^r)}\Tr\big(\pi^m_{r+k}(\pi^{-1}\pi^{p^k}_{r+k})^M\big)
\equiv \pi^{-k\phi(p^r)}\Tr(\pi_{r+k})\!\!\!\mod \pi A.$

\vs\noi
{\it For} $\os 2=:\Tr(\pi^m_{r+k})=\Tr\big((1-\zeta  _{r+k})^m\big)
=\sum^m_{j=0} (-1)^j\big(\begin{matrix}m\\j\end{matrix}\big)\Tr(\zeta 
^j_{r+k}),$ with $j=0$ term $p^k,$ so it suffices to show that
$\Tr(\zeta  ^j_{r+k})=0$ for $1\le j\le m <p^k.$  

\vs\noi
This implies that $\theta = \zeta  ^j_{r+k}$ is in $B^\pr \bs A,$
 since $j<p^k$ and $\zeta  ^{p^k}_{r+k} = \zeta  _r$ generates the
$p$-power roots of unity in $A.$ Thus the generator $\ol g$ of
$\Gal\big(\fq(B^\pr)/\fq(A)\big)$ moves $\theta  ,$ which has minimal
polynomial $X^{p^h}-1$ over $\fq(A)$ with $h> 1,$ and 
$  \us{0\le i<p^h}\sum {}^{\ol g^{i}}\theta  =0.$  

\vs\noi
{\it For} $\os 3\equiv : \pi^{-\phi(p^r)}_r p= \prod_{v\in
(\Z/p^r\Z)^\times}  \,\frac{1-\zeta  ^v_r}{1-\zeta  _r} \equiv
\prod_v v\equiv -1\!\!\!\mod \pi A.$  Now
raise to the power $k.$
\end{proof}

\vs\noi
{\sc Lemma 5B.} {\it If $a\in A,$ $b\in B,$ then $N_\fg (a+\Pi b)
\equiv a^{p^s}+\pi N_\fg (b)\!\!\!\mod \pi^{\phi(p^r)}A. $ }

\begin{proof} We adapt the proof of Lemma 4E,
with $n=1,$
using the same
$\fg$-set $\fM =$ Maps$(\Z/p^s\Z,I),$ but with the $\fg$-set map
$c:\fM \to B$ modifying $h(\mu  )$ to
$
c(\mu  ) = a^{p^s-\vert  \mu  \vert}h(\mu) 
$
 with $\vert  \mu  \vert  = \us{j\in\Z/p^s\Z}\sum \mu  (j) \le p^s.$

\vs\noi
Setting $\theta(\mu  ) =\us{j\in\Z/p^s\Z}\sum \mu  (j)g^j$ and
$h^\pr(\mu  ) = \,^{\theta(\mu  )}(\Pi a^{-1}b) =
\, ^{\theta(\mu  )} a^{-1}\,^{\theta(\mu  )}(\Pi b) = a^{-\vert  \mu 
\vert  }h(\mu  ),$ since  $^g a=a,$ implies that 
$$
\begin{aligned}
N_\fg(a+\Pi b) &= a^{p^s}
N_\fg(1+\Pi a^{-1}b) \dot= a^{p^s} \sum^s_{m=0} \,\sum_{\mu 
\in\fM^*_m} \,\Tr_{ \fg/\fg^{p^m}}\big(h^\pr(\mu  )\big)\\
&\ddot = \sum^s_{m=0}\,\sum_{\mu  \in \fM  ^*_m} \,a^{p^s-\vert  \mu
 \vert  } \Tr_{\fg/\fg^{p^m}}\big(h( \mu   )\big),
\end{aligned}
$$
with $\dot=$ by line 5 of loc.cit., and $\ddot =$ as $\Tr_{\fg/\fg
^{p^m}} $  is $A$-linear.

\vs\noi 
Continuing as in loc.cit., we have $\us{\mu  \in\fM^*_0}\sum
\,a^{p^s-\vert  \mu\vert    }\,\Tr_{\fg/\fg^{p^0}}\big(h(\mu)  \big) = a^{p^s}
h(\ul 0) +
a^0h (\ul 1) = a^{p^s} + N_\fg(\Pi b) =  
a^{p^s} + \pi N_\fg(b),$ the right side of our claim.
So it remains to show, for $1\le m\le s,$
that
$$
\Tr_{\fg/\fg^{p^m}}\big(h(\mu  )\big) \equiv 0\!\!\!\mod
\pi^{\phi(p^r)}A,\q 
\text{\rm for all}\q \mu  \in\fM^*_m.
\leqno{\rm (5B1)}
$$

\vs\noi
Factoring $\mu  $ by 
 (4E2),  inducing the factorization of
$\sum_j \mu  (j)g^{\hat j}$ two lines later, we get
$$
h(\mu  ) = ^{(\sum_{i\in\Z/p^m\Z}\mu
 _m(i)g^i_m)} \big(N_{\fg^{p^m}}(\Pi b)\big)\in ^{(\sum_i \mu 
_m(i)g^i_m)}(\pi_{r+m}B^\pr_m) = \pi^{\vert  \mu  _m\vert  }_{r+m}
B^\pr_m,
$$
with $B^\pr_m$ the subring of $B$ fixed by $\fg^{p^m}.$  Thus
$$
\Tr_{\fg/\fg^{p^m}}\big(h(\mu  )\big)\in
\,\Tr_{\fg/\fg^{p^m}}\big(\pi^{\vert  \mu  _m\vert  }_{r+m}
B^\pr_m\big) = \pi^{\Phi_m(\vert  \mu  _m\vert  )}A,
$$
with $\Phi_m(\vert  \mu  _m\vert  ) = m\phi(p^r) + \lfloor\,
\frac{\vert  \mu  _m\vert  }{p^m}\,\rfloor \ge \phi(p^r),$ as $m\ge
1.$
\end{proof}

\vs\noi
{\sc Proposition 5C.}
{\it If $n\ge 1$ and $0\le k\le s,$ write $b\in B$ as
$\sum^{p^s-1}_{i=0} a_i\Pi^i,$ with unique $a_i\in A.$  

\vs\noi
Then
$$
\pi^{-\Phi_k(n)}\Tr_{\fg/\fg^{p^k}}\big(N_{\fg^{p^k}}(\Pi^nb)\big)
\equiv (-1)^k \lambda  ^{(k)}_n (b)^{p^{s-k}}\!\!\!\mod \pi A,
$$

where $\lambda  ^{(k)}_n(b) = \sum^{n^*-n}_{i=0} a_i,$ and $n^* \ge
n$ is minimal with $n^*\equiv -1\!\!\!\mod p^k.$}

\begin{proof}
Write $c_0 =b$ and $c_i = a_i +\Pi c_{i+1},$ for $0\le i <p^s,$
with all $c_i= \os{p^s-1}{\us{j=i}\sum} \,a_j\Pi^{j-i}
\in B$ and $a_i$ as above.
  Let $B^\pr_k = A[\pi_{r+k}]$ be the fixed ring
$B^{\fg^{p^k}}.$  Then
$$
N_{\fg^{p^k}}(c_i) - \pi_{r+k} N_{\fg^{p^k}}(c_{i+1}) \equiv
a^{p^{s-k}}_i\!\!\!\mod \pi^{\phi(p^{r+k})}_{r+k} B^\pr_k
$$
by Lemma 5B for $B/B^\pr_k$ (so $r^\pr = r +k,$ $s^\pr = s-k),$ and
$\pi_{r+k} = N_{\fg^{p^k}}(\Pi),$ by (2A3).  Applying $\Tr_{\fg/\fg^{p^k}}
(\pi^{n+i}_{r+k^{-\!\!\!-}})$ to this, with $0\le i \le n^*-n,$ gives
$$
\begin{aligned}
\Tr_{\fg/\fg^{r^k}} \big(\pi^{n+i}_{r+k}N_{\fg^{p^k}}(c_i)\big)
&- \Tr_{\fg/\fg^{p^k}} \big(\pi^{n+i+1}_{r+k}
N_{\fg^{p^k}}(c_{i+1})\big)\\
&\equiv \Tr_{\fg/\fg^{p^k}}(\pi^{n+i}_{r+k})a^{p^{s-k}}_i \!\!\!\mod
\pi^{\Phi_k(n)+1}_r A
\end{aligned} \leqno{\rm (5C1)}
$$
since $\Tr_{\fg/\fg^{p^k}}(\pi^{n+i+\phi(p^{r+k})}_{r+k} B^\pr_k) =
\pi^{\Phi_k(n+i)+\phi(p^r)}_r A$ by Lemma 1B, 
as $\Phi_k(x) = k\phi(p^r)  +\lfloor \frac{x}{p^k}\rfloor,$
and $\phi(p^r)\ge 1.$ 
The definitions of $n^*, \Phi_k$ imply that $\Phi_k(n+i) = \Phi_k(n^*),$
for $0\le i\le n^*-n.$  Now Lemma 5A gives
$$
\Tr_{\fg/\fg^{p^k}}(\pi^{n+i}_{r+k}) \equiv
(-1)^k\pi^{\Phi_k(n^*)}\!\!\mod \pi^{\Phi_k(n^*)+1} A
$$ for $0\le i \le n^* - n\;\,(\le p^k -1 \le p^s -1).$  Thus summing (5C1)
over $0\le i\le n^*-n,$ and multiplying by $\pi^{-\Phi_k(n^*)},$
yields the telescoping sum, with $c_0 =b,$

\noi
$ \pi^{-\Phi_k(n^*)}\Tr_{\fg/\fg^{p^k}}\big(\pi^n_{r+k}N_{\fg^{p^k}}(b)
\big)
-\pi^{-\Phi_k(n^*)}
\Tr_{\fg/\fg^{p^k}}\big(\pi^{n^*+1}_{r+k}
N_{\fg^{p^k}}(c_{n^*-n+1})\big)
\equiv (-1)^k \sum^{n^*-n}_{i=0} a^{p^{s-k}}_i\!\!\!\mod \pi A.
$

\noi
 Here
$\Phi_k(n^*) = \Phi_k(n),$
\,$\sum^{n^*-n}_{i=0} a^{p^{s-k}}_i \equiv
\lambda  ^{(k)}_n(b)^{p^{s-k}}\!\!\!\mod \pi A,$ and

\noi
$
\pi^{-\Phi_k(n)}\Tr_{\fg/\fg^{p^k}}\big(\pi^{n^*+1}_{r+k}N_{\fg^{p^k}
} (c_{n^*-n+1})\big) \in
\pi^{-\Phi_k(n^*)}\Tr_{\fg/\fg^{p^k}}(\pi^{n^*+1}_{r+k}B^\pr_k) 
= \pi^{\Phi_k(n^*+1)-\Phi_k(n^*)}A
$
with $\Phi_k(n^*+1) - \Phi_k(n^*) = \lfloor\,
\frac{n^*+1}{p^k}\,\rfloor - \lfloor\,\frac{n^*}{p^k}\,\rfloor =1.$
Thus
$\pi^{-\Phi_k(n)}\Tr_{\fg/\fg^{p^k}}\big(N_{\fg^{p^k}}(\Pi^nb)\big)
\equiv (-1)^k \lambda  ^{(k)}_n (b)^{p^{s-k}}
\!\!\!\mod \pi A.$
\end{proof}

\noi
We next focus on the case $k\in \{w,w+1\}$ of most interest for Theorem
4G.  Since $n$ determines the weight $w$ of $n,$ we may define $\lambda  _n(b) =
\lambda  ^{(w)}_n(b)$ when $k=w,$ in the notation of Proposition 5C,
and also $\lambda  ^+_n(b) = \lambda  ^{(w+1)}_n (b)$ when $k=w+1,$ i.e.,
in the special
case $n$ is transitional of weight $w\le s-1,$ provided that we use
the special notation $n^+$ instead of $n^*,$ i.e.,
$$
\lambda  ^+_n(b) =\sum\nolimits^{n^+-n}_{i=0} a_i, \;\text{with}\;
n^+\ge n\;\text{minimal so that} \; n^+ \equiv -1\!\!\!\mod p^{w+1},
\leqno{\rm (+)}
$$
to distinguish $n^+$ from $n^*.$ Note that $n,n^*,n^+$ all have
weight $w.$  Thus Proposition~5C implies the congruences for $\lambda
 _n,\lambda  ^+_n$ below
$$
\begin{aligned}
\pi^{-\Phi(n)}\Tr_{\fg/\fg^{p^w}} \big(N_{\fg^{p^w}}(\Pi^n b)\big)
&\equiv (-1)^w \lambda  _n (b)^{p^{s-w}}\!\!\!\mod \pi A,
\;\text{and}\\
\pi^{-\Phi(n)}\Tr_{\fg/\fg^{p^{w+1}}}\big(N_{\fg
^{p^{w+1}}}(\Pi^nb)\big) &\equiv (-1)^{w+1}\lambda 
^+_n(b)^{p^{s-(w+1)}}\!\!\!\mod \pi A.
\end{aligned}
$$
 Also $\Phi(n)$ makes sense in the transitional case, by Corollary~1D. 

\vs\noi
{\sc Proposition 5D.} {\it Suppose $d=1+\Pi^n\sum_{j\in \Z/p^s\Z}
b_j u^j$ has level $n\ge 1,$ with  $\nr (d)
=1.$  Then}

\vs
\begin{description}
\item  {(i)} {\it If $n$ is ordinary, then $\lambda  _n (b_{p^w
v})\equiv 0\!\!\!\mod \pi A$ for $0\le v<p^{s-w}.$}

\item  {(ii)} {\it If $n$ is transitional, then, for $0\le v
<p^{s-(w+1)},$
$$
\lambda  ^+_n(b_{p^{w+1}v}) \equiv \sum^{p-1}_{z=0} \lambda  _n
(b_{p^w v+p^{s-1}z})^p(1+T)^z\!\!\!\mod \pi A.
$$
In particular, $\lambda  ^+_n(b_{p^{w+1}v})$ and $\{\lambda 
_n(b_{p^w v+p^{s-1}z}): 0\le z <p\}$ determine each other $\mod \pi
A.$ }
\end{description} 

\vs\noi
\begin{proof} {\it For} (i):
Multiplying the Proposition~5C congruence by $(1+T)^v$ and summing
over

\noi $0\le v < p^{s-w}$ gives the congruence
$$
(-1)^w\sum^{p^{s-w}-1}_{v=0}
\lambda  _n(b_{p^wv})^{p^{s-w}}(1+T)^v\equiv
\pi^{-\Phi(n)}\sum^{p^{s-w}-1}_{v=0}\,Tr_{\fg/\fg^{p^w}}
(N_{\fg^{p^w}}(\Pi^n b_{p^wv})(1+T)^v \dot\equiv 0 \!\!\!\mod \pi A,$$
with
$\dot\equiv$ by Theorem~4G with $n$ ordinary.  The
result follows from 1. of the proof of Corollary 4H, as $A/\pi
A$ is a commutative domain.

\noi
{\it For} (ii): Similarly using the two congruences just before
Proposition~5D imply that
$$
\begin{aligned}
&(-1)^w \sum^{p^{s-w}-1}_{V=0} \lambda 
_n(b_{p^wV})^{p^{s-w}}(1+T)^V
 + (-1)^{w+1}\sum^{p^{s-(w+1)}-1}_{v=0} \lambda  ^+_n
(b_{p^{w+1}v})^{p^{s-(w+1)}}(1+T)^v\\
&\q\equiv \pi^{-\Phi(n)}\Big(\sum_{0\le V<p^{s-w}}
Tr_{\fg/\fg^{p^w}}\big(N_{\fg^{p^w}}(\Pi^n b_{p^wV})\big)(1+T)^V \\
&\q\q\q\q\q\q+ \sum_{0<v<p^{s-(w+1)}}\Tr_{\fg/\fg^{p^{w+1}}}
\big(N_{\fg^{p^{w+1}}}(\Pi^n
b_{p^{w+1}v})\big) (1+T)^v\Big)
\dot\equiv\; 0\! \!\!\mod \pi A,
\end{aligned}
$$
with $\dot\equiv$ by Theorem 4G, with $n$ transitional. 

\vs\noi
Next
 the capital $V$ is used to enable the re-ordering of the first sum
according to $V=v+p^{s-(w+1)}z$ with the unique $v,z$ so that $0\le v
<p^{s-(w+1)}$ and $0\le z <p.$  This transforms the above
congruence, multiplied by $(-1)^w$ and
after some organizing, into
$$
\begin{aligned}
&\sum\nolimits ^{p^{s-(w+1)}-1}_{v=0} \Big(\sum\nolimits^{p-1}_{z=0}
\,\lambda
_n\big(b_{p^w(v+p^{s-(w+1)}z)})^{p^{s-w}}(1+T)^{p^{s-(w+1)}z}\Big)(1+
T)^v\\
&\q\q
- \sum\nolimits^{p^{s-(w+1)}-1}_{v=0} \lambda  ^+_n
(b_{p^{w+1}v})^{p^{s-(w+1)}}(1+T)^v \equiv 0\!\!\!\mod \pi A,
\end{aligned}
$$
which, combining the coefficients of $(1+T)^v,$ turns into
$$
\sum^{p^{s-(w+1)}-1}_{v=0}\,\Big[\sum^{p-1}_{z=0}\,\lambda 
_n\big(b_{p^w v+p^{s-1}z}\big)^p (1+T)^z - \lambda  ^+_n
\big(b_{p^{w+1}v}\big)\Big]^{p^{s-(w+1)}}(1+T)^v\equiv 0\!\!\!\mod
\pi A.
$$

\noi
Applying 1. of the proof of Corollary 4H implies our claimed
congruences; the same argument shows that the $\lambda  _n(b_{p^w
v+p^{s-1}z})^p$ are the co-ordinates of $\lambda 
^+_n(b_{p^{w+1}v})$ in the $\ol A^p$-basis $\{(1+\ol T)^z: 0\le z<p\}$
of $\ol A=A/\pi A.$
\end{proof}

\vs\noi
{\sc Lemma 5E.} (i) {\it Let $\Delta  _{[k]}$ denote the crossed
product $B^{\fg^{p^k}}$-suborder $B*\fg^{p^k} = \bigoplus_{0\le
j<p^{s-k}}Bu^{p^kj}$}
{\it  of $\Delta= \Delta  _{r,s,\IF} $ if  for $0\le k\le s-1,$ with
$u^{p^k}b=
\,^{g^{p^k}}bu^{p^k},$ for $b\in B,$ and $(u^{p^k})^{p^{s-k}}=1+T.$
Then the map $\partial_{[k]} :\Delta  _{[k]} \to \Delta  _{[k]}$
defined by
$$
\partial_{[k]}\Big(\sum_{0\le j <p^{s-k}} b_{p^kj} u^{p^kj}\Big)
=\sum_{0\le j <p^{s-k}}  jb_{p^kj} u^{p^kj}
$$
induces a derivation on $\Delta  _{[k]}/\Pi\Delta  _{[k]},$ a
commutative domain of characteristic $p.$}

\begin{description}
\item{(ii)} 
{\it If $x\in \Delta_{[k]},$ then $\Pi x \Pi^{-1} \equiv
x-\Pi^{p^{r+k}-1}\partial_{[k]}(x)\!\!\!\mod
\Pi^{p^{r+k}}\Delta_{[k]}.$}

\item{(iii)} 
{\it If $z\in \Delta^\times  _{[k]},$ then}
$ [\Pi,z]\equiv 1 -\Pi^{p^{r+k}-1}
\partial_{[k]}(z) z^{-1}\!\!\!\mod \Pi^{p^{r+k}}\Delta 
_{[k]}.
$
\end{description}

\vs\noi
\begin{proof}
{\it For} (i):
By (2A1), $\Delta/\Pi\Delta \simeq \IF\big[[\,\ol {u-1}\,]\big]$ is
a commutative domain of characteristic $p,$ for all $r,s.$  Since
$\Delta_{r+k,s-k,\IF} \to \Delta_{[k]}$ (induced by identifying the
left ``$B$'' with the right ``$B$'', and sending the left $u$ to
the right $u^{p^k}\in \Delta)$ is an isomorphism (note that $B^{\fg^{p
^k}}= \SO_{r+k}\big[[T]\big] = B[\pi_{r+k}]$ and that the field
extension $\fq(B)/\fq(B^{\fg^{p^k}}) $ has Galois group $\fg/\fg^{p^k}),$
this applies to all relevant $k.$

\noi
{\it For} (ii): Multiplying by $\Pi^{-1}$ and taking inverses in
Lemma 2G(iv) implies

$$
^{1-g^{p^k}} \Pi\equiv 1- \Pi^{p^{r+k}-1}\!\!\!\mod \Pi^{p^{r+k}}
B, \q
 \text{\rm hence }
$$
$$
\Pi u^{p^kj}\Pi^{-1} \equiv u^{p^{k}j} -
 \Pi^{p^{r+k}-1} ju^{p^kj}
\!\!\!\mod
\Pi^{p^{r+k}}\Delta  _{[k]}.
\leqno{\rm (5E1)}
$$ 
This follows from 
$$
\begin{aligned}
\Pi u^{p^kj} \Pi^{-1} &=
(\Pi u^{p^k} \Pi^{-1})^j = \big((^{1-g^{p^k}} \Pi) u^{p^k}\big)^j\\
&\equiv \big((1-\Pi^{p^{r+k}-1})u^{p^k}\big)^j
\;\dot\equiv\, (1-\Pi^{p^{r+k}-1})^j u^{p^kj}\q\q\q\q\q\q\;{} \\
&\equiv
(1-j\Pi^{p^{r+k}-1})u^{p^kj}\!\!\! \mod \Pi^{p^{r+k}}\Delta 
_{[k]}\,,
\end{aligned}
$$
with $\dot\equiv $ as $^{g^{p^k}}\Pi \equiv \Pi\!\!\!\mod
\Pi^{p^{r+k}}\Delta_{[k]},$ by loc.cit., implies  
the commutation 
$$
u^{p^k}(1-j\Pi^{p^{r+k-1}})
=\, {}^{g^{p^k}} 
(1-j\Pi^{p^{r+k -1}})u^{p^k}
= \big(1-\,^{g^{p^k}} j\Pi^{p^{r+k-1}})u^{p^k}
 \equiv (1-j\Pi^{p^{r+k-1}})u^{p^k}\!\!\!\mod
\Pi^{p^{r+k}}\Delta_{[k]}.
$$
If $x=\us j{\sum} b_{p^kj} u^{p^{k}j},$
left multiply (5E1) by  $b_{p^{k}j}$
to get
$\Pi b_{p^kj} u^{p^{k}j}\Pi^{-1} \equiv
b_{p^kj} u^{p^{k}j} -
\Pi^{p^{r+k}-1}jb_{p^kj}
u^{p^{k}j}
\!\!\!\mod
\Pi^{p^{r+k}}\Delta_{[k]}.$
Summing over $j$ gives (ii).

\vs\noi
It follows, for all $x\in \Delta  _{[k]},$ that

\vs\noi
(5E2) \q 
${\rm If} \; n\ge p^{r+k} \;{\rm and}\; \wt x = 1+\Pi^{n-(p^{r+k}-1)}
x,\; {\rm then}\; [\Pi,\wt x]\equiv 1
-\Pi^n\partial_{[k]}(x)\!\!\!\mod \Pi^{n+1}\Delta  _{[k]}.$

\vs\noi
For, setting $N= n-(p^{r+k}-1)\ge 1,$ we have $\Pi\wt x\Pi^{-1} =
1+\Pi^N\cdot \Pi x\Pi^{-1} \os 1\equiv
1+\Pi^N(x-\Pi^{p^{r+k}-1}\partial_{[k]}(x)\big) = \wt x - \Pi^n
\partial_{[k]}(x) \os 2 \equiv \big(1-\Pi^n\partial_{[k]}(x)\big)
\wt x\!\!\mod \Pi^{n+1}\Delta  _{[k]}\,,$ with $\os 1\equiv$ by
(ii), and $\os 2\equiv$ as every $c\in \Pi^n\Delta  _{[k]}$ has $c(\wt x-1)\in \Pi^n
\Delta  _{[k]} \Pi^N\Delta  _{[k]} \subseteq \Pi\Delta  _{[k]}.$
Right multiplication by $\wt x^{-1}$ yields (5E2).

\vs\noi
Finally, (5E2) implies that $[\Pi,\wt x](1+\Pi^n y) \equiv
\big(1-\Pi^n \partial_{[k]}(x)\big)(1+\Pi^n y) = 1+
\Pi^n\big(y-\partial_{[k]}(x)\big) -
\Pi^{2n}\big(\Pi^{-n}\partial_{[k]}(x)\Pi^n\big)y\!\!\!\mod
\Pi^{n+1}\Delta  _{[k]}$
and thus
(iii), by $2n\ge n+1.$
\end{proof}

\vs\noi
{\sc Proposition 5F.} {\it Suppose that $d=1+\Pi^n
\sum^{p^s-1}_{j=0} b_ju^j\in \Delta  ,$ with $n$ of weight $w$  and
$p^r-1 \le n< p^{r+s-1}.$  Then there exists a $\Pi$-commutator
product $c\in \Delta  ^\times$ so that }
$$
c^{-1}d\equiv 1+\Pi^n\sum\nolimits _{0\le j<p^{s-W}} b_{p^Wj}\,
u^{p^Wj}\!\!\!\mod \Pi^{n+1}\Delta  ,\;{\it where} \; W=\begin{cases}
{}\q w, &\;{\it if}\; n<p^{r+w}\\ w+1, &\;{\it if}\; n\ge p^{r+w}.
\end{cases}
$$
{\it Moreover, the $b_{p^Wj}$ are unchanged for $0\le j<p^{s-W}.$}

\vs\noi
\begin{proof}
We start with

\vs\noi
{\it Step} 1: If $n\ge p^{r+z-1},$ with $z\ge 0$, there exists a
$\Pi$-commutator product $c\in \Delta  ^\times$ with
$$
c^{-1}d\equiv 1 +\Pi^n \sum\nolimits _{0\le j <p^{s-z}} b_{p^zj}
u^{p^zj}\!\!\!\mod \Pi^{n+1}\Delta  ,
$$
and the $b_{p^zj}$ unchanged for $0\le j<p^{s-z}.$

\vs\noi
We show that there exist $\Pi$-commutator products $c_e\in \Delta 
^\times$ so that
$$
c^{-1}_e d \equiv 1+\Pi^n y_e\!\!\!\mod  \Pi^{n+1}\Delta  ,
$$
with $y_e =\sum^{p^{s-e}-1}_{j=0}b_{p^ej} u^{p^ej}\in \Delta 
_{[e]},$ by induction on $0\le e\le z.$  This holds for $e=0$ with
$c_0=1$ (and we're done if $z=0),$ and, for $e\ge 1,$ the induction
hypothesis gives $c_{e-1}.$  Solving $b_{p^{e-1}j} =
jb^\pr_{p^{e-1}j}$ in $B,$ for $0\le j< p^{s-(e-1)}$ with
$j\not\equiv 0 \!\!\!\mod p,$ and setting $x_{e-1} =\sum_{p\nmid j}
b^\pr_{p^{e-1}j} u^{p^{e-1}j}$ arranges that $ y_e \os{\rm{def}}= y_{e-1} -
\partial_{[e-1]}(x_{e-1}) \in \Delta  _{[e]}.$

\vs\noi
Now Lemma 
5E(ii) with $k=e-1$ applies, since $0\le k <s-1$ by
$p^{r+k}\le p^{r+z-1} \le n < p^{r+s-1},$ to provide $\wt
x_{e-1}\in \Delta  _{[e-1]}$ so that $c_e \os{\rm def}=
c_{e-1}[\Pi,\wt x_{e-1}]^{-1}$ has
$$
\begin{aligned}
c^{-1}_e d &= [\Pi, \wt x_{e-1}](c^{-1}_{e-1}d) \equiv [\Pi,\wt
x_{e-1}](1+\Pi^n y_{e-1})\\
&\equiv 1+\Pi^n \big(y_{e-1} - \partial_{[e-1]}(x_{e-1})\big)
= 1+\Pi^n y_e\!\!\!\mod \Pi^{n+1}\Delta  ,
\end{aligned}
$$
completing the induction and settling {\it Step} 1.


\vs\noi
{\it Step} 2:  Define $\kappa  = \kappa  (n) =\max\{z\ge 0: n\ge
p^{r+z-1}\}.$  Then $\kappa  =W.$

\vs\noi
Suppose first that $w\ge 1,$ so $p^{r+w-1} + p^{w-1}\le n <p^{r+w} +
p^w.$  Then $\kappa  =w$ if $n<p^{r+w},$ hence $\kappa  =w+1$ if
$n\ge p^{r+w} = p^{r+(w+1)-1}.$  When $w=0,$ then $n=p^r-1$ or
$p^r,$ hence $\kappa  =0$ or $1,$ respectively, giving {\it Step} 2.

\vs\noi
Finally, {\it Step} 2 allows us to take $z=\kappa  =W$ in {\it Step}
1, finishing the proof.
\end{proof}

\vs\noi
{\sc Lemma 5G.} 
\begin{description}
\item  {(i)} {\it Let $j\ge 1$ have $0\le \v_p(j)\le s-1$ and write $j=
p^{\v_p(j)}j^\pr.$  Then}
$$
^{g-1} \Pi^j \equiv 1+
j^\pr(1-\Pi^{p^{\v
_p(j)}})\Pi^{(p^r-1)p^{\v_p(j)}}+ \Big(\begin{matrix} j^\pr\\
2\end{matrix}\Big)\big(1-2\Pi^{p^{v_p(j)}}\big)\Pi^{2(p^r-1)p^{v_p(j)
}}
\!\!\!\mod
\Pi^{2p^{r+v_p(j)}}B.
$$
\item  {(ii)}
{\it If $j^\pr \equiv 1\!\!\!\mod p^k,$ then}
$\v_\Pi\Big(\begin{pmatrix} j^\pr\\ 2\end{pmatrix}\Big) \ge k(p-1)p^{r+s-1}.$

\item  {(iii)} {\it If $0\le q\le s-1$ and $h\ge 0,$ then}
$$
^{g^h-1} \Pi^{p^q} \equiv 1 + h(1-\Pi^{p^q})\Pi^{(p^r-1)p^q}
\!\!\!\mod\Pi^{2(p^r-1)p^q}B.
$$
\item{(iv)} {\it  If $0\le q\le s-1,$ $m\ge 1$ and $x\in \Delta 
,$ then
$$\Pi^{p^qm}x\Pi^{-p^qm} \equiv
x-m(1-\Pi^{p^q})\Pi^{(p^r-1)p^q}\partial (x)\!\!\!\mod
\Pi^{2(p^r-1)p^q}\Delta  ,$$
where $x = \os{p^s-1}{\us{h=0}\sum} b_h
u^h,$ with all $b_h\in B,$
 has $\partial (x) =\sum_h b_h hu^h.$}

\item  {(v)} {\it Notation and hypotheses as in} (i), {\it with}
 $j^\pr \equiv 1\!\!\!\mod p,$ {\it then}

\centerline{
$
^{g-1}
\Pi^j\equiv 1 +
(1-\Pi^{p^{\v_p(j)}})\Pi^{(p^r-1)p^{\v_p(j)}}
\!\!\!\mod
\Pi^{(p^r-1)p^{v_p(j)+1}}B.
$}
\end{description}

\begin{proof}
Write $v=\v_p(j)$  for brevity.  

\vs\noi
{\it For} (i): Starting from $^{g-1}\Pi \equiv
1+(1-\Pi)\Pi^{p^r-1}\!\!\!\mod pB,$ from Lemma 2G(iii), it follows
by induction on $v\ge 0,$ that
$^{g-1}\Pi^{p^v}\equiv 1+ (1-\Pi^{p^v}) \Pi^{(p^r-1)p^v}\!\!\!\mod
p$ by taking $p^{\rm th}$ powers.
Setting $v= v_p(j)$ and then 
raising both sides
to the power $j^\pr, $ by the binomial expansion, still $\!\!\!\mod
p,$ and then dropping all terms after the third, while going
$\!\!\!\mod \Pi^{2p^{r+v}},$ gives the result:
more precisely, the first and second terms are as claimed, the third
$\begin{pmatrix} j^\pr\\
2\end{pmatrix}(1-\Pi^{p^v})^2\Pi^{2(p^r-1)p^v}$ is congruent
$\!\!\!\mod \Pi^{2p^{r+v}}$ to the claimed third, and the higher
terms are congruent to $0$ (by $3(p^r-1)p^v \ge 2p^{r+v}$ for the
fourth term).
This works since
$v_\Pi(p) = \phi(p^{r+s})\ge \phi(p^{r+v+1}) =\phi(p)p^{r+v} \ge
2p^{r+v}.$

\vs\noi
{\it For} (ii): $\v_\Pi\Big(\begin{pmatrix} j^\pr\\
2\end{pmatrix}\Big) = \v_\Pi(j^\pr -1)\ge k\v_\Pi (p) =
k\phi(p^{r+s}).$
\vs\noi
{\it For} (iii): The $h=0$ case is trivial.  Fix simpler notation
$M= (p^r-1)p^q,$ $\alpha  _h = \,^{g^h-1}\Pi^{p^q},$ and $\beta  _h
= (1-\alpha  _h\Pi^{p^q})\alpha  ^{p^r-1}_h.$  Note that
$^{g^h}\Pi^{p^q} = \alpha  _h\Pi^{p^q}$ and $^{g^h}\Pi^M =
(\,^{g^h}\Pi^{p^q})^{p^r-1} = \alpha  ^{p^r-1}_h \Pi^M.$
We next prove, by full induction on $h\ge 0,$ that $\alpha  _h \equiv
1+h(1-\Pi^{p^q})\Pi^M\!\!\!\mod\Pi^{2M}B.$

\noi
The case $h=1$ is  (i), with $j=p^q,$ 
which implies
the congruence, for $h\ge 0,$
$$
\begin{aligned}
\alpha  _h &= \,^{1+g+\cdots+g^{h-1}}\alpha  _1
= \prod\nolimits^{h-1}_{i=0} \,^{g^i}\alpha  _1 
\equiv \prod\nolimits_i [1+(1- \,^{g^i}\Pi^{p^q}) \,^{g^i}\Pi^M]\\
&= \prod\nolimits_i [1+(1-\alpha  _i \Pi^{p^q})(\alpha  ^{p^r-1}_i \Pi^M)]
= \prod\nolimits_i(1+\beta  _i \Pi^M)
\equiv 1+\sum\nolimits^{h-1}_{i=0} \beta  _i\Pi^M\!\!\!\mod \Pi^{2M}.
\end{aligned}
$$

\noi
Now  $\alpha  _i \equiv 1+
i(1-\Pi^{p^q})\Pi^M$ (always $\!\!\!\mod \Pi^{2M}),$ for $0\le i\le
h-1,$ implies $\alpha  ^p_i\equiv 1 + pi(1-\Pi^{p^q})\Pi^M\dot\equiv
1,$ with $\dot\equiv$ by $v_\Pi(p) =\phi(p^{r+s}) \ge
\phi(p^{r+q+1}) \ge M;$
similarly, $\alpha  ^{p-1}_i \equiv 1+(p-1)i(1-\Pi^{p^q})\Pi^M
\equiv 1-i(1-\Pi^{p^q})\Pi^M.$ 
It follows, from 
$p^r-1 = p(p^{r-1}-1)+ (p-1),$ that $\beta  _i\equiv (1-\alpha 
_i\Pi^{p^q})\alpha  ^{p-1}_i = \alpha  ^{p-1}_i - \alpha  ^p_i
\Pi^{p^q} \equiv
1-i(1-\Pi^{p^q})\Pi^M - \Pi^{p^q},$ hence $\beta _i \Pi^M\equiv
\Pi^M - \Pi^{M+p^q} = (1-\Pi^{p^q})\Pi^M.$  Substituting this into
 $\alpha  _h \equiv 1 +\sum^{h-1}_{i=0} \beta 
_i \Pi^M\equiv 1+\sum_i(1-\Pi^{p^q})\Pi^M $ yields $\alpha  _h\equiv
1+h(1-\Pi^{p^q})\Pi^M\!\!\!\mod \Pi^{2M}.$

\vs\noi
{\it For} (iv):
Set $M= (p^r-1)p^q,$ so taking $m^{\rm{th}}$ powers in (iii) implies

\noi
\centerline{$
^{g^h-1}\Pi^{p^qm} \equiv 1+ mh(1-\Pi^{p^q})\Pi^M\!\!\!\mod
\Pi^{2M}B.$}

\noi
Now{}\q  $\Pi^{p^qm} x\Pi^{-p^qm}-x = \sum_h b_h
\Pi^{p^qm}u^h\Pi^{-p^qm}-x \; =  \;\sum_h b_h(^{1-g^h}\Pi^{p^qm}-1)u^h\;
\dot\equiv$

\centerline{
$\sum_h b_h\big(-mh(1-\Pi^{p^q})\Pi^M\big) u^h
\equiv -m(1-\Pi^{p^q}) \Pi^M\sum_h b_h hu^h\!\!\!\mod \Pi^{2M}\Delta 
,$}

\noi
with $\dot\equiv$ by 
taking inverses of both sides of
the first sentence
and introducing $\partial (x).$ 

\vs\noi
{\it For} (v): 
Writing $v=v_p(j),$ 
$M= (1-\Pi^{p^v})\Pi^{(p^r-1)p^v}$ and
$j^\pr = 1+pk$ for the congruence
${}^{g-1}\Pi^j \equiv (1+M)^{j^\pr}
\!\!\!\mod pB,$ 
which is the $j^\pr$-th power of the congruence on the second line of
the proof of
``For (i)''.
Then $(1+M)^{1+pk} \equiv
(1+M)(1+M^p)^k\!\!\!\mod pMB$ and $(1+M^p)^k\equiv 1\!\!\!\mod M^pB$
implies
$
(1+M)^{j^\pr} \equiv 1+M\!\!\!\mod M^pB
$
because $v_\Pi(pM)\ge v_\Pi(M^p),$ which is equivalent to 
$v_\Pi (p)
\ge (p-1)v_\Pi(M)$ which holds since $v_\Pi(p) = \phi(p^{r+s}) =
(p-1)p^{r+s-1} \ge (p-1)p^{r+v} \ge (p-1)(p^r-1)p^v = (p-1)v_\Pi(M).$
\end{proof}

\vs\noi
{\sc Proposition 5H.}  {\it Suppose that $d= 1+\Pi^n
\sum\limits^{p^s-1}_{j=0} b_ju^j \in \Delta,$  with $n$ of weight
$w$ and $p^{r+w} \le n < p^{r+s-1}.$  Then there exists a commutator
product $c\in \Delta^\times$ so that
$$
c^{-1}d\equiv 1 +\Pi^n \sum\limits_{0\le j
<p^{s-(w+1)}} b_{p^{w+1}j} u^{p^{w+1}j} \!\!\!\mod \Pi^{n+1}\Delta.
$$ 
Moreover the $b_{p^{w+1}j}$ are unchanged for $0\le j<p^{s-(w+1)}.$}

\vs\noi
{\sc Remark.}  Proposition 5H is the $W=w+1$ special case of
Proposition~5F.  Our purpose is to discuss $\Pi$-commutator products
for exceptional levels in analogy with Proposition~6G (for ordinary
levels).  The difficulty here is that a better understanding of
Theorem~4G is necessary to replace Proposition~4D(i).

\section{Special Levels and Type}\label{Sec6}

Inferences from $\lambda  _n(b)\!\!\!\mod \pi A$ to $b\!\!\!\mod \Pi
B$ are hard to find, so we first focus on some special cases of
particular interest.

\vs\noi
{\sc Lemma 6A.} {\it If $n$ has weight $w$ and $n \equiv
-1\!\!\!\mod p^w,$ then
$\lambda  _n(b)\equiv b\!\!\!\mod \Pi B$ for all $b\in B.$ 
Moreover,  if $n$ is transitional, then either

(i) $n=p^{r+w}-1,$ {\it which implies $n^+ = n^*,$ or}

(ii) $n=p^{r+w} + p^w-1,$ {\it which implies} $n^+ = n^*
+\phi(p^{w+1}).$

\begin{proof}
Note, from Proposition~5C,
and the definition of $n^*,$
that $n\equiv -1\!\!\!\mod p^w$ implies $n^*=n.$ 
Writing $b= \sum^{p^s-1}_{i=0} a_i\Pi^i,$ with $a_i\in A,$ implies $\lambda  _n(b)-b
\equiv \sum^{n^*-n}_{i=0} a_i - a_0 =\sum^{n^*-n}_{i=1} a_i
=0\!\!\!\mod \Pi B.$ 

\noi
When $n$ is transitional, then $p^{r+w} - p^w\le n <p^{r+w} + p^w$
implies $p^r-1 <\frac{n+1}{p^w} \le p^r+1$ with $\frac{n+1}{p^w} \in
\Z,$ hence either $n=p^{r+w}-1$ or $n= p^{r+w} + p^w -1.$ So (i)
holds since $p ^{r+w}-1 \equiv -1\!\!\!\mod p^{w+1},$ and (ii) holds
since $n^+ = p^{r+w} + p^{w+1} -1,$  both by (+).
\end{proof}

\vs\noi
{\sc Proposition 6B.}  {\it Let $d=1+ \Pi^n\sum_{0\le j < p^s}b_j
u^j,$ with $n\equiv -1\!\mod p^w$ and ordinary of weight $ w\le s-1.$ 
If $\nr(d)=1$ then $d$ is a
$\Pi$-commutator product $\!\!\!\mod \Pi^{n+1}\Delta  .$}

\vs\noi
\begin{proof}  First, Corollary 4H implies that $n\ge p^r-1,$ so $n$
ordinary implies that $w\ge 1,$ and Proposition 5F that

\vs\centerline{
$
c^{-1}d \equiv 1+\Pi^n \us{0\le v< p^{s-w}}\sum b_{p^wv} u^{p^w v}
\!\!\!\mod \Pi^{n+1}\Delta  ,
$}

\vs\noi
for a suitable $\Pi$-commutator product $c,$ because $n$ ordinary
has $W=w.$  Now $b_{p^wv}\equiv \lambda  _n(b_{p^wv}) \equiv
0\!\!\mod \Pi B,$ for $0\le v < p^{s-w},$ by Lemma 6A and
Proposition~5D(i), so  the previous sentence implies that $d\equiv
c\!\! \mod \Pi^{n+1}\Delta  .$
\end{proof} 

\noi
{\rm We next prepare for the critical level $n.$}

\vs\noi
{\sc Proposition 6C.} {\it Let $C=\IF\big[[U]\big]$ and $u=1+U.$ 
Then $\{u^j:j\in \Z/p\Z\}$ is a $C^p$-basis of $C,$ with $C^p
=\IF\big[[U^p]\big],$ and $\partial (\sum_j c_j u^j) =\sum_j jc_j
u^j$ defines a $C^p$-derivation of $C.$  Then $0\ne c =\sum_j c_j
u^j\in C$ has $c^p =c_0,$ if, and only if, $c=\partial (z) z^{-1}$
for some $z\in C^\times.$ }

\begin{proof}
Looking for $z$ in the kernel of the $C^p$-linear map $z\mapsto
cz-\partial (z),$ we first need

\vs\noi
{\it Step} 1: $\det(c-\partial) = c^p - c_0.$

\noi
For this, we follow the proof of Lemma~4B(ii) , and its notation, so work
with  a group $\fg$ of order $p$ acting
trivially on $C.$  Thus
$$
cu^{\hat j} =\sum_k c_ku^{\hat j+\hat k} = \sum_k c_k
u^{\widehat{j+k}+ph_{j,k}} = \sum_i c_{i-j}(u^p)^{h_j(i)} u^{\hat
i}, \q\text{\rm  with} \q h_j(i) =\begin{cases} 1, &\hat j>\hat i\\
0, &\hat j \le \hat i, \end{cases}
$$
noting that $u^p\in C^p$ and that $\partial (u^{\hat j}) = ju^{\hat
j}$ implies that $\det ( c-\partial) =\det[{\pmb c}_{i,j}]_{(\hat i,\hat
j)}$ with
$$
{\pmb c}_{i,j} = \begin{cases} c_{i-j}u^p, &\hat i <\hat j\\
c_0 - i. &\hat i =\hat j\\
c_{i-j}, &\hat i >\hat j
\end{cases}
$$
and ${\pmb c}(\sigma  ) =\prod_i {\pmb c}_{i,\sigma  (i)} = u^{pf_+(\sigma  )}
\prod_{i\in\CF(\sigma  )}(c_0 -i)\prod_{i\notin \CF(\sigma 
)}c_{i-\sigma  (i)}.$  

\noi
It follows that $\det (c-\partial)=
\sum_{\sigma  \in \fG_1} \sgn(\sigma  ){\pmb c}(\sigma  ),$ with $\fG_1 =
\,\Sym(\Z/p\Z).$
Letting $\fg$ act on $\fG_1$  by $^g\sigma  
=\omega\sigma  \omega^{-1},$ as in Lemma 4C,  write
$$
\det(c-\partial) = \sum_{\sigma  \in \fR_0} \,\sgn(\sigma  )\pmb c(\sigma
 ) + \sum_{\sigma  \in\fR_1} \,\sgn(\sigma  )\Tr_\fg\big(\pmb c(\sigma 
)\big)
$$
with $\fR_k$ a set of representatives of the $\fg$-orbits on $\fG_1$
with stabilizer of index $p^k.$  We have $\fR_0 =\{\omega^{-\hat j}:j\in
\Z/p\Z\}$ and $f_+(\omega^{-\hat j}) = \hat j,$ by Step~1 of the
proof of  Proposition~4F. 
It follows that
$$
\begin{aligned}
\sum_{\sigma  \in \fR_0}&\sgn(\sigma  ){\pmb c}(\sigma  )
 =\sum\nolimits _j {\pmb c}(\omega^{-\hat
j})
= {\pmb c}(1) + \sum\nolimits_{j\ne 0} u^{p\hat j}\prod\nolimits_i
c_{i-(i-j)}\\
&= \prod\nolimits_i (c_0-i) +\sum\nolimits _{j\ne 0} u^{p\hat
j}c^p_j
= c^p_0 -c _0 +\sum\nolimits _{j\ne 0} c^p_j u^{p\hat j} 
= \Big(\sum\nolimits _j c_j  u^{\hat j}\Big)^p - c_0 = c^p-c_0.
\end{aligned}
$$
It now suffices to show that $\Tr_\fg\big({\pmb c}(\sigma  )\big) =0,$ for
$\sigma  \in \fR_1.$  We use the polynomial $P_\sigma  (x) =
\prod_{i\in \CF(\sigma  )}(x-i)$ in $\IF_p [x].$  Involving it uses
$\Tr_\fg \big({\pmb c}(\sigma  )\big) = \sum_{j\in\Z/p\Z}{\pmb
c}(\omega^j \sigma 
\omega^{-j})$ and two properties, for all $j \in \Z/p\Z:$

\begin{description}
\item  {1.} $f_+(\omega^{\hat j}\sigma  \omega^{-\hat j}) =f_+ (\sigma  ),$ from 
Lemma 4C(ii), and

\item  {2.} $i\mapsto i + j$ induces bijections $\CF(\sigma  )\to
\CF(\omega^{\hat j}\sigma  \omega^{- \hat j})$ and $\CF^\pr(\sigma  ) \to
\CF^\pr(\omega ^{\hat j}\sigma  \omega^{-\hat j}),$ with $\CF^\pr = \Z/p\Z\bs
\CF,$ which
follow from $(i+j) - (\omega ^{\hat j}\sigma  \omega^{-\hat j})(i+j) = i-\sigma  (i).$
\end{description}

\vs\noi
Combining these we get
$$
\begin{aligned}
\Tr_\fg\big({\pmb c}(\sigma  )\big) &= \sum\nolimits _j u^{pf_+(\omega^{\hat
j}\sigma 
w^{-\hat j})} \prod \nolimits _{i\in \CF(\omega^{\hat j}\sigma 
\omega^{-\hat j})} (c_0 -i)
\prod\nolimits_{i\in \CF^\pr(\omega^{\hat j}\sigma  \omega^{-\hat j})} c_{i-\sigma  (i)}\\
&= u^{pf_+(\sigma  )} \sum\nolimits _j \prod\nolimits_{i\in \CF(\sigma  )} \big(c_0 -
(i+j)\big) \prod\nolimits_{i\in\CF^\pr(\sigma  )}
c_{(i+j)-(\omega^{\hat j}\sigma  \omega^{-\hat j})(i+j)}\\
&= \big(u^{pf_+(\sigma  )} \prod\nolimits_{i\in\CF^\pr(\sigma  )}
c_{i-\sigma  (i)}\big) \sum\nolimits _j P_{\sigma  }(c_0-j),
\end{aligned}
$$
so it now remains to prove $\sum_j P_\sigma  (x-j) =0.$

\noi
We next observe that $P_\sigma  (x)$ has degree $<p-1:$ for $\vert
 \CF(\sigma  )\vert  \ge p-1$ implies that $\sigma  = 1\in \fR_0,$
since $\sigma  $ can only fix $\CF^\pr(\sigma  )$ with $\vert 
\CF^\pr(\sigma  )\vert  \le 1.$  Thus $P_\sigma  (x)$ is an $\IF_p$-linear
combination of $1,x,\dots,x^{p-2},$ so it would suffice to show 
that $\sum_j (x-j)^h =0$ for $0<h<p-1,$ as the $h=0$ case is
$\us j\sum 1=0.$  Now $(x-j)^h = x^h +\sum^{h-1}_{k=1}\big(\begin{matrix} h\\
k\end{matrix}\big) (-j)^k x^{h-k} +
(-j)^h$ implies
$$
\sum_j (x-j)^h =\sum^{h-1}_{k=1} (-1)^k \begin{pmatrix} h\\
k\end{pmatrix} \Big(\sum_j j^k\Big)x^{h-k} + (-1)^h\Big(\sum_j
j^h\Big),
$$
so we need only observe that the coefficients $\sum_j j^h
=\sum_{j\ne 0} j^h =0,$ for $0<h<p-1,$ because $j\mapsto j^h$ is
then a non-trivial $\IF_p$-valued character of $\IF^\times_p.$  This
completes Step 1.

\vs\noi
{\it For} `$\Lra$':  Step 1 implies a solution $z\ne 0,$ in $C,$ of $\partial(z) =
cz.$  Since $C$ is a PID with unique prime element $U,$ we can write
$z= U^iz_1$ with $i\ge 0$ and $z_1\in C^\times.$  If $i=0,$ we're
done.

\vs\noi
If $i\ge 1,$ then $\partial (U^i) = iuU^{i-1}$ by induction
on $i,$ hence $cU^iz_1= cz = \partial (z) =\partial(U^iz_1) =
iuU^{i-1}z_1 + U^i\partial (z_1)$ implies $cU= iu +
U\partial(z_1)z^{-1}_1$ and $iu = U\big(c-\partial(z_1)z^{-1}_1\big).$
 Since $u\in C^\times,$ this implies $i\equiv 0\!\!\!\mod UC,$ hence
$iu=0$ in $C$ and $c=\partial (z_1)z^{-1}_1.$

\vs\noi
{\it For} `$\Longleftarrow$': This follows from a formula of Hochschild
(Proposition A.7.1 of [GS])
$$
\partial (c)^p = c^{p-1} \partial^{[p]} (c) -
\partial^{[p-1]}\big(c^{p-1}\partial (c)\big).
$$
Our $\partial$ has $\partial^{[p]} = \partial$
(since $\partial(u^j) = ju^j = j^pu^j = \partial^{[p]}(u^j))$
and  $C^p$-linear $\partial^{[p-1]}$ 
(as $\partial(U^p) = puU^{p-1} =0$
implies $\partial^{[p-1]}(U^p) = 0).$
  Substituting
$c=z\in C^\times,$ this becomes $\partial (z)^p = z^{p-1}\partial(z)
-\partial^{[p-1]}\big(z^{p-1}\partial(z)\big),$ hence
$$
\Big(\frac{\partial(z)}{z}\Big)^p = \frac{\partial(z)}{z} -
\partial^{[p-1]} \,\Big(\frac{\partial(z)}{z}\Big) 
$$
after multiplying by $(z^{-1})^p\in C^p,$ so
$c=\frac{\partial (z)}{z}$ has $c^p=
c-\partial^{[p-1]}(c) =\sum_j (1-j^{p-1}) c_j u^{\hat j} = c_0\,.$

 {}\q\q{} \end{proof}
\vs\noi
{\sc Theorem 6D.} {\rm (i):} {\it If 
$d= 1+\Pi^n\sum_{0\le j<p^s} b_j u^j\in \Delta  $ has 
critical level
$n= p^{r+w}-1$ with weight 
$0\le w\le s-1,$  and\,
$\nr(d)=1,$
then $d$ is congruent to a $\Pi$-commutator product $\!\!\!\mod
\Pi^{n+1}\Delta  .$}

\noi {\rm (ii):} {\it If $s=1$ then $SK_1(D)=0.$}

\vs\noi
\begin{proof} {\it For} (i): Since $n\ge p^{r+w-1}$ is transitional of weight $w,$
multiplying by a suitable $\Pi$-commutator product, from Proposition
5F, permits us to assume that
$$
d \equiv 1+ \Pi^n \sum_{0\le i <p^{s-w}} b_{p^wi}u^{p^{w_i}}\!\!\!\mod
\Pi^{n+1}\Delta  ,
$$
while keeping $\nr(d)=1$ and the $b_{p^wi}$ unchanged.

\noi
We start with the congruence
$$
\begin{aligned}
&\sum\nolimits _{0\le i <p^{s-w}}\lambda 
_n(b_{p^wi})^{p^{s-w}}(1+T)^i\\
&\q \equiv \sum\nolimits _{0\le j <p^{s-(w+1)}}\lambda 
^+_n(b_{p^{w+1}j})^{p^{s-(w+1)}}(1+T)^j \!\!\!\mod \pi A
\end{aligned} \leqno{\rm (6D1)}
$$
from (5D1), with $n=p^{r+w}-1,$ and $n^+ = n^* = n$ from Lemma
6A, hence $\lambda  ^+_n (b_{p^{w+1}j}) = \lambda  _n(b_{p^{w+1}j}),$
for $0\le j <p^{s-(w+1)},$ and $\lambda  _n(b) \equiv b\!\!\!\mod
\Pi B,$ for all $b\in B,$ by Lemma 6A.  Setting $1+T = u^{p^s}$
in (6D1) yields the congruence
$$
\sum_i b^{p^{s-w}}_{p^wi} u^{p^si} \equiv \sum_j
b^{p^{s-(w+1)}}_{p^w(pj)} u^{p^sj} \!\!\!\mod \Pi\Delta  , \;{\rm
hence}
$$
$$
\Big(\sum_i b_{p^wi} u^{p^wi}\Big)^{p\cdot p^{s-w-1}} 
\equiv \Big(\sum_j b_{p^w(pj)} u^{p^w(pj)}\Big)^{p^{s-w-1}}\!\!\!\mod
\Pi\Delta  ,
$$
when changing the $j$-sum to an $i$-sum, by setting $i=pj,$ gives the congruence
$$
\Big(\sum_i b_{p^wi}u^{p^wi}\Big)^p \equiv \sum_{i\equiv 0\!\!\!\mod
p} b_{p^wi}u ^{p^wi} \!\!\!\mod \Pi \Delta  ,
$$
by the last sentence of Lemma 5E(i).

\noi
Setting $\xi=\sum_{0\le i <p^{s-w}}b_{p^wi}u^{p^wi}\in \Delta  _{[w]}$
 and $\xi 
_{\la j\ra} =\sum_{i\mapsto j} b_{p^wi} u^{p^w(i-\hat j)}$ for $j\in\Z/p\Z,$
where $j= \hat j + p\Z$ with $0\le \hat j <p,$ then $\xi 
=\sum_{j\in\Z/p\Z}\xi  _{\la j\ra}u^{p^w{\hat j}},$ so the display and $\Delta 
_{[w]}\cap \Pi\Delta  = \Pi\Delta  _{[w]}$ imply
$$
\xi  ^p \equiv \xi  _{\la 0\ra}\!\!\!\mod \Pi\Delta  _{[w]}\,. \leqno{\rm
(6D2)}
$$
To apply Proposition 6C, we change its notation $u=1+U$ into
$u_C=1+U_C$ to avoid confusion with our $u\in\Delta  .$  We define a
ring homomorphism $\Delta  _{[w]} \to C$ by sending $u^{p^w}\mapsto
u_C$ and $B= \fO_{r+s}\big[[T]\big]\to\!\!\!\!\to\IF\big[[T]\big]\to
\IF \big[[U_C]\big] = C,$ by $T= u^{p^s}-1 = (u^{p^w})^{p^{s-w}} -1
\mapsto u^{p^{s-w}}_C - 1= U^{p^{s-w}}_C.$  This ring homomorphism
is surjective since $B$ maps onto $C^{p^{s-w}} \subseteq C^p,$ and
$C=\sum_{0\le i<p^{s-w}}C^{p^{s-w}}u^i_C.$  Moreover, the kernel of
$\Delta  _{[w]}\to C$ is $\Pi\Delta_  {[w]},$ by the first paragraph
of the proof of Lemma 2A(i) (and $\Delta  _{[w+1]}$ maps onto $C^p).$
 We have a commutative square between $\partial _{[w]}$ and
$\partial:$
$$
\begin{matrix}
\Delta  _{[w]} &\os{\partial_{[w]}}\lra &\Delta  _{[w]}\\
\\
\downarrow &&\downarrow\\
\\
C   &\os{\partial}\lra &C.\end{matrix}
$$
Now let $c$ be the image of $\xi  $ under $\Delta  _{[w]} \to C$ and
observe that each $\xi  _{(j)}$ maps into $C^p,$ since the image of
$B,$ under $\Delta  _{[w]}\to C,$ is contained in $C^p.$  It follows
that $\xi  _{(j)}\mapsto c_j,$ for all $j\in\Z/p\Z,$ since
$\{u^j_C:j\in\Z/p\Z\}$ is a $C^p$-basis of $C.$  We conclude that
$c^p=c_0,$ hence that $c=\partial (z)/z,$ with $z\in C^\times,$
follows from Proposition 6C.  Choose a pre-image $z^\pr$ of $z$
under $\Delta  _{[w]}\to C,$ and note that $z^\pr$ is a unit of
$\Delta  _{[w]}$ since $\ker (\Delta  _{[w]}\to C) \subseteq \,\rad
\Delta  _{[w]}.$

\noi
Finally, $d\equiv 1 + \Pi^{p^{r+w}-1}\xi\!\!\!\mod
\Pi^{p^{r+w}}\Delta_{[w]}  $ and Lemma 5E(iii) implies that
$$
[\Pi,z^\pr]d\equiv [\Pi,z^\pr](1+\Pi^{p^{r+w}-1}\xi) \equiv 1 +
\Pi^{p^{r+w}-1} \big(\xi  -\partial_{[w]}(z^\pr)/z^\pr\big) \equiv
1\!\!\!\mod \Pi^{p^{r+w}}\Delta_{[w]} ,
$$
since $\xi  -\partial_{[w]}(z^\pr)/z^\pr \in\,\ker(\Delta  _{[w]}
\to C)\subseteq \Pi\Delta_{[w]}.$  Thus $d\equiv
[\Pi,z^\pr]^{-1}\!\!\!\mod \Pi^{n+1}\Delta_{[w]}.$

\vs\noi
{\it For} (ii): 
Any counterexample is, by Theorem~2F,
the image of some $d\in \Delta^\times$ with $nr(d)=1,$ and $d$ has
level $p^r-1$ by the Remark ending \S4.  By Theorem~6D, $d\equiv
c\!\!\mod \Pi^{p^r}\Delta$ for some $\Pi$-commutator product $c,$
hence the class of $c^{-1}d$ is in $\ker\big(SK_1(\Delta)\to
K_1(\Delta/\pi_1\Delta)\big),$ by $\pi_1\Delta = \Pi^{p^r}\Delta.$  Thus
$c^{-1}d \in [\Delta^\times,\Delta^\times],$ by Theorem~3E, and
$d\in c [\Delta  ^\times,\Delta  ^\times]\subseteq [D^\times,
D^\times].$
\end{proof}

\vs\noi
{\sc Remark.} {\rm
  In particular, we may
asume $s\ge 2$ from now on.

\noi
It is useful, for organizing the various commutators we will use, to
have the following notions.}

\vs\noi
{\sc Definition.} {\it With $r,s,\IF$ fixed as in} {\rm \S1, }
{\it define $n\ge 1$ to have}  type $t\ge 0$ {\it
if $n=p^{r+t} + p^tm-1$ for some $m\ge 1$ with $m\not\equiv
0\!\!\!\mod p.$}

\vs\noi
{\sc Lemma 6E.} {\it Suppose that $n$ has type $t$ and weight
$w.$  Then}

\begin{description}
\item  {\rm (i)} $t$ {\it and} $m$ {\it are unique.}
\item  {\rm (ii)} $t\le w,$ {\it with equality only when} $m=1.$ 
\item  {\rm (iii)} {\it If} $n\equiv -1\!\!\mod p^w,$ {\it then $n$ is
exceptional with} $t=w.$
\end{description} }

\begin{proof}
{\it For} (i): $n+1 = p^t(p^r+m),$ with $p^r +m \not\equiv 0\!\!\mod
p,$ implies $t= v_p(n+1),$ which implies uniqueness of $m.$

\noi
{\it For} (ii): $w$ minimal with $p^{r+t} + p^tm  -1<p^{r+w} + p^w
\Llra p^t(p^r+m) \le p^w(p^r+1) \Llra 1+\frac{m-1}{p^r+1} \le
p^{w-t}.$  This implies $t=w$ if $m=1,$ and $t<w$ if $m>1.$

\noi
{\it For} (iii): $p^{r+t} + p^tm-1 \equiv -1\!\!\mod p^w$ implies
$p^t(p^r+m) \equiv 0\!\!\mod p^w,$ hence $p^t\equiv 0\!\!\mod p^w$
and $t\ge w.$  3
Comparing with (ii) gives $t= w$ and $m=1,$ hence $n=
p^{r+w} + p^w -1.$ 
\end{proof}

\vs\noi
{\sc Proposition 6F.} {\it Suppose that\,} $n,$ {\it of weight} $w,$ {\it is the level of some} $d\in \Delta
^\times$ {\it with} $nr(d)=1$ {\it which is not congruent to any}
$\Pi${\it -commutator
product} $\!\!\!\mod \Pi^{n+1}\Delta.$ {\it Then} $p^r -1 \le n < p^{r+s-1}$
{\it and} $n$ {\it has type} $t< w\le s -1,$ {\it or is exceptional
with type $t=w< s-1.$}

\vs\noi
{\it Proof:}
First, $n\ge 1$ by Lemma 4A,then $n\ge p^r-1$ by Corollary~4H,
and, second, $n<p^{r+s-1},$ by  Theorem~3E(ii), (iii),
which also implies $w\le s-1.$
  Thus
$p^r-1\le n<p^{r+s-1}.$

\vs\noi
Now $n>p^{r+w-1}$ (by $n\ge p^{r+w-1} +
p^{w-1}$ if $w\ge 1,$ and $n\ge p^r-1 > p^{r-1}$ if $w=0);$
and $p^{r+w-1} <n<p^{r+s-1}$ imply that 
$$
n- (p^{r+w-1}-1) = p^h m_0 \q\text{and} \q w\le s-1,
$$
with unique $h\ge 0$ and $1\le m_0
\not\equiv 0\!\mod p.$  There are two cases:

\vs\noi
{\it Case} 1: $h<w.$  Here $n=p^{r+h}+p^hm-1$ with $m=
p^r(p^{w-1-h}-1) + m_0\equiv m_0 \not\equiv 0\!\!\mod p$ and $m\ge
m_0\ge  1.$ So
 $n$ has
type $t=h<w\le s-1.$

\vs\noi
{\it Case} 2: $h\ge w.$  Now $n=p^{r+w-1} + p^h m_0 -1\equiv
-1\!\!\mod p^w.$ 
  So if $n$ is ordinary of weight $w$ then Proposition~6B contradicts our
$\Pi$-commutator product hypothesis on $d.$  If $n$ is transitional, then Lemma~6A
puts $n$ into two subcases.  In subcase (i), Theorem 6D  contradicts our
$\Pi$-commutator product hypothesis on $d.$  Finally, in subcase (ii), $n= p^{r+w} + p^w-1$ is
exceptional of type $t=w <s-1,$ by Corollary~1D(ii), as  $w=s-1$
contradicts
$n<p^{r+s-1}.$ \hfill $\square$

\vs\noi
{\sc Remark.} The essential levels $n$ for Proposition~6F are those with
$p^r-1\le n<p^{r+s-1}.$  The number of essential levels of type $t,$
for $0\le t\le s-2,$ is

\noi
\centerline{$ \phi (p^r) (p^{s-1-t} -1). $}

\noi
For when  $t\ge 0,$ then  $p^{r+t} + p^tm-1 <p^{r+s-1} $
implies 
$1\le m\le p^rM$ with $M= p^{s-1-t}-1$ (and $p^{r+t} + p^tm-1 >p^r-1$ always), and
the essential $n$ are in bijection with these $m$ with  $m\not\equiv
0\!\!\mod p,$ hence there are $p^rM - \lfloor \frac{p^rM}{p}\rfloor =
\phi (p^r)M$ of them.

\noi
It follows that the number of essential $n$ with no
type is $p^r(p^{s-1}-1) +1 -\os{s-2}{\us{t=0}\sum}
\,\phi(p^r)(p^{s-1-t}-1) = \phi (p^r)(s-1)+1,$ all of which appear
in the proof of Case~2 of Proposition~6F, i.e. Theorem~6D provides
$s$ levels and $\big(\phi(p^r)-1\big)(s-1)$ come from
Proposition~6B: for this it suffices to show that there are
$\phi(p^r)-1$ ordinary levels, with $n\equiv -1\!\!\!\mod p^w,$ for
each weight $w$ with $1\le w\le s-1.$

\noi
Such levels $n$ have
$p^{r+w-1} + p^{w-1}\le n <p^{r+w}-p^w$ and $n= -1+p^wN,$ hence are
in bijection with $N$ satisfying $1\le N-p^{r-1}\le \phi(p^r)-1.$ It
follows that
$n= p^{r+w-1} + p^wk-1$ with $1\le k\le\phi(p^r)-1,$ as required.

\vs\noi
{\sc Hypothesis (6Gs)}: Let $\Delta  = \Delta  _{r,s,\IF}.$  Then
every $d\in \Delta  ^\times,$ with $nr(d)=1,$  having level $n$ of type $t<s-1$ 
 is congruent to some $\Pi$-commutator product
$c\!\mod \Pi^{n+1}\Delta  .$  Note that $\Pi$-commutator product is
defined after Lemma~2H.

\vs\noi
{\sc Corollary 6G.} {\it  If $\Delta=\Delta  _{r,s,\IF}  $
satisfies Hypothesis}
(6Gs), {\it then
 $SK_1(D)=0.$ }

\begin{proof}
Any counterexample to $SK_1(D)=0$ is the image of some $d\in \Delta 
^\times$ with $nr(d)=1,$ by Theorem~2F.  If $d$ has level $n,$
then  $d\equiv c\!\!\mod \Pi^{n+1}\Delta
 $ for some $\Pi$-commutator product $c:$
 more precisely, either $d$ is congruent to
some $\Pi$-commutator product $c\!\!\mod \Pi^{n+1}\Delta  ,$ or, by
Proposition~6F, $n$ has type $t<s-1,$ hence $d$ is congruent to some
$\Pi$-commutator product $c\!\!\mod \Pi^{n+1}\Delta  $ by
Hypothesis~(6Gs).
Then $c^{-1}d\equiv 1\!\!\!\mod \Pi^{n+1}\Delta  $ (as $c\in \Delta 
^\times),$ hence $d^\pr = c^{-1}d$ has $nr(d^\pr)=1$ (by $c\in
[D^\times,D^\times ]),$ and $d^\pr$ has level $n^\pr >n,$ so
$d^\pr\equiv 1\!\!\!\mod \Pi^{n^\pr}\Delta  .$ 

\noi
It follows, by induction on $i\ge 0,$ that there exists a sequence
$\{(d_i,n_i,c_i)\}_{i\ge 0},$ starting with $(d_0,n_0,c_0) =
(d,n,c),$ so that $d_i\in \Delta  ^\times$ has $nr(d_i)=1$ and level
$n_i,c_i$ is a $\Pi$-commutator product with $c_i\equiv
d_i\!\!\!\mod \Pi^{n_i+1}\Delta  ,$ and $d_{i+1} = c^{-1}_i d_i,$
with
$n_{i+1}>n_i$ for all $i.$  This follows from the above argument
with $d_{i+1} =d^\pr_i$ and $n_{i+1} = n^\pr_i.$  

\noi
Note that $d=
(c_0c_1 \dots c_{k-1})d_k$ for all $k\ge 1.$
Eventually, the first $k$ with  $n_k\ge p^{r+s-1} $ implies that $d_k \in [\Delta 
^\times,\Delta  ^\times],$ as in the first paragraph of the proof of
Proposition~6F, when $d= (c_0c_1\dots c_{k-1})d_k\in
[\Delta^\times,\Delta^\times]\subseteq [D^\times, D^\times].$
 \end{proof}

\vs\noi
{\sc Remark.} Hypothesis (6Gs) is too strong: the first sentence in
the "proof" of Corollary~6G requires {\it every} counterexample in
$SK_1(D)$ to have a preimage in $K_1(\Delta).$  Similary, Theorem~2F
is too strong: the first sentence in its "proof" requires
essentially the same thing.  This conjunction seems incredibly
unlikely.

\vs\noi
{\sc Lemma 6H.} {\it With $n=p^{r+t} + p^tm-1$ of type $t\ge 0,$ let
$\kappa  (n)=\max\{z\ge 0: n\ge p^{r+z-1},$ and $\ell(n) =\min\{j\ge
1: m\le p^r(p^j-1)\}.$  Then $\kappa  =t+\ell.$ }

\vs\noi
\begin{proof} We have $p^{r+(t+j)-1} \le n = p^{r+t} + p^tm -1\Llra
p^{r+t+j-1} <p^t(p^r+m) \Llra p^r(p^{j-1}-1) <m.$  It follows, since
$n\ge p^{r+t},$ that $\kappa  = t+j\Llra j$ is maximal with respect
to $p^{r +(t+j)-1} \le n\Llra j$ is maximal with respect
$p^r(p^{j-1}-1)<m\Llra j$ is minimal with respect to $p^r(p^j-1)
\ge m.$  Thus
$j=\ell.$
\end{proof}

\noi
{\sc Remark.} The above $\kappa  $ appears in Step~2 of the proof of 
Proposition 5F, but is not related to the $\kappa $ in the proof of
Lemma~2A(ii)1.

\section{Commutators and Type} \label{Sec7}

\vs\noi
We fix the usual notation $\Delta  =\Delta  _{r,s,\IF}$ of \S1.  We
also extend the notation $^gx$ for $x\in B$ to all $x\in D,$ by
writing $^g x$ for $uxu^{-1},$ a convenient abuse of notation below.
The new focus of Hypothesis~(6Gs)
is continued by searching for
commutators on levels of a single type, hence fixing $n= p^{r+t} +
p^tm-1,$ with $1\le m\not\equiv 0\!\!\!\mod p$ and $t\ge 0,$ as in
Lemma~6E.

\vs\noi
{\sc Lemma 7A.}
 {(i)} {\it If $1\le m\not\equiv 0\!\!\mod p$ and $0\le t\le
s-1,$ then, for $x\in \Delta  ,$}
$$
[u,1-\Pi^{p^tm}x] \equiv 1-\Pi^{p^tm}(\,^g x-x)(1-\Pi^{p^tm}x)^{-1}-
\Pi^{p^{r+t}+p^t(m-1)} m\;^g x\!\!\!\mod \Pi^{p^{r+t}+p^tm}\Delta  .
$$

\begin{description}
\item  {(ii)} {\it If $x\in\Delta  $ and $y\in A[u]$ 
 then the condition
$$
[u, 1-\Pi^{p^tm}x]\equiv 1-\Pi^n y\!\!\!\mod \Pi^{n+1}\Delta  
$$
is equivalent to the condition}
$$
\Pi^{p^{r+t}-p^t} m\;^gx + (\,^gx-x)\equiv \Pi^{p^{r+t}-1}
y\!\!\!\mod \Pi^{p^{r+t}}\Delta  .
$$

\item{(iii)} If $x\in \Delta$ then ${}^g x\equiv x\!\!\mod
\Pi^{p^r}\Delta.$
\end{description}

\begin{proof}
{\it For} (i): $\,^{g-1}\Pi^{p^tm} \equiv 1 +
m(1-\Pi^{p^t})\Pi^{p^{r+t}-p^t}\!\!\!\mod \Pi^{2(p^{r+t}-p^t)},$ by
 Lemma~5G(i) with $j=p^tm,$ implies $^{g-1}\Pi^{p^tm} \equiv 1 +
m\Pi^{p^{r+t}-p^t}\!\!\!\mod \Pi^{p^{r+t}}B,$ by $p^r\ge 2;$ used at
$\dot\equiv$ below.  Turning to our commutator, $\!\!\!\mod
\Pi^{p^{r+t}+p^tm}\Delta,$
$$
\begin{aligned}
{}  [u,1-\Pi^{p^tm}x] &= \,^g(1-\Pi^{p^tm}x) (1-\Pi^{p^tm} x)^{-1}\\
&= \big(1-\Pi^{p^tm}\cdot \,^{g-1}(\Pi^{p^tm})\cdot\,^gx
\big)(1-\Pi^{p^tm}x)^{-1}\\
&\dot\equiv\big(1-\Pi^{p^tm}(1+m\Pi^{p^{r+t}-p^t})\, ^g
x\big)(1-\Pi^{p^tm}x)^{-1}\\
&= (1-\Pi^{p^tm}\cdot \,^gx)(1-\Pi^{p^tm}x)^{-1}
-(\Pi^{p^t(m-1)+p^{r+t}}m\,^gx)(1-\Pi^{p^tm}x)^{-1}\\
&= (1-\Pi^{p^tm}x)^{-1} -\Pi^{p^tm}\cdot \,^gx(1-\Pi^{p^tm}x)^{-1}
-\Pi^{p^{r+t}+p^t(m-1)} m\,^gx(1-\Pi^{p^tm}x)^{-1}\\
&\ddot\equiv \big(1+\Pi^{p^tm}x(1-\Pi^{p^tm}x)^{-1}\big) - \Pi^{p^tm}\cdot
\,^gx(1-\Pi^{p^tm}x)^{-1}\\
&\q - \Pi^{p^{r+t}+p^t(m-1)}m\, ^gx - \Pi^{p^{r+t}+p^t(m-1)}m\;^g x
\cdot \Pi^{p^tm}x(1-\Pi^{p^tm}x)^{-1}\\
&= 1-\Pi^{p^tm}(\,^gx -x)(1-\Pi^{p^tm}x)^{-1} -
\Pi^{p^{r+t}+p^t(m-1)}m\;^gx-0\!\!\!\mod \Pi^{p^{r+t}+p^tm}\Delta  ,
\end{aligned}
$$
with $\dot=$ since $^{g-1}\Pi^{p^tm} - (1+m\Pi^{p^{r+t}})$ in
$\Pi^{p^{r+t}}\Delta$ implies that the difference between line~2 and
line~3 is indeed in $\Pi^{p^tm}(\Pi^{p^{r+t}}\Delta) = \Pi^{p^{r+t}
+p^tm}\Delta;$ 
and with $\ddot= $ since the first entry in line~5 is the
first entry in line~6, as $y = \Pi^{p^tm}x \in \Pi\Delta$ has
$(1-y)^{-1}= 1+y(1-y)^{-1},$ and the last entry in line~5 is the sum 
 in line~7 because $z=\Pi^{p^{r+t}+p^t{(m-1)}} m\;^gx$ has
$zy = z+zy(1-y)^{-1}$
with $zy(1-y)^{-1}\in \Pi^{p^{r+t}+p^t(m-1)}\Delta  \Pi^{p^tm}\Delta 
\subseteq \Pi^{p^{r+t}+p^tm}\Delta  .$

\vs\noi
{\it For} (ii): The first condition is, by (i), equivalent to
$$
1- \Pi^{p^{r+t}+p^t(m-1)} m\;^gx -
\Pi^{p^tm}(\,^gx-x)(1-\Pi^{p^tm}x)^{-1} \equiv 1-
\Pi^{p^{r+t}+p^tm-1} y\!\!\!\mod \Pi^{p^{r+t}+p^tm}\Delta  .
$$
Adding $-1$ and cancelling $-\Pi^{p^tm}$ on both sides of this turns
it into
$$
\Pi^{p^{r+t}-p^t} m\;^g x + (\,^gx-x)(1-\Pi^{p^tm}x)^{-1}
\equiv \Pi^{p^{r+t}-1}y\!\!\!\mod \Pi^{p^{r+t}}\Delta , 
$$
and conversely.  Next, right multiplying both sides by the unit
$(1-\Pi^{p^tm}x)$ turns this into
$$
\Pi^{p^{r+t}-p^t}m\;^gx(1-\Pi^{p^tm}x) + (\,^g x-x)\equiv
\Pi^{p^{r+t}-1} y (1-\Pi^{p^tm}x)\!\!\!\mod \Pi^{p^{r+t}}\Delta  ,
$$
which becomes
$$
\Pi^{p^{r+t}-p^t} m\,^gx + (\,^g x-x)\equiv \Pi^{p^{r+t}-1}y\!\!\!\mod
\Pi^{p^{r+t}} \Delta  ,
$$
because $\Pi^{p^{r+t}-p^t}\Delta  \Pi^{p^tm}\Delta  =
\Pi^{p^{r+t}+p^t(m-1)}\Delta  $ and $\Pi^{p^{r+t}-1} \Delta 
\Pi^{p^tm}\Delta  = \Pi^{p^{r+t}+p^tm-1} \Delta  $ are both
contained in $\Pi^{p^{r+t}}\Delta  ,$ and conversely.

\vs\noi
{\it For} (iii): Write $x= \os{p^s-1}{\us{j=0}\sum} \Pi^j x_j,$ with
$x_j\in A[u],$ by Lemma~1G(i), hence ${}^gx_j = ux_ju^{-1} = x_j$
for all $j.$  Since ${}^g\Pi \equiv \Pi\!\!\mod \Pi^{p^r}\Delta,$
by Lemma~2G(iii), it follows that ${}^gx= \sum_j {}^g\Pi^j x_j
\equiv \sum_j \Pi^jx_j = x\!\!\mod \Pi^{p^r}\Delta.$
\end{proof}

\vs\noi
{\sc Corollary 7B.}
{\it If $d\in \Delta  ^\times$ has level $n= p^{r+0} + p^0 m-1$ of type
$0,$ then $d\equiv [u,1-\Pi^m y/m]\!\!\!\mod \Pi^{n+1}\Delta  ,$ for
some $x\in \Delta  .$}

\vs\noi
\begin{proof}
Now $d\equiv 1-\Pi^n y\!\!\!\mod \Pi^{n+1}\Delta  ,$ with $y\in
A[u],$ by Lemma~1G(i).  It suffices, by  Lemma~7A(ii), to find $x\in\Delta 
$ with $\Pi^{p^r-1}m\;^gx+(\,^gx-x)\equiv \Pi^{p^r-1}y\!\!\!\mod
\Pi^{p^r}\Delta  ,$ i.e., to solve $\Pi^{p^r-1} m\;^gx \equiv
\Pi^{p^r-1}y\!\!\!\mod \Pi^{p^r}\Delta $ for $x,$ by 
 Lemma~7A(iii).  This amounts to $x\equiv\,u^{-1}yu/m =
y/m\!\!\!\mod \Pi\Delta,  $ as required. 
\end{proof}

\vs\noi
{\sc Remark.} 
When $s=2,$ Corollary~7B verifies Hypothesis (6G2) so Corollary~6G
implies $SK_1(D)=0.$
 In particular, we may assume $s\ge 3$
from now on.

\vs\noi
We continue from Lemma 7A(ii), with level $n$ of type $t\ge 1.$
Fixing $y\in A[u],$ we seek $x\in \Delta  $ so that

\noi
($*$)\q\q\q\q\q\q\q\q$ \Pi^{p^{r+t}-p^t}m\;^gx +(\,^gx-x)\equiv \Pi^{p^{r+t}-1}
y\!\!\!\mod \Pi^{p^{r+t}}\Delta  .$\q\q\q\q\q\q\q\q

\vs\noi
{\sc Proposition 7C.}
\begin{description}
\item  {(i)} {\it Define $f(j)\os{\rm def}= (p^r-1)p^{\v_p(j)} + j =
(p^r-1+j^\pr)p^{\v_p(j)} = f(j^\pr)p^{\v_p(j)}$ for $j\ge 1.$  Then}
\end{description}

\noi
\centerline{${}\q
[\,^g\Pi^j - \Pi^j] \equiv j^\pr
\big(1-\Pi^{p^{\v_p(j)}})\Pi^{f(j)} +
\big(\begin{matrix}j^\pr\\
2\end{matrix}\big)(1-2\Pi^{p^{v_p(j)}}\big)\Pi^{f(j)+(p^r-1)p^{v_p(j)}}
\!\!\!\mod
\Pi^{f(j)+(p^r+1)p^{\v_p(j)}}B,
$}

\begin{description}
\item{\q} {\it which implies $\v_\Pi(\,^g\Pi^j - \Pi^j) = f(j),$
because of the $\Pi^{f(j)}$ factor on the right side.}

\item{(ii)} {\it Write $x= \sum_{0\le j <p^s} \Pi^j x_j$ with unique
$x_j\in A[u].$ Then congruence $(*)$ becomes
$$
\Pi^{p^{r+t}-p^t} \sum_{0\le i<p^t} \,^g\Pi^i mx_i +
\sum_{j\in \CL} [\,^g\Pi^j - \Pi^j]x_j \equiv
\Pi^{p^{r+t}-1} y\!\!\!\mod \Pi ^{p^{r+t}}\Delta  
$$
with $\CL =\{1\le j<p^s: \v_p(j)<t$ and $j^\pr
<p^r(p^{t-\v_p(j)}-1)\}= \underset{0\le \v<t}{\dot\bigcup}\CL_\v,$
where $\CL_\v=\{1\le j<p^s:\v_p(j) =\v$ and $j^\pr <\CV_\v\} $ with
$\CV_\v = p^r(p^{t-\v}-1).$  Note that $0\not\in \CL,$ and the $i$
sum and $j$ sum are not disjoint: $i=j$ still implies $x_i = x_j.$}

\item  {(iii)} {\it If $x= \sum_{0\le j<p^s} \Pi^j x_j,$ as above, solves the congruence
of}
(ii), {\it with $y\in A[u],$ then $x$ also solves
$$
\sum_{j\in L} [\,^g \Pi^j - \Pi^j]x_j \equiv 0\!\!\!\mod
\Pi^{p^{r+t}-p^t}\Delta  ,
$$
with $L=  \big\{1\le j<p^s: \v_p(j) <t\; 
and
\; j^\pr <
(p^r-1)\big(p^{t-\v_p(j)}-1\big)\big\}=\underset{0\le
\v<t}{\dot\bigcup} L_\v,$ where $L_\v =\{1\le j<p^s:\v_p(j)=\v$ and
$j^\pr < V_\v\}$ with
$V_v = (p^r-1)(p^{t-v}-1).$  Thus $L\subseteq \CL.$ }
\end{description}

\vs\noi
{\it Proof.}

\noi
{\it For} (i): Multiply the Lemma 5G(i) congruence by $\Pi^j,$ subtract  $\Pi^j$ 
 and keep track of $f(j).$

\noi
{\it For} (ii):  Substituting $x=\sum_j \Pi^j x_j,$ from Lemma~1G(i), into $(*)$ 
clarifies
the first sum 

\noi
$\Pi^{p^{r+t}-p^t} m\;^g x = \Pi^{p^{r+t}-p^t} \sum_j \, ^g\Pi^j
mx_j$ and the terms with $j\ge p^t$ redundant $\!\!\!\mod
\Pi^{p^{r+t}}\Delta  $ as $\v_\Pi(\Pi^{p^{r+t}-p^t}\cdot \,^g\Pi^j
m) = p^{r+t}-p^t + j\ge p^{r+t}.$  
Renaming the $j$'s in the first sum as $i$'s makes them easier to
refer to later.  Special care is needed here.

\vs\noi
Similarly, $^gx-x =\sum_j[\,^g\Pi^j - \Pi^j]x_j$ has $j=0$ term equal
to zero, eliminating this term,  and redundancy eliminates all but the claimed terms. 
Namely, $j\ge 1$ is redundant if, and only if, $f(j) \ge p^{r+t}$
leaving, by (i), 
only those $j$ with $(p^r-1+j^\pr)p^{\v_p(j)} < p^{r+t}.$
Since $0\le j^\pr-1,$ this implies $p^{r+\v_p(j)} <p^{r+t},$ hence 
$\v_p(j)<t.$  Thus $p^r-1+j^\pr \le p^{r+(t-\v_p(j))} -1,$ hence
$0<j^\pr <p^r\big(p^{t-\v_p(j)}-1\big),$ as
$j^\pr\not\equiv 0\!\!\!\mod p.$

\noi
{\it For} (iii): Every solution $x$ of $(*)$ has 
$\us{j\in\CL}\sum  [\,^g\Pi^j - \Pi^j]x_j\equiv 0\!\!\!\mod
\Pi^{p^{r+t}-p^t}\Delta ,$ so it remains to eliminate the redundant
terms, i.e., those $j\in \CL$ with $f(j)\ge (p^r-1)p^t.$  
This leaves
 $j$ with $(p^r-1 + j^\pr)p^{\v_p(j)} <(p^r-1)p^t \Llra
j^\pr< V_{\v_p(j)},$ i.e., 
$j\in L_{\v_p(j)}\,.$ 
\hfill $\square$

\vs\noi
{\sc Corollary 7D.} {\it If $d\in \Delta  ^\times $ has level
$n=p^{r+1}+p^1m-1$ of type $1,$ then}
$$
d\equiv [u,1-\Pi^{pm}x]\!\!\!\mod \Pi^{n+1}\Delta  , \;
\text{\it for some}\; x\in \Delta  .
$$

\vs\noi
\begin{proof}
Just as in the proof of Corollary 7B, it suffices to find, for each
$y\in A[u],$ an

\noi
$x=\sum_{j} \Pi^j x_j$ in $\Delta  ,$ with $x_j \in
A[u]$ for all $j,$ solving
$$
\Pi^{p^{r+1}-p}\sum\nolimits_{0\le i<p} \,^g\Pi^i mx_i +
\sum\nolimits_{j\in \CL} [\,^g\Pi^j - \Pi^j]x_j \equiv
\Pi^{p^{r+1}-1} y\!\!\!\mod \Pi^{p^{r+1}}\Delta  , 
\leqno{\rm (7D1)}
$$
with $\CL =\{1\le j <p^s:\v_p(j)=0$ and $j^\pr <p^r(p-1)\},$
as in Proposition~7C(ii).  To find such a solution $x,$ first
observe that $x$ must also solve, by Proposition~7C(iii),
$$
\sum\nolimits_{j\in L} [\,^g\Pi^j - \Pi^j]x_j \equiv 0\!\!\!\mod
\Pi^{p^{r+1}-p} \Delta   \leqno{\rm (7D2)}
$$
with $L =\{1\le j <p^s:\v_p(j) =0$ and $j^\pr <V\}$ and $V=V_0=
(p^r-1)(p-1).$

\vs\noi
There are two cases, the first of which delayed to Corollary~7D$^*,$
following Corollary~7D.  The
second case is

\vs\noi
{\it Step} 1: If $s\le r,$ and $x$ solves (7D2), then
$$
\big\lfloor \frac{p-1+j}{p^s}\big\rfloor + v_\pi(x_j) \ge
p^{r-s}(p-1), \;{\rm for\, all}\; j\in L_0.
$$

\noi
First $v_\bullet ([\;^g\Pi^j - \Pi^j]x_j) = v_\bullet (\,^g\Pi^j -
\Pi^j) + v_\bullet (x_j) \dot = v_\Pi(\,^g\Pi^j -\Pi^j) +
p^sv_\pi(x_j) \ddot = f(j) + p^sv_\pi(x_j) = p^r-1+j +
p^sv_\pi(x_j),$ for all $j\in L = L_0,$ with $\dot =$ by
Proposition~7C(i) and $\ddot =$ by Corollary~2B(ii) and
Corollary~2B(iii).
Thus $v_\bullet([\,^g\Pi^j -\Pi^j]x_j) \equiv
p^r -1 +j\!\!\!\mod p^s,$  so if $j_1,j_2$ in $L_0$ give the same
value here then $j_1\equiv j_2\!\!\!\mod p^s,$ hence $j_1=j_2$ and

\noi
\centerline{$\min_{j\in L_0}
\,\v_\bullet([\,^g\Pi^j - \Pi^j]x_j) = \v_\bullet(\sum_{j\in L_0}[\,^g\Pi^j -
\Pi^j]x_j)\ge p^{r+1} -p,$}

\noi
by (7D2), so that $p^r -1 +j+p^s\v_\pi(x_j)\ge
p^{r+1}-p$ for some, hence  all,  $j\in L_0.$ Rearranging this slightly yields
$\frac{p-1+j}{p^s}\ge p^{r-s}(p-1) - v_\pi(x_j)\in \Z,$ giving {\it
Step}~1.

\vs\noi
{\it Step} 2: This has consequences,  the first of which we 
list,  to display the needed notation.
$$
\begin{aligned}
\text{(2.1)} &\q {\rm If}\q j\le p^s-p, \q{\rm then}\q \v_\bullet\big([\,^g\Pi^j -
\Pi^j]x_j\big) \ge p^{r+1};\\
&\q {\rm if}\q 1\le i <p,\;\q
 {\rm then}\q \v_\bullet\big(\Pi^{p^{r+1}-p}\;
^g\Pi^i mx_i\big) 
\ge p^{r+1}.\\
\text{(2.2)} &\q {\rm If}\; j(k)= p^s -
p +k \q\text{with}\q 1\le k\le p-1,\\
&\q {\rm then}\q x_{j(k)} = \pi^h z_k\;{\rm with}\;
h=p^{r-s}(p-1)-1,\q {\rm and}\q z_k\in A[u].\\
\text{(2.3)} & \q 
{\rm With}\q \eta \in B\q \text{defined by}\q\pi = \Pi^{p^s}\eta 
,\q{\rm then}\q \eta  \equiv 1\!\!\!\mod \Pi^{2(p^r-1)}.
\q\q\q\q\q\q\q\q{}
\end{aligned}
$$
{\it For} (2.1): By {\it Step} 1, $j$ has $\lfloor
\frac{p-1+j}{p^s}\rfloor =0,$ so the $\v_{\bullet}$ of interest is
$p^r-1+j +\v_\bullet(x_j)
\ge p^r-1+j + p^sp^{r-s}(p-1) = p^{r+1}+j -1 \ge p^{r+1}.$

\noi
Similarly, $i$ has $i<p\le p^s-p,$ as $s\ge 2$ by the Remark
after Corollary~7B, hence $\lfloor \frac{p-1+j}{p^s} \rfloor =0.$
So the $\v_{\bullet}$ of interest is
$\ge p^{r+1} -p + i + p^s p^{r-s}(p-1) = p^{r+1} +(p^r-1)(p-1) +
(i-1)\ge p^{r+1}.$

\noi
{\it For} (2.2): Now $\lfloor\frac{p-1+j(k)}{p^s}\rfloor = 1,$ hence
$v_\pi(x_{j(k)})\ge h,$ by {\it Step} 1; 
writing $x_{j(k)} = \pi^l c,$ with $c\in A[u]$ and maximal
$l,$ implies that  $v_\pi(c)=0,$ by the
final paragraph of the proof of Corollary~
2B(iii).  Thus $l=
v_\pi(x_{j(k)}) \ge h$ and $x_{j(k)} = \pi^h(\pi^{l-h}c)$ implies
 $z_k = \pi^{l-h}c\in A[u].$

\noi
{\it For} (2.3): Putting $M= p^r-1$ and using
Lemma~5G(iii) with $q$ there $=t$ here $=0,$ we have $\eta  =\us{0\le
h<p^s}\prod\,^{g^h-1}\Pi \equiv \prod_h 
\big(1+h(1-\Pi)\Pi^M\big) \equiv 1+\sum_h h(1-\Pi)\Pi^M=1 +
\frac{p^s(p^s-1)}{2}\,(1-\Pi)\Pi^M\equiv
1\!\!\!\mod \Pi^{2M}B,$ with  $\v_\bullet
\big(\frac{p^s(p^s-1)}{2}\,(1-\Pi)\Pi^M\big) = s\v_\bullet(p)+M =
s\phi(p^{r+s}) + p^r-1 >2(p^r-1),$
completing (2.3) and {Step}~2.

\vs \noi
Next we turn to solutions $x$ of (7D1), and use Step~2 to
eliminate many of its redundant terms from (7D1).
First  (2.1) allows us to
ignore all $j$ in the $j$-sum except for $j= j(k)= p^s-p +k$ with $0
<k<p,$
and also all $i$ in the $i$-sum except $i= 0,$ leaving
$$
\Pi^{p^{r+1}-p}
mx_0 + \sum_{0<k<p} [\,^g\Pi^{j(k)} - \Pi^{j(k)}] x_{j(k)}\equiv \Pi^{p^{r+1}-1}
y \!\!\!\mod \Pi^{p^{r+1}}\Delta  .
\leqno{\rm (7D3)}
$$
Next, Proposition~7C(i) for
$j(k)\equiv k\!\!\!\mod p,$ hence $\v_p\big(j(k)\big) =0,$ gives the
starting 

\centerline{$[\,^g\Pi^{j(k)} - \Pi^{j(k)}]\equiv j(k)^\pr(1-\Pi)\Pi^{p^r-
1+j(k)}
+\Big(\begin{matrix} j(k)^\pr\\
2\end{matrix}\Big)(1-2\Pi)\Pi^{2(p^r-1)+j(k)}
\!\!\!\mod \Pi^{2p^r+j(k)}.$}

\noi
Right multiplying this by
$x_{j(k)}
=
\pi^h z_k = \Pi^{p^sh}\eta  ^hz_k,$ from (2.2) and (2.3), turns it into
$$
\begin{aligned}
{}[^g\Pi^{j(k)} - \Pi^{j(k)}] x_{j(k)}  &\equiv k(1-\Pi)\Pi^{p^{r+1}-p+k-1}\eta^h
z_k\\
&\q +\Big(\begin{matrix} j(k)^\pr\\ 2\end{matrix}\Big)
(1-2\Pi)\Pi^{p^{r+1}+(p^r-p)+(k-2)}n^hz_k \!\!\!\mod
\Pi^{p^{r+1}+p^r-p+k}
\end{aligned}
$$
since $p^s(h+1)=p^{r+1}-p^r$ implies $p^r -1+j(k)  + p^sh = p^r-1 +p^s
-p +k + p^sh = p^r-p+p^s(h+1) + (k-1) = (p^{r+1}-p) + (k-1),$ and
$2(p^r -1) + j(k) + p^sh = (p^r-1) + p^{r+1} -p + k - 1= p^{r+1} +
(p^r-p) + (k-2)\ge p^{r+1}$ unless $r=1$ (and $k=1),$  
when
 
\noi
$\big(\begin{matrix} j(k)^\pr\\ 2\end{matrix}\big)\equiv
0\!\!\!\mod p$  (with $\v_\bullet(p)\ge p^{r+1})$ implies that we
always have
$$
[^g\Pi^{j(k)} - \Pi^{j(k)}] x_{j(k)} \equiv
k(1-\Pi)\Pi^{p^{r+1}-p+k-1}\eta^h z_k\!\!\!\mod \Pi^{p^{r+1}}\Delta  .
$$
Substituting this into (7D3), for every $k,$ gives
$$
mx_0 +(1-\Pi)\eta  ^h \sum^{p-1}_{k=1} \Pi^{k-1} kz_k\equiv
\Pi^{p-1}y\!\!\!\mod \Pi^p\Delta  
\leqno{\rm (7D4)}
$$
after left cancelling the common $\Pi^{p^{r+1}-p}.$

\noi
Now, (2.3) removes $\eta ^h,$ so multiplying by $(1-\Pi)^{-1}$ and re-arranging gives
$$
\sum^{p-1}_{k=1} \Pi^{k-1} (mx_0 + kz_k) + \Pi^{p-1}(mx_0 -y) \equiv
0\!\!\!\mod \Pi^p \Delta  .
$$
It follows that $mx_0 + kz_k \equiv 0\!\mod\pi A[u]$ for $1\le
k\le p-1:$ for considering the smallest $k$ for which this is false
 mod $\Pi^k\Delta$ implies that $mx_0 + kz_k \in \Pi\Delta \cap
A[u] = \pi A[u],$ by (2B1), contradiction.  Continuing
 reveals the solution
$x_0\equiv y/m$ and
 $z_k \equiv -y/k\!\!\!\mod \pi A[u]$ for $1\le k\le p-1.$ 
\end{proof}

\vs\noi
{\sc Remark:} When $s=3,$ Corollary 7D verifies Hypothesis (6Gs) so
Corollary~6G implies $SK_1(D)= 0.$
In particular, we may assume $s\ge 4$ from now on.

\vs\noi
{\sc Remark.}
  Lemma 7A studies a single commutator $[u,
1-\Pi^{p^tm}x],$ for  levels $n$ of type $t,$ by finding a
linearized form of it $\!\!\!\!\mod \Pi^{n+1}\Delta  $ in (ii), resulting in
the congruence $(*)\!\mod \Pi^{p^{r+t}}\Delta  $ just before
Proposition~7C.  Proposition~7C(iii) is a reduction step exhibiting a
rather smaller space containing all solutions of $(*),$ which
permits a simplified form of the Proposition~7C(ii) congruence,
equivalent to $(*),$ to be found.  In the Corollary~7D study of type
$t=1,$ this simplified form (7D4)  can be completely solved.  

\vs\noi
{\sc Corollary 7D$^*.$} 
{\it If $d\in \Delta^\times$ has level $n=p^{r+1} + p^1m-1$ of
type~1,then $d\equiv [u,1-\Pi^{pm}x]\!\!\!\mod \Pi^{n+1}\Delta,$
for some $x\in \Delta.$ }

\vs\noi
\begin{proof}
This is the first case of Corollary 7D, which deals
with the second case.  We start with the identical first paragraph as
in Corollary~7D, based on (7D1) and (7D2).

\vs\noi
{\it Step} 1: If $s\ge r+1,$ and $x$ solves (7D2), then
$x_j\equiv 0\!\!\!\mod \pi A[u]$ for all $j\in L\,.$

\noi
For suppose $x$ is a counterexample,  and recall the notation of
Corollary~2B.  Then
$L(x) \os{\rm def}=
\{j\in L:\v_\bullet(x_j) =0\} \ne \emptyset:$  for every $j\in L\bs
L(x)$ is redundant in (7D2), since $\v_\bullet(x_j)>0$ implies
$\v_\bullet(x_j) \ge
p^s \ge p^{r+1}-p,$ by Corollary 2B(iii).
Now  the $j\in L(x)$ have 

\centerline{$
\v_\bullet([\,^g\Pi^j - \Pi^j]x_j)
= f(j) +0 = p^r-1+j,
$}

\noi
which are all different.  Thus
$\v_\bullet(\sum_{j\in L(x)}[\,^g\Pi^j - \Pi^j]x_j) = \min_{j\in L(x)}(p^r-1+j) =
p^r-1+j_*,$ for some $j_* \in L,$ which has 
$p^r-1+j_*<p^r-1+V= p^{r+1}-p,$  contradicting (7D2).  This settles
Step~1.

\vs\noi
Step 1\; implies that \q $\CL\bs L =\{1\le j<p^s:\v_p(j)=0$ \; and \; $V\le
j^\pr <p^r(p-1)\}\; =$

\noi
$
\{1\le j<p^s:j=j^\pr$ and $V\le j^\pr <V-1+p \}=\{V-1+k:0<k<p\},$
since $V-1+k\equiv k\!\!\!\mod p;$ the $i$ sum has only $i=0$ to
contribute (as $0\not\in \CL$ and $V-1+k\ge V>p-1).$ 
Thus (7D1) becomes
$$
\Pi^{p^{r+1}-p} mx_0 +\sum^{p-1}_{k=1} [\,^g\Pi^{V-1+k} -
\Pi^{V-1+k}] x_{V-1+k} \equiv \Pi^{p^{r+1}-1} y\!\!\!\mod
\Pi^{p^{r+1}} .
\leqno{\rm (7D3)}
$$
Now Proposition~7C(i), with $j=V-1+k\equiv k\!\!\mod p$ and $p^r-1+V =
p^{r+1}-p,$ implies $$
[^g\Pi^j - \Pi^j] \equiv j^\pr(1-\Pi)\Pi^{\wt k} +
\big(\begin{matrix} j^\pr\\ 2\end{matrix}\big)(1-2\Pi)\Pi^{\wt
k+p^r-1}
\!\!\!\mod
\Pi^{\wt k+p^r+1},
$$
with $\wt k = p^{r+1}-p +k-1.$  Here $\wt k+p^r-1\ge p^{r+1} (\Llra
p^r-p+k-2\ge 0)$ holds unless $r=1$ (and $k=1),$ when $\big(\begin{matrix}
j^\pr\\ 2\end{matrix} \big) \equiv 0\!\!\!\mod p$ (with $\v_\bullet(p) =
\phi(p^{r+s}) \ge p^{r+1})$
and $\wt k+ p^r +1>
p^{r+1},$ so we always have

\centerline{$[^g\Pi^j- \Pi^j ]\equiv k(1-\Pi)\Pi^{\wt k}
\!\!\!\mod \Pi^{p^{r+1}},$}

\noi
hence 
right multiplying by
$x_{V-1+k},$  summing over $k,$ and substituting into (7D3) gives 
$$
mx_0 + (1-\Pi)\sum_k 
k\Pi^{k-1}x_{V-1+k}\equiv \Pi^{p-1}y\!\!\!\mod \Pi^p\Delta,  
\leqno{\rm (7D4)}
$$

\noi
after
left cancelling the common $\Pi^{p^{r+1}-p}$ in $\Pi^{\wt k}.$  Next left multiplying
this congruence by the unit $(1-\Pi)^{-1}\equiv \sum^p_{k=1}
\Pi^{k-1}\!\!\!\mod \Pi^p$ and rearranging a little results in
$$
\sum^{p-1}_{k=1} \Pi^{k-1}(mx_0 + kx_{V-1+k}) + \Pi^{p-1}(mx_0-y)
\equiv 0\!\!\!\mod \Pi^p\Delta.
$$

\noi
It follows that $mx_0 + kx_{V-1+k}\equiv 0\!\!\!\mod \pi A[u]$ for
$1\le k\le p-1:$ for considering the smallest 
$k$ for which this is false  {\rm modulo} $\Pi^k\Delta  $
implies that $mx_0 + kx_{V-1+k}\in
\Pi\Delta  \cap A[u]\subseteq \pi A[u],$ by (2B1), contradiction.

\noi
 Continuing reveals that
$ x_0 \equiv y/m$ and $ x_{V-1+k} \equiv -mx_0/k \equiv -y/k \!\!\!\mod
\pi A[u],$  for 

\noi
$1\le k\le p-1,$ which
 settles the first case, since all of these moves are reversible.
\end{proof}

\vs\noi
{\sc Lemma 7E.}  {\it The valuation $v_\pi$ of Corollary}
2B(iii) {\it has} 
$$
v_\pi(x) = v^\bullet\big(N_{\fq(A[u])/\fq(A)}(x)\big),\;
\text{\it for} \;x\in A[u].
$$

\vs\noi
\begin{proof}
{\it Step 1:} Fix the valuation $v_\bullet :
D\to\Z\cup \{\infty  \}$ of $D.$  Then each subfield $L$ of $D,$
containing the centre $Z(D)$ of $D,$ has a unique valuation $w$ with
$w\vert  _{Z(D)} = v_\bullet \vert  _{Z(D)}.$  It has 
$$
w(x) =
(v_\bullet \vert 
_{Z(D)})\big(N_{L/Z(D)}(x)\big)/\,[L:Z(D)],\;\text{\rm
for all}\;
x\in L.
$$

\vs\noi
The unique existence of $w$ follows from Theorem~1.4 of ([TW],
1.2.2), and the formula for $w$ is given by Lemma~1.6 of loc.cit.. 
Note that the extension convention of [TW] is  that of
restriction of  functions (which differs from our use of
normalization).  This gives Step~1.  

\vs\noi
 We revert to our usual $\Delta  =\Delta  _{r,s,\IF}$
notation and $v_\bullet :D\to \Z\cup \{\infty  \},$ hence $Z(D) =
\fq(A).$  We apply Step~1 to $L=\fq (A[u]),$ hence $v_\bullet \vert 
_{\fq(A)} = p^sv^\bullet,$ by Corollary~2B(iv), to find that the
valuation $w$ on $\fq(A[u])$ has
$w(x)
 = N_{\fq(A[u])/q(A)}(x).$
\end{proof}

\section{Toward the 7C(iii) Congruence} \label{Sec8}

Fixing $\Delta  =\Delta  _{r,s,\IF},$ an understanding of the
solution set of the title congruence
$$
\sum_{j\in L}\, [^g\Pi^j - \Pi^j] x_j \equiv 0\!\!\!\mod \Pi^{p^{r+t}-p^t}
\Delta  ,
$$
with all $x_j\in A[u],$ was the initial obstacle in both cases of
the proof of
Corollary~7D, and $L$ becomes more
complicated as $t$ grows. Recall, from Proposition~7C(iii), that,
for $0\le \v <t,$

\vs\noi \centerline{$L_\v =\{1\le j<p^s: \v_p(j) =\v \q\text{\rm  and} \q j^\pr
<V_{\v}\} \q\text{\rm with} \q V_\v = (p^r -1)(p^{t-\v}-1).
$}

\noi

\vs\noi
This extra complexity seems to require more suitable
valuation arguments when $t\ge 2.$  
Our approach to the $L$'s is intended to be susceptible to induction
on $t\ge 2.$
The key to this is the equivalence relation on $L$ induced by the
$f:L\to\Z$ of Proposition~7C(i), which is described in Lemma~8A,
especially in (ii) and (iv).

\noi
We continue with the conventions of the first
paragraph of \S 7, and also
assume that $s\ge 4$ (by the first Remark 
after Corollary~7D), and that $t\ge 2$ from now on.

\vs\noi
{\sc Lemma 8A.} (i)
{\it $f: L\to \Z,$ defined in Prop.~7C(i), induces an
equivalence relation $\sim_f$ on }
{}\q  {\it $L,$ i.e., $j_1\sim _f j_2$ if, and
only if, $f(j_1) = f(j_2).$} 

\begin{description}
\item  {(ii)} {\it There is a  function $\tau:L\to L,$
so that $f\big(\tau  (j)\big)=f(j),$ for all $j\in L,$
defined by 
$$
\tau  (j) = \begin{cases}
(p^r-1)\big(p^{\v_p(j)} - p^{\v_p(j)-1}\big) + j, &\q\text{\it if} \q j\in L\bs
L_0\\
j, &\q\text{\it if} \q j\in L_0. 
\end{cases}
$$
This has
$$
\tau  (j)^\pr = \begin{cases}
(p^r-1)(p-1) + pj^\pr, &\q\text{\it if} \q j\in L\bs L_0\\
 j^\pr  , &\q\text{\it if} \q j\in L_0.
\end{cases}
$$
Also $\tau  (j)^\pr
\equiv 1\!\!\!\mod p,$
if $j\in L\bs L_0.$}

\item  {(iii)} {\it With $\tau^k:L\to L$ denoting the composition of
$k\,\tau$'s including $\tau  ^0=\id_L,$ then every} $j\in L\bs L_0$
{\it has}
$\tau  ^k(j)\; = (p^r-1)(p^{\v_p(j)} - p^{\v_p(j)-k}) + j$  {\it and}
$ \tau  ^k(j)^\pr = (p^r  -1)(p^k-1) + p^kj^\pr,$ \,{\it for} $ 0\le
k\le \v_p(j),$ {\it and} $\tau  ^k(j)= (p^r-1)(p^{\v _p{(j)}}-1) + j$
{\it for} $k\ge \v_p(j).$  {\it Also every} $j\in L_0$ {\it has}
$\tau  ^k(j)=j$ {\it for} $k\ge 0.$

 \item  {(iv)} {\it The $f$-equivalence class with least element $l,$ in
the ordering of $L\subseteq \Z,$ is given by
$[l]_f = \{\tau  ^k(l): 0\le k \le \v_p(l)\}.$  Note
that  $\tau  (j)>j$ and $\v_p\big(\tau  (j)\big)
<\v_p(j)$ for  $j\notin L_0\,.$
\newline Note also that $\v_p(j)\le s-1$ for all $j\in L.$}

\item  {(v)} {\it
The  $f$ of Proposition~7C(i) has $f(L) = p^r-1+L_0.$  It follows
that
$L/\!\!\!\sim\; \to L_0,$ by $[l]\mapsto f(l) - (p^r-1),$ is bijective.}  
\end{description}

\vs\noi
{\it Proof.}
(i) is clear from Proposition~7C(i).

\noi
(ii) Note that $(p^r-1)p^{\v_p(j)} + j = f(j) =
f\big(\tau  (j)\big) = (p^r-1)p^{\v_p(\tau  (j))} + \tau 
(j),$  is equivalent to
 $$
\tau  (j) = (p^r - 1)(p^{\v_p(j)} - p^{\v_p(\tau(j))}) + j, \q\text{for} \q
j\in L.
\leqno{\rm (8A0)}
$$
We may regard (8A0) as  defining  $\tau$ provided that
$\v_p\big(\tau  (j)\big)< \v_p(j),$ for $j\in L\bs L_0,$ and
$\v_p\big(\tau  (j)\big) = \v_p(j),$ for $j\in L_0.$
For $\v_p\big((p^r-1)(p^{\v_p(j)} - p^{\v_p(\tau  (j))}) + j\big) 
= \min \big(\v_p(p^{\v_p(j)} - p^{\v_p(\tau  (j)}, 
\v_p(j)\big) = \min\big(\v_p(\tau  (j)), \v_p(j)\big) 
= \v_p\big(\tau  (j)\big)$ for $j\in L\bs L_0,$ and
$\v_p(0+j) = $

\noi
$\v_p(j) = \v_p\big(\tau  (j)\big)$ for 
$j\in L_0,$ do not contradict (8A0).
 It follows that
$$
\v_p\big(\tau  (j)\big) = \begin{cases}\v_p(j)-1, &\q j\in L\bs L_0\\
\v_p(j), &\q j\in L_0
\end{cases}\leqno{\rm (8A1)}
$$
and (8A0) permits the definition of  our $\tau .$

\noi
If $j\in L \bs L_0$ then $1\le \v_p(j)<t$ and $j^\pr
<V_{\v_p(j)},$ hence $\tau (j)  
 = p^{\v_p(j)-1} [(p^r-1)(p-1) +
pj^\pr]$ 
implies our formula for $\tau  (j)^\pr $ which implies
$\tau  (j)^\pr \equiv 1\!\!\!\mod
 p.$

\noi
Continuing $\tau  (j)^\pr < (p^r-1)(p-1) + pV_{\v_p(j)} = (p^r-1)(p-1) +
p(p^r-1)(p^{t-\v_p(j)} -1) = (p^r-1)\cdot$
$(p-1 + p^{t-\v_p(j)+1} - p) = (p^r-1)(p^{t-\v_p(\tau(j))} -1) =
V_{\v_p(\tau  (j))},$ i.e. $\tau  (j)\in L_{\tau  (j)}.$  Thus $\tau 
(j)\in L,$ i.e. $\tau  :L\to L$
is well-defined.  If $j\in L_0$  this is clear.

\vs\noi
(iii)  This is clear for $k=0,$ and we 
proceed by induction for $1\le k< \v_p(j),$ starting from~(ii)
for $k=1.$  If $j\in L\bs L_0$ and $2\le k\le \v_p(j)$ then
$$
\begin{aligned}
\tau  ^k(j) &= \tau  ^{k-1}\big(\tau  (j)\big) 
= (p^r-1)(p^{v_p(\tau  (j))} - p^{\v_p(\tau 
(j))-(k-1)}) + \tau  (j)\\
&\;\dot =(p^r-1) [(p^{\v_p(j)-1} - p^{\v_p(j)-k}) 
+ (p^{\v_p(j)} - p^{\v_p(j)-1})] + j\\
&= (p^r-1)(p^{\v_p(j)} - p^{\v_p(j)-k}) + j,
\end{aligned}
$$
with $ \dot=
$ by (8A1),
which also implies the $\tau  ^k(j)^\pr$ formula for $k\le \v_p
(j).$  If $j\in L\bs L_0$ and $k\ge \v_p(j)$ then the argument is
similar.  If $j\in L_0,$ it's easier.

\vs\noi
(iv) (ii) implies `$\supseteq$'. Conversely, if $j\sim _f l,$ then
$0\le \v_p(j)\le \v_p(l),$ since the maximum value of $\v_p$ on
$[l]_f$ is $v_p(l),$ 
hence $\v_p(j) =\v_p(l)-k$ with $0\le k\le \v_p(l).$ 
So $j
= f(j) - (p^r-1)p^{\v_p(j)} 
= f(l) - (p^r-1)p^{\v_p(j)}
=(p^r-1)p^{\v_p(l)}+l-(p^r-1)p^{\v_p(j)}
= (p^r-1)(p^{\v_p(l)}- p^{\v_p(j)}) + l
= (p^r-1)(p^{\v_p(l)}- (p^{\v_p(l) - k}) + l
= \tau  ^k(l), $
by (iii).

\vs\noi
 Also $L\subseteq \{1,2,\dots, p^s-1\}$ implies $\v_p(j)\le
s-1.$

\vs\noi
(v) First,  $\tau  ^{\v_p(j)} (j) \in L_0$ for all $j\in L.$
If $j\in L\bs L_0,$ then $\tau  ^{\v_p(j)}(j) \dot=
(p^r-1)(p^{\v_p(j)}-1) + j,$ with $\dot = $ by (iii), has
$\v_p\big(\tau  ^{\v_p(j)}(j)\big) = 0$ since
$\v_p\big ((p^r-1)(p^{\v_p(j)} -1)\big) =0,$ $\v_p(j)>0,$ and
$\min\big(0,\v_p(j)\big) =0.$
If $j\in L_0,$ then $\tau  ^{\v_p(j)}(j) = j,$ and $\v_p(j)=0.$

\vs\noi
Next, turning to $f(L) = p^r-1+L_0,$ every $j\in L$ has 
$f(j) \dot =
f\big(\tau  ^{\v_p(j)}(j)\big)
= (p^r-1)p^{\v_p(\tau  ^{\v_p(j)}(j))}
+ \tau  ^{\v_p(j)}(j) 
\in p^r-1+L_0,$ with $\dot = $ by (i).  And any
$j\in L_0$ has $p^r- 1+j=f(j)\in f(L).$

\vs\noi
Finally, our map $L/\!\! \sim\,\, \to L_0$ is well-defined by (i) and $f(L) =
p^r-1+L_0,$ surjective for the same reason, and injective: if $[l],
[l^\pr]$ have $p^r-1 + f(l^\pr) = p^r-1+f(l),$ then $f(l^\pr) =
f(l),$ hence $[l^\pr] = [l].$ 
~\qed

\vs\noi
{\sc Proposition 8B.} 
  {\it Decomposing $L$ into $\sim_f$ -equivalence classes
$[l]_f,$ with $l$ chosen as the
least element of $[l]_f$ (so with
largest $\v_p(l)),$ define, for $x=\us {0\le j <p^s}\sum \Pi^jx_j$ as in
Proposition}~7C,

 \vs\noi
$ E_{[l]} (x) =\us{j\in [l]}\sum [^g\Pi^j- \Pi^j]x_j\in \Delta,$
\,
$e_{[l]}(x) = \us{j\in [l]}\sum j^\pr x_j \in A[u]
$ {\it and} $[l]_x=\{j\in[l]:\v_\bullet (x_j) =0\}.$

\vs\noi
{\it Then} 
{\rm (i)}
{\it the congruence of Proposition{\rm ~7C(iii)} becomes }
$$
\sum\nolimits_{[l]\in L/\sim} E_{[l]}(x) \equiv 0\!\!\!\mod
\Pi^{p^{r+t}-p^t}\Delta  .
$$

\begin{description}
\item  {\rm (ii)} 
$
\Pi^{-f(l)} E_{[l]}(x) \equiv e_{[l]} (x)\!\!\!\mod \Pi\Delta  .
$
{\it In particular, $\v_\bullet \big(E_{[l]}(x)\big) = f(l)$ if $e_{[l]}(x)
\not\equiv 0\!\!\!\mod \pi A[u],$ and $\v_\bullet\big(E_{[l]}(x)\big)
>f(l)$ if $e_{[l]}(x) \equiv 0\!\!\!\mod \pi A[u].$  }

\item  {\rm (iii)} {\it If\; $[l]_x
 \ne \emptyset,$ let $J=J_{l,x}$ be the 
largest $j$ in $[l]_x.$  
   Then}
$$
\Pi^{-f(l)} E_{[l]}(x) - e_{[l]}(x)\equiv -J^\pr \Pi^{p^{\v_p(J)}}
x_J +\Big(\begin{matrix} J^\pr\\ 2\end{matrix}\Big)
\Pi^{(p^r-1)p^{\v_p(J)}}x_J 
\!\!\!\mod \Pi^{p^{\v_p(J)+1}}\Delta  .
$$
{\it 
In particular, 
if $e_{[l]}(x) \equiv 0\!\!\mod \pi A[u]$
and $[l]_x \ne \emptyset$
then
$\v_\bullet\big(E_{[l]}(x)\big) =f(l) + p^{\v_p(J)}.$}

\item{\rm (iv)} {\it If $[l]_x = \emptyset$ then $x_j\equiv 0\!\!\!\mod
\pi A[u],$ for all $j\in [l],$ hence
$e_{[l]}(x)\equiv 0\!\!\!\mod \pi A[u].$ 
This implies that 
$\v_\bullet\big(E_{[l]}(x)\big) \ge f(l) + p^s.$}

\item  {(v)} {\it We have $\v_p\big(f(l) + p^{\v_p(j)}\big) = \v_p(j),$
for all $j\in [l].$
\newline And, if $\v_p(j_1) < v_p(j_2),$ then $f(j_1+j_2) = f(j_1) +
j_2.$}

\item  {(vi)} {\it If $x=\us{0\le j<p^s}\sum \Pi^j x_j$ has 
$E_{[l]}(x)\equiv 0\!\!\!\mod \Pi^{p^{r+t}-p^t}\Delta  $
(i.e., $\v_\bullet\big(E_{[l]}(x)\big)\ge p^{r+t} - p^t),$
 then $x_j\equiv 0\!\!\!\mod \pi A[u]$ for all $j\in [l].$} 
\end{description}

\noi
\begin{proof} (i) just fixes the notation according to Lemma
8A(i), (iv).

\noi
{\it For} (ii): Reading the congruence of Proposition~7C(i)
$\!\!\!\mod \Pi^{f(l)+1}  ,$ for $j\in [l],$ gives

\centerline{$
^g\Pi^j- \Pi^j \equiv  \Pi^{f(l)}j^\pr\!\!\!\mod \Pi^{f(l)+1}B .
$}

\noi
Right multiplying the $j^{\rm th}$ congruence by $x_j,$  summing
on $j\in [l],$ and left multiplying by $\Pi^{-f(l)}$
gives the claimed congruence.

\noi
From $A[u]\cap \Pi\Delta   = \pi A[u],$ by (2B1), it follows that
$e_{[l]}(x)
\equiv 0\!\!\!\mod \Pi\Delta  $ if, and only if, $e_{[l]}(x)\equiv
0\!\!\!\mod \pi A[u].$  Thus, $e_{[l]}(x)\not\equiv 0\!\!\!\mod
\Pi\Delta  $ implies $\v_\bullet\big(\Pi^{-f(l)}E_{[l]}(x)\big)=0$ hence
$\v_\bullet\big(E_{[l]}(x)\big) = f(l),$ and $e_{[l]}(x)\equiv
0\!\!\!\mod \Pi\Delta  $ implies $\Pi^{-f(l)}E_{[l]}(x)\equiv 0
\!\!\!\mod \Pi \Delta  $ hence
$\v_\bullet\big(\Pi^{-f(l)}E_{[l]}(x)\big) >0$ and
$\v_{\bullet}\big(E_{[l]}(x)\big) >f(l).$

\vs\noi
{\it For} (iii): The
congruences of Proposition~7C(i), for $j\in [l],$ and right multiplied by
$x_j,$  read
$$
\begin{aligned}
\Pi^{-f(l)}[\,^g\Pi^j &- \Pi^j] x_j - j^\pr x_j\\
&\equiv -j^\pr
\Pi^{p^{\v_p(j)}} x_j
  + \Big(\begin{matrix} j^\pr\\
2\end{matrix}\Big)(1-2\Pi^{p^{\v_p(j)}})\Pi^{(p^r-1)p^{\v_p(j)}}x_j
\!\!\!\mod \Pi^{(p^r+1)p^{\v_p(j)}} B  x_j\,
\end{aligned}
\leqno{\rm (8B1)}
$$
since $f(j) = f(l)$ for all $j\in [l];$
noting that these congruences have varying strength for different
$j,$ we 
need to work with a common modulus, for all $j\in [l],$ which must
involve the choice of $J.$  The rationale for why this should work is

\vs\noi
\begin{description}
\item  {1.} If $j<J,$ then $\v_p(j) >\v_p(J),$ by Lemma~8A(iii).
\item  {2.} If $j>J,$ then $\v_\bullet (x_j) \ge p^s,$ by
Corollary~2B(iii).
\item  {3.} $s\ge \v_p(j) +1,$ for all $j,$ by Lemma~8A(iii).
\end{description}

\noi
These are our main tools.

\begin{description}
\item  {{\it Step} 1.} $\Pi^{(p^r+1)p^{\v_p(j)}}Bx_j \subseteq
\Pi^{p^{\v_p(J)+1}}\Delta  ,$ for $j\in [l].$
\newline If $j<J,$ then $\Pi^{(p^r+1)p^{\v_p(j)}}Bx_j \overset 1\subseteq
\Pi^{(p^r+1)p^{\v_p(J)+1}}\Delta  \subseteq
\Pi^{p^{\v_p(J)+1}}\Delta  .$
\newline If $j=J,$ then $\Pi^{(p^r+1)p^{\v_p(J)}}Bx_J\subseteq
\Pi^{p^{r+\v_p(J)}}\Delta  \subseteq \Pi^{p^{v_p(J)+1}}\Delta  .$
\newline If $j> J,$ then $\Pi^{(p^r+1)p^{\v_p(j)}}Bx_j \overset
2\subseteq \Pi^{p^r+1}\Delta  \Pi^{p^s} \overset 3\subseteq
\Pi^{p^{\v_p(J)+1}}\Delta  .$

\item {\it {Step} 2.} $\Pi^{p^{\v_p(j)}}x_j$ and
$\Pi^{(p^r-1)p^{\v_p(j)}}x_j$ are in $\Pi^{p^{\v_p(J)+1}}\Delta  ,$
if $j\ne J.$ 
\newline First, $\Pi^{p^{\v_p(j)}}x_j \in \Pi^{p^{\v_p(J)+1}}
\Delta  ,$ by $1.$  if $j<J,$ and by $2., 3.$ if $j>J.$
\newline
For $j<J,\,\Pi^{(p^r-1)p^{\v_p(j)}}x_j \os 1\in
\Pi^{(p^r-1)p^{\v_p(J)+1}}\Delta  \subseteq
\Pi^{p^{\v_p(J)+1}}\Delta  ,$
\newline
For $j>J,$ \,$\Pi^{(p^r-1)p^{\v_p(j)}}x_j\os 2\in
\Pi^{(p^r-1)p^{\v_p(j)}}\Delta  \Pi^{p^s}\subseteq \Delta\Pi^{p^s}  \os
3\subseteq \Pi^{p^{\v_p(J)+1}}\Delta  .$

\item{\it {Step} 3.} $\Pi^{(p^r-1)p^{\v_p(J)}}x_J\in
\Pi^{p^{\v_p(J)+1}}\Delta  ,$ unless $r=1,$ and
$\Pi^{p^{\v_p(J)}}x_J\notin \Pi^{p^{\v_p(J)+1}}.$
\newline For $v_\bullet(\Pi^{(p^r-1)p^{\v_p(J)}}x_J) =
(p^r-1)p^{\v_p(J)} + 0 = p^{\v_p(J)+1} + (p^r-p-1)p^{\v_p(J)} \ge
p^{\v_p(J)+1}$ when $p^r-p-1\ge 0,$ i.e. $r\ge 2.$  And, $v_\bullet
(\Pi^{p^{\v_p(J)}}x_J) = p^{\v_p(J)} + 0,$ completing Step~3.
\end{description}

\vs\noi
Now, by Step 1, we can read the (8B1) congruences, for all $j\in
[l], \!\!\mod \Pi^{p^{\v_p(J)+1}}\Delta  .$  By Step~2, the right side
of these congruences for $j\ne J$ are all $0\!\!\mod
\Pi^{p^{\v_p(J)+1}}\Delta  .$ By Step~3, the right side of the
congruence for $j=J$ has
 second term (starting with
$\big(\begin{matrix} J^\pr\\ 2\end{matrix}\big))$  congruent to 
$0\!\!\mod \Pi^{p^{\v_p(J)+1}}\Delta  ,$ for $r>1;$
for $r=1,$ it is $\big(\begin{matrix} J^\pr\\
2\end{matrix}\big)\big(1-2\Pi^{p^{\v_p(J)}}\big)
\Pi^{(p-1)p^{\v_p(J)}}\equiv \big(\begin{matrix} J^\pr\\
2\end{matrix}\big) \Pi^{(p-1)p^{\v_p(J)}}\!\!\!\mod
\Pi^{p^{\v_p(J)+1}}\Delta  ,$ since $p^{\v_p(J)}+ (p-1)p^{\v_p(J)} =
p^{\v_p(J)+1}.$  Summing 
these \!\!$\mod \Pi^{p^{\v_p(J)+1}}\Delta  $ congruences, over all $j\in
[l],$ now gives the claimed congruence of (iii) (with the second
term on the right side $\equiv 0$ if $r>1).$

\vs\noi
Now if, $e_{[l]}(x)\not\equiv 0\!\!\!\mod \pi A[u],$ reading our new
congruence $\!\!\!\mod \Pi\Delta$ recovers (ii) and thus
\newline $\v_\bullet
\big(E_{[l]}(x)\big)=f(l).$ And also, in particular,
if $e_{[l]}(x)\equiv 0\!\!\!\mod \pi
A[u]$ then reading it $\!\!\!\mod \Pi^{p^{\v_p(J)}+1}\Delta$
determines
$\v_\bullet\big(E_{[l]}(x)\big),$ since $\v_\bullet(x_J) =0$
and $s> \v_p(J)$ (and $(p-1)p^{\v_p(J)}\ge p^{\v_p(J)}+1$ when
$r=1).$

\vs\noi
{\it For} (iv): Every $j\in [l]\bs [l]_x$  has
$v_\bullet(x_j)>0,$ hence $x_j \in A[u]\cap \Pi\Delta = \pi A[u].$
So $[l]_x = \emptyset$ implies
$e_{[l]}(x) \equiv 0\!\!\!\mod \pi A[u].$
Multiplying the congruence of Proposition~7C(i) by $\Pi^{-f(j)}$
gives
$$
\Pi^{-f(j)}[\,^g\Pi^j-\Pi^j]
\equiv j^\pr(1-\Pi^{p^{\v_p(j)}}) +
\Big(\begin{matrix} j^\pr\\
2\end{matrix}\Big)\big (1-2\Pi^{p^{\v_p(j)}}\big)
 \Pi^{(p^r-1)p^{\v_p(j)}}\!\!\!\mod \Pi^{(p^r+1)p^{\v_p(j)}}B,
$$
when right multiplying the $j^{\text{\rm th}}$ congruence by $x_j$ $(\equiv
0\!\!\!\mod \pi A[u])$ and summing over $j\in [l]$ (so $f(j) =f(l))$ yields
$$ 
\Pi^{-f(l)} E_{[l]}(x) 
\!\equiv\! e_{[l]}(x) -\sum_{j\in [l]} j^\pr \Pi^{p^{\v_ p(j)}}x_j
 + \sum_{j\in [l]}
\Big(\begin{matrix} j^\pr\\ 2\end{matrix}\Big) (1-2\Pi^{p^{\v_p(j)}})
 \Pi{(p^r-1)}p^{\v_p(j)} x_j\!\!\!\!\mod\!
\Pi^{(p^r+1)p^{\v_p(j)}}\pi \Delta  .
$$
In particular, we get $\Pi^{-f(l)}E_{[l]}(x) \equiv e_{[l]}
(x)\!\!\!\mod \Pi \pi\Delta,  $
hence
$$
\Pi^{-f(l)} E_{[l]} (x) \equiv e_{[l]}(x)\equiv 0\!\!\!\mod
\pi\Delta\q\text{\rm and}\q
 \v_\bullet \big(E_{[l]}(x)\big) \ge f(l)+\v_\bullet (\pi) = f(l) +
p^s.
$$

\vs
\noi
{\it For} (v): First, $f(j) + p^{\v_p(j)} = f(j^\pr)p^{\v_p(j)}+ p^{\v_p(j)} =
\big(f(j^\pr) + 1\big)p^{\v_p(j)} = (j^\pr + p^r)p^{\v_p(j)},$ for
all $j\in L.$  When $j\in [l],$ $f(l) = f(j),$ hence $f(l) +
p^{\v_p(j)} = f(j) + p^{\v_p(j)} = (j^\pr + p^r)p^{\v_p(j)}$ and
$\v_p\big(f(l) + p^{\v_p(j)}\big) = \v_p(j^\pr + p^r) + \v_p(j)
=\v_p(j).$

\noi
And $(p^r-1)p^{\v_p(j_1+j_2)} + (j_1+j_2) = \big((p^r-1)
p^{\v_p(j_1)} + j_1\big) + j_2.$

\vs\noi
{\it For} (vi): If this were false, then $[l]_x \ne \emptyset,$ hence
there is a largest $J\in [l]_x$
 which, by 
(ii) and (iii), has $\v_\bullet\big(E_{[l]}(x)\big) $ equal to $f(l)$
or $f(l) + p^{\v_p(J)}.$   Thus $\v_\bullet\big(E_{[l]}(x)\big)
\le f(l) + p^{\v_p(J)}=f(J) + p^{\v_p(J)} \os 1= 
(J^\pr + p^r)p^{\v_p(J)} \os 2\le (V_{\v_p(J)}
- 2 + p^r)p^{\v_p(J)} =
\big((p^r-1)p^{t-\v_p(J)}-1\big)p^{\v_p(J)} = (p^r-1)p^t-p^{\v_p(J)} < p^{r+t} - p^t,$
contradicting the Proposition hypothesis
$\v_\bullet\big(E_{[l[}(x)\big)\ge p^{r+t}-p^t.$  Here $\os 1=$ holds
 by the proof of (v), and $\os 2\le$ by Proposition~7C(iii),
since
$J\in L_{\v_p(J)}$  implies
$J^\pr \le V_{\v_p(J)}-1,$
which cannot be an equality by $V_{\v_p(J)}\equiv
1\!\!\!\mod p.$  
\end{proof}

\vs\noi
{\sc Lemma 8C.}
 {(i)} {\it If $l\in L$ is least in $[l]$ and $j\in [l],$ then
$\tilde j \os{\rm def}= j+p^{\v_p(j)}\in L$ is least in $[\,\tilde j\,]$ and has
 $f(\,\tilde j\,) = f(\ell) + p^{\v_p(j)},$
also $\v_p(\wt j) =\v_p(j),$
unless  $j=l$ and $l^\pr \equiv -1\!\!\!\mod p.$  We call this the
general case.}

\vs\noi
 {\it 
In the special case when $j=l,$ $\v_p(l)>0,$ and $l^\pr \equiv
-1\!\!\!\mod p,$ then
$\tilde l\os{\rm def} = \tau  (l) + p^{\v_p(l)}\in L$ is least in
$[\,\tilde l\,]$ and has $f(\,\tilde l\,) =f(l) + p^{\v_p(l)},$
also $\v_p(\wt \ell) = \v_p(\ell)-1.$} 

\vs \noi
{\it 
Finally, there is no solution $\tilde l\in L$ of
$f(\,\tilde l\,) = f(l) + p^{\v_p(l)} $ when $l=l^\pr\equiv -1\!\!\!\mod
p.$ }

\vs\noi
(ii) {\it If $e_{[l]}(x)\equiv 0\!\!\!\mod \pi A[u]$ and $[l]_x \ne
\emptyset,$ then $\v_p(l)>0.$}

\vs\noi
(iii) {\it If $d_i\in D,$ for $0\le i\le k,$ have $d_0\ne 0$ and
$\v_\bullet(d_i)> \v_\bullet(d_0)$ for $i>0,$ then

$\v_\bullet\big(\us{i\ge 0}\sum d_i) = \v_\bullet (d_0).$} 

\vsk\noi
{\it Proof} {\it For} (i): 
We first need

\vs\noi
{\it Step 1}: $v_p(j^\pr +1) =0,$ unless
 $j=l$ and
 $l^\pr \equiv -1\!\!\!\mod p.$

\vs\noi
We have $j=l+(p^r-1)(p^{\v_p(l)}-p^{\v_p(j)}),$ by Lemma~8A(v).
Now if $\v_p(j) < \v_p(l),$ then $j^\pr = lp^{-\v_p(j)} +
(p^r-1)(p^{\v_p(l)-\v_p(j)} -1) \equiv 1\!\!\!\mod p,$ hence $j^\pr
+1\equiv 2\!\!\!\mod p,$ as required.  If $\v_p(j) =\v_p(l),$ then
$j=l,$ by Lemma~8A(iv), and $l^\pr\not\equiv -1\!\!\!\mod p$
implies $j^\pr + 1 = l^\pr +1\not\equiv 0\!\mod\!p.$ This settles
Step~1.

\noi
When $v_p(j^\pr +1) =0,$  $\tilde j = (j^\pr +1)p^{\v_p(j)}$ implies $\v_p(\,\tilde j) =
\v_p(j)$ and $\tilde j^{\,^\pr} =(\tilde j)^\pr= j^\pr +1,$ so $f(\tilde
j^{\,^\pr}) =f(j^\pr
+1) = j^\pr + p^r,$ hence 
$f(\,\tilde j\,) =
f(\,\tilde j^{\,^\pr})p^{\v_p(\,\tilde j\,)} = 
(j^\pr + p^r)p^{\v_p(j)} \dot = f(l) + p^{\v_p(j)},$ with
$\dot =$ by the proof of
Proposition~8B(v). If $\tilde j$ were not least in $[\,\widetilde j\,],$ then
$\tilde j = \tau  (i)$ with $i\in [\,\tilde j\,],$ by Lemma~8A(iv),
implies $\tilde j^{\,^\pr} = \tau  (i)^\pr \equiv 1\!\!\!\mod p$ by Lemma~8A(ii),
hence that $f(\,\tilde j^{\,^\pr}\,)= f\big(\tau  (i)^\pr\big) =
(p^r-1)p^{\v_p(\tau  (i)^\pr)} + \tau  (i)^\pr \equiv (p^r-1)p^0 +1
\equiv 0\!\!\!\mod p,$
 contradicting $\v_p\big(f(\,\tilde
j^{\,^\pr}\,)\big)= v_p(j^\pr + p^r) = 0.$  Finally, $\tilde j\in L$
because $\v_p(\tilde j) = \v_p(j)<t$ and $j^\pr <V_{\v_p(j)}-1,$ since
$V_{\v_p(j)} -1\equiv 0\!\!\!\mod p$ implies $\tilde j^{\,^\pr} = j^\pr +1
<V_{\v_p(j)} = V_{\v_p(\tilde j)}.$

\noi
In the special case, 
$\tilde
l = \tau  (l) + p^{\v_p(l)}$ solves $f(\,\tilde l\,) = f(l) +
p^{\v_p(l)}$ since $\v_p\big(\tau  (l)\big)= \v_p(l) -1$ implies
$f\big(\tau  (l) + p^{\v_p(l)}\big) =
f\big((\tau  (l)
^\pr +p)
p^{\v_p(l)-1}\big) 
\dot= f\big(\tau  (l)^\pr + p\big)p^{\v_p(l)-1}
\dot= \big(p^r-1+\tau  (l)^\pr + p\big)p^{\v_p(l)-1} 
\ddot = p(l^\pr+p^r) p^{\v_p(l)-1}
= (l^\pr + p^r)p^{\v_p(l)}
\dot= f(l) + p^{\v_p(l)},$
with
$\dot=$ by Proposition~7C(i)  and \"= by 
(8A2).
And $\tilde l$ is least in
$[\,\tilde l\,]$ since $\tau  (l) + p^{\v_p(l)} =\tau  (i),$
with $i\in [\,\tilde l\,],$ has
$p^{v_p(\tau(l))}\big(\tau(l)^\pr +p\big) = p^{v_p(i)-1}
\tau(i)^\pr,$ with $\tau  (l)^\pr$ and $\tau(i)^\pr$ both $\equiv
1\!\!\!\mod p,$ by Lemma~8A(ii)
and its proof,
hence $v_p(l) = v_p(i)$ and
$\tau(l)^\pr +p = \tau(i)^\pr.$  Thus 
 $i^\pr = l^\pr +1\equiv
0\!\!\!\mod p,$ contradiction.  Finally, $\wt l\in L$ because
$\v_p(\,\tilde l\,) =\v_p(l) -1<t-1$ and 
 $\tilde l^{\,\pr}=\tau(l)^\pr +p = (p^r-1)(p-1)+ pl^\pr +p \dot<
(p^r-1)(p-1) + p\big(V_{\v_p(l)}-1\big) + p= (p^r-1)(p-1) +
pV_{\v_p(l)} = V_{\v_p(l)-1} = V_{\v_p(\tilde l)},$ with $\dot<$ by
$l^\pr <V_{\v_p(l)}-1 \equiv 0\!\!\!\mod p.$

\noi
Finally, this argument fails when $\v_p(l)=0,$ since $l\in L_0.$
But $f(\,\tilde l^{\,^\pr}\,)p^{\v_p(\,\tilde l\,)} = f(\,\tilde
l\,) = f(l) + p^{\v_p(l)} = f(l^\pr) + 1 = l^\pr + p^r$ implies
$\v_p(\,\tilde l\,)=0,$ hence $f(\,\tilde l^{\,^\pr}\,) = l^\pr + p^r $
and $\tilde l^{\,^\pr} = l^\pr + 1\equiv 0\!\!\!\mod p,$ contradiction.

\vs\noi
{\it For} (ii): If $\v_p(l) =0$ then
$\vert  [l]\vert  =\v_p(l)+1=1,$ by Lemma~8A(ii), hence
$l^\pr x_l = e_{[l]}(x)\equiv
0\!\!\!\mod \pi A[u].$  This implies that  $x_l\equiv 0\!\!\!\mod \pi
A[u],$ contradicting $[l]_x\ne \emptyset.$  

\vs\noi
{\it For} (iii): Taking $c_i=d^{-1}_0 d_i$ for all $i,$ we have
$c_0=1$ and $\v_\bullet(c_i)>0$ for $i>0.$  Then $1+\us{i>0}\sum
c_i$ is not in the maximal ideal $\Pi\Delta  _{\bullet}$ of the
valuation ring $\Delta  _\bullet,$ in the notation of
Corollary~2B(i), hence $\v_\bullet\big(1+\us{i>0}\sum c_i\big) =0$ and
$\v_\bullet\big(\us{i\ge 0}\sum d_i\big) =
\v_\bullet\big(d_0(1+\us{i>0}\sum c_i)\big) = v_\bullet (d_0).$ 
{\hfill  $\square$}

\vsk\noi 
{\sc Lemma 8D.} {\it Notation as in Proposition}~7C, {\it with level
$n$ of type $t\ge 1,$
we have
$$
\us{j\in\CL\bs L}\sum \big[^g\Pi^j - \Pi^j\big]x_j =\sum_{0\le \v<t}
\Big(\us{0<k<p^{t-\v}}{\sum\,^{\!\!^*}}
\big[^g\Pi^{p^\v(V_\v-1+k)} -
\Pi^{p^\v(V_\v-1+k)}\big]x_{p^\v(V_\v-1+k)}\Big),
\leqno{\rm (i)}
$$ 
where $\sum^*$ means $\v_p(k)=0, and $ $\v_\Pi(\,^g
\Pi^{p^\v(V_\v-1+k)}-\Pi^{p^v(V_\v-1+k)}) = (p^{r+t}-p^t) + (k-1)p^\v.$}

\vs\noi
(ii) {\it  If $j\in\CL\bs L$ then $j\ge p^t. $}

\vs\noi
(iii) $\v_\Pi(p) \ge p^{r+t}(p-1). $

\vs\noi
(iv) Write $j\in \CL_\v\bs L_\v$ as $j= p^\v(V_\v -1+k),$ for
$0\le \v <t$ and  $ 0<k<p^{t-\v}$ as in (i), and 
abbreviate $\big\lceil\frac{p^{t-v}}{p^r-1}\big\rceil$ by $[v\rangle,$ to get
the congruence
$$
\Pi^{-(p^{r+t}-p^t)}[\,^g\Pi^j - \Pi^j] \equiv \sum_{0<i\le [v\rangle } 
\Big(\begin{matrix} V_v -1+k\\ i \end{matrix}\Big)
(1-\Pi^{p^v})^i\; \Pi^{p^v[(p^r-1)(i-i) +(k-1)]}\!\!\!\mod \Pi^{p^t}
B.
$$

\vs
\begin{proof}
{\it For} (i): Since $L=\us{0\le \v<t}{\dot\bigcup} L_{\v},$ $\CL =\us{0\le
\v<t}{\dot\bigcup} \CL_\v$ with $L_\v\subseteq  \CL_\v,$ by
Proposition~7C(ii) and (iii), we have $\CL\bs L =\us{0\le
\v<t}{\dot\bigcup} (\CL_\v\bs L_\v)$
with
$
\CL_\v\bs L_\v
=\{1\le j<p^s:\v_p(j)=\v \;
\text{\rm
and}\;
V_\v\le j^\pr <\CV_\v\}.
$

\noi
Thus $j\in \CL_\v\bs L_\v$ has $j=p^\v j^\pr,$ so
defining $k$ by $j^\pr = V_\v -1+k,$ hence
  $V_\v \le V_\v-1 +k <\CV_\v$ and
$1\le k < \CV_\v -(V_\v-1) = 1+ p^{t-\v}-1 = p^{t-\v},$ i.e.
$0<k<p^{t-\v},$ and $k\equiv j^\pr\!\!\!\mod p$ implies $\v_p(k)=0.$

\noi
The $\Pi$-valuation of [\q] is, by Proposition 7C(i), equal to
$$
\begin{aligned}
f\big(p^\v(V_\v -1+k)\big) &= f(V_\v -1+k)p^\v
 =\big(p^r-1+(p^r-1)(p^{t-\v}-1)+k-1\big)p^\v\\
& =\big((p^r-1)p^{t-\v}+k-1\big)p^\v
= (p^r-1)p^t + (k-1)p^\v.
\end{aligned}
$$

\noi
{\it For} (ii): Suppose that $j\in\CL_\v \bs L_{\v},$ with $0\le \v
<t,$ has $j<p^t.$  Then $j= p^\v (V_\v-1+k),$ with $0<k<p^{t-\v}$ by
(i), implies $p^t >j \ge p^\v V_{\v},$ hence $p^{t-\v} >V_{\v} =
(p^r-1)(p^{t-\v}-1),$ by Proposition~7C(iii).  
So
$p^{t- \v}-1\ge (p^r-1)(p^{t-\v}-1)$ and $p\le p^r \le 2,$
 contradiction.

\vs\noi
{\it For} (iii):
$\v_\Pi(p) = \phi(p^{r+s})= p^{r+s-1}(p-1) \os 1\ge p^{r+w}(p-1) \os 2\ge
p^{r+t}(p-1),$ with $\os 1\ge$ by the first paragraph of the proof of
Theorem~4G and $\os 2\ge$ by Lemma~6E(ii).

\vs\noi
(iv) Now $\v_p(j) = \v$ hence $^{g-1}\Pi \equiv 1+
(1-\Pi)\Pi^{p^r-1}\!\!\!\mod pB,$ from Lemma~2G(iii), and
$^{g-1}\Pi^{p^\v} \equiv 1 +
(1-\Pi^{p^\v})\Pi^{(p^r-1)p^\v}\!\!\!\mod pB,$ by induction on
$\v\ge 0.$  Raising this to the $j^{\rm'th}$ power, in binomial
expansion, 
gives $^{g-1}\Pi^j \equiv \us{0\le i\le j^\pr}\sum
\begin{pmatrix} j^\pr\\ i\end{pmatrix} (1-\Pi^{p^\v})^i
\Pi^{i(p^r-1)p^{\v}}\!\!\!\mod pB,$ and multiplying this congruence
by $\Pi^j$ and subtracting $\Pi^j$ from both sides gives 
$$
\big[^g\Pi^j-\Pi^j\Big] \equiv \sum_{0<i\le j^\pr} \begin{pmatrix}
j^\pr\\ i\end{pmatrix} (1-\Pi^{p^\v})^i\,
\Pi^{p^\v(p^r-1)i+j}\!\!\!\mod p\Pi^jB.
$$
Observing that $j=(p^r-1)(p^t-p^\v) +p^\v(k-1)$
has $j-(p^{r+t}-p^t) = -p^\v(p^r-1) + p^\v(k-1) = p^\v(k-p^r),$ and,
by (iii), that
$\v_\Pi\big(p\Pi^{j-(p^{r+t}-p^t)}\big)
= \v_\Pi(p) + 
\big(j-(p^{r+t}-p^t)\big) \ge p^{r+t}(p-1) + p^\v(k-p^t)
= p^{r+t}(p-2) + \big(p^{r+t}-p^{r+\v}\big) + p^\v k
\ge p^{r+t},$
 may cancel  $\Pi^{p^{r+t}-p^t}$ throughout the above congruence to get
$$
\Pi^{-(p^{r+t}-p^t)} \big[^g\Pi^j-\Pi^j\big]
\equiv \sum_{0<i\le j^\pr} \begin{pmatrix} j^\pr\\ i\end{pmatrix} 
(1-\Pi^{p^\v})^i\, \Pi^{p^\v\big((p^r-1)(i-1)+(k-1)\big)}\!\!\!\mod
\Pi^{p^t} B.
$$

\noi
Note that the $i^{\rm th}$ term of the sum can be ignored provided
$p^\v\big((p^r-1)(i-1)+(k-1)\big)\ge p^t$ for all
$k\Longleftrightarrow (p^r-1)(i-1) \ge p^{t-\v} \Longleftrightarrow i
\ge 1+ \Big\lceil\frac{p^{t-\v}}{p^r-1}\Big\rceil,$ yielding the shorter sum
claimed since $[v\rangle \le j^\pr,$ for which it suffices to show that
$
\frac{p^{t-v}}{p^r-1} \le V_v+ k-1$ for all $k
\Longleftrightarrow\,\frac{p^{t-v}}{p^r-1} \le V_v = 
 (p^r-1)(p^{t-v}-1), $ i.e.
$(p^r-1)^2(1-\frac{1}{p^{t-v}})\ge 1.$ Since $r\ge 1$ and $t-v\ge 1,$
this follows from $(p-1)^2(1-1/p) \ge 1.$
\end{proof}

\section{Toward the 7C(ii) Congruence} \label{Sec9}

\vs\noi
Next we take first steps generalizing Corollary~
7D, for $t=1,$ to
search for solutions of congruence~7C(ii), for $t\ge 2,$ by using
the methods of \S8.  The hypothetical solutions of these congruences
$\!\!\!\mod \Pi^{p^{r+t}}\Delta$ also solve congruence 7C(iii)
$\!\!\!\mod\!\!\!\!\;\;\Pi^{p^{r+t}-p^t}\Delta,$ with relatively simpler structure,
which have been studied in Proposition~8B.

\vs\noi 
{\sc Theorem 9A.} 
{\it  Suppose that $s\ge r+t.$ }

\noi
(i) {\it If $ x=\us{0\le j<p^s}\sum \Pi^j x_j,$ with all $x_j\in
A[u],$
has $\sum_{[l]\in L/\sim} E_{[l]}(x)\equiv 0 \!\!\!\mod
\Pi^{p^{r+t}-p^t}\Delta,$ then $x_j\equiv 0\!\!\!\mod \pi A[u]$ for all
$j\in L.$}

\noi
{\rm(ii)} 
{\it If $x=\os {p^s-1}{\us{j=0}\sum} \Pi^j x_j$ solves the
Proposition~{\rm 7C(iii)} congruence then $x$ solves the
Proposition~{\rm 7C(ii)} congruence if, and only if, $x$ solves the
congruence}
$$
\begin{aligned}
mx_0 &+\sum_{0\le v<t} \Big(\us{0<k<p^{t-v}}{\sum{}^*}
\big(\sum_{0<i\le [v)}
\big(\begin{matrix} V_\v-1+k\\ i\end{matrix} \big)
\big(1-\Pi^{p^v}\big)^i\,
\Pi^{p^\v[(p^r-1)(i-1)+(k-1)]}\big)x_{p^\v(V_\v-1+k)}\Big)\\
& \equiv \Pi^{p^t-1} y\!\!\!\mod \Pi^{p^t}\Delta.
\end{aligned}
$$

\vs\noi
{\it Proof.} {\it For} (i): Let $\mu  _x =\,\text{min}\big\{\v_\bullet
(E_{[l]}(x)\big) :[l]\in L/\sim\big\} \le \infty  $ and $M_x
=\big\{[l] \in L/\sim: \v_\bullet\big(E_{[l]}(x)\big)=\mu  _x\big\}.$ 
Note that $M_x\ne \emptyset.$  Also note that  $x_j$ with $j\in L$
occurs only in $E_{[j]}(x).$ 

\vs\noi
{\it Step} 1. There is a unique $[l_0]\in M_x$ with $f(l_0) <f(l)$
for all $[l]\in M_x$ with $[l] \ne [l_0].$  Moreover, the least $l_0$
in $[l_0]$ is unique in $L.$

\noi
$M_x$ contains an
 $[l_0]$ with $f(l_0)\le f(l)$ in $\Z$ for all
$[l]\in M_x.$ This $[l_0]$ must
have $f(l_0) < f(l)$ for all $[l]\in M_x$ with $[l]\ne [l_0],$
because $f(l_0) = f(l)$ 
implies $l_0\sim l,$ by Lemma~8A(i), hence $[l_0] = [l],$
contradiction.
The rest
of Step~1 follows from Lemma 8A(iv).

\noi
The value of $e_{[l_0]}(x)\!\!\mod \pi A[u]$ in Proposition 
8B(ii) creates two basic cases. 
The first is

\vs\noi
{\it Step} 2.  If $e_{[l_0]}(x)\not\equiv 0\!\!\!\mod \pi A[u],$
then $\mu_x = f(l_0)$ and $\v_\bullet\big(E_{[l]}(x)\big) >\mu  _x$ for
all $[l]\ne [l_0]$ in $L/\!\!\sim.$ 
It follows that $x_j\equiv 0\!\!\!\mod \pi A[u]$ for all $j\in L.$

\noi
Now Proposition~8B(ii) and 
Step~1 imply that $\v_\bullet \big(E_{[l_0]}(x)\big)=f(l_0)$ and 
$\v_\bullet\big(E_{[l_0]}(x)\big) = \mu_x,$ since $[l_0]\in M_x,$
hence
$f(l_0) = \mu  _x.$
In particular, we have $\v_\bullet 
\big(E_{[l]} (x)\big) \ge f(l)>f(l_0) = \mu  _x$ 
for $[l]\in M_x\bs \{[l_0]\},$ 
 and $\v_\bullet\big(E_{[l]}(x)\big)> \mu_x$ for $[l]
\notin M_x,$ proving
the first assertion of Step~2.

\noi 
The first assertion checks the hypothesis of Lemma 8C(iii)
which implies the  equality

\vs
\centerline{$
v_\bullet\big(E_{[l_0]}(x) \big)
= \v_\bullet\big(\us {[l]\in L /\sim}\sum\,E_{[l]}(x)\big) \ge p^{r+t} -
p^t,
$}

\vs\noi
with the inequality by the hypothesis of (i), with $\mu_x=f(l_0).$
This
 implies that
$\v_\bullet\big(E_{[\ell]}(x)\big)\ge \v_\bullet\big(E_{[l_0]}(x)\big) \ge
p^{r+t} - p^t$ for all $[l]\in L/\sim,$ hence Proposition~8B(vi) completes
Step~2.

\noi
The remaining basic case $e_{[l_0]}(x)\equiv 0\!\!\!\mod \pi A[u]$ has two
subcases by Proposition~8B(iii).
The first is

\noi
{\it Step} 3.  If $e_{[l_0]}(x)\equiv 0\!\!\!\mod \pi A[u]$ and
$[l_0]_x\ne \emptyset$ then
$\mu  _x =f (l_0) + p^{\v_p(J_0)},$ with $J_0$ largest in
$[l_0]_x,$
and $\v_\bullet\big(E_{[l]}
(x)\big) >\mu  _x$
for all $[l] \ne [l_0]$ in $L/\sim.$  It follows that
$x_j \equiv 0\!\!\!\mod \pi A[u]$ for
all $j\in L.$

\noi Since the unique $[l_0]$ of Step~1  has $e_{[l_0]}(x)\equiv
0\!\!\!\mod \pi A[u],$ and $\mu  _x = \v_\bullet
\big(E_{[l_0]}(x)\big),$
by $[l_0]\in M_x,$ it has $\mu  _x =
f(l_0) + p^{\v_p(J_0)},$ by
Proposition~8(iii).
 Since $\v_\bullet \big(E_{[l]}(x)\big) >\mu  _x$ for
$[l] \not\in M_x,$
 we investigate the $[l_1]\in M_x\bs \{[l_0]\}.$

\vs\noi
{\it Step} 3.1.
There is no $[l_1]\in M_x\bs \{[l_0]\}$ with
$[l_1]_x\ne \emptyset.$ 

\noi
If such an $[l_1]$ did exist, then 
 Proposition 8B(iii) finds the  largest $J_1$ in $[l_1]_x$ with

\noi
\centerline{$f(l_1) + p^{\v_p(J_1)} =
\v_\bullet\big(E_{[l_1]}(x)\big) = \mu  _x =
\v_\bullet\big(E_{[l_0]}(x)\big) =
f(l_0)
+ p^{\v_p(J_0)}.$}

\noi
Comparing $\v_p$ of the first and last terms above by 
 Proposition~8B(v) implies that $\v_p(J_1) =
\v_p(J_0),$  hence  $f(l_1) = f(l_0).$ Thus Lemma~8A(i) implies $[l_1] =
[l_0],$ contradiction.

\vs\noi
{\it Step} 3.2.  There is no $[l_1]\in M_x\bs \{[l_0]\}$ with
$[l_1]_x=\emptyset.$

\noi
If such an $[l_1]$ did exist, then we would have
$$
f(l_1) + p^{s-1} \os 1> f(l_0)
+ p^{\v_p(J_0)} \os 2= \mu  _x \os
3= \v_\bullet\big(E_{[l_1]} (x)\big)\os 4\ge f(l_1) + p^s,
$$
with $\os 1 >$ by Step~1 and Lemma~8A(iv), $\os 2=$ by the start of
the proof of Step~3, $\os 3=$ by $[l_1] \in M_x,$ and $\os 4\ge$ by
Proposition 8B(iv).  This  implies the contradiction
$p^{s-1}\ge p^s.$

\noi
Combining Step~3.1 and Step~3.2
shows that $M_x =
\{[l_0]\},$ which implies the first assertion of Step~3.  The first
assertion checks the hypothesis of Lemma~8C(iii), which
implies the equality
$$
\v_\bullet\big(E_{[l_0]}(x)\big) = v_\bullet \Big(\sum_{[l]\in
L/\sim} E_{[l]}(x)\Big) \ge p^{r+t} - p^t,
$$
with the inequality by the hypothesis of Theorem~9A(i), now with 
$\mu  _x = f(l_0)+ p^{\v_p(J_0)}.$  Step~3 is complete.

\vs\noi
{\it Step} 4.  If $e_{[l_0]}(x)\equiv 0\!\!\!\mod \pi A[u]$ and
$[l_0]_x=\emptyset$
then 
$x_j \equiv 0\!\!\!\mod \pi A[u]$ for all
$j\in L.$

\noi
Suppose that this is false, so there exists $[l_1]$ with 
$x_{l_1} \not\equiv 0\!\!\!\mod \pi A[u],$ i.e. $[l_1]_x\ne
\emptyset.$ 
 Thus
$$
\v_\bullet \big(E_{[l_1]}(x)\big)  \ge \mu  _x = 
\v_\bullet\big(E_{[l_0]}(x)\big) 
\os 1\ge f(l_0) + p^s \ge 0+p^{r+t} \ge p^{r+t} - p^t,
$$
with $\os 1\ge$
by Proposition~8B(iv).
This implies $\v_\bullet \big(E_{[l_1]} (x)\big) \ge p^{r+t}-p^t$
hence $x_j\equiv 0\!\!\!\mod \pi A[u]$ for all $j\in [l_1],$
by Proposition~8B(vi).
In particular, $x_{l_1} \equiv 0\!\!\!\mod \pi A[u],$ contradiction.

\vs\noi
{\it For} (ii):If $x$ solves the Proposition~7C(iii) congruence then
$\us{[l]\in L/\sim}\sum E_{[l]}(x)\equiv 0\!\!\!\mod
\Pi^{p^{r+t}-p^t}\Delta,$ by Proposition~8B(i), hence $x_j
\equiv0\!\!\!\mod \pi A[u]$ for all $j\in L,$ by (i).

\vs\noi
Thus $x$
solves the Proposition~7C(ii) congruence 
$$
\Pi^{p^{r+t}-p^t} 
\us{0\le i<p^t}\sum{}^{g}\Pi^{i}mx_i 
+\us{j\in\CL\bs L}\sum [^g\Pi^j-\Pi^j]x_j \equiv 
 \Pi^{p^{r+t}-1} y\!\!\!\mod \Pi^{p^{r+t}}\Delta,
$$
since \; $x_j\equiv 0\!\!\!\mod \pi A[u],$ \, with \; $v_\Pi(\pi) = s\ge r+t$
\; for all\; $j\in L,$ by (i).
{}\; Note 
that the

\noi
$i$-sum consistency with redundancy of Proposition~7C(ii)
and the $j\in \CL\bs L$ of
Lemma~8D(ii) combine to reduce the $i$-sum to the $i=0$ term (as
$0\notin \CL),$ 
when
  cancellation by $\Pi^{p^{r+t}-p^t}$ implies that $x$ solves
$$
mx_0 
 +\sum_{j\in \CL\bs L} \Pi^{-(p^{r+t}-p^t)}[^g \Pi^j-\Pi^j]x_j
\equiv \Pi^{p^t-1}y\!\!\!\mod \Pi^{p^t}\Delta.
$$
Now Lemma~8D(i)  and Lemma~8D(iv) complete the proof.
\hfill $\square$

\vs\noi
{\sc Corollary 9B.} {\it  If $s\ge r+t,$ then the congruence of
Theorem}~9A(ii) {\it would have a solution. }

\vs\noi
Since the preparation of Proposition~8B, Lemma~8C and Lemma~8D also
goes smoothly, we even hoped to complete the special case $s=2.$
Corollary 9B apparently does not extend to $s<r+t:$ see Lemma~10C.

\section{More Commutators} \label{Sec10}

\vs\noi
This section records some commutators which may eventually become
useful if the inductive approach of \S9 turns out to be inadequate.

\vsk\noi 
{\sc Lemma 10A.}
(i) {\it If $0\ne \beta  \in \,\rad(B),$ then $u-\beta  \in
\Delta  ^\times,$ $1+T-N_\fg (\beta  ) \in A^\times,$ and}
$$
(u-\beta  ) ^{-1} = \big(1+T - N_\fg(\beta  )\big)^{-1}
\sum^{p^s-1}_{i=0} \,^{\sigma  _i}\beta   \, u^i,
$$
\begin{description}
\item{\q}
{\it where $\sigma  _i = \sum_{i<k<p^s} g^k\in \Z [\fg]$ for $0\le i\le 
p^s-1.$}

\item  {(ii)} {\it $\wt\beta  =\sum^{p^s-1}_{i=1}\, ^{\sigma  _i} \beta 
u^i$ has $\wt \beta$ and $\tilde \beta    + \,^{\sigma  _0}\beta  $ in $\Delta  ^\times,$
if $\beta  \ne 0.$ }

\item  {(iii)} {\it For $x\in \Delta  ,$ $(u-\beta  )x(u-\beta  )^{-1}
=\,^gx +\theta  _\beta  (x)$ with $\theta  _\beta  (x) =
\frac{(
\,^gx\beta  -\beta  x)(\wt\beta  +\,^{\sigma  _0}\beta 
)}{1+T-N_\fg(\beta  )}\,.$ }
\end{description}

\vs\noi
{\it Proof.} {\it For} (i): First, $u-\beta  =1+ (u-1-\beta  )\in 1
+ \,{\rm rad}(\Delta  ),$ by Lemma 2A; similarly
$1+\big(T-N_\fg(\beta  )\big)\in 1+\,{\rm rad}(A),$ as $N_\fg
(\beta  ) =\, ^{1+\sigma_0}\beta      \in\, {\rm
rad}(B)\cap A= \,{\rm rad}(A)$ (since $A\hookrightarrow B\to B/{\rm rad}(B)$
has kernel ${\rm rad}(A)).$  Now $B^\times$ is a multiplicative
$\Z[\fg]$-module, and 
$$
g\sigma  _{i-1}=1 +\sigma  _i, {\rm for}\; 1\le i\le p^s-1
\leqno{\rm (10A1)}
$$
(for $g\sum_{i-1<k<p^s}g^k=\sum_{i<k<p^s}\, g^k + g^{p^s} = 1+\sigma 
_i).$  Thus
$$
\begin{aligned}
(u-\beta  ) &\sum\nolimits ^{p^s-1}_{i=0} \,^{\sigma  _i}\beta  
\,u^i = \sum\nolimits_i\; ^{g\sigma_i }\beta  \,u^{i+1}
-\sum\nolimits  _i \,^{1+\sigma  _i}\beta  \,u^i\\
&\q =\sum\nolimits^{p^s}_{i=1} \,^{g\sigma  _{i-1}}\beta  \, u^i -
\sum\nolimits^{p^s-1}_{i=0} \,^{1+\sigma  _i}\beta  \,u^i
= {g\sigma  _{p^s-1}}\beta  u^{p^s} - ^{1+\sigma  _0} \beta  \,
u^0\\
&\q = \,^{g\cdot 0}\beta  (1+T)-N_\fg(\beta  )
= 1+T-N_\fg(\beta  ).
\end{aligned}
$$
Note that $\beta  =0$ is not allowed.

\vs\noi
{\it For} (ii): The $i=p^s -1$ term of $\wt\beta  $ is $u^{p^s-1}=
(1+T)u^{-1};$ the other terms are in $\beta  \Delta  \subseteq
\,{\rm rad}(B)\Delta  \subseteq \,{\rm rad}(\Delta  ).$

\vs\noi
{\it For} (iii): Writing $e=1+T-N_{\fg}(\beta  ),$ we have
$$
\begin{aligned}
&(u-\beta  )x(u-\beta  )^{-1}
= (\,^gxu-\beta  x)\frac 1e\,\sum\nolimits ^{p^s-1}_{i=0} \,^{\sigma 
_i}\beta  \,u^i\\
&\q= \frac 1e \Big[\,^gx\sum\nolimits ^{p^s}_{i=1}\,^{g\sigma 
_{i-1}}\beta  u^i-\beta  x\,\sum\nolimits ^{p^s-1}_{i=0} \,^{\sigma 
_i}\beta  u^i\Big]
= \frac 1e \Big[\,^gx\big(\beta  \wt\beta   + (1+T)\big) -\beta 
x(\wt\beta  +\,^{\sigma  _0}\beta  )\Big]\\
&\q= \frac 1e \Big[\,^gx(1+T) + (\,^g x\beta  -\beta  x)\wt \beta  -\beta 
x\;^{\sigma  _0}\beta  \Big]\\
&\q= \frac 1e \Big[\,^gx \big(1+T-N_{\fg}(\beta  )\big) + (\,^gx\beta 
-\beta  x)\wt\beta   + (\,^gx\beta  -\beta  x)\,^{\sigma  _0}\beta 
\Big]
= \,^gx + \theta  _\beta  (x). \q\q\q\q \square
\end{aligned}
$$

\noi
{\sc Remark.} Since the above units of $\Delta  $ are all in
$1+\,\text{rad}\,\Delta  ,$ note that $\Delta 
^\times/\,1+\,\text{rad}\,\Delta  \simeq \IF^\times,$ by Lemma~2A(i), with
$\IF^\times$ liftable to $\Delta  ^\times$ by central roots of unity in the Witt
vectors.

\vs\noi
{\sc Proposition 10B.}
{\it Let $0\ne \beta  \in \,{\rm rad}(B)$ and $n=p^{r+t} +
p^tm -1$ with $1\le t\le s-1$ and $p\!\!\not\vert m\ge 1.$}

\noi
 {(i)}
$\begin{aligned} \\
\,[u-\beta  ,1&-\Pi^{p^tm}x] \equiv 
1-\Pi^{p^tm} (\,^g x-x  )(1-\Pi^{p^tm})^{-1}-\Pi^{n+1-p^t}m^g x
\\
 & - \Pi^{n+1-p^t} m\;^gx \;\beta(u-\beta 
)^{-1}-\Pi^{p^tm}\theta_\beta  (x)\big(1-\Pi^{p^tm}x\big)^{-1}
\!\!\!\mod
\Pi^{n+1}\Delta  .
\end{aligned}
$

\vs\noi
  {(ii)} {\it If $x\in \Delta  ,$ $y\in A[u],$ then $[u-\beta 
,1-\Pi^{p^tm}x] \equiv 1-\Pi^ny\!\!\!\mod \Pi^{n+1}\Delta  $ if, and
only if,} 
$$
\Pi^{p^{r+t}-p^t} m\,^g x\;u(u-\beta  )^{-1} + \big(\,^g x-x +\theta
_\beta  (x)\big)\equiv \Pi^{p^{r+t}-1} y\!\!\!\mod
\Pi^{p^{r+t}}\Delta  .
$$

\noi
  {(iii)} $[u-\beta  , 1-\Pi^{p^tm}x]\equiv [u,1-\Pi^{p^tm}x] -
\Pi^{n+1-p^t}m\,^gx\beta  (u-\beta  )^{-1} - \Pi^{p^tm}\theta_\beta 
(x)\big(1-\Pi^{p^tm}x\big)^{-1}$

$\!\!\mod \Pi^{n+1}\Delta  .$

\vs\noi
\begin{proof}
{\it For} (i): The details are much like Lemma~7A(i).

\vs\noi
{\it Step 1.} 
$
[u-\beta  ,1-\Pi^{p^tm}x]$
$$
\equiv [1-\Pi^{p^tm}(u-\beta  )x(u-\beta
 )^{-1} -\Pi^{n+1-p^t}m\;^gx\; u(u-\beta  )^{-1} ]
 (1-\Pi^{p^tm}x)^{-1} \!\!\!\mod \Pi^{n+1}\Delta  .
$$

\noi
For  
$ (u-\beta  )(1-\Pi^{p^tm}x) (u-\beta  )^{-1}
= 1-(u-\beta )(\Pi^{p^tm}x)(u-\beta  )^{-1} $
$$
\begin{aligned}
= &1-(\,^g\Pi^{p^tm} u-\Pi^{p^tm}\beta  ) x(u-\beta  )^{-1}\\
= &1-\Pi^{p^tm}\big((\,^{g-1}\Pi^{p^tm}-1)u + (u-\beta  )\big)x(u-\beta  )^{-1}\\
= &1-\Pi^{p^tm}(u-\beta  )x(u-\beta  )^{-1}
-\Pi^{p^tm}(\,^{g-1}\Pi^{p^tm}-1)\,^gx\; u(u-\beta  )^{-1}\\
\dot\equiv &1-\Pi^{p^tm}(u-\beta  )x(u-\beta  )^{-1} -
m(1-\Pi^{p^t})\Pi^{n+1-p^t}\;^gx\;u (u-\beta  )^{-1}\\
\equiv &1-\Pi^{p^tm}(u-\beta  )x(u-\beta  )^{-1} -\Pi^{n+1-p^t}
m\;^gx \;u(u-\beta  )^{-1} \!\!\!\mod \Pi^{n+1}\Delta  
\end{aligned}
$$
with $\dot\equiv$ by Lemma~5G(i), $p^tm + (p^r-1)p^t = n+1-p^t,$
and $p^t + (n+1-p^t)\ge n+1.$  Right
multiplying by the unit $(1-\Pi^{p^tm}x)^{-1}$ completes {\it Step}
1.

\vs\noi
We next start from the $\!\!\!\mod \Pi^{n+1}\Delta  $ congruence of
{\it Step} 1 by replacing $(u-\beta  )x(u-\beta  )^{-1}$ via 
Lemma 10A(iii), using $(1-\Pi^{p^tm}x)^{-1} = 1+
\Pi^{p^tm}x(1-\Pi^{p^tm}x)^{-1}$ with the first and last terms
inside $[\q],$ and splitting the middle one to continue (always
$\!\!\!\mod \Pi^{n+1}\Delta  ):$
$$
\begin{aligned}
\big(1&+\Pi^{p^tm}x (1-\Pi^{p^tm}x)^{-1}\big) - \Pi^{p^tm}\;
^gx(1-\Pi^{p^tm}x) ^{-1} 
 - \Pi^{p^tm}\theta  _\beta 
(x)(1-\Pi^{p^tm}x)^{-1}\\
&\q\q - \Pi^{n+1-p^t}m\;^gx\,u(u-\beta 
)^{-1}\big(1+\Pi^{p^tm}x(1-\Pi^{p^tm} x)^{-1}\big)\\
&\equiv 1-\Pi^{p^tm}(\,^g x-x)(1-\Pi^{p^tm}x)^{-1}
\\
&\dot\equiv -
\Pi^{p^tm}\theta  _\beta  (x) (1-\Pi^{p^tm}x)^{-1}
 -
\Pi^{n+1-p^t}m\, ^gx\big(1+\beta   (u-\beta  )^{-1}\big)+0\\
& \ddot\equiv 1- \Pi^{p^tm}(\,^g
x-x)(1-\Pi^{p^tm}x)^{-1}-\Pi^{n+1-p^t} m^g x\\
&\q\q-
\Pi^{n+1-p^t} m\;^gx\,\beta  (u-\beta  )^{-1}
 - \Pi^{p^tm} \theta  _\beta  (x)(1-\Pi^{p^tm}x)^{-1}
\end{aligned}
$$
with $\dot\equiv$ by $u(u-\beta ) ^{-1}= 1+\beta  (u-\beta  )^{-1}$
and
`0' before $\ddot\equiv$ by $\Pi^{n+1-p^t}\Delta  \Pi^{p^tm} =
\Pi^{n+1+p^t(m-1)}\Delta   \subseteq \Pi^{n+1}\Delta  .$

\vs\noi
{\it  For} (ii): This amounts to making the same moves as for 
Lemma~7A(ii), for the same reasons i.e., add $-1,$ cancel
$-\Pi^{p^tm},$ right multiply by $(1-\Pi^{p^tm}).$

\vs\noi
{\it For} (iii): The first three terms on the right side of (i) agree
with the right side of Lemma~7A(i) because $p^{r+t} +p^t(m-1) = n+1
-  p^t.$
\end{proof}

\noi
{\sc Remark:} 
So each
type $t$ has a family of commutators, parametrized by
$\beta  .$

\vs
\noi
What was to be Theorem~9C has become

\noi
{\sc Lemma 10C.}
{\it Suppose that $s<r+t.$}

\noi
{\it If 
$x= \us{0\le j <p^s}\sum \Pi^jx_j,$ 
with all $x_j\in A[u],$  
has 
$\sum_{[l]\in L/\sim} E_{[l]}(x) \equiv
0\!\!\!\mod \Pi^{p^{r+t}-p^t}\Delta  ,$ then
\newline $x_j\equiv 0\!\!\!\mod \pi A[u]$ for all $j\in L.$ }

\vs\noi
{\it Proof.}  Let $\mu  _x =\,\text{min}\big\{\v_\bullet
(E_{[l]}(x)\big) :[l]\in L/\sim\big\} \le \infty  $ and $M_x
=\big\{[l] \in L/\sim: \v_\bullet\big(E_{[l]}(x)\big)=\mu  _x\big\}.$ 
Note that $M_x\ne \emptyset.$  Also note that  $x_j$ with $j\in L$
occurs only in $E_{[j]}(x).$ 

\vs\noi
{\it Step} 1. There is a unique $[l_0]\in M_x$ with $f(l_0) <f(l)$
for all $[l]\in M_x$ with $[l] \ne [l_0].$  Moreover, the least $l_0$
in $[l_0]$ is unique in $L.$

\noi
$M_x$ contains an
 $[l_0]$ with $f(l_0)\le f(l)$ in $\Z$ for all
$[l]\in M_x.$ This $[l_0]$ must
have $f(l_0) < f(l)$ for all $[l]\in M_x$ with $[l]\ne [l_0],$
because $f(l_0) = f(l)$ 
implies $l_0\sim l,$ by Lemma~8A(i), hence $[l_0] = [l],$
contradiction.
The rest
of Step~1 follows from Lemma 8A(iv).

\noi
The value of $e_{[l_0]}(x)\!\!\mod \pi A[u]$ in Proposition 
8B(ii) creates two basic cases. 
The 
first is

\vs\noi
{\it Step} 2.  If $e_{[l_0]}(x)\not\equiv 0\!\!\!\mod \pi A[u],$
then $\mu_x = f(l_0)$ and $\v_\bullet\big(E_{[l]}(x)\big) >\mu  _x$ for
all $[l]\ne [l_0]$ in $L/\!\!\sim.$ 
It follows that $x_j\equiv 0\!\!\!\mod \pi A[u]$ for all $j\in L.$

\noi
Now Proposition~8B(ii) and 
Step~1 imply that $\v_\bullet \big(E_{[l_0]}(x)\big)=f(l_0)$ and 
$\v_\bullet\big(E_{[l_0]}(x)\big) = \mu_x,$ since $[l_0]\in M_x,$
hence
$f(l_0) = \mu  _x.$
So 
we have $\v_\bullet 
\big(E_{[l]} (x)\big) \ge f(l)>f(l_0) = \mu  _x$ 
for $[l]\in M_x\bs \{[l_0]\},$ 
 and $\v_\bullet\big(E_{[l]}(x)\big)> \mu_x$ for $[l]
\notin M_x,$ proving
the first assertion of Step~2.

\noi 
The first assertion checks the hypothesis of Proposition 8C(iii)
which implies the  equality

\vs
\centerline{$
\v_\bullet\big(E_{[l_0]}(x) \big)
= \v_\bullet\big(\us {[l]\in L /\sim}\sum\,E_{[l]}(x)\big) \ge p^{r+t} -
p^t,
$}

\vs\noi
with the inequality by the Proposition hypothesis.
This
 implies
$\v_\bullet\big(E_{[\ell]}(x)\big)\ge \v_\bullet\big(E_{[l_0]}(x)\big) \ge
p^{r+t} - p^t$ for all $[l]\in L\sim,$ hence Proposition~8B[vi] completes
Step~2.

\noi
The remaining basic case $e_{[l_0]}(x)\equiv 0\!\!\!\mod \pi A[u]$ has two
subcases by Proposition~8B(iii).
The first is

\noi
{\it Step} 3.  If $e_{[l_0]}(x)\equiv 0\!\!\!\mod \pi A[u]$ and
$[l_0]_x\ne \emptyset$ then
$\mu  _x =f (l_0) + p^{\v_p(J_0)}$ and $\v_\bullet\big(E_{[l]}
(x)\big) >\mu  _x$
for all $[l] \ne [l_0]$ in $L/\sim.$  It follows that
$x_j \equiv 0\!\!\!\mod \pi A[u]$ for
all $j\in L.$

\noi Since the unique $[l_0]$ of Step~1 also has $e_{[l_0]}(x)\equiv
0\!\!\!\mod \pi A[l_0],$ it  has $\mu  _x = \v_\bullet \big(E_{[l_0]}(x)\big) =
f(l_0) + p^{\v_p(J_0)}$ with   $J_0$ largest in $[l_0]_x,$ by
Proposition~8(iii).
 Since $\v_\bullet \big(E_{[l]}(x)\big) >\mu  _x$ for
$[l] \not\in M_x,$
 it remains to investigate the $[l_1]\in M_x\bs \{[l_0]\}.$

\vs\noi
{\it Step} 3.1.
There is no $[l_1]\in M_x\bs \{[l_0]\}$ with
$[l_1]_x\ne \emptyset.$ 

\noi
If such an $[l_1]$ did exist, then 
 Proposition 8B(iii) finds the  largest $J_1$ in $[l_1]_x$ with

\noi
\centerline{$f(l_1) + p^{\v_p(J_1)} =
\v_\bullet\big(E_{[l_1]}(x)\big) = \mu  _x =
\v_\bullet\big(E_{[l_0]}(x)\big) =
f(l_0)
+ p^{\v_p(J_0)}.$}

\noi
Comparing $\v_p$ of the first and last terms above by 
 Proposition~8B(v) implies that $\v_p(J_1) =
\v_p(J_0),$  hence  $f(l_1) = f(l_0).$ Thus Lemma~8A(i) implies $[l_1] =
[l_0],$ contradiction.

\vs\noi
{\it Step} 3.2.  There is no $[l_1]\in M_x\bs \{[l_0]\}$ with
$[l_1]_x=\emptyset.$

\noi
If such an $[l_1]$ did exist, then we have
$$
f(l_1) + p^{s-1} \os 1> f(l_0)
+ p^{\v_p(J_0)} \os 2= \mu  _x \os
3= \v_\bullet\big(E_{[l_1]} (x)\big)\os 4\ge f(l_1) + p^s,
$$
with $\os 1 >$ by Step~1 and Lemma~8A(iv), $\os 2=$ by the start of
the proof of Step~3, $\os 3=$ by $[l_1] \in M_x,$ and $\os 4\ge$ by
Proposition 8B(iv).  This  implies the contradiction
$p^{s-1}\ge p^s.$

\noi
Combining Step~3.1 and Step~3.2 proves the first assertion of Step~3.

\noi
The
 argument of the last paragraph of the
proof of Step~2, now with $\mu  _x = f(l_0) + p^{\v_p(J_0)},$  completes Step~3.

\vs\noi
{\it Step} 4.  If $e_{[l_0]}(x)\equiv 0\!\!\!\mod \pi A[u]$ and
$[l_0]_x=\emptyset$
then 
$x_j \equiv 0\!\!\!\mod \pi A[u]$ for all
$j\in L.$ 

\noi
This is apparently false, e.g.: 
$
 \os 1\ge f(l_0) + p^s \ge 0+p^{r+t}
$
 \hfill $ \square$

\newpage
\vsk\noi
\section*{Appendix A}

{\rm The following lemma repeats [RW2, Lemma 2] and provides a more
detailed proof.}

\vs\noi
{\sc Lemma A1.} {\it Let $\chi  $ be a $\Qbar_p\,^c$-valued
character of $G_\infty  $ with open kernel and $F/\Qbar_p$ a finite
field extension so that $\chi  $ has a realization $M_\chi  $ ove
the ring $\fo$ of integers of $F.$  Moreover, let $A$ be a finitely
generated left $\Z_p\big[[G_\infty  ]\big]$-module.  Then
$\Hom_{\fo [H]}(M_\chi  ,\fo \otimes _{\Z_p} A)$ is a finitely
generated $\fo \big[[\Gamma  _k[\big]$-module such that $(\gamma  _k
f)(m) = \gamma  \cdot f(\gamma  ^{-1}m)$ for $\gamma  \!\!\!\mod H =
\gamma  _k\,.$}

\vs\noi
\begin{proof} We show that the $\Gamma  _k$-action above extends to
an $\fo\big[[\Gamma  _k]\big]$-action on
$$
\SM_\chi   =\,\Hom_{\fo [H]}(M_\chi  ,\fo \otimes_{\Z_p} A).
$$
Choose a central open $\Gamma  $ in the kernel of $\chi  ,$ with
image $\ol\Gamma  $ in $\Gamma  _k.$  Then $(\ol \gamma  f)(m) =
\gamma  \cdot f(m),$ for $\gamma    \in \Gamma  ,$ $f\in \SM_\chi  $
$m\in M_\chi  \,.$  Since $A$ is a $\Z_p\big[[G_\infty 
]\big]$-module $\SM_\chi  $ is an $\fo\big[[\Gamma  ]\big]$-module
and finite generation follows.  The lemma will now result from
$$
\fo\big[[\Gamma  _k]\big] = \bigoplus_{0\le i<n} \,\fo\big[[\,\ol
\Gamma  \,]\big]\gamma ^i _k\,,
$$
with $n = [\Gamma  _k:\ol \Gamma  \,]$ and $\gamma  _k$ a chosen
generator of $\Gamma  _k\,.$

\noi
First, observe that if $\SM_\chi  $ has the desired $\fo\big[[\Gamma
 _k]\big]$-module structure, then every $x\in \fo\big[[\Gamma 
_k]\big]$ has $x=\sum_{0\le i<n} \ol x_i\gamma  ^i$ with unique
$x_i\in \fo\big[[\Gamma  ]\big],$ and
$$
(xf)(m) =\sum\nolimits_{0\le i <n}(x_i\gamma  ^i)\cdot f(\gamma  ^{-i}m),
\;{\rm for all}\; m\in M_\chi  \,,
\leqno{\rm (*)}
$$
where $\gamma  $ is the unique generator of $\Gamma  $ with
$\ol\gamma  =\gamma  _k\,.$  For
$$
\begin{aligned}
(xf)(m) &=\sum_i \big((\,\ol x_i\gamma  ^i_k)f\big) (m)
= \sum_i\big(\,\ol x_i\cdot (\gamma  ^i_kf)\big)(m)\\
&= \sum_i x_i\cdot \big(\gamma  ^i\cdot f(\gamma  ^{-i}m)\big)
= \sum_i (x_i\gamma  ^i)\cdot f(\gamma  ^{-i}m).
\end{aligned}
$$
Conversely, we show that property $(*
)$ defines an $\fo\big[[\Gamma  _k]\big]$-structure on $\SM_\chi 
\,,$ i.e., that (i) $(x+y)f = xf+ yf,$ (ii)~$x(f+f^\pr) = xf +
xf^\pr,$ and (iii)~$(xy)f = x(yf),$ hold for all $x,y\in
\fo\big[[\Gamma  _k]\big]$ and $f,f^\pr\in \SM_\chi  .$  As (i),
(ii) are straightforward, we concentrate on (iii).  Given $x=
\sum_{0\le i< n} \ol x  _i\gamma  ^i_k,$ $y=\sum_{0\le j<n} \ol
y_j\gamma  ^j_k$ in $\fo\big[[\Gamma  _k]\big],$ we first observe
that $\gamma  $ central in $G_\infty  $ implies
$$
\Big(\sum_i x_i\gamma  ^i\Big)\Big(\sum_j y_j\gamma  ^j\Big)
= \sum_{(i,j)} x_iy_j\gamma  ^{i+j} =\sum_{0\le t<n} z_t\gamma  ^t
+\sum_{0\le t<n} z_{t+n}\gamma  ^{t+n},
$$
with $z_t=\sum_{\os{(i,j)}{i+j=t}}x_iy_j,\,z_{t+n}= \sum_{\os
{(i,j)}{i+j=t+n}}x_iy_j,$ which equals $\sum_{0\le
t<n}(z_t + z_{t+n}\gamma  ^n)\gamma  ^t$ with $z_t + z_{t+n}\gamma 
^n \in \fo\big[[\Gamma  ]\big],$ since $\gamma  ^n_k\in \ol \Gamma 
.$  This implies that $xy = \sum_t (\,\ol z_t+\ol z_{t+n}\gamma 
^n_k)\gamma  ^t_k\,,$ hence $\big((xy)f\big)(m) = \sum_t(z_t +
z_{t+n}\gamma  ^n\big)\gamma  ^t\cdot f(\gamma  ^{-t}m).$  On the
other hand,
$$
\begin{aligned}
\big(x(yf)\big)(m) &= \sum\nolimits_i(x_i\gamma  ^i)\cdot
(yf)(\gamma  ^{-i}m)\\
&=\sum\nolimits_i(x_i\gamma  ^i)\cdot\big(\sum\nolimits_j(y_j\gamma 
^j)\cdot f\big(\gamma  ^{-j}\gamma  ^{-i}m)\big)
= \sum\nolimits_{(i,j)}(x_iy_j\gamma  ^{i+j})\cdot f\big(\gamma 
^{-(i+j)}m\big)\\
&=\sum\nolimits _t(z_t\gamma  ^t)\cdot f(\gamma  ^{-t}m)
+\sum\nolimits _t (z_{t+n}\gamma  ^{t+n})\cdot f(\gamma 
^{-(t+n)}m)\\
&\dot = \sum\nolimits_t(z_t + z_{t+n}\gamma  ^n)\gamma  ^t\cdot
f(\gamma  ^{-t}m),
\end{aligned}
$$
with $\dot=$ by $\gamma  ^n\in \Gamma  \subseteq \ker(\chi  ).$ 
Thus $\SM_\chi  $ is an $\fo\big[[\Gamma  _k]\big]$-module.

\noi
In fact, this is the unique $\fo\big[[\Gamma  _k]\big]$-module
structure on $\SM_\chi  $ with the prescribed
$\Gamma  _k$-module structure because $(*)$ applies to every choice
of $\Gamma
 $ and $\gamma  _k.$
\end{proof}
\section*{Appendix B}

\vs\noi
{\sc Proposition B1.} {\it The centre fields of the Wedderburn components
of $\SQ G$ all have cohomological dimensions $3.$ }

\vs\noi
\begin{proof}  Choose $\Gamma  \simeq \Z_p$ central open in $G,$
hence $\dim_{\Qbar\Gamma  }\SQ G=[G:\Gamma  ],$ and each Wedderburn
component of $\SQ G$ is $\SQ G\cdot e$ with $e$ a primitive central
idempotent of $\SQ G.$  It follows that $\SQ \Gamma  \cdot
e\subseteq \,\cent(\SQ G\cdot e)$ with $\dim_{\SQ \Gamma  \cdot
e}\SQ G\cdot e$ finite, hence the centre field of $\SQ G\cdot e$ is a
finite extension of $\SQ \Gamma  \cdot e\simeq \SQ \Gamma  ,$ so has
the same cohomological dimension as $\SQ \Gamma  ,$ by [Se2,
Proposition 14(ii), p.17], which is $3$  (see [Lau2].
\end{proof}

\vs\noi
{\sc Proposition B2.} {\it If $d\in D^\times$ has $d\Delta 
d^{-1}=\Delta  ,$ then $d\in \Pi^{\Z}\Delta ^\times \fq(A)^\times.$}

\vs\noi
\begin{proof}
Let $S$ be the multiplicative set $\{\pi^n:n\ge 0\}$ in $A,$ hence 
disc$(S^{-1}\Delta  / S^{-1}A)= S^{-1}A$ by 1. of the proof of
Lemma~2A(ii), and $S^{-1}A$ is a PID by [Ma, Theorem 20.8 and
Theorem 20.1].  We next establish

\begin{description}
\item  {1.} Picent$(S^{-1}\Delta  )=1$
\item  {2.} $(S^{-1}\Delta  )^\times = \Pi^{\Z}\Delta  ^\times.$
\end{description}

\vs\noi
{\it For} 1: We apply [Re, Theorem 37.28] to the maximal
$S^{-1}A$-order $S^{-1}\Delta  ,$ observing first that
Picent$\big((S^{-1}\Delta  )_P\big)=1$ for all maximal ideals $P$ of
$S^{-1}A,$ by Step 2 of the proof of loc.cit.. It follows, from
37.29, that Picent$(S^{-1}\Delta  ) =$ Picent$(S^{-1}A) =1,$ by
Corollary 37.24.

\vs\noi
{\it For} 2: We first observe that $A_\bullet\cap S^{-1}A=A,$ by
Theorem~4.25 of loc.cit., which implies that $\Delta  _\bullet
\cap S^{-1}\Delta  =\Delta  ,$ because an $A$-basis of $\Delta  $ is
also an $A_\bullet$-basis of $\Delta  _\bullet$ and an
$S^{-1}A$-basis of $S^{-1}\Delta  .$  Also $\Pi^{-1}\in
(S^{-1}\Delta  )^\times,$ since $\pi\in \Pi B$ implies $\Pi^{-1}\in
\pi^{-1}B\subseteq S^{-1}\Delta  .$  Thus, if $d\in (S^{-1}\Delta 
)^\times,$ then $\Pi^{-v(d)}d\in \Delta  ^\times_\bullet,$ in the
notation of Corollary 2B,hence $\Pi^{-v(d)}d\in \Delta 
^\times_\bullet \cap (S^{-1}\Delta  )^\times \subseteq (\Delta 
_\bullet \cap S^{-1}\Delta  )^\times =\Delta  ^\times,$ settling $2.$

\vs\noi
Finally, from $d\Delta  d^{-1}=\Delta  ,$ it follows that $d\cdot
S^{-1}\Delta  \cdot d^{-1}=S^{-1}\Delta  ,$ hence $d\in
(S^{-1}\Delta  )^\times \fq(A)^\times,$ by 1. and Theorem (37.25) of
loc.cit. (as $\omega^\pr$ is injective), with $(S^{-1}\Delta  )^\times
\fq(A)^\times = \Pi^{\Z}\Delta  ^\times \fq(A)^\times,$ by 2.  Note
that $\Delta  ^\times \fq(A)^\times$ is a normal subgroup of
$\Pi^{\Z}\Delta  ^\times \fq(A)^\times$ with index $p^s,$ since
$\Pi^{p^s}\in B^\times \pi\subseteq \Delta  ^\times \fq(A)^\times.$
\end{proof}

\section*{Appendix C}

\vs
\noi
This joint paper
is unfinished: it is now a quarry which can be
better understood in a somewhat different order.  The Abstract and
Introduction serve as a short informal sketch of the motivation and
the goals of this quarry, and all references are internal.
The main unfinished aspect is that there may be errors of a few
kinds,
because the structure must be determined exactly.

\vsk\noi
\begin{description}
\item{a:}
\S1 starts with the formal statement Theorem~1A of the $SK_1(QG)=0$
goal to the vanishing of $SK_1(D_{r,s,\IF})$ for the division
algebras $D_{r,s,\IF}$ that are described above, based on [Lau1] and
[Lau2], and to an example $D_{r,s,\IF}$ to analyze.  
This quickly leads to the definition of $\Phi$ and of
weight (either ordinary or transitional), and the development until Lemma~1F.  
\newline
However the role of $A[u]$ in Lemma~7A in \S7 requires knowledge of
Lemma~1G.  This phenomenon will reappear later, especially in
Lemma~2A and Corollary~2B.

\vs\noi
\item{b:}
\S7 concerns the oldest collection of commutators $[u,1-\Pi^{p^tm} x]$ 
organized by type $t\ge 0$ in Lemma~7A.  This already handles
Corollary~7B, so reducing to $(*)$ when $t\ge 1.$
\newline
Now Proposition~7C divides $(*)$ into two parts (ii) and (iii), of
different kinds, which can be used to decide whether
$d\in \Delta^\times$ with $nr(d)=1$ and level $n$ of type $t$ is a
$\Pi$-commutator $\!\!\mod~\Pi^{n+1}\Delta.  $
\newline
An example of this, of type 1, in Corollary~7D,  and a valuation
$v_\pi,$ in Corollary~2B(iii), completes \S7.     

\vs\noi
\item{c:}
\S2 is about some of the many properties of the $A$-order $\Delta$
on $D,$ the valuations on its maximal commutative $A$-suborders,
discriminants, homological algebra,
and
the Galois action on them.
\newline
Further, there are studies of several topics which are needed only
later, in the spirit of Lemma~1G in \S1.  In particular the Remark
after Theorem~2F needs attention.

\vs\noi
\S4 links the reduced norm of $d\in \Delta  ^\times$ to the level of
$d.$  The main result is Theorem~4G; in particular, if $d\in \Delta 
^\times$ has $nr(d) = 1$ then $d\in 1+\Pi^{p^r-1}\Delta  .$  

\noi\item{d:}
\S5 is technical preparation for \S6 (and beyond), which, by
Proposition~5C, expresses the $A$-coefficients $\mod
\pi^{\Phi(n)+1}A$ (in Theorem~4C) by structure revealing $\lambda 
_n^{(k)}(b)$'s. Then Proposition~5D shows that $\lambda 
_n(b_{p^w\v})
\equiv 0\mod \pi A$ for $0\le \v < p^{s-w},$ whenever $n$ is
ordinary of weight $w.$  However, if $n$ is transitional of weight
$w,.$  then $\lambda  ^+_n$ is more complicated; fortunately,
transitional $n$ are rare.

\noi
\S3 uses the Hom description of $K_1(\Lambda_\SO G)$ and the
Mayer-Victoris sequence to prove that $SK_1(\Delta) \to
K_1(\Delta/\pi_1\Delta)$ is injective to deduce, in Theorem~3E(ii),
that if $d\in 1+\Pi^{p^{r+s-1}}\Delta$ has $nr(d)=1$ then $d\in
[\Delta^\times, \Delta^\times].$

\noi\item{e:} \S6 benefits from the above preparation, which can
deal with ordinary levels, by Proposition~6B, and critical levels,
by Theorem~6D.  This leads to the concept of the type $t$ of level
$n,$ defined by $n=p^{r+t} + p^tm -1,$ for some $m\ge 1$ with
$m\not\equiv 0\!\!\!\mod p.$  Then Proposition~6F shows that the
remaining levels (not yet dealt with) have $p^r-1\le n<p^{r+s-1}$
and type $t\le$ the weight $w$ of $n.$  The following Remark
even provides a list of all of the remaining levels, and formulates
a hypothesis, which Corollary~6G and its Remark would find decisive.

\S8 takes the first step toward the program of Proposition~7C, in
Lemma~8A, by using the function $f:L\to\Z$ of Proposition~7C(i). 
The equivalence classes induced by $f$ are then described as the
orbits of a special function $\tau :L\to L,$ leading toward a
bijection $L/_\sim \to L_0.$  The next step Proposition~8b, interprets
the congruence of Proposition~7C(iii)\! as 
$\us{[\ell]\in L/\sim}\sum E_{[\ell]} (x) \equiv 0\!\!\!\mod
\Pi^{p^{r+t}-p^t}\Delta$ with 
$$
E_{[\ell]}(x)\! =\! \us{j\in [\ell]}\sum \big[\,^g\Pi^j - \Pi^j]x_j
\subseteq \Delta \q\text{and} \q e_{[\ell]} (x) = {\us{j\in[\ell]}\sum} j^\pr
x_j \subseteq A[u],
$$
 which is studied further, notably in (iii) and
(vi).  The rest of \S8 prepares for \S9.

\noi\item{f:} \S9 We take first steps generalizing Corollary~7D for
$t=1,$ to search for solutions of congruence 7C(ii) for $t\ge 2,$
using the methods of \S8.  Solutions for the congruence of
Theorem~9A(ii) seem accessible but need more structure.

\S10 concerns itself with two topics, the first with families of
commutators of higher type $t,$ and the second, in Lemma~10C, about
the failure of Theorem~9C when $s<r+t.$
\end{description}

\noi
We can improve on Lemma~1C(ii) by

\vs\noi
{\sc Lemma C1}$^{\#}$(ii).
{\it If $\Phi_k(n) = \Phi_{k+1}(n)$ then }
$n =
\left(
\begin{matrix} p^k(p^r-1), & \text{\it if}\q n_k=p-1.\\
p^{r-1}, &\!\!\!\!\!\!\!\!\!\text{\it if} \q
n_0=0.\end{matrix}\right.
$
{\it provided that $\ell(n) >k.$}

\vs\noi  \begin{proof} 
 The $p$-adic expansion of $n\ge 0$ has
$n=\sum\limits^{\ell}_{i=0}\! n_ip^i,$ with $0\!\le\!n_i\!<\!p$ for all $i,$
and $ \ell(n)=
\max\{i:n_i>0\}.$  Then $\sum\limits_{0\le i<k}
n_ip^i\le \sum\limits_i (p-1)p^i = (p-1) \,\frac{p^k-1}{p-1} = p^k-1
<p^k,$ for $0 \le k\le \ell =\ell(n),$ hence
$$
\Big\lfloor \frac {n}{p^k}\Big\rfloor = \Big\lfloor \frac{\us{0\le i
<k}\Sigma n_ip^i +\us{k\le i\le \ell}\Sigma\, n_i p^i}{p^k}\Big\rfloor = 0
+\sum_{k\le i\le \ell} n_i p^{i-k}, \q\text{\rm and}
$$
$$
\Big\lfloor \frac{n}{p^{k+1}}\Big\rfloor =\sum_{k+1\le i\le \ell}
n_i p^{i-(k+1)} \q \text{\rm similarly}.
$$

\vs\noi
This implies $\lfloor \frac{n}{p^k}\rfloor -\lfloor
\frac{n}{p^{k+1}}\rfloor \:
=\;\sum\limits_{k\le i\le\ell} n_i p^{i-k}
-\sum\limits_{k+1\le i\le\ell} n_ip^{i-(k+1)}
=n_k\;+\;\sum\limits_{k+1\le i\le \ell} n_i p^{i-k} \!
-\!\sum\limits_{k<i\le \ell} n_ip^{i-(k+1)}
=\! n_k + p\sum\limits_{k<i\le \ell} n_ip^{i-(k+1)}
-\sum\limits_{k<i\le \ell} n_i p^{i-(k+1)}
= n_k + (p-1)\sum\limits_{k<i\le \ell} n_i p^{i-(k+1)}.$

\vs\noi
Since $\Phi_k(n) =\Phi_{k+1}(n)$ implies 
$\lfloor \frac{n}{p^k}\rfloor - \lfloor
\frac{n}{p^{k+1}}\rfloor 
= \phi(p^r) =(p-1)p^{r-1},$ by Lemma~1B(i),
we have
$n_k + (p-1)\sum\limits_{k<i\le\ell} n_ip^{i-(k+1)} 
= (p-1)p^{r-1},$  hence
$$
n_k = (p-1)\big(p^{r-1} -\sum_{k<i\le \ell} n_i p^{i-(k+1)}\big).
$$
Thus $p-1$ divides $n_k$ so there are two possibilities:

\vs\noi
{\it For} 3: $n_k = p-1$ and $p^{r-1} -1 = \us{k<i\le \ell}\sum n_i
p^{i-(k+1)}.$

\vs\noi
Comparing this to $p^{r-1} -1 =
{\sum\limits^{r-2}_{j=0}}(p-1)p^j,$ the uniqueness of $p$-adic
expansion (of $p^{r-1}-1)$ tells us that 
 $n_i = p-1$ when $k<i\le \ell$ and
$\ell- (k+1) = r-2,$ hence $\ell = r+k-1.$  Combining with $n_k =
p-1$ this becomes $n_i = p-1$ for $k\le i \le r+k-1,$ and $n_i=0$
when $0\le i<k.$
Thus $n={\sum\limits^\ell_{i=0}} n_ip^i 
= (p-1) {\sum\limits^{r+k-1}_{i=k}} p^i 
= (p-1) \Big({\sum\limits^{r+k-1}_{i=0}} p^i 
- {\sum\limits^{k-1}_{i-0}} p^i\Big)
= (p-1) \Big(\frac{p^{r+k}-1}{p-1} - \frac{p^k-1}{p-1}\Big)
= p^{k+r} - p^r = p^k(p^r-1).$

\vs\noi
{\it For} 4: $n_0 =0$ and $p^{r-1} =\us{0<i\le \ell}\sum n_i
p^{i-(k+1)}.$

\vs\noi
Now $p^{r-1}$ has $p$-adic expansion $p^{r-1} = 1p^{r-1}$ hence $n_i
=0$ for $i<\ell,$  so $\ell = r-1.$  Thus $p^{r-1} = 1p^{r-1}.$
\end{proof}

\vs\noi
{\sc Remark} 1: This result is simpler than Lemma~1C(ii), 
and could influence the weight concept.  It also contradicts
Lemma~1C(ii) at $n=0.$

\vs\noi
{\sc Remark} 2: Keeping track of quotation data is essential when it
changes: for example, Lemma~5E(ii) quotes Lemma~2G(iv), or
Proposition~7C(i) quotes Lemma~5G(i) quotes Lemma~2G(iv).  This is
as stated right now, but of course could change.

\vs\noi
{\sc Remark} 3: This consists of the results which form the
structure.  Much of it is correct, and more of it a consequence of
that, by some of it mysterious, e.g. Proposition~5D(ii).  But, for
example, some important material, e.g. Proposition~5D(i), appeals to
lost details in the proof from Proposition~5C and the paragraph
on the way to Proposition~5D(i).  And everything in between.

\vs\noi
{\sc Remark} 4. Dubrovin valuation rings appeared already  in Corollary~7D
and 7D$^*$, and will certainly have a bigger role. But the influence
of H.~Brungs,  a long term colleague, has been lost on me.  
There is no time now to remedy this.

\vskip.5in\noi
Department of Mathematics, University of Alberta, Edmonton, AB,
Canada \q T6G 2G1

\noi
E-mail address: weissa@ualberta.ca
\end{document}